\documentclass{memo-l}

\usepackage{array}
\usepackage{float}
\usepackage{longtable} \setlongtables
\usepackage{hhline}
\usepackage{graphics}
\usepackage{color}
\usepackage{amssymb}
\usepackage[matrix,arrow,curve]{xy}
\usepackage{epsfig}
\usepackage{euscript}
\usepackage{xspace}
\usepackage{ulem} 

\numberwithin{equation}{section}
\newtheoremstyle{my}{1.5em}{0.5em}{\em}{}{\sc}{.}{0.5em}{}
\newtheoremstyle{mydef}{1.5em}{0.5em}{}{}{\sc}{.}{0.5em}{}

\theoremstyle{my}
\newtheorem{thm}{Theorem}[section]
\newtheorem{theorem}[thm]{Theorem}
\newtheorem{cor}[thm]{Corollary}
\newtheorem{corollary}[thm]{Corollary}
\newtheorem{lemma}[thm]{Lemma}

\newtheorem{prop}[thm]{Proposition}
\newtheorem{proposition}[thm]{Proposition}
\newtheorem{assumptions}[thm]{Assumptions}
\newtheorem{assumption}[thm]{Assumption}

\newtheorem{defn}[thm]{Definition}

\newtheorem{notation}[thm]{Notation}

\newtheorem{remark}[thm]{Remark}

\newtheorem{example}[thm]{Example}

\newcommand{\acknowledgments}{{\em Acknowledgments.} }

\newcommand{\FIGURE}[3]{%
\begin{figure}[#3]
\begin{center}
\epsfig{file=#2} \\ \caption{\label{#1}}
\end{center}
\end{figure}}

{

\begin{enumerate}
\itemsep1em \leftmargin2em
}{\end{enumerate}}

{\begin{list}{{\rm(\alph{enumi}) }}{\usecounter{enumi}

\leftmargin0cm \labelsep0cm \rightmargin0cm \parsep1em
\listparindent0em \itemsep-1em \topsep1em \parskip1em
\setlength{\labelwidth}{\fill}}} {\parskip0em \end{list}}

{

\begin{enumerate}}{\end{enumerate}}

\newenvironment{romanlist}%
{

\begin{enumerate} \parsep0cm
\itemsep0cm \parskip0cm}{\end{enumerate}}

{

\begin{enumerate}}{\end{enumerate}}

\newcommand{\R}{\mathbb{R}}
\newcommand{\Z}{\mathbb{Z}}
\newcommand{\Q}{\mathbb{Q}}
\newcommand{\C}{\mathbb{C}}
\newcommand{\N}{\mathbb{N}}
\newcommand{\half}{{\textstyle\frac{1}{2}}}
\newcommand{\quarter}{{\textstyle\frac{1}{4}}}

\newcommand{\iso}{\cong}           
\newcommand{\htp}{\simeq}          
\newcommand{\smooth}{C^\infty}
\newcommand{\CP}[1]{\C {\mathrm P}^{#1}}
\newcommand{\RP}[1]{\R {\mathrm P}^{#1}}
\newcommand{\leftsc}{\langle}
\newcommand{\rightsc}{\rangle}
\newcommand{\Rgeq}{\R^{\scriptscriptstyle \geq 0}}
\newcommand{\Rleq}{\R^{\scriptscriptstyle \leq 0}}
\newcommand{\suchthat}{\; : \;}

\newcommand{\id}{\mathrm{id}}

\newcommand{\im}{\mathrm{im}}

\renewcommand{\o}{\omega}
\renewcommand{\O}{\Omega}
\newcommand{\Aut}{\mathit{Aut}}
\newcommand{\Diff}{\mathit{Diff}}

\newcommand{\Fuk}{{\mathcal F}}

\newcommand{\Cat}{{\mathcal C}}
\newcommand{\A}{\mathcal A}
\newcommand{\B}{\mathcal B}
\newcommand{\T}{\mathcal T}
\newcommand{\F}{\mathcal F}
\newcommand{\G}{\mathcal G}

\newcommand{\OO}{\mathcal O}
\renewcommand{\AA}{\mathfrak A}
\newcommand{\GG}{\mathfrak G}

\newcommand{\g}{\mathfrak g}
\newcommand{\semidirect}{\rtimes}
\newcommand{\inner}{\llcorner\,}

\newcommand{\UU}{\mathfrak U}
\newcommand{\SHom}{\underline{\mathit{Hom}}}
\newcommand{\cech}[1]{\check{#1}}
\newcommand{\I}{\mathcal I}
\renewcommand{\SS}{\mathcal S}

\newcommand{\lcm}{l.c.m.}
\newcommand{\LL}{{\mathfrak L}}
\newcommand{\orient}{\mathrm{or}}

\newcommand{\alphagr}{\tilde{\alpha}}
\newcommand{\Lgr}{\tilde{L}}
\newcommand{\phigr}{\tilde{\phi}}
\newcommand{\psigr}{\tilde{\psi}}
\newcommand{\hgr}{\tilde{h}}
\newcommand{\Autgr}{\widetilde{\Aut}}
\newcommand{\taugr}{\tilde{\tau}}
\renewcommand{\Diff}{\mathit{Diff}}

\renewcommand{\O}{\mathcal{O}}
\newcommand{\K}{\mathcal{K}}

\newcommand{\red}{\mathrm{red}}

\newcommand{\Spin}{\$}
\newcommand{\br}{{\scriptscriptstyle \flat}}
\newcommand{\rbr}{{\scriptscriptstyle \sharp}}
\newcommand{\mbr}{{\scriptscriptstyle \natural}}

\newcommand{\trans}{\pitchfork}
\newcommand{\CC}{\mathfrak{C}}
\newcommand{\MM}{\mathfrak{M}}
\newcommand{\Fol}{\mathcal{F}}
\newcommand{\Dat}{\mathcal{D}}

\newcommand{\QQ}{\mathcal Q}

\newcommand{\One}{I}
\newcommand{\twist}[1]{\langle #1 \rangle}

\begin{document}
\frontmatter

\title{Homological mirror symmetry\\ for the quartic surface}


\author{Paul Seidel}
\address{MIT Department of Mathematics, 77 Massachussetts Ave, Cambridge, MA 02139, USA}
\thanks{Partially supported by NSF grant DMS-1005288 and a Simons Investigator Grant from the Simons Foundation.}

\date{}




\maketitle

\tableofcontents

\begin{abstract}
We prove Kontsevich's form of the mirror symmetry conjecture for (on the symplectic geometry side) a quartic surface in $\C P^3$.
\end{abstract}

\mainmatter

\section{Introduction}

This paper deals with a special case of Kontsevich's ``homological
mirror symmetry'' conjecture \cite{kontsevich94}. Before formulating
the statement, we need to introduce the relevant
coefficient rings.

\begin{notation}
Let $\Lambda_\N = \C[[q]]$ be the ring of formal power series in one
variable. We denote by $\Lambda_\Z$ its quotient field, obtained by
formally inverting $q$. By adjoining roots $q^{1/d}$ of all orders to
$\Lambda_\Z$, one obtains its algebraic closure $\Lambda_\Q$. Elements of
this field are formal series
\begin{equation}
 f(q) = \sum_m a_m q^m
\end{equation}
where $m$ runs over all numbers in $(1/d)\Z \subset \Q$ for some $d
\geq 1$ which depends on $f$, and $a_m \in \C$ vanishes for all
sufficiently negative $m$. Geometrically, one should think of
$\Lambda_\N$ as functions on a small (actually infinitesimally small) disc. Then $\Lambda_\Z$ corresponds to the punctured disc (with the origin removed), and $\Lambda_\Q$ to its universal cover. Consider the semigroup $\mathit{End}(\Lambda_\N)$ of $q$-adically continuous $\C$-algebra
endomorphisms of $\Lambda_\N$. More concretely, these are
substitutions
\begin{equation}
 \psi^* : q \longmapsto \psi(q)
\end{equation}
with $\psi \in q \C[[q]]$. For each nonzero $\psi$, $\psi^*$ extends to an
endomorphism of the quotient field. The $q$-adically continuous
Galois group of $\Lambda_\Z/\C$ is the group of invertibles
$\mathit{End}(\Lambda_\N)^\times$, consisting of those $\psi$ with $\psi'(0)
\neq 0$. Moreover, the Galois group of $\Lambda_\Q/\Lambda_\Z$ is the
profinite group $\hat{\Z} = \mathit{Hom}(\Q/\Z,\C^*)$, whose topological
generator $\hat{1}$ takes $q^m$ to $e^{2\pi i m} q^m$.
\end{notation}

On the symplectic side of mirror symmetry, take any smooth quartic
surface $X_0 \subset \CP{3}$, with its standard symplectic structure.
To this we associate a triangulated category linear over
$\Lambda_\Q$, the split-closed derived Fukaya category
$D^\pi\Fuk(X_0)$, defined using Lagrangian submanifolds of $X_0$ and
pseudo-holomorphic curves with boundary on them (as well as a formal closure process). The use of formal
power series to take into account the area of pseudo-holomorphic
curves is a familiar device (the relevant coefficient rings usually go under the name ``Novikov rings'', see e.g.\ \cite{hofer-salamon95}; what we
are using here is a particularly simple instance). $\Lambda_\Q$ appears because we allow
only a certain class of Lagrangian submanifolds, namely the rational
ones (there is also another restriction, vanishing of Maslov classes,
which is responsible for making the Fukaya category $\Z$-graded).

On the complex side, we start with the quartic surface in ${\mathbb
P}^3_{\Lambda_\Q}$ defined by
\begin{equation} \label{eq:fermat}
y_0 y_1 y_2 y_3 + q(y_0^4 + y_1^4 + y_2^4 + y_3^4) = 0.
\end{equation}
Note that here $q$ is a ``constant'' (an element of the ground
field). The group $\Gamma_{16} = \{ [\mathrm{diag}(\alpha_0,\alpha_1,\alpha_2,\alpha_3)] \suchthat
\alpha_k^4 = 1, \; \alpha_0\alpha_1\alpha_2\alpha_3 = 1 \} \subset \mathit{PSL}(V)$, $\Gamma_{16} \iso
\Z/4 \times \Z/4$, acts on this surface in the obvious way. We denote
by $Z_q^*$ the unique minimal resolution of the quotient orbifold.
Note that $\Lambda_\Q$ is an algebraically closed field of characteristic
zero, so the standard theory of algebraic surfaces applies, including
minimal resolutions. The associated category is the bounded derived
category of coherent sheaves, $D^b\mathit{Coh}(Z_q^*)$, which is also a
triangulated category linear over $\Lambda_\Q$.

\begin{remark}
In fact, $Z_q^*$ can be defined over the smaller field $\Lambda_\Z$. This is easy to see if one argues
``by hand'' as follows. The action of $\Gamma_{16}$ on \eqref{eq:fermat}
has $24$ points with nontrivial isotropy; their isotropy subgroups are isomorphic to $\Z/4$. Each of these points is defined over $\Lambda_\Z$, hence
can be blown up over that field. Take the blowup, and divide it by the action of the subgroup $\Z/2 \times \Z/2 =
\Gamma_{16} \cap \{(\pm 1,\pm 1,\pm 1, \pm1)\}$. The quotient is a regular algebraic variety
over $\Lambda_\Z$, and carries an action of the group $\Gamma_{16}/(\Z_2 \times \Z_2) \iso \Z_2 \times \Z_2$,
with $24$ points that have nontrivial isotropy $\Z/2$. These points are again defined over $\Lambda_\Z$, and one
repeats the previous process to get an algebraic surface which is a $\Lambda_\Z$-model of $Z_q^*$. (There are
also other possible strategies, such as using a toric resolution of ${\mathbb P}^3_{\Lambda_\Z}/\Gamma_{16}$ and
realizing the desired algebraic surface inside that).
\end{remark}

\begin{theorem} \label{th:main}
There is a $\psi \in \mathit{End}(\Lambda_\N)^\times$ and an equivalence of
triangulated categories,
\begin{equation}
D^\pi\Fuk(X_0) \iso \hat{\psi}^*D^b\mathit{Coh}(Z_q^*).
\end{equation}
\end{theorem}

Here $\hat\psi$ is a lift of $\psi$ to an automorphism of
$\Lambda_\Q$, which we use to change the $\Lambda_\Q$-module
structure of $D^b\mathit{Coh}(Z_q^*)$, by letting $\hat\psi^*f$ act instead of
$f$. The theorem says that the outcome of this ``reparametrization'' is
equivalent to $D^\pi\Fuk(X_0)$. The choice of lift $\hat\psi$ is
irrelevant because $D^b\mathit{Coh}(Z_q^*)$ carries a $\hat{\Z}$-action (and
so does $D^\pi\Fuk(X_0)$ by \cite{fukaya03}).

The function $\psi$ which occurs in the theorem is unknown at present
(or rather, it is not known to agree with the standard
``mirror map'' obtained from the Picard-Fuchs equation for $Z_q^*$,
see \cite{nagura-sugiyama95} or the more recent \cite{hartmann11}). Nevertheless, certain consequences can already be
drawn from the result as given. Let $Aut(X_0)$ be the group of symplectic automorphisms. This admits a natural central
extension, the graded symplectic automorphism group \cite{seidel99}
\begin{equation} 
1 \rightarrow \Z \longrightarrow \Autgr(X_0) \longrightarrow Aut(X_0) \rightarrow 1.
\end{equation}
Let $\mathcal{M}^*$ be the classifying space (or moduli stack, depending on the reader's
preference) for smooth $K3$ surfaces equipped with an ample
cohomology class $A$ whose square is $A \cdot A = 4$, and with a
choice of nonzero holomorphic two-form. By applying Moser's theorem, one sees that $\mathcal{M}^*$ carries a fibre bundle
with structure group $Aut(X_0)$, and the presence of the holomorphic two-forms on the fibres provides a natural lift of this
to $\widetilde{Aut}(X_0)$. In particular, there is an associated parallel transport map
\begin{equation} \label{eq:mono-1}
\pi_1(\mathcal{M}^*) \longrightarrow \pi_0(\Autgr(X_0)).
\end{equation}
On the other hand, graded symplectic automorphisms act on the Fukaya category. Denote by 
$Aut(D^\pi Fuk(X_0)/\Lambda_\Q)$ the group of isomorphism classes of exact (with respect to the triangulated structure) 
$\Lambda_\Q$-linear autoequivalences. We use the same notation on the other side of mirror symmetry. Then, the outcome is a map
\begin{equation} \label{eq:mono-2}
\pi_0(\widetilde{Aut}(X_0)) \longrightarrow Aut(D^\pi Fuk(X_0)/\Lambda_\Q) \iso Aut(D^b\mathit{Coh}(Z_q^*)/\Lambda_\Q).
\end{equation}
(If one is willing to divide the autoequivalence groups by the subgroups of even degree shifts, one can drop the adjective ``graded'' and 
define such a map on $\pi_0(\Aut(X_0))$ itself). Finally, combining \eqref{eq:mono-1} and \eqref{eq:mono-2} yields a homomorphism
\begin{equation} \label{eq:alg-monodromy}
\pi_1(\mathcal{M}^*) \longrightarrow Aut(D^b\mathit{Coh}(Z_q^*)/\Lambda_\Q).
\end{equation}
This is interesting because the objects on both sides are algebro-geometric, but the construction itself necessarily seems to 
go through symplectic geometry and mirror symmetry. The existence of \eqref{eq:alg-monodromy} had been postulated by
several people, starting with Kontsevich himself, and motivated a large amount of work on autoequivalences of derived
categories, see e.g.\ \cite{seidel-thomas99, horja99, aspinwall01,hosono04}; however, the kernel and cokernel of that map remain unknown. Another direct implication
of Theorem \ref{th:main} is the equality of certain simpler
invariants, such as the Grothendieck groups
\begin{equation}
K_0(D^\pi\Fuk(X_0)) \iso K_0(D^b \mathit{Coh}(Z_q^*)) = K_0(Z_q^*).
\end{equation}
Every object of the Fukaya category (in particular, every oriented
Lagrangian sphere in $X_0$) determines a class in
$K_0(D^\pi\Fuk(X_0))$, and $\pi_0(Aut(X_0))$ acts linearly on this
group. Recall that Grothendieck's theorem says that $K_0(Z_q^*) \otimes \Q \iso CH_*(Z_q^*)
\otimes \Q$, and that Mumford's work shows that $CH_0(Z_q^*)$ is quite an interesting object.

Having stated our result, it is time to review briefly some
work on homological mirror symmetry preceding this paper. The first case to be considered was that of elliptic curves, where Polishchuk and Zaslow
\cite{polishchuk-zaslow98} proved a result roughly similar to Theorem
\ref{th:main}, with an explicit mirror map $\psi$ but not taking into
account the triangulated structure (although subsequent work on Massey
products \cite{polishchuk98} largely made up for this deficiency).
Many beautiful and deep results on abelian varieties have been obtained
\cite{fukaya02b, kontsevich-soibelman00, fukaya02, nishinou03}, giving partial descriptions of their Fukaya categories. A common thread of all these papers has been to establish an explicit
bijection between objects on both sides (Lagrangian submanifolds and
sheaves). The prevailing view is that this correspondence should be
set up as a Fourier-Mukai type transform associated to a
Strominger-Yau-Zaslow torus fibration. This approach was pursued in
\cite{kontsevich-soibelman00}, in the situation where the SYZ
fibration has no singularities, and where the Lagrangian submanifolds
are supposed to be transverse to the fibres. Both restrictions are
serious ones from the point of symplectic topology, and there has been some
work on how to remove them, see in particular \cite{fukaya02b}.

Actually, this entire circle of ideas is of little importance for our proof (although the SYZ conjecture does serve as a source of intuition in some places). We rely instead on techniques from symplectic topology and homological algebra, which are maybe closer to the spirit of Kontsevich's original paper. The computation on the algebraic geometry side involves a classical idea of Beilinson, together with a Massey product computation which the author learned from \cite{douglas-govindarajan-jayaraman-tomasiello01}. On the symplectic side, where the overwhelming majority of the work is invested, there are three steps, namely (slightly simplified, and in reverse logical order):
\begin{itemize} \itemsep1em \parskip1em

\item
First, we prove that $D^\pi\Fuk(X_0)$ can be entirely reconstructed
from the full $A_\infty$-subcategory of $\Fuk(X_0)$ consisting of a
particular set of 64 Lagrangian two-spheres, which are vanishing
cycles for the standard Fermat pencil (the choice of these cycles is
best explained by physics considerations, compare for instance
\cite{hori-iqbal-vafa00}). There is a general algebraic part to this, based on the notions of (spherical) twist functors and split-generators: we formulate this in Lemma \ref{th:twist-generators}, which generalizes \cite[Corollary 5.8]{seidel04}. The translation to geometry depends on the identification of Dehn twist and twist functors (at least on the level of the action of objects), which is our Proposition \ref{th:exact-sequence-2}, generalizing \cite[Corollary 17.17]{seidel04}. The additional geometric ingredient is a ``negativity'' property of the monodromy of the Fermat pencil around the ``large complex structure limit'', which is the subject of Section \ref{sec:negativity}. Corollary \ref{th:generates-everything} summarizes the outcome.

\item
By construction, all 64 vanishing cycles lie in an affine
Zariski-open subset $M_0 \subset X_0$. Following a proposal from
\cite{seidel02}, the relation between the Fukaya categories
$\Fuk(M_0)$ and $\Fuk(X_0)$ can be formulated in terms of the
deformation theory of $A_\infty$-structures (this recourse to
abstract deformation theory is the reason why we cannot determine
$\psi$ explicitly). The algebraic part of the deformation theory is set up in Section \ref{sec:deformations}, in particular Lemma \ref{th:versal-3} (in the end, we will find it convenient to use a variation on the theory which exploits additional finite symmetries present in our case, see Lemma \ref{th:recognize-q64-2}). For the geometric application, one has to introduce a ``relative Fukaya category'' which interpolates between the Fukaya categories of $M_0$ and $X_0$. This is done in Section \ref{subsec:relative-fukaya}. Then, a more specific argument shows that the resulting deformation of $\Fuk(M_0)$ is nontrivial to first order in the parameter $q$ (which implies that $\psi'(0) \neq 0$). Very roughly speaking, our argument for nontriviality depends on having non-isomorphic objects of $\Fuk(M_0)$ which survive to $\Fuk(X_0)$ and become isomorphic there. The precise algebraic criterion is given in Lemmas \ref{th:product-turns-on} and \ref{th:product-turns-on-2}; the geometric strategy is developed in Sections \ref{subsec:pss} and \ref{subsec:skip-divisor}, see in particular Corollary \ref{th:wall-crossing-2}; and its specific application to our example is carried out in Section \ref{subsec:nontrivial-geometric}.

\item
Having reduced our problem to a finite collection of Lagrangian
two-spheres in the affine four-manifold $M_0$, we apply a general
dimensional induction argument from \cite{seidel04} to compute the
relevant full $A_\infty$-subcategory of $\Fuk(M_0)$. Even though that method yields a complete computation of the $A_\infty$-structure, it is not practical to carry out all of it, so we use the general classification theory of $A_\infty$-structures (a close relative of the previous deformation theory) to reduce the necessary information to knowing the Floer cohomology groups, their product structure, and a single higher order product. The algebraic idea is introduced in Lemma \ref{th:versal-1}, and then adapted to the specific application in Lemma \ref{th:recognize-q4}. The actual Floer cohomology computation is based on Picard-Lefschetz theory one dimension down, which means fibering our $M_0$ by a pencil of affine algebraic curves. The elementary geometric theory is explained in Section \ref{subsec:braid-monodromy}, and the relation with Fukaya categories in Section \ref{subsec:induction}. In practice, the first step is to draw the actual vanishing cycles (as circles on the fibre, which is a topological surface). This is done by a braid monodromy method, based on explicitly solving the equation for the branch curve \eqref{eq:elimination-theory}. We then count polygons bounding the vanishing cycles by hand (see the figures and tables in Section \ref{subsec:computation}). The outcome is plugged into a purely algebraic computation involving twisted complexes \eqref{eq:s-objects}, which we carry out using a computer. The reason why this step in the proof is less transparent is that, unlike the previous two, it has no known meaning in terms of mirror symmetry.
\end{itemize}

Overall, the organisation of the paper is as follows. Section \ref{sec:generators} contains the necessary background material from homological algebra. Section \ref{sec:deformations} is devoted to algebraic deformation theory. Section \ref{sec:group-actions} concerns finite group actions on algebras, in particular exterior algebras, which are the specific algebraic structures relevant for our argument. Section \ref{sec:sheaves} contains the rather modest algebro-geometric background we need. Section \ref{sec:geometry} covers elementary symplectic geometry, including Picard-Lefschetz theory. Section \ref{sec:negativity} is more specifically about gradings and their application to monodromy maps. Section \ref{sec:fukaya} explains the basics of Fukaya categories in the affine, projective and relative sense, and their relationship. Section \ref{sec:induction} describes the dimensional induction machinery from \cite{seidel04}. The final part (Sections
\ref{sec:4-64}--\ref{sec:computation}) contains the computations which are really specific to the case of the quartic surface.

\acknowledgments This paper could not have been written without the
assistance of Michael Douglas, who showed me the ``quiver
presentation'' of the derived category of coherent sheaves which
appears as the algebras $Q_4$ and $Q_{64}$ in this paper, and also
the crucial Massey product formula (Lemma \ref{th:massey}) from
\cite{douglas-govindarajan-jayaraman-tomasiello01}. Maxim Kontsevich
helped me in many ways, one of them being the Hochschild cohomology
computations for semidirect products (Lemma \ref{th:twisted-hkr}).
Simon Donaldson introduced me to the notion of matching cycle, on
which the main dimensional induction procedure is based. To all of
them, thanks! ETH Z\"urich invited me to give a lecture course on
Fukaya categories, and provided a responsive audience on which parts
of this paper were tested. The following computer algebra packages were
used: {\sc singular} for elimination theory, {\sc mathematica} for
visualization, and the Python {\sc sympy} module for matrix computations.

{\em Comments on the revised version.} Because of the need for completing \cite{seidel04} and for other reasons, this paper has languished in preprint form for years (2003-2011). I have not tried to update the introduction or references to reflect more recent developments, since at least part of those developments (Sheridan's recent work, for instance) were based on this paper, which puts one in danger of creating a ``closed loop'' of references. In the body of the paper, I have made one technical simplification (a better treatment of finite group actions on Fukaya categories: Section \ref{subsec:coverings}) and removed one mistake (the original proof that the first order deformation of the affine Fukaya category is nontrivial was incorrect, and has been replaced by the new Sections \ref{subsec:skip-divisor} and \ref{subsec:nontrivial-geometric}). Elsewhere, the exposition has been slightly smoothed. Finally, parts of the computer code have been updated (and made accessible on the author's homepage). 

\section{$A_\infty$-categories\label{sec:generators}}

This section is designed to serve as quick reference for the basic
algebraic terminology, for the benefit of readers with a geometry
background. The main construction is that of the bounded derived
category of an $A_\infty$-category, following Kontsevich
\cite{kontsevich94}. We then consider the idempotent (Karoubi, or
split-closed) completion, and the corresponding notion of
split-generators. A brief look at twist functors and exceptional
collections completes the tour. Some references for more in-depth
treatments of various matters touched on here are \cite{kontsevich98,
keller99, lefevre, rudakov90, fooo, fukaya00, seidel04}. The exposition
here is closest to \cite{seidel04}, even though that is obviously not
the earliest source.

\subsection{}
A non-unital $A_\infty$-category $\Cat$ consists of a set (all our
categories are small) of objects $Ob\,\Cat$, together with a graded
$\C$-vector space $\mathit{hom}_\Cat(X_0,X_1)$ for any two objects, and
multilinear maps
\begin{equation} \label{eq:mu}
 \mu^d_\Cat: \mathit{hom}_\Cat(X_{d-1},X_d) \otimes \dots \otimes
 \mathit{hom}_\Cat(X_0,X_1) \longrightarrow \mathit{hom}_\Cat(X_0,X_d)[2-d]
\end{equation}
for all $d \geq 1$ and all $(d+1)$-tuples of objects
$(X_0,\dots,X_d)$, satisfying the $A_\infty$-associativity equations
\begin{equation} \label{eq:ainfty}
 \sum_{e,i} (-1)^\ast \mu_\Cat^{d-e+1}(a_d,\dots,
 a_{i+e+1},\mu_\Cat^e(a_{i+e},\dots,a_{i+1}),a_i,\dots,a_1) = 0,
\end{equation}
where $\ast = |a_1|+\dots+|a_i|-i$ (the decreasing numbering of the
$a_k$ is the result of trying to bring together various standard
conventions). The underlying (non-unital) cohomological category
$H(\Cat)$ has the same objects as $\Cat$, with
\begin{equation}
\mathit{Hom}^*_{H(\Cat)}(X_0,X_1) = H^*(\mathit{hom}_{\Cat}(X_0,X_1),\mu^1_\Cat)
\end{equation}
and composition induced by the cochain level map $a_2a_1 =
(-1)^{|a_1|}\mu^2_\Cat(a_2,a_1)$. There is also the subcategory
$H^0(\Cat) \subset H(\Cat)$ which has the same objects, but retains only the morphisms of degree zero. One says
that $\Cat$ is {\em c-unital} (cohomologically unital) if $H(\Cat)$
has identity morphisms. It is true, but not entirely obvious, that
this is a good class of $A_\infty$-categories to work with: the
proofs of many properties rely on the fact that c-unitality turns out
to be equivalent to the less flexible but more visibly well-behaved
notion of strict unitality \cite{lefevre,seidel04}. {\em From now on, all
$A_\infty$-categories are assumed to be c-unital}. Similarly, when
considering $A_\infty$-functors $\F: \Cat \rightarrow {\mathcal D}$,
we always require that $H(\F): H(\Cat) \rightarrow H({\mathcal D})$
be unital (carries identity morphisms to identity morphisms). $\F$ is
called a quasi-equivalence if $H(\F)$ is an equivalence. One can
prove that for any quasi-equivalence, there is a $\G: {\mathcal D} \rightarrow \Cat$ such
that $H(\G)$ is an inverse equivalence to $H(\F)$.

\begin{remark}
A linear category with one object is the same as an associative
algebra, and similarly, an $A_\infty$-category with one object is an
$A_\infty$-algebra. Slightly more generally, given a category $C$
with only finitely many objects $X_1,\dots,X_m$, one can form the
{\em total morphism algebra}
\begin{equation}
 A = \bigoplus_{j,k} \mathit{Hom}_C(X_j,X_k)
\end{equation}
which is linear over the semisimple ring $R_m = \C \oplus \dots
\oplus \C$ ($m$ summands). Again, the same holds for
$A_\infty$-categories versus $R_m$-linear $A_\infty$-algebras, and we
will use that frequently. Our notation will be that $e_1,\dots,e_m$
are the orthogonal idempotents in $R_m$, so $e_k A e_j =
\mathit{Hom}_C(X_j,X_k)$.
\end{remark}

\begin{remark} \label{th:top-free}
For concreteness, we take the ground field to be $\C$ for our exposition of the basic algebraic theory, but any other field would do just as well (later on though, we will make use of deformation theoretic ideas which work best in characteristic $0$). Defining $A_\infty$-categories over a ring is a more difficult issue, and the only elementary solution is apparently to assume that all the $hom$ spaces are projective modules. Besides $R_m$ (where this problem does not really arise because of semi-simplicity), the other ring of interest for us is $\Lambda_\N$, which requires special treatment because of its
topological nature. Recall that for a $\Lambda_\N$-module $V_q$ which is complete in the $q$-adic topology, the following conditions are equivalent: (1) $V_q$ is torsion-free, which means that multiplication by $q$ is an injective endomorphism; (2) $V_q$ is topologically free, which means that projection $V_q \rightarrow V$ to the vector space $V = V_q \otimes_{\Lambda_q} \C = V_q/qV_q$ can be lifted to an isomorphism $V_q \rightarrow V \hat\otimes \Lambda_\N$. Of course, for finitely generated modules completeness is automatic, and topological freeness is the same as ordinary freeness. A {\em topological $A_\infty$-category} $\Cat_q$ over $\Lambda_\N$ consists of a set of objects, together with graded $\Lambda_\N$-modules $\mathit{hom}_{\Cat_q}(X_0,X_1)$ whose homogeneous pieces are topologically free (the compositions are automatically continuous).
%
\end{remark}

\subsection{}
The twisted complex construction embeds $\Cat$ into a larger
$A_\infty$-category which is closed under (appropriately understood)
mapping cones. This can be broken into two steps:
\begin{itemize}
\item
The additive enlargement $\Sigma\Cat$ has objects which are formal
sums
\begin{equation}
X = \bigoplus_{f \in F} X_f[\sigma_f],
\end{equation}
where $F$ is some finite set, $X_f \in Ob\,\Cat$, and $\sigma_f \in
\Z$. The morphisms between any two objects are
\begin{equation}
 \mathit{hom}_{\Sigma\Cat}\Big(\bigoplus_{f \in F} X_f[\sigma_f],
 \bigoplus_{g \in G} Y_g[\tau_g]\Big) =
 \bigoplus_{f,g} \mathit{hom}_{\Cat}(X_f,Y_g)[\tau_g-\sigma_f],
\end{equation}
and the compositions are defined using those of $\Cat$ and the
obvious ``matrix multiplication'' rule, with additional signs put in
as follows: for $a_1 \in \mathit{hom}_{\Sigma\Cat}(X_0[\sigma_0],X_1[\sigma_1]), \dots, a_d
\in \mathit{hom}_{\Sigma\Cat}(X_{d-1}[\sigma_{d-1}],X_d[\sigma_d])$,
\begin{equation}
 \mu^d_{\Sigma\Cat}(a_d,\dots,a_1) = (-1)^{\sigma_0}
 \mu^d_\Cat(a_d,\dots,a_1).
\end{equation}

\item
A twisted complex is a pair $(C,\delta)$ consisting of $C \in
Ob\,\Sigma\Cat$ and $\delta \in \mathit{hom}^1_{\Sigma\Cat}(C,C)$, with the ``strict upper
triangularity'' property that the indexing set $F$ for $C =
\bigoplus_{f \in F} C_f$ can be ordered in such a way that all
components $\delta_{fg}$ with $f \geq g$ are zero, and subject to the
generalized Maurer-Cartan equation
\begin{equation} \label{eq:generalized-maurer-cartan}
 \mu^1_{\Sigma\Cat}(\delta) + \mu^2_{\Sigma\Cat}(\delta,\delta) + \dots = 0.
\end{equation}
Upper triangularity says that this is a finite sum, hence the
condition makes sense. We define an $A_\infty$-category $\mathit{Tw}\Cat$ of
which the $(C,\delta)$ are the objects. The spaces of morphisms are
the same as for $\Sigma\Cat$, but composition involves the $\delta$s:
for $a_k \in \mathit{hom}_{\mathit{Tw}\Cat}((C_{k-1},\delta_{k-1}),(C_k,\delta_k))$,
\begin{equation} \label{eq:insert-delta-everywhere}
 \mu^d_{\mathit{Tw}\Cat}(a_d,\dots,a_1) = \sum_{j_0,\dots,j_d \geq 0}
 \mu^d_{\Sigma\Cat}(\overbrace{\delta_d,\dots,\delta_d}^{j_d},a_d,
 \overbrace{\delta_{d-1},\dots,\delta_{d-1}}^{j_{d-1}},a_{d-1},\dots).
\end{equation}
\end{itemize}
The c-unitality of $\Cat$ implies that of $\mathit{Tw}\Cat$. The underlying
cohomological category in degree zero is (somewhat improperly, since there is no localization involved) called the
{\em bounded derived category} $D^b(\Cat) = H^0(\mathit{Tw}\Cat)$. Given two
twisted complexes $(C_k,\delta_k)$, $k = 0,1$, and a morphism $a \in
\mathit{hom}^0_{\mathit{Tw}\Cat}(C_0,C_1)$ such that $\mu^1_{Tw\Cat}(a) = 0$, one can define the
mapping cone as usual:
\begin{equation} \label{eq:mapping-cone-tw}
 \mathit{Cone}(a) = \Big\{C_0 \stackrel{a}{\rightarrow} C_1\Big\} \stackrel{\text{def}}{=}
 \Big(C_0[1] \oplus C_1, \begin{pmatrix} \delta_0 & 0 \\ a &
 \delta_1 \end{pmatrix}\Big)
 \in Ob\,\mathit{Tw}\Cat.
\end{equation}
The isomorphism class of this in the derived category depends only on
the cohomology class $\alpha = [a] \in \mathit{Hom}_{D^b(\Cat)}(C_0,C_1)$,
hence we will sometimes write sloppily $\mathit{Cone}(\alpha)$. Declare {\em
exact triangles} in the derived category to be those which are
isomorphic to the standard mapping cone triangles $C_0 \rightarrow
C_1 \rightarrow \mathit{Cone}(a) \rightarrow C_0[1]$. This makes $D^b(\Cat)$
into a triangulated category. It contains $H^0(\Cat)$ as a full
subcategory, and the objects of that subcategory {\em generate}
$D^b(\Cat)$ in the usual triangulated sense, that is, any object of
$D^b(\Cat)$ can be obtained from objects of $H^0(\Cat)$ by successive
shifts, mapping cones, and isomorphism. This does not mean that
on the cohomology level, $H^*(\Cat)$ determines $D^b(\Cat)$.

\begin{example}
Take an $A_\infty$-category $\Cat$ with three objects $X_k$ ($k = 1,2,3$), and
where the only nonzero morphism spaces are
\begin{equation}
\begin{aligned}
\mathit{hom}^0_\Cat(X_k,X_k) \iso \C, \\
\mathit{hom}_\Cat^0(X_1,X_2) \iso \C, \\
\mathit{hom}_\Cat^0(X_2,X_3) \iso \C, \\
\mathit{hom}_\Cat^1(X_3,X_1) \iso \C.
\end{aligned}
\end{equation}
The differential $\mu^1_\Cat$ then necessarily vanishes, and the only nontrivial $\mu^2_\Cat$
multiplications are those involving identity endomorphisms. Hence, $H^*(\Cat)$ is unique up to
isomorphism. However, there are different choices for the higher order operations, and those affect 
$D^b(\Cat)$. Specifically, there is an $A_\infty$-structure with a nontrivial $\mu^3_\Cat$ for which $X_3$ becomes isomorphic to the
cone $\{X_1 \rightarrow X_2\}$
\cite[Section 3g]{seidel04}; in that case, $D^b(\Cat)$ is the well-known derived category of
representations of the $(A_2)$ quiver, and in particular the $X_k$ are the only indecomposable
objects up to a shift. On the other hand, one can consider the formal $A_\infty$-structure, with vanishing
higher order products. For that $A_\infty$-structure, the twisted complex schematically
denoted by
\begin{equation}
C = \{X_1 \rightarrow X_2 \rightarrow X_3\}
\end{equation}
is an indecomposable object of $D^b(\Cat)$, not isomorphic to the shift of any $X_k$.
\end{example}

The idempotent completion $C^\pi$ of a category $C$ has objects
which are pairs $(X,p)$, where $X$ is an object of the original
category and $p \in \mathit{Hom}_C(X,X)$ an idempotent endomorphism. One thinks
of this as the ``direct summand $\mathit{im}(p)$'', and defines
$\mathit{Hom}_{C^\pi}((X_0,p_0),(X_1,p_1)) = p_1 \mathit{Hom}_C(X_0,X_1) p_0$. It is proved in
\cite{balmer-schlichting01} that the idempotent completion of a
triangulated category is again triangulated, in a unique way which
makes the embedding $C \rightarrow C^\pi$ into an exact functor. The
result of applying this construction to the bounded derived category
of an $A_\infty$-category is called the {\em split-closed derived
category} $D^\pi(\Cat)$. Suppose that we have a set of objects in an
idempotent complete triangulated category; these are called {\em
split-generators} if one can obtain any object from them by
successive shifts, mapping cones, splitting off direct summands, and
isomorphism.

Here are some easy but useful facts about derived categories of $A_\infty$-categories:

\begin{lemma}
If two $A_\infty$-categories are quasi-equivalent, their derived
categories $D^b$ and $D^\pi$ are equivalent triangulated categories.
\end{lemma}

\proof Any $A_\infty$-functor $\F: \Cat \rightarrow {\mathcal D}$
admits a canonical extension $\mathit{Tw}\F: \mathit{Tw}\Cat \rightarrow \mathit{Tw}{\mathcal
D}$, and $D^b(\F) = H^0(\mathit{Tw}\F): D^b(\Cat) \rightarrow D^b({\mathcal D})$
is exact (roughly speaking, since it takes cones to cones). Suppose that $\F$ is a
quasi-equivalence, and consider the full subcategory of $D^b(\Cat)$
consisting of those objects $X$ with the following property:
\begin{quote}
For all $Y \in \mathit{Ob}\, \Cat$, the map
\begin{equation}
D^b(\F): \mathit{Hom}^*_{D^b(\Cat)}(X,Y) \longrightarrow
\mathit{Hom}^*_{D^b({\mathcal D})}(D^b\F(X),\F(Y))
\end{equation}
is an isomorphism.
\end{quote}
By assumption this subcategory contains $H^0(\Cat)$,
and it is a standard exercise to show that it is triangulated, so it
is the whole of $D^b(\Cat)$. Another step of the same kind proves that
$D^b(\F)$ is full and faithful. Its image is a triangulated subcategory
of $D^b({\mathcal D})$ which, up to isomorphisms, contains all
objects of $H^0({\mathcal D})$. Since these are generators for
$D^b({\mathcal D})$, we conclude that $D^b(\F)$ is an equivalence.
The corresponding statement for $D^\pi$ follows immediately. \qed

\begin{lemma} \label{th:generating-subcategory}
Let $\Cat' \subset \Cat$ be a full $A_\infty$-subcategory such that the
set of all objects of $\Cat'$ split-generates $D^\pi(\Cat)$. Then
$D^\pi(\Cat') \iso D^\pi(\Cat)$ as triangulated categories.
\end{lemma}

\proof Inclusion induces a full and faithful functor $D^\pi(\Cat')
\rightarrow D^\pi(\Cat)$, the image of which contains the objects of
$H^0(\Cat')$ and is closed under shifts, cones, and idempotent splittings.
\qed

Combining the two previous Lemmas, we have:

\begin{lemma} \label{th:derived-equivalence}
Let $\Cat$, ${\mathcal D}$ be two $A_\infty$-categories. Suppose that
they contain full $A_\infty$-subcategories $\Cat' \subset \Cat,
{\mathcal D}' \subset {\mathcal D}$ which are quasi-equivalent to
each other. If the objects of $\Cat',{\mathcal D}'$ split-generate
$D^\pi(\Cat)$ respectively $D^\pi({\mathcal D})$, then these two
derived categories are equivalent. \qed
\end{lemma}

\subsection{}
Let $C$ be a $\C$-linear triangulated category, and $X,Y$ two objects
such that the graded vector space $\mathit{Hom}_C^*(X,Y) = \bigoplus_{k \in \Z} \mathit{Hom}_C(X,Y[k])$ is
finite-dimensional. Then one can form two new objects $T_X(Y)$,
$T_Y^\vee(X)$ which fit into exact triangles
\begin{equation}
\begin{aligned}
 & T_X(Y)[-1] \rightarrow \mathit{Hom}_C^*(X,Y) \otimes X \xrightarrow{ev} Y \rightarrow T_X(Y), \\
 & T_Y^\vee(X) \rightarrow X \xrightarrow{ev^\vee} \mathit{Hom}^*_C(X,Y)^\vee \otimes Y \rightarrow
 T_Y^\vee(X)[1].
\end{aligned}
\end{equation}
The notation means that we choose a homogeneous basis $\{f_k\}$ of
$\mathit{Hom}^*_C(X,Y)$, set $\mathit{Hom}^*_C(X,Y) \otimes X = \bigoplus_k X[-|f_k|]$,
and take the evaluation map $ev = \bigoplus_k f_k[-|f_k|]$. The dual map
$ev^\vee$ is defined in the same way, using the dual basis $\{f_k^\vee\}$. It is easy to see that the choice of basis is essentially irrelevant. As a consequence, these algebraic {\em twist} and {\em untwist} operations are well-defined up to non-unique isomorphism (and the isomorphism type of the new objects depends only on those of $X$ and $Y$).

The construction can be formulated in a particularly simple and natural way in the case where
$C = D^b(\Cat)$ is the bounded derived category of an $A_\infty$-category $\Cat$ such that: (1) the cochain level morphism space $\mathit{hom}_{\mathit{Tw}\Cat}(X,Y)$ is itself finite-dimensional; and (2) $X$ has a strict cochain level identity morphism $e_X$. One can then form
the tensor product $\mathit{hom}_{\mathit{Tw}\Cat}(X,Y) \otimes X \in Ob\,\mathit{Tw}\Cat$, which is a
finite direct sum of shifted copies of $X$ with a differential
$\mu^1_{\mathit{Tw}\Cat} \otimes e_X$ (added to that of $X$ itself, if $X$ is a twisted complex). The evaluation map has a canonical cochain representative, and $T_X(Y)$ can be defined as its mapping cone \eqref{eq:mapping-cone-tw}:
\begin{equation} \label{eq:explicit-twist}
 T_X(Y) = \{\mathit{hom}^*_{\mathit{Tw}\Cat}(X,Y) \otimes X \longrightarrow Y\}.
\end{equation}
In fact, if the finite-dimensionality condition holds for all $Y$,
one can realize $T_X$ by an $A_\infty$-endofunctor of $\mathit{Tw}\Cat$,
inducing an exact functor on $D^b(\Cat)$. Similarly, if (1) holds and
(2') $Y$ has a strict cochain level identity morphism $e_Y$, a definition of the untwist as
$T_Y^\vee(X) = \{X \rightarrow \mathit{hom}^*_{\mathit{Tw}\Cat}(X,Y)^\vee \otimes Y\}[-1]$ works.
Fukaya categories are not in general strictly unital, so the previous
observation does not apply directly to them. This can be remedied \cite[Section 5]{seidel04},
but instead of entering into a full discussion, we just record the
technical statement that will be used later on:

\begin{lemma} \label{th:recognize-twist}
Let $\Cat$ be a c-unital $A_\infty$-category whose $\mathit{hom}_{\Cat}$ spaces are
finite-dimen\-si\-on\-al, and $X,Y,Z$ objects of it. Any pair
\begin{equation} \label{eq:hk}
 k \in \mathit{hom}^0_\Cat(Y,Z), \quad
 h : \mathit{hom}_\Cat(X,Y) \longrightarrow \mathit{hom}_\Cat(X,Z)[-1]
\end{equation}
satisfying
\begin{equation} \label{eq:coupled-morphism}
 \mu^1_\Cat(k) = 0, \quad
 \mu^1_\Cat h(a) - h(\mu^1_\Cat a) + \mu^2_\Cat(k,a) = 0 \text{ for $a \in \mathit{hom}_\Cat(X,Y)$},
\end{equation}
determines a morphism $T_X(Y) \rightarrow Z$ in $D^b(\Cat)$. This is an isomorphism if for all
objects $W$ of $\Cat$, the following complex of vector spaces is
acyclic:
\begin{equation} \label{eq:total-cone}
\begin{split}
 & (\mathit{hom}_\Cat(X,Y) \otimes \mathit{hom}_\Cat(W,X))[2] \oplus \mathit{hom}_\Cat(W,Y)[1] \oplus \mathit{hom}_\Cat(W,Z), \\
 & \partial(a_2 \otimes a_1,b,c) =
 \big(\mu^1_\Cat(a_2) \otimes a_1 +
 (-1)^{|a_2|-1} a_2 \otimes \mu^1_\Cat(a_1),
 \mu^1_\Cat(b) \\ & \qquad + \mu^2_\Cat(a_2,a_1),
 -\mu^1_\Cat(c) + \mu^2_\Cat(k,b) + \mu^3_\Cat(k,a_2,a_1) +
 \mu^2_\Cat(h(a_2),a_1)\big).
\end{split}
\end{equation}
\end{lemma}

\proof Suppose first that $\Cat$ is strictly unital. Then $(k,h)$ do
indeed define a ($\mu^1_\Cat$-closed) cochain level morphism from the object
$T_X(Y)$ as defined in \eqref{eq:explicit-twist} to $Z$. Composition
with this morphism is a map of complexes of vector spaces
\begin{equation}
 \mathit{hom}_{Tw\Cat}(W,T_X(Y)) \longrightarrow \mathit{hom}_{\Cat}(W,Z)
\end{equation}
and \eqref{eq:total-cone} is the mapping cone of it, so the
assumption that this is acyclic means that our morphism induces
isomorphisms
\begin{equation}
\mathit{Hom}^*_{D^b(\Cat)}(W,T_X(Y)) \iso \mathit{Hom}^*_{D^b(\Cat)}(W,Z)
\end{equation}
for all $W \in Ob\,\Cat$. A straightforward long exact sequence
argument extends this to objects $W$ in the derived category, and
then the Yoneda Lemma implies that $T_X(Y) \iso Z$.

There are two possibilities for extending this argument to the c-unital case. One can use $A_\infty$-modules instead of twisted complexes, as in \cite{seidel04}. Alternatively, one could pass from the given category to a quasi-isomorphic strictly unital one (in which case one has to show that the morphisms \eqref{eq:hk} carry over). We omit the details. \qed

\begin{lemma} \label{th:twist-generators}
Let $C$ be a $\C$-linear idempotent closed triangulated category,
such that $\mathit{Hom}_C^*(X,Y)$ is finite-dimensional for all $X,Y$. Take a
collection of objects $X_1,\dots,X_m$ and another object $Y$, and
consider $Y_d = (T_{X_1} \dots
T_{X_m})^d(Y)$ for $d \geq 0$. Suppose there is a $d$ such that $\mathit{Hom}_C(Y,Y_d) = 0$
(that is, there are no nontrivial morphisms {\em of degree zero}). Then $Y$
lies in the idempotent closed triangulated subcategory of $C$
split-generated by $X_1,\dots,X_m$.
\end{lemma}

\proof By definition, we have exact triangles
\begin{equation} \label{eq:taut-triangle}
\xymatrix{
 {T_{X_{i+1}} \dots T_{X_m}(Y_j)} \ar[rr] &&
 {T_{X_i}\dots T_{X_m}(Y_j)} \ar^-{[1]}[dl] \\
 & {\hspace{-5em}
 \mathit{Hom}^*_C(X_i,T_{X_{i+1}}\dots T_{X_m}(Y_j)) \otimes X_i
 \hspace{-5em} } \ar^-{ev}[ul] &
}
\end{equation}
The octahedral axiom allows us to conflate these, for all $i$ and all
$j < d$, into a single exact triangle
\begin{equation} \label{eq:new-triangle}
\xymatrix{
 Y \ar[rr] && Y_d \ar^-{[1]}[dl] \\
 & R_d \ar[ul] &}
\end{equation}
where $R_d$ is an object built by taking repeated mapping cones from
the objects in the same positions in \eqref{eq:taut-triangle}, hence
which lies in the triangulated subcategory generated by
$X_1,\dots,X_m$. The assumption says that the $\rightarrow$ in
\eqref{eq:new-triangle} must be zero, hence $R_d \iso Y \oplus
Y_d[-1]$. \qed

Our applications will be along the following lines: if
$C = D^\pi(\Cat)$ and there are fixed $X_1,\dots,X_m \in
\mathit{Ob}\, \Cat$ such that any $Y \in \mathit{Ob}\, \Cat$ has the property
required above, then the $X_k$ are split-generators for
$D^\pi(\Cat)$.

\subsection{}
Let $C^\rightarrow$ be a $\C$-linear graded category with finitely
many objects, ordered in a fixed way as $\{X_1,\dots,X_m\}$. We say that
$C^\rightarrow$ is {\em directed} if
\begin{equation} \label{eq:directedness}
 \mathit{Hom}^*_{C^\rightarrow}(X_i,X_j) = \begin{cases}
 0 & i>j, \\ \C \cdot id_{X_i} & i = j, \\
 \text{finite-dimensional} & i < j.
 \end{cases}
\end{equation}
If $C$ is a category with finite-dimensional
$\mathit{Hom}^*$'s, and $\{X_1,\dots,X_m\}$ an ordered finite collection of
nonzero objects in it, we can define the associated directed
subcategory to be the subcategory $C^\rightarrow \subset C$ with objects $\{X_1,\dots,X_m\}$, which is directed and satisfies
$\mathit{Hom}^*_{C^\rightarrow}(X_i,X_j) = \mathit{Hom}^*_C(X_i,X_j)$ for $i<j$. This
is a full subcategory iff the $X_i$ form an {\em exceptional
collection}, which means that their morphisms in $C$ already satisfy
the conditions of \eqref{eq:directedness}.

Similarly, let $\Cat^\rightarrow$ be an $A_\infty$-category with
objects $\{X_1,\dots,X_m\}$. We call it directed if the chain level $hom$-spaces
satisfy the analogous conditions as in \eqref{eq:directedness}, where the
generators $e_i$ of $\mathit{hom}_{\Cat^\rightarrow}(X_i,X_i)$ are assumed to be strict units. Given a finite
collection of objects $\{X_1,\dots,X_m\}$ in an arbitrary
$A_\infty$-category $\Cat$ with finite-dimensional $hom$'s, one can
define an associated directed $A_\infty$-category $\Cat^\rightarrow$
by requiring that $\mathit{hom}_{\Cat^\rightarrow}(X_i,X_j) =
\mathit{hom}_\Cat(X_i,X_j)$ for $i<j$, and
$\mu^d_{\Cat^\rightarrow}(a_d,\dots,a_1) =
\mu^d_{\Cat}(a_d,\dots,a_1)$ for all $a_k \in
\mathit{hom}_\Cat(X_{i_{k-1}},X_{i_k})$ with $i_0 < \dots < i_d$. This is not
really an $A_\infty$-subcategory of $\Cat$ (unless that happens to
have strict units), but after some more work one can still construct
an $A_\infty$-functor $\Cat^\rightarrow \rightarrow \Cat$ which is the identity on the spaces
$\mathit{hom}_{\Cat^\rightarrow}(X_i,X_j)$, $i<j$. We will tolerate a
slight abuse of terminology, and call $\Cat^\rightarrow$ the {\em
directed $A_\infty$-subcategory} associated to $\{X_1,\dots,X_m\}$.

The main interest in exceptional collections and directed categories
comes from the notion of mutations. Occasional glimpses of this
theory will be visible later in the paper, in particular in Section
\ref{sec:induction}, but it is not strictly necessary for our
purpose, so we will just refer the interested reader to
\cite{rudakov90,kontsevich98,seidel00}.

\section{Deformation theory\label{sec:deformations}}

We need to put together some material (well-known in the right
circles) from the classification theory of $A_\infty$-structures.
Using the Homological Perturbation Lemma
\cite{gugenheim-lambe-stasheff90}, the ``moduli space'' of
quasi-isomorphism types of $A_\infty$-algebras with fixed cohomology
algebra can be written as the quotient of an infinite-dimensional
space of solutions to the $A_\infty$-equations, by the action of an
infinite-dimensional ``gauge'' group. There is a canonical base
point, represented by the trivial $A_\infty$-structure (the one where
the higher order composition maps are all zero); $A_\infty$-algebras quasi-isomorphic to the trivial one are called formal. The ``Zariski
tangent space'' at the base point is a product of suitable Hochschild
cohomology groups of the cohomology algebra. Moreover, the ``moduli
space'' carries a $\C$-action which contracts it to the base point.
As noted by Kadeishvili \cite{kadeishvili88}, this implies that the
``tangent space'' already carries a lot of information about the
``moduli space''. A more straightforward argument along the same
lines applies to the classification of infinitesimal deformations of
a fixed possibly nontrivial $A_\infty$-structure, using a suitable
generalization of Hochschild cohomology \cite{getzler-jones90}. We
will recall these and some similar ideas, aiming in each case for
uniqueness results with easily verifiable criteria, tailored to the
specific computations later on. The proofs are exercises in abstract
deformation theory \`a la Deligne \cite{goldman-millson88}, and we
postpone them to the end.

\subsection{\label{subsec:classification}}
Let $A$ be a graded algebra. Consider the set $\AA(A)$ of those
$A_\infty$-algebras $\A$ whose underlying graded vector space is $A$,
and where the first two composition maps are
\begin{equation} \label{eq:mu1zero}
 \mu^1_\A = 0, \quad
 \mu^2_\A(a_2,a_1) = (-1)^{|a_1|} a_2a_1.
\end{equation}
This means that the cohomology algebra is $H(\A) = A$ itself. $\A,\A'
\in \AA(A)$ are considered equivalent if there is an
$A_\infty$-homomorphism $\G: \A \rightarrow \A'$ whose underlying
linear map is $\G^1 = id_A$. After writing out the relevant
equations, of which the first one is
\begin{equation} \label{eq:pullback}
\begin{split}
 \mu^3_{\A'}(a_3,a_2,a_1)
 & = \mu^3_\A(a_3,a_2,a_1) \\ &
 + \G^2(a_3,\mu^2_\A(a_2,a_1)) + (-1)^{|a_1|-1} \G^2(\mu^2_\A(a_3,a_2),a_1) \\
 & - \mu^2_\A(a_3,\G^2(a_2,a_1))
 - \mu^2_\A(\G^2(a_3,a_2),a_1),
\end{split}
\end{equation}
one sees that $\A$ and $\G$ determine $\A'$ completely, and that
there are no constraints on $\G$. In other words, one can define a
group $\GG(A)$ of {\em gauge transformations} which are arbitrary
sequences of multilinear maps $\{\G^1 = id_A, \G^2: A \otimes A
\rightarrow A[-1], \dots\}$, and then \eqref{eq:pullback} and its
successors define an action of $\GG(A)$ on $\AA(A)$, whose quotient
is the set of equivalence classes (this becomes entirely transparent
in a geometric formulation, where it is the action of formal
noncommutative diffeomorphisms on integrable odd vector fields).
Additionally, there is a canonical action of the multiplicative
semigroup $\C$ on $\AA(A)$, by multiplying $\mu^d$ with
$\epsilon^{d-2}$, which we denote by $\A \mapsto \epsilon^*\A$. For
$\epsilon \neq 0$, it is still true that $\A$ is isomorphic to
$\epsilon^*\A$, even though not by a gauge transformation: the
isomorphism is multiplication by $\epsilon^k$ on the degree $k$
component. The $\C$-action obviously descends to $\AA(A)/\GG(A)$.

The main justification for imposing conditions \eqref{eq:mu1zero} is
the Homological Perturbation Lemma, which says the following. Let
$\B$ be an $A_\infty$-algebra with an isomorphism (of graded algebras) $F: A \rightarrow
H(\B)$. Then one can find an $A_\infty$-structure $\A$ on $A$
satisfying \eqref{eq:mu1zero}, and an $A_\infty$-homomorphism $\F: \A
\rightarrow \B$ such that $H(\F) = F$. This establishes a bijection
\begin{equation}
 \frac{\AA(A)}{\GG(A)} \iso \frac{
 \{\text{pairs $(\B,F: A \rightarrow H(\B))$}\}
 }
 {
 \{\text{$A_\infty$-quasi-isomorphisms compatible with $F$}\}
 }.
\end{equation}

\begin{remark} \label{th:explicit-hpt}
The construction of $\A$ is explicit: take maps
\begin{equation} \label{eq:hpt-data}
 \pi: \B \rightarrow A, \quad
 \lambda: A \rightarrow \B, \quad
 h: \B \rightarrow \B[-1]
\end{equation}
such that $\pi \mu^1_\B = 0$, $\mu^1_\B \lambda = 0$ (with the map
induced by $\lambda$ on cohomology being $F$), $\pi \lambda = id_A$,
$\lambda \pi = id_\B + \mu^1_\B h + h \mu^1_\B$. Then $\mu_\A^d$ is a
sum over planar trees with $d$ leaves and one root. The contribution
from each tree involves only the $\mu^k_\B$ for $2 \leq k \leq d$, and
the auxiliary data \eqref{eq:hpt-data} (this particular formulation
is taken from \cite{kontsevich-soibelman00}, see also \cite{markl04}). For instance, writing
$b_k = \lambda(a_k)$ one has
\begin{align} \label{eq:mu4-hpt}
 & \mu^4_\A(a_4,a_3,a_2,a_1) =
 \pi \mu^4_\B(b_4,b_3,b_2,b_1) + {} \\ \notag
 & \qquad \pi \mu^2_\B(h\mu^3_\B(b_4,b_3,b_2),b_1)
 + \pi \mu^2_\B(b_4,h\mu^3_\B(b_3,b_2,b_1)) + {} \\ \notag
 & \qquad \pi \mu^3_\B(h\mu^2_\B(b_4,b_3),b_2,b_1)
 + \pi \mu^3_\B(b_4,h\mu^2_\B(b_3,b_2),b_1) + {} \displaybreak[0] \\ \notag
 & \qquad \pi \mu^3_\B(b_4,b_3,h\mu^2_\B(b_2,b_1))
 + \pi \mu^2_\B(h\mu^2_\B(b_4,b_3),h\mu^2_\B(b_2,b_1)) + {} \\ \notag
 & \qquad \pi \mu^2_\B(h\mu^2_\B(h\mu^2_\B(b_4,b_3),b_2),b_1)
 + \pi \mu^2_\B(h\mu^2_\B(b_4,h\mu^2_\B(b_3,b_2)),b_1) + {} \\ \notag
 & \qquad \pi \mu^2_\B(b_4,h\mu^2_\B(h\mu^2_\B(b_3,b_2),b_1))
 + \pi \mu^2_\B(b_4,h\mu^2_\B(b_3,h\mu^2_\B(b_2,b_1))).
\end{align}
There is nothing particularly memorable about this formula, but $\mu^4_\A$ happens
to be the higher order product that is crucial in our specific computation 
later on.
\end{remark}

The Hochschild cohomology of the graded algebra $A$ is a bigraded vector space. In our
slightly unusual notation, elements of $\mathit{HH}^{s+t}(A,A)^t$ are
represented by cocycles $\tau \in \mathit{CC}^{s+t}(A,A)^t$, which are
multilinear maps $\tau: A^{\otimes s} \rightarrow A$ of degree $t$.
Suppose that for $\A \in \AA(A)$, the compositions $\mu_\A^s$ vanish
for $2 < s < d$. Then $\mu^d_\A$ is a Hochschild cocycle; its class
\begin{equation} \label{eq:md-def-class}
o_\A^d = [\mu^d_\A] \in \mathit{HH}^2(A,A)^{2-d},
\end{equation}
the order $d$ {\em obstruction class} of $\A$, forms the obstruction
to making $\mu^d_\A$ trivial by a gauge transformation which is
itself equal to the identity to all orders $<d-1$. In particular, if
all the obstruction groups $\mathit{HH}^2(A,A)^{2-s}$, $s > 2$, are zero,
$\AA(A)/\GG(A)$ reduces to a point, which means that each
$A_\infty$-algebra structure on $A$ can be trivialized by a gauge
transformation (this is the {\em intrinsically formal} situation from
\cite{halperin-stasheff79, kadeishvili88}). A slightly more general
result allows the total obstruction group to be one-dimensional:

\begin{lemma} \label{th:versal-1}
Suppose that there is some $d>2$ such that
\begin{equation}
\begin{cases}
\mathit{HH}^2(A,A)^{2-d} \iso \C, \\ \mathit{HH}^2(A,A)^{2-s} = 0 & \text{for all $s > 2$,
$s \neq d$}.
\end{cases}
\end{equation}
Suppose also that there exists an $\A \in \AA(A)$ such that $\mu_\A^s
= 0$ for $2 < s < d$, and whose obstruction class $o_\A^d$ is
nonzero. Then any $\A' \in \AA(A)$ is equivalent to $\epsilon^*\A$
for some $\epsilon \in \C$, which means that the $\C$-orbit of $[\A]$ is
the whole of $\AA(A)/\GG(A)$.
\end{lemma}

We say that $\A$ is {\em versal}. We will now consider another deformation problem, which is different but related.
%

\subsection{}
Let $\A$ be an $A_\infty$-algebra. A {\em one-parameter deformation}
of $\A$ is a topological $A_\infty$-algebra $\A_q$ over $\Lambda_\N$,
whose underlying graded vector space is $\A \hat\otimes \Lambda_\N$,
and whose composition maps reduce to those of $\A$ if one sets $q =
0$. More concretely,
\begin{equation}
\mu_{\A_q}^d = \mu_{\A}^d + q \mu_{\A_q,1}^d + q^2 \mu_{\A_q,2}^d +
\dots =\mu_{\A}^d + O(q),
\end{equation}
where each coefficient $\mu_{\A_q,k}^d$ is a $\C$-linear map $\A^{\otimes d}
\rightarrow \A[2-d]$. Two one-parameter deformations are equivalent
if there is a $\Lambda_\N$-linear $A_\infty$-homomorphism $\G_q$ between them of the form
\begin{equation} \label{eq:zero-phi}
\G_q^d = \begin{cases} id + O(q) & d = 1, \\ O(q) & d>1. \end{cases}
\end{equation}
Let $\AA_q(\A)$ be the set of all one-parameter deformations, and
$\GG_q(\A)$ the group of sequences of multilinear maps $\{\G^1_q,
\G^2_q,\dots\}$ on $\A \hat{\otimes} \Lambda_\N$ which satisfy
\eqref{eq:zero-phi}. Then $\AA_q(\A)/\GG_q(\A)$ is the set of
equivalence classes. The semigroup $\mathit{End}(\Lambda_\N)$ acts on
$\AA_q(\A)$ by reparametrizations $q \mapsto \psi(q)$; we denote the
action by $\A_q \mapsto \psi^*\A_q$. In this context, versality
means that the $\mathit{End}(\Lambda_\N)$-orbit of some equivalence
class $[\A_q]$ is the whole of $\AA_q(\A)/\GG_q(\A)$.

\begin{remark}
We will also occasionally use terminology from finite order
deformation theory. An order $d$ infinitesimal deformation of $\A$
would be an $A_\infty$-algebra $\A_q$ as before, except that the
ground ring is now $\C[[q]]/q^{d+1}$. The notion of equivalence is
obvious. Of course, a one-parameter deformation can be truncated to
any desired finite order in $q$.
\end{remark}

Hochschild cohomology generalizes to $A_\infty$-algebras
\cite{getzler-jones90} as a graded (no longer bigraded) space
$\mathit{HH}^*(\A,\A)$. Consider 
\begin{equation} \label{eq:cc-prod}
\mathit{CC}^d(\A,\A) = \prod_{s+t = d}
\mathit{Hom}(\A^{\otimes s},\A[t])
\end{equation}
with the Gerstenhaber product
\begin{equation} \label{eq:gerstenhaber}
 (\sigma \circ \tau)^d(a_d,\dots,a_1) =
 \sum_{i,j} (-1)^\sharp \sigma^{d-j+1}(a_d,\dots,
 \tau^j(a_{i+j},\dots,a_{i+1}),\dots,a_1),
\end{equation}
where $\sharp = (|\tau|-1)(|a_1|+\dots+|a_i|-i)$, and the
corresponding Lie bracket $[\sigma,\tau] = \sigma \circ \tau -
(-1)^{(|\sigma|-1)(|\tau|-1)} \tau \circ \sigma$. As $\mu_\A^*$
itself lies in $\mathit{CC}^2(\A,\A)$ and satisfies $\mu_\A^* \circ \mu_\A^* =
0$ by definition, one can define a differential $\delta\tau =
[\mu_\A^*,\tau]$, and this is the $A_\infty$ version of the
Hochschild differential. The (complete decreasing) filtration of $\mathit{CC}^*(\A,\A)$ by length, which means by $s$ in \eqref{eq:cc-prod}, yields a spectral
sequence with
\begin{equation}
E_2^{s,t} = \mathit{HH}^{s+t}(A,A)^t
\end{equation}
for $A = H(\A)$, and this converges to $\mathit{HH}^*(\A,\A)$ under suitable
technical assumptions (for instance, if the nonzero terms in the
$E_2$ plane are concentrated in a region cut out by $t \leq c_1$, $t
\geq c_2 - c_3 s$ for $c_1,c_2$ arbitrary and $c_3 < 1$). If
$\mu_{\A}^d$, $d \geq 3$, is the first nonzero higher order
composition, then the first potentially nontrivial differential in
the spectral sequence is
\begin{equation} \label{eq:bracket-differential}
 \delta_{d-1} = [o_{\A}^d,\cdot]:
 E_2^{s,t} \longrightarrow E_2^{s+d-1,t+2-d}
\end{equation}
where $o_{\A}^d$ is \eqref{eq:md-def-class}. The
spectral sequence is a quasi-isomorphism invariant, so one can
study it under the assumption \eqref{eq:mu1zero}. Then, either $\A$
is equivalent to the trivial $A_\infty$-structure on $A$ and the
spectral sequence degenerates; or else, after using gauge
transformations to make as many higher order terms zero as possible,
one encounters a nontrivial obstruction class $o_{\A}^d$ for some
$d>2$. In that case, taking the Euler element $\tau \in \mathit{HH}^1(A,A)^0$
given by the derivation which is $k$ times the identity on $A^k$, one
finds that
\begin{equation} \label{eq:euler-element}
 [o_\A^d,\tau] = (d-2) o_\A^d.
\end{equation}
By \eqref{eq:bracket-differential}, $\delta_{d-1}$ is nonzero and
kills $o_\A^d \in E_{d-1}^{d,2-d} = \mathit{HH}^2(A,A)^{2-d}$.

Having come this far by straightforward adaptation of the previously
used notions, we encounter a slight snag, which is that the
deformation problem governed by $\mathit{HH}^*(\A,\A)$ reaches beyond the
class of $A_\infty$-algebras into that of ``curved'' or
``obstructed'' ones (with a $\mu^0_{\A_q}$ term of order $O(q)$), and the correspondingly ``extended'' notion of $A_\infty$-homomorphisms (with a $\G^0_q$ term of order $O(q)$). 

\begin{example}
Start with an $A_\infty$-algebra $\A$, and let $\delta \in \A^1[[q]]$ be an element of order $O(q)$ which solves the generalized Maurer-Cartan equation \eqref{eq:generalized-maurer-cartan} in $\A[[q]]$. Consider on one hand the trivial (constant) $A_\infty$-deformation $\A$, and on the other hand the $A_\infty$-deformation $\A_q$ obtained by inserting arbitrarily many copies of $\delta$ into the $A_\infty$-operations, just as in \eqref{eq:insert-delta-everywhere}.
In general, these two deformations are not equivalent in the ordinary sense (they may not even have isomorphic cohomology algebras). However, they are related by an extended equivalence, which is defined by setting $\G^0_q = \delta$, $\G^1_q = \mathit{id}$, and $\G^d_q = 0$ for $d \geq 2$.
\end{example}

Even though this extension of the $A_\infty$-deformation problem is natural, it is not always desirable for applications. One can get rid of the additional parts as follows. Let $\A$ be an $A_\infty$-algebra with $\mu^1_{\A} = 0$, and denote by $A$ the associated graded algebra. The {\em truncated Hochschild cohomology} is a graded vector space $\mathit{HH}^*(\A,\A)^{\scriptscriptstyle \leq 0}$, built from the subcomplex $\mathit{CC}^*(\A,\A)^{\scriptscriptstyle \leq 0}$ of
Hochschild cochains $\tau: \A^{\otimes s} \rightarrow \A[t]$ with $t
\leq 0$. There is a spectral sequence leading to it as before, where
now
\begin{equation} \label{eq:truncated-ss}
 E_2^{s,t} = \begin{cases}
 \mathit{HH}^{s+t}(A,A)^t & \text{$t \leq 0$}, \\
 0 & \text{otherwise.}
 \end{cases}
\end{equation}

All remarks made above carry over to this modification, including
\eqref{eq:euler-element} because the Euler element lies in
\eqref{eq:truncated-ss}. For the application to deformation theory,
consider the subset $\AA_q(\A)^{\scriptscriptstyle \leq 0} \subset
\AA_q(\A)$ of those one-parameter deformations $\A_q$ during which
the differential remains trivial, $\mu^1_{\A_q} = 0$. This still
carries actions of $\GG_q(\A)$ and of $\mathit{End}(\Lambda_\N)$. The first order
(in $q$) term of a deformation of this kind defines a deformation
class
\begin{equation}
m_{\A_q}^{\scriptscriptstyle \leq 0} = [\mu_{\A_q,1}^*] \in
\mathit{HH}^2(\A,\A)^{\scriptscriptstyle \leq 0},
\end{equation}
which is the obstruction to transforming $\A_q$ into a deformation with leading
order $q^2$ (or in other words, making the associated first order
deformation trivial). The desired versality result is:

\begin{lemma} \label{th:versal-3}
Suppose that $\mathit{HH}^2(\A,\A)^{\scriptscriptstyle \leq 0} \iso \C$. Let
$\A_q$ be a one-parameter deformation with $\mu^1_{\A_q} = 0$, whose
deformation class is nonzero. Then $\A_q$ is versal in
$\AA_q(\A)^{\scriptscriptstyle \leq 0}/\GG_q(\A)$, which means that any other one-parameter deformation of $\A$ is equivalent to $\psi^*\A_q$ for some $\psi \in \mathit{End}(\Lambda_\N)$.
\end{lemma}

\begin{corollary} \label{th:versal-3b}
Suppose that $\A$ satisfies the conditions from Lemma \ref{th:versal-3}, and let $\A_q$ be a versal deformation as described there. Let $\B$ be an $A_\infty$-algebra quasi-isomorphic to $\A$, and $\B_q$ a one-parameter deformation of it such that $H(\B_q)$ is torsion-free. Then $\B_q$ is quasi-isomorphic to some reparametrization $\psi^*\A_q$, where $\psi \in \mathit{End}(\Lambda_\N)$.
\end{corollary}

\proof
Split $\B$ into a copy of $H(\B)$ and an acyclic complex. The Perturbation Lemma then constructs a minimal $A_\infty$-structure on $H(\B)$ quasi-isomorphic to $\B$. The same machinery yields an $A_\infty$-deformation of that minimal structure quasi-isomorphic to $\B_q$. With this in mind, we can assume from now on that the original $\B$ was already minimal. 

In that case, if $\mu^1_{\B_q}$ was nonzero, then $H(\B_q)$ would have $q$-torsion; so, the torsion-freeness assumption implies that $\mu^1_{\B_q} = 0$. Now $\B$ is isomorphic to $\A$, and hence $\B_q$ will be isomorphic to a deformation $\A_q$ of the kind considered in Lemma \ref{th:versal-3}.
\qed

\begin{remark}
When $H(\A)$ is concentrated in even degrees, the torsion-freeness of $H(\B_q)$ is automatic. This follows directly from the argument just given: one passes to the minimal case, and then $\mu^1_{\B_q}$ must vanish for degree reasons.
\end{remark}

\subsection{}
The abstract framework of one-parameter deformation theory \cite{goldman-millson88} is as follows. Let $\g$ be a dg Lie algebra, with differential $\delta$. We consider the space $\AA_q(\g)$ of elements $\alpha_q \in \g^1 \hat{\otimes}\, q\Lambda_\N$ which satisfy the 
 equation
\begin{equation} \label{eq:classical-mc}
 \delta\alpha_q + \half [\alpha_q,\alpha_q] = 0.
\end{equation}
In particular, if we expand $\alpha_q = \alpha_{q,1}q + \alpha_{q,2} q^2 + \cdots$, the leading order term defines a deformation class $[\alpha_{q,1}] \in H^1(\g)$. An infinitesimal abstract gauge
transformation $\gamma_q \in \g^0 \hat{\otimes} q\Lambda_\N$ defines
a formal vector field $X_q$ on $\AA_q(\g)$, $X_q(\alpha_q) =
\delta\gamma_q + [\alpha_q,\gamma_q]$. These vector fields
exponentiate to an action of a group $\GG_q(\g)$ on $\AA_q(\g)$,
which we denote by $\exp(\gamma_q)(\alpha_q)$. The
deformation class is an invariant of this action. As usual, there is
also an action of $\mathit{End}(\Lambda_\N)$ on $\AA_q(\g)$ by
reparametrizations. The abstract formality and versality theorems are:

\begin{lemma} \label{th:formal-abstract}
If $H^1(\g) = 0$, $\AA_q(\g)/\GG_q(\g)$ is a point.
\end{lemma}

\begin{lemma} \label{th:versal-abstract}
Suppose that $H^1(\g) \iso \C$. Then, every $\alpha_q \in \AA_q(\g)$
whose deformation class is nonzero is versal. As before, this means 
that any other element of $\AA_q(\g)/\GG_q(\g)$ can be obtained from $\alpha_q$ by
applying a reparametrization in $End(\Lambda_\N)$.
\end{lemma}

These are both standard results, but for the convenience of the
reader we reproduce the conventional induction-by-order proof of the
second one. 

{\sc Proof of Lemma \ref{th:versal-abstract}.}
Take an arbitrary $\beta_q \in \AA_q(\g)$. Suppose that
for some $d \geq 1$, we have $\psi^{(d)} \in \mathit{End}(\Lambda_\N)$ and
$\beta^{(d)}_q \in \AA_q(\g)$ gauge equivalent to $\beta_q$, such that
\begin{equation} \label{eq:discrepancy}
 \beta^{(d)}_q - (\psi^{(d)})^*\alpha_q \stackrel{\text{def}}{=} \epsilon_q \in O(q^d).
\end{equation}
Because $\delta\epsilon_q = \half [(\psi^{(d)})^*\alpha_q+\beta^{(d)}_q,\epsilon_q] \in O(q^{d+1})$, the
leading order term $\epsilon_{q,d} \in \g^1$ is closed. By assumption, one can therefore write $\beta^{(d)}_q = (\psi^{(d)})^*\alpha_q + \epsilon_q = (\psi^{(d)})^* \alpha_q + c q^d \alpha_{q,1} + q^d \delta b
+ O(q^{d+1})$ for some $c \in \C$, $b \in \g^0$. With
\begin{equation} \label{eq:new-approximation}
 \psi^{(d+1)}(q) = \psi^d(q) + c q^d, \quad
 \beta^{(d+1)}_q = \exp(-q^d b)(\beta^{(d)}_q)
\end{equation}
one finds that \eqref{eq:discrepancy} holds for $d+1$ as well. Since
the parameter change and gauge transformation in
\eqref{eq:new-approximation} are equal to the identity up to
$O(q^d)$, the process converges to $(\psi^{(\infty)})^*\alpha_q =
\beta^{(\infty)}_q$, with $\beta^{(\infty)}_q$ gauge equivalent to $\beta_q$.
\qed

When applied to $\g_\A = \mathit{CC}^{*+1}(\A,\A)^{\scriptscriptstyle \leq
0}$, this yields Lemma \ref{th:versal-3}. To get Lemma \ref{th:versal-1}, one needs to modify the framework a little. The traditional way to do that is to consider dg Lie algebras $\g$ which are pro-nilpotent. By this we mean that $\g$ comes with a complete decreasing filtration $\g = F^1\g \subset F^2\g \subset \cdots$, such that
\begin{equation}
d(F^u\g) \subset F^u g, \quad [F^u\g,F^v\g] \subset F^{u+v}\g.
\end{equation}
There is an associated group, defined in terms of formal exponentation, and this will act on the set of solutions of the Maurer-Cartan equation in $\g$. Use of the filtration replaces that of the parameter $q$ in ensuring convergence of order-by-order arguments, such as the proof of Lemma \ref{th:versal-abstract}. To obtain Lemma \ref{th:versal-1} in this way, one sets $\g[-1] = \prod_{t<0} \mathit{CC}^*(A,A)^t$, with $F^u\g$ the subspace where $t \leq -u$. There is also an alternative (but ultimately equivalent) viewpoint, which consists in introducing an auxiliary formal parameter. For that, let's think of the bigrading on the Hochschild complex of $A$ as giving a $\C^*$-action, which has weight $t$ on $\mathit{CC}^{s+t}(A,A)^t$. Combine this with the $\C^*$-action on $\Lambda_\N$ which has weight $d$ on $\C q^d$, so as to obtain an action on $\mathit{CC}^*(A,A)[1] \hat\otimes \Lambda_\N$. Given an $A_\infty$-structure $\A$ on $A$ satisfying \eqref{eq:mu1zero}, one can define a one-parameter deformation $\A_q$ of the trivial $A_\infty$-structure by $\mu^d_{\A_q} = q^{d-2}\mu^d_\A$, and this corresponds to a $\C^*$-invariant element of $\mathit{CC}^*(A,A)[1] \hat\otimes q\Lambda_\N$. The only possible changes of parameter in this context are the $\C^*$-equivariant (which means homogeneous) ones $q \mapsto \epsilon q$, $\epsilon \in \C$. An equivariant version of the argument from Lemma \ref{th:versal-abstract} can then be used to prove Lemma \ref{th:versal-1}.

\subsection{}
The previous discussion generalizes effortlessly to
$A_\infty$-categories. For expository reasons, start with a category $C$
with finitely many objects $\{X_1,\dots,X_m\}$; and let
$A$ be its total morphism algebra, which is linear over $R_m$.
Consider Hochschild cochains that preserve this structure, which
means that they are $R_m$-bimodule maps $A \otimes_{R_m} \dots
\otimes_{R_m} A \rightarrow A$. By writing this out, one sees that
the relevant complex can be generalized to arbitrary categories,
leading to a Hochschild cohomology theory $\mathit{HH}^*(C,C)$. The same
applies to Hochschild cohomology of $A_\infty$-categories, and to
the other notions of deformation theory.

\begin{remark}
For any category $C$ with finitely many objects, the Hochschild cohomology $\mathit{HH}^*(C,C)$ of the category is the same as that of the total
morphism algebra $A = \bigoplus_{i,j} \mathit{Hom}_C(X_i,X_j)$ seen as an algebra over $\C$.
In other words, the inclusion of the subcomplex of $R_m$-multilinear
Hochschild cochains into the whole cochain complex is
a quasi-iso\-mor\-phism. The reason is that the relative bar resolution
\begin{equation}
\cdots \longrightarrow A \otimes_{R_m} A \otimes_{R_m} \otimes A
\longrightarrow A \otimes_{R_m} A \longrightarrow A
\end{equation}
is also a resolution by projective $A$-bimodules in the non-relative
sense (which means over $\C$). By quasi-isomorphism invariance of abstract deformation
theories, this shows that any $A_\infty$-structure $\A$ on $A$
satisfying \eqref{eq:mu1zero} comes, up to gauge equivalence, from an
$A_\infty$-category structure $\Cat$ on $C$. The previous observation
combined with the standard spectral sequence argument also proves
that for any $A_\infty$-category $\Cat$ with finitely many objects and total morphism
$A_\infty$-algebra $\A$, we have
\begin{equation}
\mathit{HH}^*(\Cat,\Cat) \iso \mathit{HH}^*(\A,\A),
\end{equation}
and similarly for the truncated version (when that is defined, meaning in the case of vanishing differential). This implies that the one-parameter deformation theory of $\A$ as an $A_\infty$-algebra is
the same as that of $\Cat$ as an $A_\infty$-category. On occasion,
these facts will allow us to be somewhat sloppy about the distinction
between categories and their total morphism algebras.
\end{remark}

We will need an elementary criterion for non-triviality of the Hochschild cohomology classes associated to deformations of $A_\infty$-categories. Let $\Cat$ be an $A_\infty$-category. Suppose that it contains objects $Y_0,Y_1,Y_2$ with the following properties:
\begin{itemize}
\itemsep1em
\item For all $i \leq j$, the space $H^0(\mathit{hom}_{\Cat}(Y_i,Y_j))$ is one-dimensional.
\item The product
\begin{equation}
H^0(\mathit{hom}_{\Cat}(Y_1,Y_2)) \otimes H^0(\mathit{hom}_{\Cat}(Y_0,Y_1)) \longrightarrow
H^0(\mathit{hom}_{\Cat}(Y_0,Y_2))
\end{equation}
vanishes.
\end{itemize}

\begin{lemma} \label{th:product-turns-on}
Suppose that $\Cat_q$ is a deformation of $\Cat$, such that the map induced by $\mu^2_{\Cat_q}$,
\begin{equation} \label{eq:q2-prod}
\xymatrix{
H^0(\mathit{hom}_{\Cat_q}(Y_1,Y_2) \otimes_{\Lambda_\N} \C[q]/q^2) \otimes 
H^0(\mathit{hom}_{\Cat_q}(Y_0,Y_1) \otimes_{\Lambda_\N} \C[q]/q^2) \ar[d] \\
H^0(\mathit{hom}_{\Cat_q}(Y_0,Y_2) \otimes_{\Lambda_\N} \C[q]/q^2)}
\end{equation}
is nonzero. Then, the class in $\mathit{HH}^2(\Cat,\Cat)$ describing this deformation to first order is nontrivial.
\end{lemma}

\proof
Without loss of generality, we may assume that $\Cat$ is minimal, since that can be achieved by a quasi-isomorphism (more precisely, one should say that any deformation of $\Cat$ carries over to its minimal model, in such a way that the ``truncated'' cohomology groups appearing in \eqref{eq:q2-prod} are preserved; the first statement is abstract deformation theory, and the second is an easy Five Lemma argument). We may further assume that $\Cat$ is strictly unital. Take generators $x_2 \in \mathit{hom}_{\Cat}^0(Y_1,Y_2)$, $x_1 \in \mathit{hom}_{\Cat}^0(Y_0,Y_1)$. We necessarily have
\begin{equation} \label{eq:mu1-xk}
\mu^1_{\Cat_q}(x_k) = O(q^2),
\end{equation}
and $H^0(\mathit{hom}_{\Cat_q}(Y_{k-1},Y_k) \otimes_{\Lambda_\N} \C[q]/q^2)$ is generated by $[x_k]$; otherwise, that cohomology would be zero, in contradiction with the assumption on \eqref{eq:q2-prod}. Suppose that the result is false, meaning that the first order term of the deformation is the boundary of a Hochschild cochain $\gamma$. Taking minimality and \eqref{eq:mu1-xk} into account, one finds that
\begin{equation}
\begin{aligned}
& \textstyle\frac{1}{q}\big(\mu^2_{\Cat_q}(x_2,x_1) - \mu^2_{\Cat}(x_2,x_1)\big) \\ & = -\mu^3_{\Cat}(\gamma^0,x_2,x_1) -\mu^3_{\Cat}(x_2,\gamma^0,x_1) - \mu^3_{\Cat}(x_2,x_1,\gamma^0) \\ & \qquad - \mu^2_{\Cat}(\gamma^1(x_2),x_1) -\mu^2_{\Cat}(x_2,\gamma^1(x_1))
 \\ & \qquad + \textstyle \gamma^1(\mu^2_{\Cat}(x_2,x_1)) + O(q).
\end{aligned}
\end{equation}
By assumption $\mu^2_{\Cat}(x_2,x_1) = 0$, and since $\gamma^1(x_k)$ is necessarily a multiple of $x_k$, the other two $\mu^2_{\Cat}$ terms also vanish. Similarly, three different $\gamma^0 \in \mathit{hom}_{\Cat}^0(Y_k,Y_k)$ appear, but each is a multiple of the identity, hence by strict unitality the $\mu^3_{\Cat}$ terms vanish. This shows that $\mu^2_{\Cat_q}(x_2,x_1) \in O(q^2)$, which is a contradiction.
\qed

To make the connection with Section \ref{sec:generators}, we need to relate the deformation theory of an $A_\infty$-algebra with that of its derived category. Suppose that we have a finite set of split-generators of $\Cat$, whose total endomorphism $A_\infty$-algebra is $\A$. Take a deformation $\Cat_q$, and the associated deformation $\A_q$ of $\A$.  By the general derived invariance of Hochschild cohomology, the restriction map
\begin{equation}
\mathit{HH}^*(\Cat,\Cat) \longrightarrow \mathit{HH}^*(\A,\A)
\end{equation}
is an isomorphism. Hence, if $\Cat_q$ gives rise to a nontrivial class in $\mathit{HH}^2(\Cat,\Cat)$, then the same applies to the class of $\A_q$ in $\mathit{HH}^2(\A,\A)$. Assume from now on that $H(\A_q)$ is torsion-free. Without essential loss of generality, one can assume that $\A$ is minimal (has vanishing differential); and thanks to our assumption, the same applies to the deformation $\A_q$ (compare the proof of Corollary \ref{th:versal-3b}). Hence, the associated deformation class lies in $\mathit{HH}(\A,\A)^{\scriptscriptstyle \leq 0}$. 

\begin{lemma} \label{th:product-turns-on-2}
Suppose that $H(\A_q)$ is torsion-free, and that the deformation $\Cat_q$ gives rise to a nontrivial class in $\mathit{HH}^2(\Cat,\Cat)$. Suppose also that each split-generator $X_i$ in the given set satisfies $H^1(\mathit{hom}_{\A}(X_i,X_i)) = 0$. Then the deformation class of $\A_q$ in $\mathit{HH}^2(\A,\A)^{\scriptscriptstyle \leq 0}$ is nontrivial.
\end{lemma}

\proof
By definition, we have a long exact sequence
\begin{equation}
\cdots \rightarrow H^1(\mathit{CC}^*(\A,\A)^{\scriptscriptstyle > 0}) \longrightarrow
\mathit{HH}^2(\A,\A)^{\scriptscriptstyle \leq 0} \longrightarrow \mathit{HH}^2(\A,\A) \rightarrow \cdots
\end{equation}
where $\mathit{CC}^*(\A,\A)^{\scriptscriptstyle > 0} = \mathit{CC}^*(\A,\A)/\mathit{CC}^*(\A,\A)^{\scriptscriptstyle \leq 0}$ is the quotient complex of Hoch\-schild cochains $\A^{\otimes s} \rightarrow \A[t]$ with $t \geq 1$. These contribute in total degree $s+t$, hence $\mathit{CC}^1(\A,\A)^{\scriptscriptstyle > 0} = \bigoplus_i \mathit{hom}_\A^1(X_i,X_i)$, which vanishes by assumption. Hence, the Hochschild cohomology class associated to our deformation in $\mathit{HH}^2(\A,\A)^{\scriptscriptstyle \leq 0}$ is necessarily nontrivial.
\qed

\section{Group actions\label{sec:group-actions}}

The semidirect product $\Lambda V \semidirect \Gamma$ of an exterior algebra and a finite abelian group appears naturally in the context of homological mirror symmetry for Calabi-Yau hypersurfaces in projective spaces. This is a consequence of Serre duality and Beilinson's
description of $D^b\mathit{Coh}(\CP{n})$ \cite{beilinson78}. The Hochschild
cohomology of such algebras can be computed by a version of the
classical Hochschild-Kostant-Rosenberg (HKR) theorem
\cite{hochschild-kostant-rosenberg62}, and this will later form the
basis for applying the deformation theory from the previous section.
Besides discussing that theorem, we also use the opportunity to assemble
some elementary facts about semidirect products, both of algebras and
$A_\infty$-algebras. The exposition here is by no means as general as
possible; readers interested in learning more could start with \cite{montgomery92}
(note that what we call the ``semidirect product'' $A \semidirect \Gamma$ also appears in the literature as ``skew group ring'' $A \ast \Gamma$ or ``smash product'' $A \# \Gamma$).

\subsection{\label{subsec:actions-generalities}}
Let $A$ be a graded (unital, as usual) algebra carrying an action $\rho$ of a finite group $\Gamma$. The semidirect product $A \semidirect \Gamma$ (or $A \semidirect_\rho \Gamma$ if one wants to remember the action) is the tensor product $A \otimes \C \Gamma$ with the twisted multiplication
\begin{equation}
 (a_2 \otimes \gamma_2) (a_1 \otimes \gamma_1) =
 a_2 \rho(\gamma_2)(a_1) \otimes \gamma_2\gamma_1.
\end{equation}
The purpose of this construction is to turn the automorphisms $\rho(\Gamma) \subset \mathit{Aut}(A)$ into inner ones: in fact, if $e \in \Gamma$ is the unit then $(1 \otimes \gamma) (a \otimes e) (1 \otimes \gamma)^{-1} = \rho(\gamma)(a) \otimes e$. Therefore, it is not surprising that $A \semidirect \Gamma$ depends
only on the ``outer part'' of the group action, in the following sense: suppose that $\tilde{\rho}$ is another action related to the previous one by
\begin{equation} \label{eq:conjugate-action}
 \tilde{\rho}(\gamma) = \iota(\gamma) \rho(\gamma) \iota(\gamma)^{-1},
\end{equation}
where $\iota: \Gamma \rightarrow A^\times$ is a map (landing in the multiplicative group of degree $0$ invertible elements) which satisfies $\iota(\gamma_2\gamma_1) = \iota(\gamma_2) \rho(\gamma_2)(\iota(\gamma_1))$. Then $A \semidirect_\rho \Gamma$ and $A \semidirect_{\tilde\rho} \Gamma$ are isomorphic by $a \otimes \gamma \mapsto a \iota(\gamma)^{-1} \otimes \gamma$. This kind of consideration appears naturally in connection with the following construction: suppose that we have a short exact sequence of groups
\begin{equation}
0 \rightarrow \Gamma_0 \longrightarrow \Gamma_1 \longrightarrow \Gamma_1/\Gamma_0 \rightarrow 0
\end{equation}
together with a splitting $\theta: \Gamma_1/\Gamma_0 \rightarrow \Gamma_1$, so that $\Gamma_1 \iso
\Gamma_0 \semidirect (\Gamma_1/\Gamma_0)$ (here, $\semidirect$ means the ordinary semidirect product of groups \cite[p.\ 68]{suzuki}; recall that, unless the groups are abelian, the existence of a splitting does not mean that $\Gamma_1$ is actually the product of $\Gamma_0$ and $\Gamma_1/\Gamma_0$). If $\Gamma_1$ acts on $A$, we then have
\begin{equation} \label{eq:in-stages}
 A \semidirect \Gamma_1
 \iso (A \semidirect \Gamma_0) \semidirect (\Gamma_1/\Gamma_0).
\end{equation}
More precisely, let $\rho_{\Gamma_1}$ be the given action, and $\rho_{\Gamma_0}$ its restriction. Through the splitting $\theta$, we also have an action $\rho_{\Gamma_1/\Gamma_0}$ of the quotient on $A \semidirect \Gamma_0$:
\begin{equation} \label{eq:quotient-action}
 \rho_{\Gamma_1/\Gamma_0}([\gamma_2])(a_1 \otimes \gamma_1) =
 \rho_{\Gamma_1}(\theta([\gamma_2]))(a_1) \otimes
 \theta([\gamma_2]) \gamma_1 \theta([\gamma_2])^{-1},
\end{equation}
and these are the three actions involved in \eqref{eq:in-stages}. Different choices of $\theta$ lead to actions \eqref{eq:quotient-action} which differ by a map $\iota: \Gamma_1/\Gamma_0 \rightarrow (A \semidirect \Gamma_0)^\times$ as in \eqref{eq:conjugate-action}. The discussion there shows that the right hand side of \eqref{eq:in-stages} is independent of $\theta$, as indeed it ought to be in order for the isomorphism to hold.

From now on, we will assume for simplicity that $\Gamma$ is abelian. The group ring $\C \Gamma$ has a basis of orthogonal idempotents, which are the normalized characters $e_\chi = |\Gamma|^{-1}\chi$. Now $A \semidirect \Gamma$ is a (unital) algebra over $\C\Gamma$, hence can be thought of as a category with objects $X_\chi$ corresponding to characters, and morphism spaces
\begin{equation}
 \mathit{Hom}_{A \semidirect \Gamma}(X_{\chi_1},X_{\chi_2}) = e_{\chi_2} (A \semidirect \Gamma)
 e_{\chi_1}.
\end{equation}

By construction, $A \semidirect \Gamma$ has a splitting into pieces $A \otimes \C \gamma$, $\gamma \in \Gamma$, which is such that the product of the $\gamma_2$-piece and $\gamma_1$-piece lands in the $\gamma_2\gamma_1$-piece. Since we are assuming that $\Gamma$ is abelian, one can think of this splitting as being given by an action $\rho_{\Gamma^*}$ of the dual group $\Gamma^* = \mathit{Hom}(\Gamma,\C^*)$. Explicitly, 
\begin{equation}
\rho_{\Gamma^*}(\chi)(a \otimes \gamma) = \chi(\gamma)a \otimes \gamma. 
\end{equation}


\subsection{}
Because we are working over a characteristic zero field, taking $\Gamma$-invariants is an exact functor. In particular, if $\Gamma$ acts on a graded algebra $A$, $\mathit{HH}^*(A,A)^\Gamma$ is the cohomology of the subcomplex of $\mathit{CC}^*(A,A)^\Gamma$ of equivariant Hochschild cochains. There is a canonical chain map
\begin{equation} \label{eq:invariant-include-2}
\begin{aligned}
& \mathit{CC}^*(A,A)^\Gamma \longrightarrow \mathit{CC}^*(A \semidirect \Gamma, A \semidirect \Gamma)^{\Gamma^*}, \;\;
\tau \longmapsto (\tau \semidirect \Gamma), \text{ where} \\
& (\tau \semidirect \Gamma)(a_s \otimes \gamma_s,\dots,a_1 \otimes \gamma_1) = \\
& \qquad \qquad = \tau(a_s,\rho(\gamma_s)(a_{s-1}),\rho(\gamma_s
 \gamma_{s-1})(a_{s-2}),\dots) \otimes \gamma_s \cdots \gamma_1.
\end{aligned}
\end{equation}

\begin{lemma} \label{th:invariant-include-2}
\eqref{eq:invariant-include-2} is a quasi-isomorphism (of dg Lie algebras, if one includes the Lie bracket on $\mathit{CC}^{*+1}$ into our consideration). \qed
\end{lemma}

The map \eqref{eq:invariant-include-2} has an obvious deformation-theoretic meaning: $\Gamma$-equivariant $A_\infty$-deformations of $A$ induce deformations of the semidirect product, which are in their turn equivariant for the action of $\Gamma^*$. Lemma \ref{th:invariant-include-2} says that the two deformation problems are equivalent. This should be thought of as a deformation-theoretic version of results about duality of group actions \cite{blattner-montgomery85}. Even though that is the most natural explanation, we also want to outline an approach to the injectivity of the cohomology level map induced by \eqref{eq:invariant-include-2}, which is formulated in more classical terms. If $P$ is an $A$-bimodule equivariant with respect to $\Gamma$, then $P \semidirect \Gamma$ is an $A \semidirect \Gamma$-bimodule. For any $A \semidirect \Gamma$-bimodule $Q$, this satisfies
\begin{equation}
 \mathit{Hom}_{A \semidirect \Gamma-A \semidirect \Gamma}(P \semidirect \Gamma, Q)
 \iso \mathit{Hom}_{A-A}(P,Q)^\Gamma.
\end{equation}
Here, the $\Gamma$-action on the right hand side is obtained from that on $P$ and the diagonal action on $Q$, which is $q \mapsto (1 \otimes \gamma)q(1 \otimes \gamma^{-1})$. In particular, if $P$ is projective then so is $P \semidirect \Gamma$. Applying this observation to a projective $A$-bimodule resolution $P^* \rightarrow A$ which carries a compatible $\Gamma$-action, one obtains
\begin{equation} \label{eq:hh-gamma}
\begin{aligned}
 \mathit{HH}^*(A \semidirect \Gamma, A \semidirect \Gamma)
 & = \mathit{Ext}_{A \semidirect \Gamma - A \semidirect \Gamma}^*(A \semidirect \Gamma,
 A \semidirect \Gamma) \\
 & = \mathit{Hom}_{A \semidirect \Gamma - A \semidirect \Gamma}^*(P^* \semidirect \Gamma,A
 \semidirect \Gamma) \\
 & = \mathit{Hom}^*_{A-A}(P^*, A \semidirect \Gamma)^\Gamma.
\end{aligned}
\end{equation}
The splitting of $A \semidirect \Gamma$ into summands $A \otimes \C \gamma$, $\gamma \in \Gamma$, is compatible with the $A$-bimodule structure. Since $\Gamma$ is abelian, each summand is invariant under the diagonal action, hence we get an induced splitting of \eqref{eq:hh-gamma}. In particular, $\gamma = e$ contributes $\mathit{Hom}^*_{A-A}(P^*,A)^\Gamma = \mathit{HH}^*(A,A)^\Gamma$, which is precisely the image of \eqref{eq:invariant-include-2}. More generally, each summand $A \otimes \C \gamma$ is isomorphic to the graph bimodule $\mathit{Graph}(\rho(\gamma))$, which is a copy of the diagonal bimodule with the module structure on one side twisted by the automorphism $\rho(\gamma)$. Therefore,
\begin{equation} \label{eq:summands}
 \mathit{HH}^*(A \semidirect \Gamma,A \semidirect \Gamma) =
 \bigoplus_{\gamma} \mathit{Ext}^*_{A-A}\big(A,\mathit{Graph}(\rho(\gamma))\big)^\Gamma.
\end{equation}

The HKR theorem describes the Hochschild cohomology of polynomial
algebras (and more general smooth commutative algebras). There is a
generalization to graded algebras \cite[Proposition 5.4.6]{loday},
including exterior algebras as a special case, for which the
statement is as follows. Take $A = \Lambda V$ for some
finite-dimensional complex vector space $V$, with the usual grading.
Let $\mathit{Sym} V^\vee$ be the symmetric algebra of the dual space. Then there
is an explicit isomorphism
\begin{equation} \label{eq:classical-hkr}
 \phi_{\mathit{HKR}}: \mathit{Sym}^sV^\vee \otimes \Lambda^{s+t}V
 \longrightarrow \mathit{HH}^{s+t}(A,A)^t.
\end{equation}
For instance, the cocycle $[\tau] = \phi_{\mathit{HKR}}(p \otimes 1)$
corresponding to a homogeneous polynomial $p$ of degree $s$ is (in
terms of a basis $v_j$ and dual basis $v_j^\vee$)
\begin{equation}
 \tau(a_s, \dots, a_1) = \frac{1}{s!} \sum_{j_s,\dots,j_1}
 \frac{\partial p}{\partial v_{j_1} \dots \partial v_{j_s}}
 (v_{j_s}^\vee \inner a_s) \wedge \dots \wedge (v_{j_1}^\vee \inner a_1);
\end{equation}
and conversely, if $\tau \in \mathit{CC}^0(A,A)^{-s}$ is any cocycle, we
recover a homogeneous polynomial $p \otimes 1 =
\phi_{\mathit{HKR}}^{-1}([\tau])$ by setting
\begin{equation} \label{eq:reverse-hkr}
 p(v) = \tau(v,\dots,v).
\end{equation}
$\phi_{\mathit{HKR}}$ is compatible with the actions of $\mathit{GL}(V)$ on both sides.
Moreover, it carries the canonical bracket on $\mathit{HH}^{*+1}(A,A)$ to the
Schouten bracket on $\mathit{Sym} V^\vee \otimes \Lambda V$, which in particular
satisfies
\begin{equation} \label{eq:schouten}
 [v_{i_1}^\vee \dots v_{i_r}^\vee \otimes 1,
 v_{i_{r+1}}^\vee \dots v_{i_s}^\vee \otimes a] =
 \sum_{k \leq r}
 v_{i_1}^\vee \dots \widehat{v_{i_k}^\vee} \dots v_{i_s}
 \otimes v_{i_k}^\vee \inner a
\end{equation}
for $a \in \Lambda V$. The standard proof of the HKR formula is via
the Koszul bimodule resolution $P^* \rightarrow A$, which in our case
consists of $P^s = A \otimes \mathit{Sym}^s V \otimes A[-s]$ with the
differential
\begin{multline}
 \qquad \partial(1 \otimes v_{i_1} \dots v_{i_s} \otimes 1) =
 \sum_k v_{i_k} \otimes v_{i_1} \dots \widehat{v_{i_k}} \dots v_{i_s}
 \otimes 1 + \\
 +(-1)^s 1 \otimes
 v_{i_1} \dots \widehat{v_{i_k}} \dots v_{i_s} \otimes v_{i_k}. \qquad
\end{multline}
Because $A$ is supercommutative, the induced differential on
$\mathit{Hom}(P^*,A)$ vanishes, which shows that an abstract isomorphism
\eqref{eq:classical-hkr} exists; of course, the explicit formula and
other properties have to be verified separately, by exhibiting a map
between the Koszul and bar resolutions.

We now look at \eqref{eq:summands} in the special case where $A$ is
an exterior algebra, and $\rho$ the action of a finite abelian group
$\Gamma \subset \mathit{GL}(V)$. Taking the Koszul resolution as $P^*$, one
sees that the $\gamma$-summand of $\mathit{HH}^*(A \semidirect \Gamma,A
\semidirect \Gamma)$ can be obtained by computing the cohomology of
the complex $\mathit{Hom}_{\C}(SV,\Lambda V) \iso \mathit{Sym} V^\vee \otimes \Lambda V$
with the differential
\begin{equation}
 (\delta_\gamma\eta)(v_{i_1} \dots v_{i_s}) =
 (-1)^s \sum_k \eta(v_{i_1} \dots \widehat{v_{i_k}} \dots
 v_{i_s}) \wedge (\gamma(v_{i_k}) - v_{i_k}),
\end{equation}
and then restricting to the $\Gamma$-invariant part of that. To carry
out the computation, we decompose $V$ into eigenspaces of $\gamma$,
and consider the corresponding tensor product decomposition of $\mathit{Sym}
V^\vee \otimes \Lambda V$. On the fixed part $V^\gamma$, the
differential vanishes as before; the remaining components can be
identified with classical Koszul complexes, which have
one-dimensional cohomology. After chasing the bidegrees and
$\Gamma$-action through the argument, this amounts to a proof of the
following:

\begin{prop} \label{th:twisted-hkr}
The Hochschild cohomology of $A = \Lambda V \semidirect \Gamma$, for
$\Gamma$ a finite abelian subgroup of $\mathit{GL}(V)$, is
\begin{multline}
 \mathit{HH}^{s+t}(A, A)^t \iso \\ \iso
 \bigoplus_{\gamma \in \Gamma} \Big( \mathit{Sym}^s(V^\gamma)^\vee \otimes
 \Lambda^{s+t-\mathrm{codim}(V^\gamma)}(V^\gamma) \otimes
 \Lambda^{\mathrm{codim}(V^\gamma)}(V/V^\gamma)
 \Big)^\Gamma. \qed
\end{multline}
\end{prop}

\subsection{\label{subsec:ainfty-actions}}
A $\Gamma$-action on an $A_\infty$-algebra $\A$ is a linear action on the underlying graded vector space, such that all the $\mu_\A^d$ are equivariant (this may seem quite naive as a definition, but it will turn out to be sufficient for our purposes). The semidirect product $A_\infty$-algebra $\A \semidirect \Gamma$ is then defined as $\A \otimes \C\Gamma$ with the compositions
\begin{multline}
 \qquad \mu_{\A \semidirect \Gamma}^d(a_d \otimes \gamma_d, \dots, a_1 \otimes \gamma_1) = \\ =
 \mu_\A^d(a_d, \rho(\gamma_d)(a_{d-1}), \rho(\gamma_d\gamma_{d-1})(a_{d-2}), \dots)
 \otimes \gamma_d \cdots \gamma_1. \qquad
\end{multline}
Most of the discussion of semidirect products above applies to the $A_\infty$-case without significant changes (the only exception is the part involving $A^\times$; but one can at least allow ``invertible elements'' $\A^\times = \C^* e$ which are nonzero multiples of the identity, provided that $\A$ is strictly unital).

In the same way, one can consider $\Gamma$-actions on $A_\infty$-categories (if there are only finitely many objects, these are the same as $R_m$-linear $\Gamma$-actions on $R_m$-linear $A_\infty$-algebras). We will need the following weak ``transfer'' property:

\begin{lemma} \label{th:transfer-group-action}
Let $\Dat$ be an $A_\infty$-category with an action of a finite group $\Gamma$, and $D = H(\Dat)$ the associated cohomology level category. Suppose that we have another graded linear category $C$ with a $\Gamma$-action, and an equivalence $F: D \rightarrow C$ which is equivariant. Then there is an $A_\infty$-category $\Cat$ with $H(\Cat)$ isomorphic (not just equivalent) to $C$, which carries an action of $\Gamma$, and an equivariant $A_\infty$-functor $\F: \Dat \rightarrow \Cat$ such that $H(\F) = F$.
\end{lemma}

\proof We first remind the reader of the argument in the non-equivariant case. Without loss of generality, one can assume that $\Dat$ has vanishing differential (by applying the Perturbation Lemma). Consider the bigraded complex
\begin{equation} \label{eq:z-complex}
Z^{s+t,t} = \mathit{CC}^{s+t-1}(D,C)^t \oplus \mathit{CC}^{s+t}(C,C)^t.
\end{equation}
Here, the first piece is Hochschild cochains of $D$ with coefficients in the $D$-bimodule $C$ (formed by using $F$ on both sides). The differential $\delta_Z$ combines the two Hochschild differentials with the pullback map $\mathit{CC}^*(C,C) \rightarrow \mathit{CC}^*(D,C)$. Since $F$ is an equivalence, the pullback map is a quasi-isomorphism, hence $Z$ is acyclic. One now constructs $\F$ and $\mu_{\Cat}$ simultaneously by solving an obstruction problem in $Z$. Start with the trivial ansatz where $\F$ has only the linear term $\F^1 = F$, and $\mu^d_{\Cat}$ similarly is nontrivial only for $d = 2$, where it reproduces the composition of morphisms in $C$. The failure of this ansatz to satisfy the required equation, at lowest order, is an element of $Z^{3,-1}$ consisting of
\begin{equation} \label{eq:failure}
\begin{aligned}
& F \circ \mu^3_\Dat \in \mathit{CC}^2(D,C)^{-1}, \\ 
& 0 \in \mathit{CC}^3(C,C)^{-1}. 
\end{aligned}
\end{equation}
Because $F$ is a functor and $\Dat$ is an $A_\infty$-category, \eqref{eq:failure} is a cocycle in \eqref{eq:z-complex}. Acyclicity shows that there is a bounding cochain in $Z^{2,-1}$. By the definition of the differential on $Z$, the components
\begin{equation}
\begin{aligned}
& \F^2 \in \mathit{CC}^1(D,C)^{-1}, \\
& \mu^3_\Cat \in \mathit{CC}^2(C,C)^{-1}
\end{aligned}
\end{equation}
of that bounding cochain must satisfy
\begin{equation}
\begin{aligned}
& F(\mu^3_\Dat(d_3,d_2,d_1)) = (-1)^{|d_2|+|d_1|-1} F(d_3)\F^2(d_2,d_1) + (-1)^{|d_1|} \F^2(d_3,d_2)F(d_1) \\ & \qquad 
+ \mu^3_\Cat(F(d_3),F(d_2),F(d_1)) \\ & \qquad + (-1)^{|d_1|+|d_2|} \F^2(d_3d_2,d_1) + (-1)^{|d_1|-1} \F^2(d_3,d_2d_1),
\\
& 0 = (-1)^{|c_1|+|c_2|+|c_3|} \mu^3_\Cat(c_4c_3,c_2,c_1) + (-1)^{|c_1|+|c_2|-1} \mu^3_\Cat(c_4,c_3c_2,c_1) \\
& \qquad + (-1)^{|c_1|} \mu^3_\Cat(c_4,c_3,c_2c_1) \\ & \qquad
+ (-1)^{|c_1|+|c_2|+|c_3|-1} c_4\mu^3_\Cat(c_3,c_2,c_1) - \mu^3_\Cat(c_4,c_3,c_2)c_1.
\end{aligned}
\end{equation}
As already indicated in the notation, that bounding cochain contains the next term $\F^2$ of the $A_\infty$-functor, together with the next term $\mu^3_\Cat$ of the $A_\infty$-structure on $\Cat$. The subsequent steps are parallel, in that one constructs $(\F^{d-1},\mu_\Cat^d)$ at the same time order-by-order in $d$. In each step, the equation that these must satisfy can be written in the form $\delta_Z(\F^{d-1},\mu_\Cat^d) = \cdots$, where the right hand side depends only on the initial data of the problem and the previously chosen $\F^2,\dots,\F^{d-2}$, $\mu_\Cat^3,\dots,\mu_\Cat^{d-1}$. One shows that this right hand side is a cocycle, and then uses acyclicity to solve the equation.

In the presence of a group action, one can apply the Perturbation Lemma to $\Dat$ equivariantly (this, and other parts of the proof, rely crucially on the fact that our ground field has characteristic $0$). Then, the same argument as before applies to the invariant part $Z^\Gamma$ of \eqref{eq:z-complex}, yielding an equivariant $A_\infty$-structure and $A_\infty$-functor. \qed

As a final element of this circle of ideas, suppose that we have an $A_\infty$-category $\Cat$ with an action of $\Gamma$ which is {\em trivial on objects}. One can then introduce an $A_\infty$-category $\mathit{Tw}_{\Gamma} \Cat$ of equivariant twisted complexes, which has a $\Gamma$-action with the same property (and if one forgets that action, $\mathit{Tw}_{\Gamma} \Cat$ becomes equivalent to a full subcategory of the ordinary category of twisted complexes $\mathit{Tw} \Cat$). The definition is similar to the standard one for twisted complexes, with the following modifications:
\begin{itemize}
\itemsep1em
\item The equivariant version $\Sigma_\Gamma\Cat$ of the additive completion
has objects which are formal direct sums
\begin{equation} \label{eq:bigoplus-f}
 \bigoplus_{f} (X_f \otimes V_f)[\sigma_f],
\end{equation}
where the $V_f$ are finite-dimensional representations of $\Gamma$. The morphism spaces are correspondingly
\begin{multline}
 \qquad \mathit{hom}_{\Sigma_\Gamma\Cat}\Big(\bigoplus_f (X_f \otimes V_f)[\sigma_f],
 \bigoplus_g (Y_g \otimes W_g)[\tau_g]\Big) = \\ =
 \bigoplus_{f,g} \mathit{hom}_{\Cat}(X_f,Y_g)[\tau_g-\sigma_f] \otimes \mathit{Hom}_\C(V_f,W_g), \qquad
\end{multline}
and they carry the obvious tensor product $\Gamma$-actions.

\item In the definition of an equivariant twisted complex, the differential
$\delta$ has to belong to the $\Gamma$-invariant part of $\mathit{hom}^1_{\Sigma_\Gamma\Cat}$.
\end{itemize}
The ``equivariant derived category'' $D^b_\Gamma(\Cat) = H^0(\mathit{Tw}_\Gamma\,\Cat)$ carries, in addition to the usual shift functor, endofunctors $- \otimes Z$ for any (graded finite-dimensional) representation $Z$ of $\Gamma$. In particular, we have an action of $\Gamma^*$ by tensoring with characters.

\begin{remark} \label{th:twist-is-equivariant}
$D^b_\Gamma(\Cat)$ is {\em not} triangulated, since mapping cones exist only for the $\Gamma$-invariant morphisms (whereas we are allowing non-equivariant ones in the category as well). Fortunately, the evaluation maps $\mathit{ev}$, $\mathit{ev}^\vee$ used in the definition of the twist and untwist operations are invariant, so $T_X(Y)$, $T_Y^\vee(X)$ are well-defined in $D^b_\Gamma(\Cat)$ (at least assuming that the necessary finite-dimensionality conditions hold). In particular, the theory of exceptional collections and mutation works well in the equivariant context.
\end{remark}

\section{Coherent sheaves\label{sec:sheaves}}

To apply the theory from Section \ref{sec:deformations} to the
derived category of coherent sheaves, one needs to introduce an
underlying differential graded category. Concretely, this means
choosing a class of cochain complexes which compute $Ext$ groups,
with cochain maps realizing the Yoneda product. Out of the many
possibilities offered by the literature, we take the lowest-tech path
via {\v C}ech cohomology, which is easy to work with but slightly
more difficult to relate to the usual derived category (a few readers
may prefer to skip that argument, and to take the {\v C}ech version
as definition of the derived category; in which case some of the
properties discussed later on have to be taken for granted).

\subsection{}
Let $Y$ be a noetherian scheme over an algebraically closed field
($\C$ or $\Lambda_\Q$, in the cases of interest here). Fix a finite
affine open cover $\UU$. We introduce a differential graded category
$\cech{\SS}(Y)$ whose objects are locally free coherent sheaves
(algebraic vector bundles), and where the morphism spaces are {\v
C}ech cochain complexes with values in $\SHom$ sheaves,
\begin{equation}
\mathit{hom}_{\cech{\SS}(Y)}(E_0,E_1) = \cech{C}^*(\UU,\SHom(E_0,E_1)).
\end{equation}
Composition is the shuffle product combined with the local tensor product maps $\SHom(E_1,E_2)
\otimes_{\OO_Y} \SHom(E_0,E_1) \rightarrow \SHom(E_0,E_2)$. One can think of
any dg category as an $A_\infty$-category with vanishing compositions
of order $\geq 3$, and therefore define a triangulated category
$D^b\cech{\SS}(Y)$ in the way explained in Section
\ref{sec:generators}, as the cohomological category of twisted
complexes (we emphasize that there is no ``inverting of
quasi-isomorphisms'' involved in this process). By refining covers
and using standard results about {\v C}ech cohomology, one sees that
$\cech{\SS}(Y)$ is independent of $\UU$ up to quasi-isomorphism,
hence that $D^b\cech{\SS}(Y)$ is well-defined up to equivalence of
triangulated categories.

\begin{lemma} \label{th:derived-is-derived}
$D^b\cech{\SS}(Y)$ is equivalent to the full triangulated subcategory
of the usual bounded derived category of coherent sheaves,
$D^b\mathit{Coh}(Y)$, which is generated by locally free sheaves. If $Y$ is a
regular projective variety, this is the whole of $D^b\mathit{Coh}(Y)$.
\end{lemma}

\proof The usual approach to dg structures on derived categories of
coherent sheaves is through injective resolutions, see for instance
\cite[Section 2]{seidel-thomas99}. In our current language, this can
be formulated as follows. For each locally free coherent sheaf $E$,
fix a resolution $I^*_E$ by quasi-coherent injective sheaves.
Consider a differential graded category $\I(Y)$ with the same objects
as $\cech{\SS}(Y)$, but where $\mathit{hom}_{\I(Y)}(E_0,E_1) =
\mathit{hom}_{\O_Y}(I^*_{E_0},I^*_{E_1})$ is the space of sheaf homomorphisms $I^*_{E_0}
\rightarrow I^*_{E_1}$ of all degrees, with the obvious differential. If
we then define $D^b\I(Y) = H^0(\mathit{Tw}\,\I(Y))$ through twisted complexes,
it is equivalent to the subcategory of $D^b\mathit{Coh}(Y)$ appearing in the
Lemma. To see that, one associates to each object of $\mathit{Tw}\,\I(Y)$ its
total complex; this defines an exact, fully faithful functor from
$D^b\I(Y)$ to the chain homotopy category ${\mathcal K}^+(Y)$ of
bounded below complexes of injective quasi-coherent sheaves. The
image of this functor is the triangulated subcategory generated by
resolutions of locally free sheaves, and by standard properties of
derived categories, this is equivalent to the triangulated
subcategory of $D^b\mathit{Coh}(Y)$ generated by such sheaves. If $Y$ is
smooth and projective, all coherent sheaves have finite locally free
resolutions, so we get the whole of $D^b\mathit{Coh}(Y)$.

To mediate between {\v C}ech and injective resolutions, we use a
triangular matrix construction. Given two locally free sheaves $E_0$
and $E_1$, one can consider the {\v C}ech complex
$\cech{C}^*(\SHom(E_0,I_{E_1}^*))$ with coefficients in the complex of
sheaves $\SHom(E_0,I_{E_1}^*) \iso E_0^\vee \otimes_{\O_Y} I_{E_1}^*$. There are
canonical cochain maps
\begin{equation} \label{eq:t-maps}
 \mathit{hom}_{\O_Y}(I_{E_0}^*,I_{E_1}^*) \longrightarrow
 \cech{C}^*(\SHom(E_0,I_{E_1}^*)) \longleftarrow \cech{C}^*(\SHom(E_0,E_1)).
\end{equation}
The first one is defined by composing a sheaf homomorphism $I_{E_0}^*
\rightarrow I_{E_1}^*$ with the resolution map $E_0 \rightarrow I_{E_0}^*$, and then restricting
the resulting section of $\SHom(E_0,I_{E_1}^*)$ to open subsets in $\UU$.
The second one is just the map induced from ${E_1} \rightarrow I_{E_1}^*$.
Both maps are quasi-isomorphisms: for the first map, note that the
{\v C}ech complex computes the hypercohomology of $E_0^\vee \otimes
I_{E_1}^*$, which on the other hand is the ordinary cohomology of
$\mathit{hom}_{\O_Y}(E_0^*,I_{E_1}^*)$ because each $E_0^\vee \otimes I_{E_1}^k$ is
injective. For the second map one uses a standard spectral sequence
comparison argument, together with the fact that the cohomology of
$\mathit{hom}_{\O_U}(E_0|U, I_{E_1}^*|U)$ for each affine open subset $U$ is
$\mathit{Hom}_{\O_U}(E_0|U,{E_1}|U)$. Along the same lines, there are canonical
homomorphisms
\begin{equation} \label{eq:t-bistructure}
\begin{split}
 & \mathit{hom}_{\O_Y}(I_{E_1}^*,I_{E_2}^*) \otimes \cech{C}^*(\SHom(E_0,I_{E_1}^*))
 \longrightarrow \cech{C}^*(\SHom(E_0,I_{E_2}^*)), \\
 & \cech{C}^*(\SHom(E_1,I_{E_2}^*)) \otimes \cech{C}^*(\SHom(E_0,E_1))
 \longrightarrow \cech{C}^*(\SHom(E_0,I_{E_2}^*)).
\end{split}
\end{equation}
We introduce another dg category $\T$, still keeping the same objects
as before, and where the $\mathit{hom}$s are the direct sums of three
components
\begin{equation}
 \mathit{hom}_\T(E_0,E_1) =
 \begin{pmatrix}
 \mathit{hom}_{\O_Y}(I_{E_0}^*,I_{E_1}^*) & \cech{C}^*(\SHom(E_0,I_{E_1}^*))[1] \\
 0 & \cech{C}^*(\SHom(E_0,E_1)).
 \end{pmatrix}
\end{equation}
The differential consists of the given differentials on each summand
together with the chain homomorphisms \eqref{eq:t-maps}. Similarly,
composition in $\T$ combines the usual compositions of morphisms in
$\I(Y)$ and $\cech{\SS}(Y)$ with the maps \eqref{eq:t-bistructure}.
By definition, this comes with dg functors $\T \rightarrow \I(Y)$,
$\T \rightarrow \cech{\SS}(Y)$, both of which are quasi-isomorphisms,
and therefore one has induced exact equivalences $D^b\cech{\SS}(Y)
\iso D^b\T \iso D^b\I(Y)$. \qed

\begin{remark} \label{th:equivariant-sheaves}
At least if one supposes that the ground field has characteristic
zero, the entire argument carries over to the derived category of
equivariant coherent sheaves with respect to the action of a finite
group $\Gamma$ on $Y$. First, by taking intersections over
$\Gamma$-orbits one can find an invariant affine open cover. One then
introduces the dg category $\cech{S}_\Gamma(Y)$ whose objects are
locally free sheaves with $\Gamma$-actions, and whose morphism spaces
are the $\Gamma$-invariant summands of the {\v C}ech complexes. In
the regular projective case,
\begin{equation}
 D^b\cech{\SS}_\Gamma(Y) \iso D^b\mathit{Coh}_\Gamma(Y),
\end{equation}
where $\mathit{Coh}_\Gamma(Y)$ is the abelian category of $\Gamma$-equivariant
coherent sheaves. We will now make the connection between this
observation and the material from Section \ref{sec:group-actions}.
If $E$ is a coherent
sheaf with a $\Gamma$-action, then so is $E \otimes W$ for any
finite-dimensional $\Gamma$-module $W$. Suppose for simplicity that
$\Gamma$ is abelian, and take equivariant locally free sheaves
$E_1,\dots,E_m$. Then the full dg subcategory $\Cat \subset
\cech{\SS}(Y)$ with objects $E_k$ carries a $\Gamma$-action, and the
semidirect product $\Cat \semidirect \Gamma$ can be identified with
the full dg subcategory of $\cech{\SS}_\Gamma(Y)$ whose objects are
$E_k \otimes W_\chi$ for all $k$ and all one-dimensional
representations $W_\chi$.
\end{remark}

\subsection{\label{subsec:formal-schemes}}
Consider a Noetherian scheme $Y_q$ which is projective and flat over $\Lambda_\N$. Geometrically, this is a formal one-parameter family of schemes. Starting with this, one can carry out the following constructions:
\begin{itemize} \itemsep1em
\item Setting $q = 0$ yields the {\em special fibre} $Y = Y_q \times_{\Lambda_\N} \C$, which is a scheme over $\C$.
\item Inverting $q$ and passing to the algebraic closure yields the {\em general fibre} $Y_q^* = Y_q \times_{\Lambda_\N} \Lambda_\Q$, which is a scheme over $\Lambda_\Q$.
\item Finally, {\em formal completion} with respect to the ideal $(q)$ yields a formal scheme $\hat{Y}_q$. This is in a sense intermediate between $Y$ and $Y_q$: the topological space underlying $\hat{Y}_q$ agrees with that of $Y$, and it carries a sheaf of complete $\Lambda_\N$-algebras whose specialization to $q = 0$ recovers the structure sheaf of $Y$. On the other hand, Grothendieck's formal GAGA theorem (see for instance \cite[Theorem 3.2.1]{harbater03}) says that the categories of coherent sheaves on $Y_q$ and $\hat{Y}_q$ are equivalent.
\end{itemize}

There are corresponding constructions on the level of differential graded categories. Take a finite cover $\UU_q$ of $Y_q$ by affine open subsets. By the same procedure as before, one can define a differential graded category $\cech{\SS}(Y_q)$ whose objects are locally free sheaves over $Y_q$, and whose morphism spaces are flat $\Lambda_\N$-modules. The derived category $D^b\cech{\SS}(Y_q)$ is equivalent to a full triangulated subcategory of $D^b\mathit{Coh}(Y_q)$. Let $\UU_q^*$ be the associated cover of the general fibre $Y_q^*$. One can use that to define an analogous dg category $\cech{\SS}(Y_q^*)$ over $\Lambda_\Q$, whose derived category describes a part of $D^b\mathit{Coh}(Y_q^*)$. Moreover, this comes with a full and faithful dg functor
\begin{equation} \label{eq:invert-sheaf}
\cech{\SS}(Y_q) \otimes_{\Lambda_\N} \Lambda_\Q \longrightarrow \cech{\SS}(Y_q^*).
\end{equation}
On the other hand, one has an associated cover $\UU$ of the special fibre $Y$. Using that cover and the ringed space structure from $\hat{Y}_q$, we define another dg category $\cech{\SS}(\hat{Y}_q)$. The morphism spaces in this category are complete torsion-free modules over $\Lambda_\N$. It is the target of a completion dg functor
\begin{equation} \label{eq:formal-gaga}
\cech{\SS}(Y_q) \longrightarrow \cech{\SS}(\hat{Y}_q),
\end{equation}
which is cohomologically full and faithful by (the easier part of) the formal GAGA theorem. Finally, we have a full and faithful dg functor
\begin{equation} \label{eq:specialization}
\cech{\SS}(Y_q) \otimes_{\Lambda_\N} \C \longrightarrow \cech{\SS}(Y_0).
\end{equation}

This puts us in a situation to apply the formal deformation theory from Section \ref{sec:deformations}. To see explicitly how this works, fix a finite collection of locally free sheaves $E_{q,1},\dots,E_{q,m}$ on $Y_q$. The key object is the full dg subcategory of $\cech{\SS}(\hat{Y}_q)$ whose objects are the completions $\hat{E}_{q,1},\dots,\hat{E}_{q,m}$ (this is a topological dg category in the sense of Remark \ref{th:top-free}). Specialization to $q = 0$ recovers the full dg subcategory of $\cech{\SS}(Y)$ whose objects are the restrictions $E_1,\dots,E_m$ to the special fibre. On the other hand, taking the tensor product with $\Lambda_\Q$ yields a dg category which, by a combination of \eqref{eq:formal-gaga} and \eqref{eq:invert-sheaf}, is quasi-isomorphic to the full dg subcategory of $\cech{\SS}(Y_q^*)$ whose objects are the restrictions $E^*_{q,1},\dots,E^*_{q,m}$ of our sheaves to the generic fibre.


\subsection{}
We end our discussion by listing two properties of $D^b\mathit{Coh}(Y)$ for a
smooth projective variety over a field, which are well-known to
specialists (I am indebted to Drinfeld for suggesting simplified
proofs and giving references). The first result is classical, and
explains why passage to $D^\pi$ is unnecessary.

\begin{lemma}
$D^b\mathit{Coh}(Y)$ is idempotent closed.
\end{lemma}

One way to prove this is to consider the chain homotopy category $\K^+(Y)$, which is idempotent closed by an infinite direct sum trick \cite[Proposition 1.6.8]{neeman} or \cite[Lemma 2.4.8.1]{levine}.
Having done that, one recognizes that the direct summand of any object of $D^b\mathit{Coh}(Y) \subset \K^+(Y)$ again lies in $D^b\mathit{Coh}(Y)$, because its cohomology sheaves are coherent and almost all zero \cite[Proposition 2.4]{seidel-thomas99}. Alternatively, one can use the stronger result from \cite{bondal-vandenbergh02} which says that the bounded derived category is {\em saturated}.

\begin{lemma} \label{th:kontsevich}
If $\iota: Y \hookrightarrow {\mathbb P}^N$ is a projective
embedding, and $F_1,\dots,F_m$ are split-generators for
$D^b\mathit{Coh}({\mathbb P}^N)$, their derived restrictions $E_k =
L\iota^*(F_k)$ are split-gene\-ra\-tors for $D^b\mathit{Coh}(Y)$.
\end{lemma}

This is due to Kontsevich, and the proof follows the same pattern as
Lemma \ref{th:twist-generators}. By assumption, on ${\mathbb P}^N$
one can express $\O(m)$ for any $m$ as direct summand of a complex
built from the $F_k$. Hence, on $Y$ one can express $\O(m)$ in the
same way using $E_k$. It therefore suffices to show that all the
$\O(m)$ together are split-generators for $D^b\mathit{Coh}(Y)$. Any locally
free coherent sheaf $E$ on $Y$ has a finite left resolution of the
form
\begin{equation}
 0 \rightarrow E' \rightarrow \O(m_l)^{\oplus r_l} \rightarrow
 \cdots \rightarrow \O(m_1)^{\oplus r_1} \rightarrow E \rightarrow 0
\end{equation}
for some $E'$, and one can make $l$ arbitrarily large. For $l >
\mathrm{dim}\,Y$ we have $\mathit{Ext}^l(E,E') = H^l(Y,E^\vee \otimes E') = 0$, hence the
resulting exact triangle
\begin{equation}
 \xymatrix{
 \!\!\!\! {\big\{ \O(m_l)^{\oplus r_l} \rightarrow \dots \rightarrow
 \O(m_1)^{\oplus r_1} \big\}} \ar[r] & {E} \ar[d] \\
 & {E'[l]} \ar@/^1pc/[ul]^-{[1]}
}
\end{equation}
splits in such a way that $E \oplus E'[l-1]$ is isomorphic to the top
left complex. Once one has split-generated all locally free sheaves
from the $\O(m)$, the rest follows as in Lemma
\ref{th:derived-is-derived}.

\section{Symplectic terminology\label{sec:geometry}}

The purpose of this section is to assemble the necessary nuts and
bolts from elementary symplectic geometry. The origin of the main
notions, to the author's best knowledge, is as follows. Rational
Lagrangian submanifolds have a long history in quantization (going
back to Einstein \cite{einstein}); the first place where they occur in connection
with Floer homology is probably Fukaya's paper \cite{fukaya93}.
Lagrangian submanifolds with zero Maslov class also originated in
(Maslov's) quantization; the essentially equivalent notion of grading
was introduced by Kontsevich \cite{kontsevich94}. The symplectic
viewpoint on Lefschetz pencils and higher-dimensional linear systems
emerged gradually from the work of Arnol'd \cite{arnold95}, Donaldson
\cite{donaldson98,donaldson02}, Gompf \cite{gompf-stipsicz}, and
Auroux \cite{auroux99,auroux00}. In particular, matching cycles were
invented by Donaldson (unpublished).

\subsection{}
By a {\em projective K{\"a}hler manifold} we mean a compact complex
manifold $X$ carrying a unitary holomorphic line bundle $o_X$, which
is positive in the sense that the curvature of the corresponding
connection $\nabla_X$ defines a K{\"a}hler form
\begin{equation}
\o_X = \frac{i}{2\pi} F_{\nabla_X}.
\end{equation}
From a symplectic viewpoint, this is an integrality condition on the
symplectic class, which leads one to make corresponding restrictions on the
other objects which live on $X$. A Lagrangian submanifold $L \subset
X$ is called {\em rational} if the monodromy of the flat connection
$\nabla_X|L$ consists of roots of unity. This means that one can choose a
covariantly constant multisection $\lambda_L$ of $o_X|L$ which is of
unit length everywhere (as a multisection, this is the multivalued
$d_L$-th root of a section $\lambda_L^{d_L}$ of $o_X^{\otimes
d_L}|L$, for some $d_L \geq 1$).

\begin{remark} \label{th:deform-rational}
Rationality is invariant under exact Lagrangian isotopy. More precisely, take such an isotopy $(L_r)$, and denote by $\Lambda = \bigcup_r \{r\} \times L_r \subset [0,1] \times X$ its domain. The curvature of the pullback connection on $o_X|\Lambda$ can be written as $\frac{2\pi}{i} \omega_X|\Lambda = \frac{2\pi}{i} dH \wedge dr$ for some function $H$ (that is a time-dependent Hamiltonian function generating the isotopy). Subtracting $\frac{2\pi}{i} H \, dr$ makes the pullback connection flat, allowing one to transport the given $\lambda_{L_0}$ to multisections $\lambda_{L_r}$ for all $r$.
\end{remark}

Rationality takes its name from the effect on the periods of action
functionals. Take two rational Lagrangian submanifolds $L_0,L_1$. For
each of them, choose a covariantly constant unit length
$d_{L_k}$-fold multisection $\lambda_{L_k}$. The {\em mod $\Q$
action} $\bar{A}(x)$ of a point $x \in L_0 \cap L_1$ is defined by
\begin{equation} \label{eq:modq-action}
e^{2\pi i \bar{A}(x)} = \lambda_{L_1}(x)/\lambda_{L_0}(x).
\end{equation}
It is unique up to adding a number in $\lcm(d_{L_0},d_{L_1})^{-1} \Z
\subset \Q$, whence the name. Now suppose that we have two
intersection points $x_0,x_1$ and a connecting disc between them:
this is a smooth map $u: \R \times [0,1] \rightarrow X$ with $u(\R
\times \{k\}) \subset L_k$, such that $u(s,\cdot)$ converges to the
constant paths at $x_0,x_1$ as $s \rightarrow -\infty,+\infty$. By
definition of curvature,
\begin{equation} \label{eq:modq-energy}
\int_{\R \times [0,1]} u^*\o_X \in \bar{A}(x_0) - \bar{A}(x_1) +
\lcm(d_{L_0},d_{L_1})^{-1}\Z.
\end{equation}

Similarly, a symplectic automorphism $\phi$ of $X$ is rational if the
monodromy of the induced flat connection on $\mathit{Hom}(o_X,\phi^*o_X)$
consists of roots of unity; and this condition is invariant under Hamiltonian
isotopy. We write $\mathit{Aut}(X)$ for the group of such automorphisms, which
acts in the obvious way on the set of rational Lagrangian
submanifolds.

\subsection{}
Equally important for us will be the relative situation, where in
addition to the previous data we have a holomorphic section
$\sigma_{X,\infty}$ of $o_X$ such that $X_\infty =
\sigma_{X,\infty}^{-1}(0)$ is a divisor with normal crossings, along
which $\sigma_{X,\infty}$ vanishes with multiplicity one, and such
that each irreducible component $C \subset X_\infty$ is smooth (no
self-intersections). We call $M = X \setminus X_\infty$ an {\em
affine K{\"a}hler manifold}. One can use
$\sigma_{X,\infty}/\|\sigma_{X,\infty}\|$ to trivialize $o_X|M$, and
by writing the connection with respect to this trivialization as
$\nabla_X|M = d - 2\pi i \theta_M$, one gets a one-form $\theta_M$ which
satisfies $d\theta_M = \o_M = \o_X|M$, so that $M$ is an exact
symplectic manifold in the usual sense of the word.

A compact Lagrangian submanifold $L \subset M$ is called {\em exact} if $\theta_M|L = dK_L$ for some function $K_L$. Note that if this is the case, then
\begin{equation} \label{eq:k-lambda}
\lambda_L = \exp(2 \pi i K_L )
\frac{\sigma_{X,\infty}}{\|\sigma_{X,\infty}\|}
\end{equation}
is a covariantly constant section of $o_X|L$, so that must be the
trivial flat line bundle. Conversely, suppose that $L$ is some
Lagrangian submanifold of $M$, such that $o_X|L$ is trivial and
admits a covariantly constant section which can be written in the
form \eqref{eq:k-lambda} for some $K_L$. Then $dK_L = \theta_M|L$,
hence $L$ is exact. We have proved:

\begin{lemma} \label{th:check-exactness}
$L$ is exact iff $o_X|L$ is trivial, and its nonzero covariantly
constant sections lie in the same homotopy class of nowhere zero
sections as $\sigma_{X,\infty}|L$. \qed
\end{lemma}

\begin{remark} \label{th:deform-exact}
As in Remark \ref{th:deform-rational}, suppose that $(L_r)$ is an exact Lagrangian isotopy, with total space $\Lambda \subset [0,1] \times M$, and write $\omega_M|\Lambda = dH \wedge dr$. Then
\begin{equation} \label{eq:theta-lambda}
\theta_M|\Lambda - H \, dr \in \Omega^1(\Lambda)
\end{equation}
is closed. Given a function $K_{L_0}$ such that $dK_{L_0} = \theta_M|L_0$, there is a unique primitive $K$ for \eqref{eq:theta-lambda} which restricts to $K_{L_0}$ on the slice $r = 0$, and at the other end gives a function $K_{L_1}$ satisfying $dK_{L_1} = \theta_M|L_1$, in particular proving that $L_1$ is again exact.
\end{remark}

Let $L_0$, $L_1$ be exact Lagrangian submanifolds of $M$, and choose
functions $K_{L_0}$, $K_{L_1}$. One can then define the action of an
intersection point $x \in L_0 \cap L_1$ to be
\begin{equation} \label{eq:action}
A(x) = K_{L_1}(x) - K_{L_0}(x) \in \R.
\end{equation}
For any two such points $x_0,x_1$ and any connecting disc $u$ in $X$,
the intersection number $u \cdot X_\infty$ is well-defined, and
\begin{equation} \label{eq:energy}
\int_{\R \times [0,1]} u^*\o_X = A(x_0) - A(x_1) + (u \cdot X_\infty).
\end{equation}
In particular, if we allow only connecting discs which lie in $M$, the action gives an a priori bound on their symplectic areas.

Let $\phi$ be a symplectic automorphism of $M$, which extends
smoothly to $X$ in such a way that $X_\infty$ remains pointwise fixed
(we denote the extension equally by $\phi$). We say that $\phi$ exact
if $\phi^*\theta_M - \theta_M = dK_\phi$ for some function $K_\phi$
on $M$. This condition is invariant under Hamiltonian isotopies
keeping $X_\infty$ fixed. Exact symplectic automorphisms form a group
$\mathit{Aut}(M)$, which acts naturally on the set of exact Lagrangian
submanifolds. The analogue of Lemma \ref{th:check-exactness} is

\begin{lemma} \label{th:check-exactness-2}
$\phi$ is exact iff the following conditions hold. First, $\mathit{Hom}(o_X,\phi^*o_X)$ is the trivial flat bundle. Secondly, if we take a nonzero covariantly constant section of that bundle and restrict it to $M$, the result lies in the same homotopy class of nowhere
zero sections as $\phi^*\sigma_{X,\infty}/\sigma_{X,\infty}$. \qed
\end{lemma}

\begin{remark} \label{th:h1}
The conditions for exactness which appear in Lemma \ref{th:check-exactness-2} are
automatically satisfied if $H_1(X) = 0$. This is obvious for the
first one, and for the second one the argument goes as follows: take the
meridian $\gamma$ in $M$ around a component of $X_\infty$, and bound
it by a small disc going through $X_\infty$ once. The homotopy class
of nowhere zero sections of $\mathit{Hom}(o_X,\phi^*o_X)|\gamma$ containing the covariantly constant ones is characterized by the fact that they extend to nowhere zero sections
over the disc. On the other hand, $\sigma_{X,\infty}|\gamma$ extends
to a section with one simple zero on the disc, and so does
$\phi^*\sigma_{X,\infty}|\gamma$. Hence,
$(\phi^*\sigma_{X,\infty}/\sigma_{X,\infty})|\gamma$ lies in the same
homotopy class of nowhere zero sections as the covariantly constant
ones. Finally, note that the meridians generate $H_1(M)$.
\end{remark}

%
The relation between the affine and projective notions is clear:
Lemma \ref{th:check-exactness} shows that an exact $L \subset M$ is a
rational Lagrangian submanifold of $X$, and Lemma
\ref{th:check-exactness-2} does the same for automorphisms. Given a
pair of exact Lagrangian submanifolds $L_0,L_1$, the multivalued
action \eqref{eq:modq-action} is a reduction of the single-valued one
\eqref{eq:action}, and \eqref{eq:energy} is a refinement of
\eqref{eq:modq-energy}.

\subsection{}
Let $M$ be an affine K{\"a}hler manifold, coming from $X$ and
$\sigma_{X,\infty}$ as always. Suppose that we have a second
holomorphic section $\sigma_{X,0}$ of $o_X$, linearly independent from
$\sigma_{X,\infty}$. Let $\{X_z\}_{z \in \CP{1}}$ be the pencil of
hypersurfaces generated by these two sections. This means that we set
\begin{equation}
\sigma_{X,z} = \sigma_{X,0} - z\sigma_{X,\infty} \quad \text{ and }
\quad X_z = \sigma_{X,z}^{-1}(0)
\end{equation}
for $z \in \C$. Assume that $X_0 = \sigma_{X,0}^{-1}(0)$ is smooth in
a neighbourhood of the base locus $X_{0,\infty} = X_0 \cap X_\infty$,
and that it intersects each stratum of $X_\infty$ transversally; then
the same will hold for each $X_z$, $z \in \C$. Every nonsingular
$X_z$, $z \in \C$, is a projective K{\"a}hler manifold with ample
line bundle $o_{X_z} = o_X|X_z$, and it carries a preferred
holomorphic section $\sigma_{X_z} = \sigma_{X,\infty}|X_z$ whose zero
set is $X_{0,\infty}$. By assumption, this is a divisor with normal
crossings, so that $M_z = X_z \setminus X_{0,\infty}$ is an affine
K{\"a}hler manifold. One can also see the $M_z$ as the regular fibres
of the holomorphic function
\begin{equation}
\pi_M = \sigma_{X,0}/\sigma_{X,\infty} :  M \longrightarrow \C.
\end{equation}
We denote by $\mathit{Critv}(\pi_M)$ the critical value set, which is finite
by Bertini. Given $z,z' \in \C \setminus \mathit{Critv}(\pi_M)$, we denote by
$\mathit{Iso}(M_z,M_{z'})$ the space of symplectic isomorphisms $\phi: M_z
\rightarrow M_{z'}$ which extend to $X_z \rightarrow X_{z'}$ in such
a way that $X_0 \cap X_\infty$ remains pointwise fixed, and such that
$\phi^*\theta_{M_{z'}} - \theta_{M_z}$ is exact. These spaces form a
groupoid under composition, and for $z = z'$ they reduce
to the groups $\mathit{Aut}(M_z)$ defined before. Take a smooth path $c: [0,l]
\rightarrow \C \setminus \mathit{Critv}(\pi_M)$. The $c'(t)$ can be lifted in
a unique way to vector fields $Z_t \in \smooth(TM|M_{c(t)})$ which
are horizontal (orthogonal to $\mathit{ker}(D\pi_M)$ with respect to the
symplectic form). By integrating these one gets an exact {\em
symplectic parallel transport} map
\begin{equation}
h_c \in \mathit{Iso}(M_{c(0)},M_{c(l)}).
\end{equation}
To show that $h_c$ is well-defined, one embeds $M$ into the graph
$\hat{X} = \{(z,x) \in \C \times X \suchthat x \in X_z\}$. The
projection $\pi_{\hat{X}}: \hat{X} \rightarrow \C$ is a properification (fibrewise compactification)
of $\pi_M$. It is a standard fact that parallel transport with
respect to the product symplectic form $\o_{\hat{X}} = (\o_X +
\o_{\C})|\hat{X}$ gives a symplectic isomorphism $X_{c(0)} \rightarrow
X_{c(l)}$. This extends the original $h_c$, because adding the
pullback of $\o_{\C}$ to $\o_M$ does not change the vector fields
$Z_t$. Moreover, at a point $(z,x) \in \C \times X_{0,\infty}$ we
have
\begin{equation}
 T\hat{X}_{(z,x)} = \C \times \{ \xi \in TX_x \suchthat
 D\sigma_{X,0}(x)\xi - z\, D\sigma_{X,\infty}(x)\xi = 0 \},
\end{equation}
and since this splitting is $\o_{\hat{X}}$-orthogonal, the parallel
transport vector field will lie in $\C \times \{0\}$. This implies
that the extension of $h_c$ to $X_{c(0)}$ leaves $X_{0,\infty}$
pointwise fixed. Finally, the exactness of $h_c^*\theta_{M_{c(l)}} -
\theta_{M_{c(0)}}$ follows from the Cartan formula
$
 (L_{Z_t} \theta_M) \,|\, M_{c(t)} =
 d(i_{Z_t}\theta_M) \,|\, M_{c(t)} +
 i_{Z_t}\o_M \,|\, M_{c(t)}
$, where the last term vanishes because of the horizontality of
$Z_t$.

\subsection{}
We say that $\{X_z\}$ is a {\em quasi-Lefschetz pencil} if $\pi_M$
is nondegenerate in the sense of Picard-Lefschetz theory, meaning
that it has only nondegenerate critical points, no two of which lie
in the same fibre. This is the same as saying that the set of
critical points $\mathit{Crit}(\pi_M)$ is regular ($D\pi_M$ is transverse to
the zero-section) and $\pi_M: \mathit{Crit}(\pi_M) \rightarrow \mathit{Critv}(\pi_M)$
is bijective. Lefschetz pencils are the special case where $X_\infty$
is smooth.

An {\em embedded vanishing path} for a quasi-Lefschetz pencil is an
embedded path $c: [0,1] \rightarrow \C$ with $c^{-1}(\mathit{Critv}(\pi_M)) =
\{1\}$. By considering the limiting behaviour of parallel transport
along $c|[0,t]$ as $t \rightarrow 1$, one defines an embedded
Lagrangian disc $\Delta_c \subset M$, the {\em Lefschetz thimble},
whose boundary is a Lagrangian sphere $V_c \subset M_{c(0)}$, the
associated {\em vanishing cycle}. The Picard-Lefschetz formula
\cite[Proposition 1.15]{seidel01} says that the monodromy around a
loop $\gamma$ which doubles around $c$ in positive sense is the Dehn
twist associated to the vanishing cycle, up to isotopy inside the
group of exact symplectic automorphisms:
\begin{equation}
h_\gamma \htp \tau_{V_c} \in Aut(M_{c(0)}).
\end{equation}
To be precise, $V_c$ carries a small additional piece of structure,
called a framing in \cite{seidel00,seidel01}, and the definition of
Dehn twist uses that too. However, framings contain no information if
the vanishing cycles are of dimension $\leq 3$, so omitting them will
not matter for our applications.

Now choose a base point $z_* \in \C \setminus \mathit{Critv}(\pi_M)$ and base
path $c_*: [0,\infty) \rightarrow \C \setminus \mathit{Critv}(\pi_M)$, meaning
an embedded path with $c_*(0) = z_*$ which eventually goes off to infinity in a straight half-line (there is a $\zeta \in \C^*$ such that $c_*(t) = \zeta
t$ for all $t \gg 0$). By a {\em distinguished basis} of embedded
vanishing paths, we mean an ordered collection of such paths
$(c_1,\dots,c_r)$, $r = |\mathit{Critv}(\pi_M)|$, starting at $z_*$ and with
the following properties: the tangent directions $\Rgeq c_k'(0)$ are
pairwise distinct and clockwise ordered, and different $c_k$ do not
intersect except at $z_*$. Moreover, $\Rgeq c_*'(0)$ should lie
strictly between $\Rgeq c_r'(0)$ and $\Rgeq c_1'(0)$, and no $c_k$
may intersect $c_*$ except at $z_*$. The composition of the
corresponding loops $\gamma_k$ is freely homotopic in $\C \setminus
\mathit{Critv}(\pi_M)$ to a circle $\gamma_\infty(t) = R\zeta e^{it}$ of some
large radius $R \gg 0$, and hence the product of all the associated
Dehn twists is the monodromy around that circle. More precisely,
taking into account the change of base point, one has
\begin{equation} \label{eq:large-small-loop}
 \tau_{V_{c_1}} \dots \tau_{V_{c_r}} \htp
 (h_{c_*|[0,R]})^{-1} h_{\gamma_\infty} h_{c_*|[0,R]} \in Aut(M_{z_*}).
\end{equation}
This becomes more familiar for actual Lefschetz pencils, where one
can extend the graph over $\CP{1}$, in such a way that the fibre at infinity is regular.
Proceeding as before, one can define symplectic parallel transport
for paths in $\CP{1} \setminus \mathit{Critv}(\pi_M)$. In particular, since
$\gamma_\infty$ can be contracted by passing over $\infty$, its
monodromy is isotopic to the identity, so that
\eqref{eq:large-small-loop} simplifies to the following relation
between Dehn twists in the symplectic mapping class group
$\pi_0(\mathit{Aut}(M_{z_*}))$:
\begin{equation} \label{eq:relation-in-mcg}
\tau_{V_{c_1}} \circ \dots \circ \tau_{V_{c_r}} \htp id.
\end{equation}

\begin{remark} \label{th:multiple-fibres}
It is no problem (except terminologically) to allow several
nondegenerate critical points to lie in the same fibre. One should
then choose a vanishing path for every critical point. Two paths
$c_k,c_l$ going to the same critical value should either coincide or
else satisfy $c_k(t) \neq c_l(t)$ for $t \in (0;1)$, and $\Rgeq
c_k'(0) \neq \Rgeq c_l'(0)$. In the first case, the associated
vanishing cycles in $M_{z_*}$ will be disjoint.
\end{remark}

\subsection{}
Let $\{X_z\}$ be a quasi-Lefschetz pencil. Consider a smooth embedded
path $d: [-1,1] \rightarrow \C$ such that $d^{-1}(\mathit{Critv}(\pi_M)) =
\{-1,1\}$. Split it into a pair of vanishing paths $c_\pm: [0,1] \rightarrow \C$, given by $c_-(t) = d(-t)$,
$c_+(t) = d(t)$. Consider the associated vanishing
cycles
\begin{equation}
V_{c_-},V_{c_+} \subset M_{d(0)}.
\end{equation}
We say that $d$ is a {\em matching path} if $V_{c_-}$, $V_{c_+}$ are
isotopic as exact Lagrangian spheres (actually, there is an
additional condition concerning the framings, which we do not
spell out here since it is vacuous in the dimensions that we will be
interested in). If $d$ is a matching path, one can glue together the
Lefschetz thimbles $\Delta_{c_-}$, $\Delta_{c_+}$ to obtain a
Lagrangian sphere in the total space, the {\em matching cycle}
$\Sigma_d \subset M$. In the naive situation where the two vanishing
cycles coincide, this would simply be the set-theoretic union of the
Lefschetz thimbles. In general, the construction involves a choice of
isotopy between the cycles, and $\Sigma_d$ might theoretically depend
on that choice. Note however that this can never happen if the fibres
are Riemann surfaces with boundary, since there, the space of exact
Lagrangian submanifolds (circles) having a fixed isotopy class is
contractible.

\subsection{\label{subsec:braid-monodromy}}
Let $\sigma_{X,0},\sigma_{X,\infty}$ be two sections of $o_X$ which
generate a quasi-Lefschetz pencil $\{X_z\}$. Suppose that we have yet
another section $\sigma_{X,0}'$, linearly independent from the
previous ones. We require that its zero set should be smooth in a
neighbourhood of the base locus $X_{0,\infty}$, and should intersect each
stratum of $X_{0,\infty}$ transversally. Consider the associated map
\begin{equation} \label{eq:bm-map}
M \xrightarrow{b_M = \left(\frac{\sigma_{X,0}}{\sigma_{X,\infty}},
\frac{\sigma_{X,0}'}{\sigma_{X,\infty}}\right)} \C^2,
\end{equation}
its critical point set $\mathit{Crit}(b_M)$ and critical value set $C =
\mathit{Critv}(b_M)$. Concerning this, we make the following additional
assumptions:
\begin{itemize}
\item
$\mathit{Crit}(b_M)$ is regular, and $b_M: \mathit{Crit}(b_M) \rightarrow C$ is an
embedding away from finitely many points.
\end{itemize}
Regularity means that $Db_M \in \Gamma(M,\mathit{Hom}_\C(TM,\C^2))$ is never zero,
and is transverse to the subset of rank one linear maps. The second
statement implies that for generic $z \in \C \setminus \mathit{Critv}(\pi_M)$,
the curve $\mathit{Crit}(b_M)$ is transverse to $b_M^{-1}(\{z\} \times \C)$,
and the projection $b_M|\mathit{Crit}(b_M) \times b_M^{-1}(\{z\} \times \C)
\rightarrow C \times (\{z\} \times \C)$ is a bijection. In this case,
the sections $(\sigma_{X,0}'|X_z,\sigma_{X,\infty}|X_z)$ generate a
quasi-Lefschetz pencil on $X_z$, and the critical value set of the
associated holomorphic function $q_{M_z}: M_z \rightarrow \C$ is
precisely
\begin{equation} \label{eq:z-critical}
\mathit{Critv}(q_{M_z}) = C \cap (\{z\} \times \C).
\end{equation}
There is a finite subset of $\C \setminus \mathit{Critv}(\pi_M)$ where this
fails, and we denote it by $\mathit{Fakev}(\pi_M)$. These ``fake critical
values'' will play a role similar to the real ones, see
\eqref{eq:braid-monodromy} below.
\begin{itemize}
\item
Let $x$ be a critical point of $\pi_M$, $\pi_M(x) = z$. Then the
restriction of $D^2\pi_M$ to the complex codimension one subspace
$ker(D_xb_M)$ is nondegenerate. Moreover, all critical points of
$q_{M_z}\,|\,(M_z \setminus \{x\})$ are nondegenerate, and their
$q_{M_z}$-values are pairwise distinct and different from
$q_{M_z}(x)$.
\end{itemize}
As a consequence, $b_M: \mathit{Crit}(b_M) \rightarrow C$ is an embedding on
$\mathit{Crit}(b_M) \cap M_z$, and moreover the projection $C \subset \C^2
\rightarrow \C$ to the first variable has an ordinary (double) branch
point at $b_M(x)$, and is a local isomorphism at all other points of
$C \cap (\{z\} \times \C)$.

If the two conditions stated above are satisfied, we say that
$\sigma_{X,0}'$ is a {\em generic auxiliary section} of $o_X$.
Practically, the situation is that we have a quasi-Lefschetz pencil
on $X$, almost every regular fibre $X_z$ of which again admits a
quasi-Lefschetz pencil. Moreover, we have a good understanding of how
these pencils degenerate as $X_z$ becomes singular. We will now
explain the consequences of this for the vanishing cycles of $\pi_M$.
For that, suppose that our base point $z_* \notin \mathit{Fakev}(\pi_M)$.
Alongside the ordinary monodromy of the pencil,
\begin{align}
 \notag \pi_1(\C \setminus \mathit{Critv}(\pi_M),z_*) & \longrightarrow
 \pi_0(\mathit{Aut}(M_{z_*})), \\
\intertext{there is a ``relative'' or ``braid'' monodromy
homomorphism}
 \label{eq:braid-monodromy}
 \pi_1(\C \setminus (\mathit{Critv}(\pi_M) \cup \mathit{Fakev}(\pi_M)),z_*)
 & \longrightarrow
 \pi_0(\Diff^c(\C,\mathit{Critv}(q_{M_{z_*}}))),
\end{align}
which describes how \eqref{eq:z-critical} moves in
$\C$ as $z$ changes. Here, $\Diff^c(\C,\mathit{Critv}(q_{M_{z_*}}))$ is the group of compactly supported diffeomorphisms preserving the finite set $\mathit{Critv}(q_{M_{z_*}})$, hence its $\pi_0$ is isomorphic to the appropriate braid group. Take an embedded vanishing path for $\pi_M$, $c:
[0,1] \rightarrow \C$, with the additional property that $c([0,1])
\cap Fakev(\pi_M) = \emptyset$. Let $\gamma$ be a loop doubling
around $c$ in positive sense. As we know, the monodromy around
$\gamma$ is the Dehn twist along $V_c$. The corresponding statement
for \eqref{eq:braid-monodromy} is that there is an embedded path $d:
[-1,1] \rightarrow \C$, with $d^{-1}(Critv(q_{M_{z_*}})) = \{\pm
1\}$, such that the image of $\gamma$ under the relative monodromy is
the half-twist along $d$. The relation between the two is as follows,
see \cite[Section 16]{seidel04} for the proof:

\begin{prop} \label{th:vanishing-matching}
$d$ is a matching path for $q_{M_{z_*}}$, and the vanishing cycle
$V_c \subset M_{z_*}$ is, up to Lagrangian isotopy, a matching cycle
$\Sigma_d$ for that path. \qed
\end{prop}

Figure \ref{fig:braid-monodromy} summarizes the situation schematically.
\begin{figure}
\begin{centering}
\begin{picture}(0,0)%
\includegraphics{braidmonodromy.pstex}%
\end{picture}%
\setlength{\unitlength}{3355sp}%
\begingroup\makeatletter\ifx\SetFigFont\undefined%
\gdef\SetFigFont#1#2#3#4#5{%
  \reset@font\fontsize{#1}{#2pt}%
  \fontfamily{#3}\fontseries{#4}\fontshape{#5}%
  \selectfont}%
\fi\endgroup%
\begin{picture}(5466,4294)(436,-3809)
\put(451,-2836){\makebox(0,0)[lb]{\smash{{\SetFigFont{10}{12.0}{\rmdefault}{\mddefault}{\updefault}{\color[rgb]{0,0,0}point of $\pi_M$}%
}}}}
\put(3151,-3586){\makebox(0,0)[lb]{\smash{{\SetFigFont{10}{12.0}{\rmdefault}{\mddefault}{\updefault}{\color[rgb]{0,0,0}``downstairs'' branch curve $C$}%
}}}}
\put(3301,314){\makebox(0,0)[lb]{\smash{{\SetFigFont{10}{12.0}{\rmdefault}{\mddefault}{\updefault}{\color[rgb]{0,0,0}$M_z$}%
}}}}
\put(4276,-811){\makebox(0,0)[lb]{\smash{{\SetFigFont{10}{12.0}{\rmdefault}{\mddefault}{\updefault}{\color[rgb]{0,0,0}$M$}%
}}}}
\put(1951,-3736){\makebox(0,0)[lb]{\smash{{\SetFigFont{10}{12.0}{\rmdefault}{\mddefault}{\updefault}{\color[rgb]{0,0,0}$\{z\} \times \C$}%
}}}}
\put(4201,-3091){\makebox(0,0)[lb]{\smash{{\SetFigFont{10}{12.0}{\rmdefault}{\mddefault}{\updefault}{\color[rgb]{0,0,0}$\C^2$}%
}}}}
\put(2660,-2588){\makebox(0,0)[lb]{\smash{{\SetFigFont{10}{12.0}{\rmdefault}{\mddefault}{\updefault}{\color[rgb]{0,0,0}$b_M$}%
}}}}
\put(3826,-2333){\makebox(0,0)[lb]{\smash{{\SetFigFont{10}{12.0}{\rmdefault}{\mddefault}{\updefault}{\color[rgb]{0,0,0}which is also the matching}%
}}}}
\put(3826,-2558){\makebox(0,0)[lb]{\smash{{\SetFigFont{10}{12.0}{\rmdefault}{\mddefault}{\updefault}{\color[rgb]{0,0,0}cycle of a path in $\{z\} \times \C$}%
}}}}
\put(3826,-2108){\makebox(0,0)[lb]{\smash{{\SetFigFont{10}{12.0}{\rmdefault}{\mddefault}{\updefault}{\color[rgb]{0,0,0}vanishing cycle in $M_z$,}%
}}}}
\put(451,-2611){\makebox(0,0)[lb]{\smash{{\SetFigFont{10}{12.0}{\rmdefault}{\mddefault}{\updefault}{\color[rgb]{0,0,0}image of a critical}%
}}}}
\end{picture}%
\caption{\label{fig:braid-monodromy}}
\end{centering}
\end{figure}

\begin{example}
Suppose that $\mathrm{dim}_\C(M) = 2$, so that the $M_z$ are
affine algebraic curves. Suppose first that $z \notin \mathit{Critv}(\pi_M) 
\cup \mathit{Fakev}(\pi_M)$. Then $M_z$ is smooth, and the map
\begin{equation} \label{eq:generic-cover}
\sigma_{X,0}'/\sigma_{X,\infty}: M_z \longrightarrow \C 
\end{equation}
is a generic branched covering (which means that the monodromy around
each branch point in $\C$ is a transposition). As $z$ approaches a point 
of $\mathit{Critv}(\pi_M)$, $M_z$ degenerates to a curve with an ordinary (nodal) 
singularity, and the effect on \eqref{eq:generic-cover} is that two branch 
points in $\C$ come together. This degeneration gives rise to a vanishing cycle (a
simple closed curve) in each nearby smooth fibre $M_z$, which under 
\eqref{eq:generic-cover} maps to an path in $\C$ connecting the two collapsing branch
points. On the other hand, for a point $z \in \mathit{Fakev}(\pi_M)$, the fibre $M_z$ remains
smooth, but the branched covering degenerates to a non-generic one.
\end{example}

\subsection{}
We temporarily digress from K{\"a}hler to general symplectic
geometry. Let $(M^{2n},\o_M)$ be any connected symplectic manifold
with $c_1(M) = 0$. Take an almost complex structure $J_M$ which tames
$\o_M$, and a complex volume form (a nowhere vanishing $(n,0)$-form
with respect to $J_M$) $\eta_M$. Let $\LL_M \rightarrow M$ be the
bundle of (unoriented) Lagrangian Grassmannians. The phase (also sometimes called squared phase) function
associated to $\eta_M$,
\begin{equation}
\alpha_M: \LL_M \longrightarrow S^1,
\end{equation}
is defined as follows. Given $x \in M$ and $\Lambda \in \LL_{M,x}$,
take any (not necessarily orthonormal) basis $e_1,\dots,e_n$ of
$\Lambda$, and set
\begin{equation} \label{eq:squared-phase}
\alpha_M(\Lambda) = \frac{\eta_M(e_1 \wedge \dots \wedge
e_n)^2}{|\eta_M(e_1 \wedge \dots \wedge e_n)|^2}.
\end{equation}
We distinguish the special case where $J_M$ is $\o_M$-compatible by
calling $\alpha_M$ a {\em classical phase function}. These are the
most commonly used ones, and they are better behaved. General
phase functions will be needed later for technical reasons.

Given $\alpha_M$, one can associate to a Lagrangian submanifold $L
\subset M$ or a symplectic automorphism $\phi: M \rightarrow M$ a
phase function
\begin{equation}
\begin{cases}
 & \!\!\! \alpha_L: L \longrightarrow S^1, \quad
 \alpha_L(x) = \alpha_M(TL_x), \\
 & \!\!\! \alpha_\phi: \LL_M \longrightarrow S^1, \quad
 \alpha_\phi(\Lambda) = \alpha_M(D\phi(\Lambda))/\alpha_M(\Lambda).
\end{cases}
\end{equation}
A grading of $L$ or $\phi$ is a real-valued lift $\alphagr_L$ resp.\
$\alphagr_\phi$, where the convention is that $\R$ covers $S^1$ by $a
\mapsto \exp(2\pi i a)$. Pairs $\Lgr = (L,\alphagr_L)$, $\phigr =
(\phi,\alphagr_\phi)$ are called graded Lagrangian submanifolds
resp.\ graded symplectic automorphisms. The latter form a group, with
a natural action on the set of graded Lagrangian submanifolds. The
composition law and action are defined by
\begin{equation}
\label{eq:compose-gradings}
\begin{cases}
 & \!\!\! \alphagr_{\psi \circ \phi} = \alphagr_\psi \circ D\phi +
 \alphagr_\phi, \\
 & \!\!\! \alphagr_{\phi(L)} = \alphagr_L \circ \phi^{-1} + \alphagr_\phi
 \circ D\phi^{-1}.
\end{cases}
\end{equation}
In particular, for $k \in \Z$ we have the $k$-fold shift $[k] = (\phi
= \id_M,\alphagr_{\phi} = -k)$, which reduces the grading of each
Lagrangian submanifold by $k$. All the ``graded'' notions actually
depend only on the homotopy class of $\eta_M^2$ as a trivialization
of $K_M^2$, or equivalently on the cohomology class $\mu_M =
[\alpha_M] \in H^1(\LL_M)$, called a global Maslov class in
\cite{seidel99}. Namely, suppose that we have two cohomologous phase
functions, which can be written as $\alpha_M' = \alpha_M e^{2\pi i
\chi}$ for some $\chi: \LL_M \rightarrow \R$. Gradings with respect
to the two are related in an obvious way,
\begin{equation} \label{eq:new-phase}
\begin{cases}
 \!\!\! & \alphagr_L'(x) = \alphagr_L(x) + \chi(TL_x), \\
 \!\!\! & \alphagr_\phi'(\Lambda) = \alphagr_\phi(\Lambda) +
 \chi(D\phi(\Lambda)) - \chi(\Lambda).
\end{cases}
\end{equation}
Note that $\chi$ is unique only up to an integer constant. This is
irrelevant for $\alphagr_\phi'$; it affects $\alphagr_L'$ but not in
a really important way, since changing $\chi$ shifts the grading of
all Lagrangians simultaneously.

\begin{remark} \label{th:or}
There is a canonical square root $\alpha_M^{1/2}: \LL^{\orient}
\rightarrow S^1$ of $\alpha_M$ on the bundle $\LL^{\orient}$ of
oriented Lagrangian Grassmannians, defined by removing the squares
from \eqref{eq:squared-phase} and using positively oriented bases
$e_1,\dots,e_n$. As a consequence, any graded Lagrangian submanifold carries a
preferred orientation, characterized by
\begin{equation}
\exp(\pi i \alphagr_L(x)) = \alpha_M^{1/2}(TL_x).
\end{equation}
\end{remark}

\subsection{}
There is an obvious notion of graded symplectic fibre bundle, which
is a locally trivial symplectic fibration $p: E \rightarrow B$ with
an almost complex structure $J_{E_b}$ and complex volume form
$\eta_{E_b}$ on each fibre $E_b$. The monodromy maps of such a
fibration are naturally graded: given a loop $\gamma: [0,l]
\rightarrow B$ and the corresponding parallel transport maps $\phi_t
= h_{\gamma|[0,t]}: E_{\gamma(0)} \rightarrow E_{\gamma(t)}$, one
takes the unique smooth family of functions $\tilde{a}_t$ on the
Lagrangian Grassmannian bundle of $E_{\gamma(0)}$ with $\tilde{a}_0 =
0$ and
\begin{equation}
\exp(2\pi i \tilde{a}_t) = (\alpha_{E_{\gamma(t)}} \circ
D\phi_t)/\alpha_{E_{\gamma(0)}}.
\end{equation}
The graded symplectic monodromy is $\hgr_{\gamma} =
(\phi_l,\alphagr_{\phi_l} = \tilde{a}_l)$. An example which is
relevant for our purpose is where we have a trivial fibre bundle $E =
\C^* \times M$ over $B = \C^*$ with a varying complex volume form
$\eta_{E_z} = z^\mu \eta_M$. The ``invisible singularity'' of this
family at $z = 0$ leads to a monodromy which, for a circle $\gamma$
turning around the origin in {\em clockwise} direction, is the shift
\begin{equation} \label{eq:shift-monodromy}
\hgr_\gamma = [2\mu].
\end{equation}

\subsection{}
Gradings make it possible to assign absolute Maslov indices to
Lagrangian intersection points. This depends on the following notion
from symplectic linear algebra: a smooth path of Lagrangian subspaces
$\Lambda_t$ ($0 \leq t \leq 1$) in a fixed symplectic vector space is
{\em crossingless} if it satisfies $\Lambda_0 \cap \Lambda_t =
\Lambda_0 \cap \Lambda_1$ for all $t>0$, and the crossing form
\cite{robbin-salamon93} on $\Lambda_0/\Lambda_0 \cap \Lambda_1$
associated to $(d\Lambda_t/dt)_{t = 0}$ is negative definite.
Crossingless paths with given endpoints exist and are unique up to
homotopy, see again \cite{robbin-salamon93}.

\begin{lemma} \label{th:angle-bound}
Suppose that $\alpha_M$ is classical. Let $\Lambda_t$ be a
crossingless path in $TM_x$ for some $x$. Choose $\tilde{a}: [0,1]
\rightarrow \R$ such that $\exp(2\pi i \tilde{a}_t) =
\alpha_M(\Lambda_t)$. Then
\begin{equation}
\tilde{a}_1 - \tilde{a}_0 \in (-n,0].
\end{equation}
\end{lemma}

\proof We can choose an isomorphism $\Phi: TM_x \rightarrow \C^n$ in such a way
that the complex structure and symplectic form become standard,
$\Phi(\Lambda_0) = \R^n$ and $\Phi(\Lambda_1) = e^{i\pi c_1} \R
\times \dots \times e^{i\pi c_n}\R$ with $c_k \in (-1,0]$. For the
obvious crossingless path $\Lambda_t = e^{i\pi t c_1}\R \times \dots
\times e^{i\pi t c_n}\R$ we have $\tilde{a}_1 - \tilde{a}_0 =
c_1+\dots+c_n \in (-n,0]$. Because of homotopy uniqueness, the same
holds for any other crossingless path. \qed

Take two graded Lagrangian submanifolds $\Lgr_0,\Lgr_1$, and a point
$x \in L_0 \cap L_1$. Choose a crossingless path from $\Lambda_0 =
TL_0$ to $\Lambda_1 = TL_1$, and take $\tilde{a}_t \in \R$ such that
$\exp(2\pi i \tilde{a}_t) = \alpha_M(\Lambda_t)$. The absolute Maslov
index is defined as
\begin{equation} \label{eq:absolute-index}
I(x) = (\alphagr_{L_1}(x) - \tilde{a}_1) - (\alphagr_{L_0}(x) -
\tilde{a}_0) \in \Z.
\end{equation}
Lemma \ref{th:angle-bound} implies that this is approximately the
difference of the phases:

\begin{lemma} \label{th:approximate-index}
If $\alpha_M$ is classical, $I(x) - \alphagr_{L_1}(x) +
\alphagr_{L_0}(x) \in [0,n)$. \qed
\end{lemma}

This was pointed out to the author by Joyce, but it is known to many other
people in mathematics and physics, see for instance the proof of
\cite[Theorem 4.3]{thomas-yau01}. In fact, ideas of this kind can be
traced back at least to \cite{salamon-zehnder92}. Another easy
consequence of Lemma \ref{th:angle-bound}, this time for symplectic
automorphisms, will be useful later:

\begin{lemma} \label{th:graded-variation}
Suppose that $\alpha_M$ is classical. Let $\phigr =
(\phi,\alphagr_\phi)$ be a graded symplectic automorphism of $M$. For
any two Lagrangian subspaces $\Lambda_0,\Lambda_1 \in \LL_{M,x}$ at
the same point,
\begin{equation}
|\alphagr_\phi(\Lambda_1)-\alphagr_\phi(\Lambda_0)| < n.
\end{equation}
\end{lemma}

\proof Choose a crossingless path $\Lambda_t$ joining our two subspaces. Lemma
\ref{th:angle-bound} shows that as we go along this path,
$\alpha_M(\Lambda_t)$ changes by an angle in $(-2\pi n,0]$, and so
does $\alpha_M(D\phi(\Lambda_t))$. The change in
$\alpha_M(D\phi(\Lambda_t))/\alpha_M(\Lambda_t)$ is the difference
between the two, hence lies in $(-2\pi n;2\pi n)$. \qed

\subsection{}
We now return to the case of an affine K{\"a}hler manifold $M = X
\setminus X_\infty$. We say that $M$ is {\em affine Calabi-Yau} if $X$
comes with a preferred isomorphism of holomorphic line bundles,
\begin{equation} \label{eq:beta}
\beta_X : \K_X \stackrel{\iso}{\longrightarrow} o_X^{\otimes -m_X}
\end{equation}
for some $m_X \in \Z$ (called the monotonicity index). On $M$ we have
a preferred trivialization of $o_X$ given by $\sigma_{X,\infty}|M$, hence
through $\beta_X$ a holomorphic volume form $\eta_M$. In the special
case where $m_X = 0$, $X$ itself is Calabi-Yau, and $\eta_M$ extends
to a holomorphic volume form $\eta_X$ on it. We denote by $\Autgr(M)$
the group of exact symplectic automorphisms of $M$ equipped with a
grading. Similarly, if $m_X = 0$, we write $\Autgr(X)$ for the graded
rational automorphisms of $X$.

\begin{remark} \label{th:extend-grading}
In the $m_X = 0$ case, there is a canonical embedding
\begin{equation}
\Autgr(M) \longrightarrow \Autgr(X).
\end{equation}
We already know that every $\phi \in Aut(M)$ extends to a rational
symplectic automorphism of $X$. There are no obstructions to
extending gradings over subsets of codimension $\geq 2$, so the
grading of $\phi$ extends uniquely from $M$ to $X$.
\end{remark}

If $\{X_z \subset X\}$ is a quasi-Lefschetz pencil, we have for $z
\in \C \setminus Critv(\pi_M)$ a canonical isomorphism $\beta_{X_z}:
\K_{X_z} \iso (\K_X \otimes o_X)|X_z \iso o_X^{\otimes -m_X+1}$, so
the $M_z$ are again affine Calabi-Yaus, with monotonicity index
$m_{X_z} = m_X-1$. More explicitly, the complex volume form is
obtained as a quotient $\eta_M/dz$, which means that it is the unique
solution of
\begin{equation} \label{eq:division}
dz \wedge \eta_{M_z} = \eta_M.
\end{equation}
The first consequence of this is that the monodromy maps along loops
$\gamma$ in $\C \setminus \mathit{Critv}(\pi_M)$ have canonical gradings,
$\hgr_\gamma \in \Autgr(M_{\gamma(0)})$. If $c$ is a vanishing path,
the associated vanishing cycle $V_c$ admits a grading; this is
trivial if its dimension is $\geq 2$, and otherwise follows from the
fact that it is bounded by a Lagrangian disc in $M$. Hence, the
associated Dehn twist $\tau_{V_c}$ has a canonical grading
$\taugr_{V_c}$, which is zero outside a neighbourhood of $V_c$
itself. Taking the loop $\gamma$ which doubles $c$, we have the
following graded version of the Picard-Lefschetz theorem:
\begin{equation}
\hgr_\gamma \htp \taugr_{V_c} \in \Autgr(M_{c(0)}).
\end{equation}
The analogue of \eqref{eq:relation-in-mcg} is

\begin{lemma} \label{th:shift-factor}
Suppose that $\{X_z\}$ is a Lefschetz pencil, and $c_1,\dots,c_r$ a
distinguished basis of vanishing paths. Then there is a graded isotopy
\begin{equation}
\taugr_{V_{c_1}} \circ \dots \circ \taugr_{V_{c_r}} \htp [4-2 m_X]
\in \Autgr(M_{z_*}).
\end{equation}
\end{lemma}

\proof What we need to show is that the graded monodromy around a
large circle $\gamma_\infty$ is $[4-2 m_X]$. For simplicity, change
variables to $\zeta = z^{-1}$, so that $\gamma_\infty$ becomes a
small circle going clockwise around $\zeta = 0$. Recall that
$\eta_{M_z}$ is obtained by dividing $\eta_M$, which has a zero or
pole of order $-m_X$ along $X_\infty$, by $dz = d\zeta/\zeta^2$,
which has order $-2$. Hence, after extending the graph of our Lefschetz pencil over infinity, one has a locally trivial symplectic fibration near $\zeta = 0$, with a family of complex
volume forms on the fibres that grows like $\zeta^{2-m_X}$. Now we
are in the same situation as in \eqref{eq:shift-monodromy}. \qed

\section{Monodromy and negativity\label{sec:negativity}}

Continuing the discussion above, suppose that $m_X = 1$, so that the fibres $X_z$ are Calabi-Yau manifolds. In the Lefschetz pencil case, Lemma \ref{th:shift-factor} tells us that the ``large complex structure limit monodromy'' is a downward shift. We will need a generalization of this statement to the case when $X_\infty$ has normal crossings. This has a single purpose, namely to provide Proposition \ref{th:negativity}, which enters into the proof of Corollary \ref{th:generates-everything}; readers willing to take that Proposition for granted may want to skip most of this section.

\subsection{}
Before entering into the details, some motivation may be appropriate. The intuition underlying our approach comes from both the homological and geometric (SYZ) aspects of mirror symmetry. In the context of homological mirror symmetry, we have a conjectured general formula for the autoequivalence which is mirror to the large complex structure limit monodromy, which goes back to Kontsevich \cite[Equation (5)]{aspinwall01}. For simplicity, let's consider $Y_q^* \subset P = \mathbb{P}(\Lambda_\Q^4)$, the quartic surface over $\Lambda_\Q$ defined by \eqref{eq:fermat}, and work with $\Gamma_{16}$-equivariant coherent sheaves on it (this is derived equivalent to working with coherent sheaves on $Z_q^*$, see Lemma \ref{th:kapranov-vasserot00} below for further discussion). Then, the mirror autoequivalence is given by tensoring with a line bundle followed by a shift:
\begin{equation} \label{eq:phi-monodromy}
\Phi(E) = E \otimes {\OO}_P(-4)[2].
\end{equation}
In particular, for any pair of objects $(E_0,E_1)$ and any integer $k$, there is a $d > 0$ such that all morphisms from $E_0$ to $\Phi^d(E_1)$ are of degree $<k$, because the effect of the shift will predominate as we iterate.

To explain the geometric side, let's look at a toy example one dimension down, namely the pencil of elliptic curves
\begin{equation}
X_z = \{ z \cdot x_0x_1x_2 + \textstyle\frac{1}{3}(x_0^3 + x_1^3 + x_2^3) = 0\} \subset \C P^2.
\end{equation}
$X_\infty = \{x_0x_1x_2 = 0\}$ is a closed chain of three rational curves. Hence, the monodromy around it is the product of three inverse Dehn twists along parallel curves on a torus, combined with  a shift $[2]$ (the shift comes from the same consideration as in Lemma \ref{th:shift-factor}). If we identify the generic fibre $X_z$ with the standard symplectic torus $T^2 = \R^2/\Z^2$ in an appropriate way, the monodromy is isotopic to the linear map
\begin{equation} \label{eq:matrix}
\begin{pmatrix} 1 & 0 \\ 3 & 1 \end{pmatrix} \in \mathit{SL}_2(\Z).
\end{equation}
Now look at graded Lagrangian submanifolds on our torus. The map \eqref{eq:matrix} (with its choice of grading coming from writing it as a product of three inverse Dehn twists) preserves the grading of the circles $\{p = \mathit{const}\}$, which are the fibres of the SYZ fibration. Generally, it raises gradings by $<1$, and the same uniform bound holds for all its iterates. In contrast, the shift adds $-2$ to all the gradings, and this number grows linearly as we iterate. This ``negativity'' property has implications for Maslov indices of Lagrangian intersections, which mirror the property of \eqref{eq:phi-monodromy} mentioned previously.

In our application, there is no overall description of the monodromy as simple as \eqref{eq:matrix} (Mark Gross has pointed out that one can potentially obtain such a description by hyperkaehler rotation, under which the SYZ fibration becomes an elliptic fibration, for which holomorphic automorphisms can be constructed easily; however, it has not been proved that this description is indeed correct). Instead, we will look at the behaviour in various local models around points of $X_\infty$. The local behaviour near smooth points of $X_\infty$ is the same as in the Lefschetz pencil case. Near singular points of $X_\infty$ which do not lie on $X_0$ (hence are not contained in the base locus of our pencil), the local picture is the same as that for a degeneration to a normal crossing: its structure is classical and not too different from that of Dehn twists in one dimension (a simplified local model is described by Assumptions \ref{as:1}). The most complicated local behaviour (Assumptions \ref{as:2}) occurs near the singular points of $X_\infty$ which also lie on $X_0$. In the elliptic curve case there were no such points, for dimension reason; in the $K3$ surface case there are finitely many, and that still makes it relatively easy to deal with them, since we can adjust the K{\"a}hler form to be standard locally near those points.

\subsection{}
We start with some generalities. Let $M^{2n}$ be a compact connected
symplectic manifold (possibly with boundary), with a tame almost
complex structure $J_M$, complex volume form $\eta_M$, and associated
phase function $\alpha_M$.

\begin{defn} \label{th:negativity-definition}
A graded symplectic automorphism $\phigr$ of $M$ is {\em negative} if
there is a $d_0>0$ such that
\begin{equation}
\alphagr_{\phi^{d_0}}(\Lambda) < 0 \quad \text{ for all $\Lambda \in
\LL_M$.}
\end{equation}
\end{defn}

Here, the grading of $\phi^{d_0}$ is defined by the composition rule
\eqref{eq:compose-gradings}. The same rule shows that if negativity
holds, the grading of the iterates goes to $-\infty$. Namely, the
grading of $\phi^{d_0}$ is $\leq -\epsilon$ for some $\epsilon>0$,
and the grading of $\phi$ is $\leq C$ for some $C$ (which could be positive). By writing $d =
k d_0 + l$ with $0 \leq l \leq d_0-1$ one sees that
\begin{equation} \label{eq:growth}
 \alphagr_{\phi^d}(\Lambda) \leq C l - k\epsilon <
 Cd_0 - \epsilon [d/d_0],
\end{equation}
as claimed. A closely related property is that the negativity of a
given graded symplectic automorphism depends only on the global
Maslov class $\mu_M$, and not on the choice of $\alpha_M$ within that
class. For this take two cohomologous phase functions $\alpha_M' =
e^{2 \pi i \chi} \alpha_M$. Since $\chi$ is bounded, it follows from
\eqref{eq:new-phase} that the gradings of $\phi^d$ with respect to
$\alpha_M,\alpha_M'$ differ by a uniformly bounded amount, which in
view of \eqref{eq:growth} gives the desired property.

\begin{lemma} \label{th:maslov-growth}
Let $\phigr$ be a negative graded symplectic automorphism, and
$\Lgr_0,\Lgr_1$ two graded Lagrangian submanifolds. For any $k \in
\Z$ one can find a $d>0$ and a small exact perturbation $L_1'$ of
$L_1$ such that the intersection $L_0 \cap \phi^d(L_1')$ is transverse, and each intersection point $x$ has absolute Maslov index
\begin{equation}
I(x) < k.
\end{equation}
\end{lemma}

\proof Since negativity is independent of the choice of phase
function, we may assume that that function is classical. Suppose that
the grading of $L_1$ is bounded above by $C_1$, and that of $L_0$
bounded below by $C_0$. Choose some $d$ such that $C_1 - C_0 + C d_0
- \epsilon[d/d_0] \leq k-n$. From \eqref{eq:growth} and
\eqref{eq:compose-gradings} it follows that
\begin{equation} \label{eq:image-d}
 \alphagr_{\phi^d(L_1)}(x) - \alphagr_{L_0}(x)
 < k-n
\end{equation}
at every intersection point $x$. This estimate continues to hold if
one perturbs $L_1$ slightly, to make the intersections transverse.
Now apply Lemma \ref{th:approximate-index}. \qed

Negativity is well-behaved with respect to abelian symplectic
reduction (this fact is related to the quotient construction of
special Lagrangian submanifolds, see for instance \cite{gross01b}).
The situation is as follows: $M = M^{2n}$ carries a free Hamiltonian
$T^k$-action, with moment map $\mu: M \rightarrow \R^k$. We have an
almost complex structure $J_M$ which is $\o_M$-compatible and
$T^k$-invariant, and a $T^k$-invariant complex volume form $\eta_M$,
with associated (classical and $T^k$-invariant) phase function
$\alpha_M$. Denote by $K_1,\dots,K_k$ the Killing fields of the
action. Let $\phi$ be a symplectic automorphism which preserves
$\mu$, hence is $T^k$-equivariant; and $\alphagr_\phi$ a grading of
it (necessarily $T^k$-invariant). For every moment value $r$ we have
the symplectic quotient $M^{\red,r} = \mu^{-1}(r)/T^k$. This carries
the reduced symplectic form $\o_{M^{\red,r}}$ and a compatible almost
complex structure $J_{M^{\red,r}}$, the latter of which is obtained
by taking the complexified tangent spaces along the orbits, and
identifying their orthogonal complement with the tangent spaces of
$M^{\red,r}$. Moreover, $\phi$ induces symplectic automorphisms
$\phi^{\red,r}$. In our situation, we also have a reduced complex
volume form $\eta_{M^{\red,r}}$, defined by restricting
$
i_{K_1} \dots i_{K_r} \eta_M
$
to the same orthogonal complement as before. Projection identifies
Lagrangian subspaces $\Lambda \subset TM$ containing the tangent
space to the orbit with Lagrangian subspaces $\Lambda^{\red} \subset
TM^{\red,r}$. The phase functions defined by $\eta_M$,
$\eta_{M^{\red,r}}$ are related by
$\alpha_{M^{\red,r}}(\Lambda^{\red}) = \alpha_M(\Lambda)$. In
particular, we can define a grading of the reduced symplectic
automorphisms by
\begin{equation} \label{eq:lift-lambda}
\alphagr_{\phi^{\red,r}}(\Lambda^{\red}) = \alphagr_\phi(\Lambda).
\end{equation}

\begin{lemma} \label{th:negative-reduction}
Suppose that the $\phigr^{\red,r}$ are all negative. Then $\phigr$ is
negative.
\end{lemma}

\proof Choose some $r$. By assumption there is a ${d_0}$ such that
the grading of $(\phigr^{\red,r})^{d_0}$ is negative, and then the
same thing holds for neighbouring $r$. By a covering argument (remember the general assumption that $M$ is compact), one
sees that there is a $d_0$ such that the gradings of the $d_0$-th
iterates of all reduced maps are negative. Applying \eqref{eq:growth}
shows that for some large $d$, the grading of each
$(\phigr^{\red,r})^d$ will be $<-n$ everywhere. By
\eqref{eq:lift-lambda}, this means that at each point of $M$ there
are some Lagrangian tangent subspaces whose $\phigr^d$-grading is
$<-n$. From this and Lemma \ref{th:graded-variation} one gets the
desired conclusion. \qed

\subsection{}
We will now look at the local models for the ``large complex
structure limit'' monodromy, first under the rather unrealistic
assumption that both the symplectic and complex structures are
standard, which allows one to integrate the monodromy vector field
explicitly.

\begin{assumptions} \label{as:1}
Fix $n \geq 2$ and $2 \leq k \leq n$. We take $Y = \C^n$ with the
standard complex structure and symplectic form $\o_Y$, and the map
\begin{equation}
p: \C^n \longrightarrow \C, \quad y \longmapsto y_1 y_2 \dots y_k.
\end{equation}
$Y$ should carry the complex volume form (away from the zero fibre of $p$)
\begin{equation}
\frac{dy_1}{y_1} \wedge \dots \wedge \frac{dy_k}{y_k} \wedge dy_{k+1}
\wedge \dots \wedge dy_n.
\end{equation}
\end{assumptions}

The clockwise monodromy vector field $Z$ is the unique horizontal lift of the
rotational field $-i\zeta \partial_\zeta$. Explicitly,
\begin{equation}
 Z = -i
 \frac{\left(\frac{1}{\bar{y}_1}, \dots, \frac{1}{\bar{y}_k},0,\dots,0\right)}
 {\frac{1}{|y_1|^2} + \dots + \frac{1}{|y_k|^2}},
\end{equation}
In particular, since each $|y_j|^2$ is invariant, the flow $\psi_t$
of $Z$ is well-defined for all times, in spite of the noncompactness
of the fibres. As in our discussion of Lefschetz pencils, we define
complex volume forms $\eta_{Y_z}$ on the fibres by (formally)
dividing $\eta_Y$ by some complex one-form on the base, which in this
case will be $d\zeta/\zeta^2$. As a consequence, if at some point $y
\in Y_\zeta = p^{-1}(\zeta)$ we have a Lagrangian subspace $\Lambda
\subset TY_y$ with $Dp(\Lambda) = a\R$, $a \in S^1$, then the phase
of the vertical part $\Lambda^v = \Lambda \cap \mathit{ker}(Dp)$ inside
$Y_\zeta$ is related to the phase of $\Lambda$ by
\begin{equation} \label{eq:divide-degree}
\alpha_{Y_\zeta}(\Lambda^v) = \frac{\zeta^4}{a^2 |\zeta|^4}
\alpha_Y(\Lambda).
\end{equation}

For each point $y$ with $p(y) = \zeta \neq 0$, consider the
Lagrangian tangent subspace $\Lambda_y = \R i y_1 \oplus \dots \oplus
\R i y_k \oplus \R^{n-k} \subset TY_y$. This satisfies
$\alpha_Y(\Lambda_y) = (-1)^k$ and $Dp(\Lambda_y) = i\zeta \R$, hence
\eqref{eq:divide-degree} says that
\begin{equation} \label{eq:lambda-shift}
\alpha_{Y_z}(\Lambda_y^v) = \frac{\zeta^2}{|\zeta|^2} (-1)^{k-1}.
\end{equation}
$\psi_t$ takes the $\Lambda_y$ to each other, and rotates the base
coordinate $\zeta$, hence it also maps the $\Lambda_y^v$ to each
other. Let $h_\zeta = \psi_{2\pi}|Y_\zeta: Y_\zeta \rightarrow Y_\zeta$ be the
monodromy on the fibre over $\zeta$. By definition and
\eqref{eq:lambda-shift}, the canonical grading of this map satisfies
\begin{equation} \label{eq:integrate-shift}
\begin{aligned}
 \alphagr_{h_{\zeta}}(\Lambda_y^v)
 & = \frac{1}{2\pi} \int_0^{2 \pi}
 \frac{d}{dt} \arg(\alpha_{Y_z}(D\psi_t(\Lambda_y^v))) \, dt
 \\ & = \frac{1}{2\pi} \int_0^{2 \pi}
 \frac{d}{dt} \arg( (-1)^{k-1} e^{-2it}) \, dt = -2.
\end{aligned}
\end{equation}
If one passes to a sufficiently large iterate such as $h_\zeta^n$,
the grading on $\Lambda_y^v$ is $-2n$, hence by Lemma
\ref{th:graded-variation} the grading of any $\Lambda$ is $<-n-1$. In
a slight abuse of terminology (since $Y_\zeta$ is not compact), we
conclude that the monodromy is a negative graded symplectic
automorphism, in a way which is uniform for all $\zeta$.

\subsection{}
A different local model appears near the singular points of the base
locus of a quasi-pencil. In Ruan's work on the SYZ conjecture, see in
particular \cite{ruan02b}, these are responsible for the
singularities of the Lagrangian torus fibration. For closely related
reasons, the negativity of the monodromy breaks down at such points.

\begin{assumptions} \label{as:2}
Fix $n \geq 3$ and $3 \leq k \leq n$. We again take $Y = \C^n$ with
the standard complex and symplectic forms, but now with the rational
map
\begin{equation}
p: Y \dashrightarrow \C, \qquad p(y) = \frac{y_2 \dots y_k}{y_1},
\end{equation}
in the presence of the complex volume form (away from the closure of
the zero fibre)
\begin{equation}
\eta_Y = dy_1 \wedge \frac{dy_2}{y_2} \wedge \dots \wedge
\frac{dy_k}{y_k} \wedge dy_{k+1} \wedge \dots \wedge dy_n.
\end{equation}
\end{assumptions}

Denote by $Y_\zeta$ the closure of $p^{-1}(\zeta)$. For $\zeta \neq 0$,
each such fibre can be identified with $\C^{n-1}$ by projecting to
$(y_2,\dots,y_n)$. We will use this identification freely, hoping
that this does not cause too much confusion. For instance, the
symplectic form on $Y_\zeta \iso \C^{n-1}$ is given by
\begin{equation} \label{eq:zeta-symplectic-form}
 \o_{Y_\zeta} = \frac{i}{2} \partial\bar\partial
 \Big(\frac{|y_2\dots y_k|^2}{|\zeta|^2}
 + |y_2|^2 + \dots + |y_n|^2 \Big)
\end{equation}
As in our previous model, we write $Z$ for the clockwise monodromy
vector field (defined by passing to the graph of $p$). In view of our
identification, this becomes a vector field on
$\C^* \times \C^{n-1}$, given in coordinates $(\zeta,y_2,\dots,y_n)$ by
\begin{equation}
Z = \Big(-i\zeta\partial_\zeta, -i
\frac{\left(\frac{1}{\bar{y}_2},\dots,\frac{1}{\bar{y}_k},0,\dots,0\right)}
{\frac{|\zeta|^2}{|y_2|^2 \cdots |y_k|^2} + \frac{1}{|y_2|^2} + \dots + \frac{1}{|y_k|^2}}\Big).
\end{equation}
Again, this preserves $|y_2|^2, \dots, |y_n|^2$, so that its flow $\psi_t$ is well-defined for all times.

Taking $\eta_Y$ on the total space and $d\zeta/\zeta^2$ on the base,
one finds that (up to a sign which is irrelevant for us) the induced volume form $\eta_{Y_\zeta}$ is the
standard form on $\C^{n-1}$. At each point $y \in Y_\zeta$ with $y_2
\dots y_k \neq 0$, consider the Lagrangian subspace $\Lambda_y^v = \R
i y_2 \oplus \dots \oplus \R i y_k \oplus \R^{n-k}$. These subspaces
are preserved by the flow $\psi_t$, and a computation similar to
\eqref{eq:integrate-shift} shows that the canonical grading of the
monodromy map $h_\zeta = \psi_{2\pi}|Y_\zeta$ satisfies
\begin{equation}
 \alphagr_{h_\zeta}(\Lambda_y^v) = - 2
 \frac{1}{1 + \frac{1}{\frac{|y_1|^2}{|y_2|^2} +
 \cdots + \frac{|y_1|^2}{|y_k|^2}}}.
\end{equation}
This is obviously negative, but of course we also need to consider
the grading at the points where $\Lambda_y^v$ is not defined. Suppose
for simplicity that the sizes of the first $k-1$ coordinates of $y
\in Y_\zeta$ are ordered, $|y_2| \leq |y_3| \leq \dots |y_{k-1}|$.
Then
\begin{equation}
 \alphagr_{h_\zeta}(\Lambda_y^v) \leq -2
 \frac{1}{1 + \frac{|\zeta|^2}{|y_3|^{2k-4}}},
\end{equation}
which means that the grading is bounded above by something which
depends on the size of the second smallest of the $(y_2,\dots,y_k)$. As a
consequence, if one has a sequence
\begin{equation}
 y^{(m)} \rightarrow y^{(\infty)}, \quad
 y_2^{(m)} \dots y_k^{(m)} \neq 0 \text{ for all $m$,}
\end{equation}
and such that only one of the $(y_2^{(\infty)}, \dots,
y_k^{(\infty)})$ vanishes, the
$\alphagr_{h_\zeta}(\Lambda^v_{y^{(m)}})$ are bounded above by a
negative constant, hence cannot go to zero. Roughly speaking, this
means that $h_\zeta$ is negative away from the subset where at least
two of the $y_2,\dots,y_k$ vanish. In the special case $k = n = 3$,
this subset is a point $0 \in Y_{\zeta}$; one can then apply
\eqref{eq:compose-gradings} and Lemma \ref{th:graded-variation} much
as in our previous local model, to get the following statement:

\begin{lemma} \label{th:bad-points}
Suppose that $k = n = 3$. Let $W \subset Y_\zeta$ be some open ball
around the origin. Then $h_\zeta(W) = W$, and there is a $d>0$ such
that the grading of $h_\zeta^d|(Y_\zeta \setminus W)$ is negative.
\qed
\end{lemma}

\begin{remark}
As mentioned before, ideas related to the SYZ conjecture suggest
that the monodromy should be a fibrewise translation in a SLAG
fibration. From this point of view, the negativity of the monodromy
is based on thinking that the phases of the tangent spaces to the
fibration should be shifted by $-2$. In the situation of Assumptions
\ref{as:1}, this is indeed the case: $L_r \iso T^{k-1} \times
\R^{n-k} \subset Y_\zeta$, given by
\begin{equation}
 L_r = \{|y_1|^2 = r_1, \dots, |y_k|^2 = r_k,\, \im(y_{k+1}) = r_{k+1}, \,
 \dots, \, \im(y_n) = r_n\}
\end{equation}
with $r_1 \dots r_k = |\zeta|^2$, are special Lagrangian with respect
to $\eta_{Y_\zeta}$. The monodromy translates each of them, and the
shift is computed in \eqref{eq:integrate-shift}. For Assumptions
\ref{as:2}, there is a standard SLAG fibration with singularities on
$Y_\zeta \iso \C^{n-1}$ (with respect to
\eqref{eq:zeta-symplectic-form} and the standard complex volume
form), see e.g.\ \cite{gross01b}:
\begin{multline} \label{eq:syz}
 \qquad L_r = \{
 |y_3|^2 - |y_2|^2 =
 r_1, \, \dots \, |y_k|^2 - |y_2|^2 = r_{k-2}, \\
 \im(y_2\dots y_k) = r_{k-1}, \,
 \im(y_{k+1}) = r_k, \dots, \im(y_n) = r_{n-1}\}. \qquad
\end{multline}
However, the monodromy preserves this only asymptotically far away
from the singular locus, and hence the grading is only asymptotically
equal to $-2$. To see more precisely what happens, we apply a
symplectic reduction procedure, which means mapping $Y_\zeta
\rightarrow \C$ by sending $y$ to $w = y_2y_3 \dots y_k$. The fibres
of the SLAG fibration lie over horizontal lines in the base $\C$. The
clockwise monodromy also fibres over a diffeomorphism of $\C$, which
rotates each point $w$ by some $|w|$-dependent angle. For $|w| \gg
0$, this is almost a full clockwise rotation, hence it approximately
preserves the horizontal lines, but near $w = 0$ the angle of
rotation undergoes a full change from $0$ to $-2\pi$.
\end{remark}

\subsection{}
We now relax our previous assumptions on the complex and symplectic
structures, for the first local model. We choose to work in Darboux
coordinates, and with a complex structure and holomorphic function
which are not standard, but which still preserve a $T^k$ symmetry.

\begin{assumptions} \label{as:1prime}
Fix $n \geq 2$ and $2 \leq k \leq n$.
\begin{itemize} \itemsep1em
\item
$Y \subset \C^n$ is an open ball around the origin, carrying the
standard symplectic form $\o_Y$ and the standard diagonal
$T^k$-action
\begin{equation}
\rho_s(y) = (e^{is_1}y_1,\dots,e^{is_k}y_k,y_{k+1},\dots,y_n),
\end{equation}
with moment map $\mu: Y \rightarrow (\Rgeq)^k$. For any regular
moment value $r$, the quotient $Y^{\red,r}$ can be identified with an
open subset of $\C^{n-k}$, with the standard symplectic form
$\o_{Y^{\red,r}}$. We denote points in the reduced spaces by
$y^{\red} = (y_{k+1},\dots,y_n)$.

\item
$J_Y$ is a complex structure which is tamed by $\o_Y$. At the point
$y = 0$ (but not necessarily elsewhere), it is $\o_Y$-compatible and
$T^k$-invariant.

\item
$p: Y \rightarrow \C$ is a $J_Y$-holomorphic function with the
following properties:
\begin{romanlist}
\item \label{item:t-sym}
$p(\rho_s(y)) = e^{i(s_1+\dots+s_k)}p(y)$.
\item \label{item:nondeg}
$\partial_{y_1}\dots \partial_{y_k} p$ is nonzero at the point $y =
0$.
\end{romanlist}

\item
$\eta_Y$ is a $J_Y$-complex volume form on $Y \setminus p^{-1}(0)$,
such that $p(y)\eta_Y$ extends to a smooth form on $Y$, which is
nonzero at $y = 0$. Let $\alpha_Y$ be the corresponding phase
function.
\end{itemize}
\end{assumptions}

\begin{lemma} \label{th:product-function}
One can write
\begin{equation}
p(y) = 2^{-k/2} y_1 y_2\dots y_k\, q(\half |y_1|^2,\dots,\half
|y_k|^2,y_{k+1},\dots,y_n)
\end{equation}
for some smooth function $q$ satisfying $q(0) \neq 0$.
\end{lemma}

\proof Consider the Taylor expansion of $p$ around $y = 0$. From
property \ref{item:t-sym} of $p$ one sees that each monomial which
occurs with nonzero coefficient in that expansion contains one of the
factors $y_l$ for each $1 \leq l \leq k$. As a consequence,
$p(y)/(y_1\dots y_k)$ extends smoothly over $y = 0$. The same
argument can be applied to points where only certain of the
coordinates $y_1 \dots y_k$ vanish. Again applying \ref{item:t-sym},
one sees that $p(y)/(y_1\dots y_k)$ is $T^k$-invariant, which means
that it can be written in the form stated above. Finally, $q(0) \neq
0$ is property \ref{item:nondeg}. \qed

\begin{lemma} \label{th:shrink}
After making $Y$ smaller if necessary, one has
\begin{align*}
 p(y) = 0 & \Longleftrightarrow  y_l = 0
 \text{ for some } l = 1,\dots,k; \\
 Dp(y) = 0 & \Longleftrightarrow y_l = 0
 \text{ for at least two } l = 1,\dots,k.
\end{align*}
\end{lemma}

\proof Suppose that $q(y) \neq 0$ on the whole of $Y$. Then the first
statement is obvious, and so is the $\Leftarrow$ implication in the
second one. Conversely, suppose that $y$ is a critical point.
Property \ref{item:t-sym} implies that $p(y) = 0$, so that $y_l = 0$
for some $l = 1,\dots,k$. From 0 = $Dp(y) = y_1 \dots \hat{y}_l \dots
y_k \, q(y)\, dy_l$ one sees that another of the $y_1,\dots, y_k$ has
to vanish too. \qed

We will assume from now on that the conclusions of Lemma
\ref{th:shrink} hold. Consider the function $H(y) = -\half |p(y)|^2$,
its Hamiltonian vector field $X$ and flow $\phi_t$. Since $H$ is
$T^k$-invariant, the flow (where defined) is equivariant and
preserves the level sets of the moment map. For every regular moment
value $r$, we can consider the induced function
\begin{equation} \label{eq:h-red}
H^{\red,r}(y^\red) = -\half r_1 \dots r_k \,
|q(r_1,\dots,r_k,y_{k+1},\dots,y_n)|^2,
\end{equation}
on $Y^{\red,r}$, its vector field $X^{\red,r}$ and flow
$\phi_t^{\red,r}$. If the moment value is small, so are
\eqref{eq:h-red} and all its $y^\red$-derivatives, which means that
the reduced flow moves very slowly. We need a quantitative version of
this:

\begin{lemma} \label{th:small-flow}
Given $\epsilon_1,\epsilon_2>0$, there are
$\delta_1,\delta_2,\delta_3>0$ such that for all
\begin{equation}
\begin{cases}
 r \in (\R^{\scriptscriptstyle>0})^k &
 \text{with $\|r\| < \delta_1$}, \\
 y^\red \in Y^{\red,r} &
 \text{with $\|y^\red\| < \delta_2$}, \\
 t \in \R &
 \text{with $|t| < \delta_3 r_1^{-1} \dots r_k^{-1}$},
\end{cases}
\end{equation}
the following holds: $\phi_t^{\red,r}$ is well-defined near $y^\red$,
and satisfies
\begin{align}
 \label{eq:c0} &
 \|\phi_t^{\red,r}(y^{\red}) \| < \epsilon_1, \\
 \label{eq:c1} &
 \|(D\phi_t^{\red,r})_{y^{\red}} - \mathit{Id} \| < \epsilon_2.
\end{align}
\end{lemma}

\proof By \eqref{eq:h-red} there are $\delta_1>0$,
$\delta_2 \in (0;\min(\epsilon_1/2,\epsilon_2/2))$ and a $C>0$, such
that
\begin{equation}
\|X^{\red,r}_{y^\red}\|, \; \|DX^{\red,r}_{y^{\red}}\| \leq C
r_1\dots r_k
\end{equation}
for each $\|r\| < \delta_1$, $\|y^{\red}\| < 2\delta_2$ (we use the
Hilbert-Schmidt norm for matrices $DX^{\red,r}$, normalized to
$\|\mathit{Id}\| = 1$). We claim that one can take $\delta_3 =
\delta_2/C $. The first estimate \eqref{eq:c0} is clear: if one
starts at a point $y^\red$ of norm less than $\delta_2$, and moves
with speed at most $C r_1 \dots r_k$ for a time less than $\delta_2
C^{-1} r_1^{-1} \dots r_k^{-1}$, the endpoint will be of norm less
than $2\delta_2 < \epsilon_1$. Next, by considering the linearisation of the flow $\phi_t^{\mathit{red},r}$ along an orbit, we get
an exponential growth bound for the time interval $t$ relevant to us:
$\|D\phi_t^{\red,r}\| \leq \exp(C r_1 \dots r_k t) \leq
\exp(\delta_2)$. One may assume that $\exp(\delta_2) \leq 2$. By
plugging this back into the linearized equation, one gets the linear
estimate $\|D\phi_t^{\red,r}-\mathit{Id}\| \leq 2 C r_1 \dots r_k t
\leq 2 \delta_2 < \epsilon_2$, from which \eqref{eq:c1} follows. \qed

We now need to deal with the fact that $J_Y$ and $\eta_Y$ are only
approximately $T^k$-invariant. Let $J_Y'$ be the constant complex
structure on $Y$ which agrees with $J_Y$ at $y = 0$. This is
$\o_Y$-compatible and $T^k$-invariant, so
\begin{equation}
J_Y' = \begin{pmatrix} i \cdot 1_k & 0 \\ 0 & J_{Y^{\red}}'
\end{pmatrix}.
\end{equation}
$J_{Y^{\red}}'$ is the induced complex structure on $Y^{\red,r}$ for
all $r$. Similarly, we consider the $J_Y'$-complex volume form
$\eta'_Y$ which is obtained by taking the value of $p(y)\eta_Y$ at
the point $y = 0$, and dividing that by $y_1\dots y_k$. This can be
written as
\begin{equation}
\eta_Y' = i^k \frac{dy_1}{y_1} \wedge \dots \wedge \frac{dy_k}{y_k}
\wedge \eta_{Y^\red}',
\end{equation}
hence is $T^k$-invariant. $\eta_{Y^\red}'$, which is the induced
complex volume form on the reduced spaces, is again constant and
$r$-independent. Let $\alpha_Y'$ be the phase function associated to
$(J_Y',\eta_Y')$ (unlike $\alpha_Y$, this is a classical phase
function), and $\alpha_{Y^{\red,r}}'$ the induced phase functions on
the quotients.

\begin{lemma} \label{th:phase-difference}
Given $\epsilon_3>0$, there is $\delta_4>0$ such that for all $\|y\|
< \delta_4$ with $p(y) \neq 0$, and all $\Lambda \in \LL_{Y,y}$,
\begin{equation}
\Big|\frac{1}{2\pi}\arg(\alpha_Y'(\Lambda)/\alpha_Y(\Lambda))\Big| <
\epsilon_3
\end{equation}
(by this we mean that there is a branch of the argument with that
property).
\end{lemma}

\proof The quotient $\alpha_Y'(\Lambda)/\alpha_Y(\Lambda)$ does not
change if we multiply both $\eta_Y$, $\eta_Y'$ with $y_1\dots y_k$.
But by construction,
\begin{align*}
 & y_1 \dots y_k \eta_Y = 2^{-k/2} q(\half |y_1|^2,\dots, \half |y_k|^2,
 y_{k+1},\dots,y_n)^{-1} \sigma_y, \\
 & y_1 \dots y_k \eta_Y' = 2^{-k/2} q(0)^{-1} \sigma_0
\end{align*}
where $\sigma$ is some smooth $n$-form with $\sigma_0 \neq 0$, so the
associated phases also become close as $y \rightarrow 0$. \qed

$\phi_t$ has a canonical grading $\alphagr_{\phi_t}$ with respect to
$\alpha_Y$, obtained by starting with the trivial grading for
$\phi_0$ and extending continuously. Denote by $\alphagr_{\phi_t}'$
the grading obtained in the same way, but with respect to
$\alpha_Y'$. There are corresponding gradings
$\alphagr_{\phi_t^{\red,r}}'$ on the reduced spaces, with respect to
$\alpha_{Y^{\red,r}}'$.

\begin{lemma} \label{th:small-phase}
Given $\epsilon_4>0$, one can find $\delta_5,\delta_6,\delta_7>0$
such that for all
\begin{equation}
\begin{cases}
 r \in (\R^{\scriptscriptstyle>0})^k &
 \text{with $\|r\| < \delta_5$}, \\
 y^\red \in Y^{\red,r} &
 \text{with $\|y^\red\| < \delta_6$}, \\
 \Lambda^\red \in \LL_{Y^{\red,r},y^\red} & \text{and} \\
 t \in \R &
 \text{with $|t| < \delta_7 r_1^{-1} \dots r_k^{-1}$},
\end{cases}
\end{equation}
$\phi_t^{\red,r}$ is well-defined near $y^\red$, and its grading is
$|\alphagr_{\phi_t^{\red,r}}'(\Lambda^\red)| < \epsilon_4$.
\end{lemma}

\proof This is an easy consequence of Lemma \ref{th:small-flow}: one
can achieve that the derivative $(D\phi_t^{\red,r})_{y^\red}$ remains
very close to the identity in the whole range of $t$'s that we are
considering. Since the volume form $\eta_{Y^{\red,r}}$ is constant
and independent of $r$, it follows that the grading remains small.
\qed

\begin{lemma} \label{th:phi-grading-estimate}
Given $\epsilon_5>0$, one can find $\delta_8,\delta_9,\delta_{10} >
0$ such that for all
\begin{equation}
\begin{cases}
 y \in Y & \text{with $0 < |p(y)| < \delta_8$ and $\|y\| < \delta_9$,} \\
 \Lambda \in \LL_{Y,y}, & \text{and} \\
 t \in \R & \text{with $|t| < \delta_{10} |p(y)|^{-2}$,}
\end{cases}
\end{equation}
$\phi_t$ is well-defined near $y$, and its grading satisfies
$
|\alphagr_{\phi_t}(\Lambda)| < n + \epsilon_5.
$
\end{lemma}

\proof Apply Lemma \ref{th:small-flow} to the corresponding point
$y^\red$ in the quotient. In particular, one can choose
$\delta_8,\delta_9,\delta_{10}$ in such a way that the reduced flow
exists near $y^\red$ for $|t| < \delta_{10} |p(y)|^{-2}$. Since the
fibres of the quotient map are compact (tori), it follows that
$\phi_t$ is well-defined near $y$ in the same range of $t$. We now
apply the same basic reasoning as in Lemma
\ref{th:negative-reduction}. After possibly making
$\delta_8,\delta_9$ smaller, Lemma \ref{th:small-phase} ensures that
there is a $\Lambda \in \LL_{Y,y}$ such that
$|\alphagr_{\phi^t}'(\Lambda)| < \epsilon_5/2$. By Lemma
\ref{th:graded-variation}, the grading is $< n + \epsilon_5/2$ for
all other Lagrangian subspaces at the same point. Finally, one uses
Lemma \ref{th:phase-difference} to pass from the grading associated
to $\eta_Y'$ to that for $\eta_Y$, again making the $\delta$s smaller
so that the phase difference becomes $< \epsilon_5/2$. \qed

The next step is to relate $\phi_t$ to the monodromy. As usual, we
take the vector field $-i\zeta \partial_\zeta$ on $\C^*$, and lift it
in the unique way to a horizontal vector field $Z$ on $Y \setminus
p^{-1}(0)$. This is well-defined because $\o_Y$ tames $J_Y$, hence
the regular fibres $Y_\zeta$, $\zeta \neq 0$, are symplectic
submanifolds. We claim that
\begin{equation} \label{eq:f}
Z = f X
\end{equation}
for some positive function $f$ on $Y \setminus p^{-1}(0)$. First,
because $H$ is constant on each fibre, the associated vector field
$X$ is horizontal. Second, because $\o_Y(X,Z) = 0$ and the space of
horizontal vectors $TY^h$ is two-dimensional, it is true that $X =
f^{-1}Z$ for some $f$. A little more thought shows that
\begin{equation} \label{eq:explicit-psi}
f = \frac{\o_Y|TY^h}{p^*\o_\C|TY^h} > 0,
\end{equation}
where $TY^h$ is the $\o_Y$-horizontal subspace (which is
two-dimensional, so any two two-forms on it differ by a scalar). Let
$(\psi_s)$ be the flow of $Z$, so that $h_\zeta =
\psi_{2\pi}|Y_\zeta$ is the clockwise monodromy (of course, like
$\phi_s$ this flow is only partially defined). As a consequence of
\eqref{eq:f}, the unparametrized flow lines of $\phi_t$ and $\psi_t$
agree, so
\begin{equation} \label{eq:gt}
\psi_t(y) = \phi_{g_t(y)}(y), \quad \text{where }\;\; g_t(y) =
\int_0^t f(\psi_\tau(y)) d\tau.
\end{equation}

\begin{lemma} \label{th:minimum-speed}
For every $d>0$, $\epsilon_6>0$ there are $\delta_{11},\delta_{12}
> 0$ such that the following holds. For all $0 < |\zeta| < \delta_{11}$ and
all $y \in Y_\zeta$ with $\|y\| < \delta_{12}$, the $d$-fold
monodromy $h^d_\zeta$ is well-defined, and \eqref{eq:gt} holds for $t
= 2\pi d$ with
\begin{equation}
g_{2\pi d}(y) \leq \frac{\epsilon_6}{|\zeta|^2}.
\end{equation}
\end{lemma}

\proof Using property \ref{item:t-sym} of $p$ and its
$J_Y$-holomorphicity, one obtains that $Dp_y(J_Y(-iy_1,0,\dots,0)) =
p(y)$ and hence $dH_y(J_Y(-iy_1,0,\dots,0)) = |p(y)|^2 = |\zeta|^2$.
It follows that for small $y$, $\|X_y\| \geq c_1 |\zeta|^2/\|y\|$ for
some $c_1>0$. Hence there is a $C_2>0$ such that
\begin{equation}
 C_2^{-1} |\zeta|^4/\|y\|^2 \leq \o_Y(X,J_YX) = f^{-2} \o_Y(Z,J_YZ).
\end{equation}
On the other hand, $\o_Y(Z,J_YZ) = |\zeta|^2 f$ by definition of $Z$
and \eqref{eq:explicit-psi}, so finally,
\begin{equation}
 f(y) \leq C_2 \frac{\|y\|^2}{|\zeta|^2}.
\end{equation}
Integrating gives a similar inequality $g_{2\pi d}(y) \leq C_3
\|y\|^2/|\zeta|^2$. Clearly, by restricting to smaller $y$ one can
make the right hand side less than $\epsilon_6/|\zeta|^2$ for any
desired $\epsilon_6$ (this, by the way, is the crucial estimate in
our whole computation). \qed

We use $J_Y$-complex volume forms $\eta_{Y_\zeta}$ on the fibres
obtained in the usual way from $\eta_Y$, so that
\eqref{eq:divide-degree} holds. Take $d>0$, $y \in Y_\zeta$, and a
Lagrangian subspace $\Lambda \subset TY_y$ with $Dp(\Lambda) =
i\zeta\R$. This is necessarily of the form $\Lambda = \Lambda^v
\oplus \R Z_y$ for some Lagrangian subspace $\Lambda^v \subset
T_y(Y_\zeta)$. By \eqref{eq:f} and \eqref{eq:gt} one has
\begin{equation} \label{eq:relation-between-diff}
 D(\psi_t)_y = D(\phi_{g_t})_y +
 X_{\psi_t(y)} \otimes dg_t =
 D(\phi_{g_t})_y + Z_{\psi_t(y)} \otimes \frac{dg_t}{f(\psi_t(y))}
\end{equation}
Since $\Lambda$ contains $Z_y$, and $\psi_t$ is the clockwise
parallel transport flow, $D(\psi_t)_y(\Lambda)$ contains
$Z_{\psi_t(y)}$. From this and \eqref{eq:relation-between-diff} it
follows that
\begin{equation}
 D\psi_t(\Lambda) = D(\phi_{g_t})(\Lambda), \qquad
 D(\psi_t|Y_\zeta)(\Lambda^v) = D(\phi_{g_t})(\Lambda) \cap \mathit{ker}(Dp),
\end{equation}
and hence by \eqref{eq:divide-degree} that
$\alpha_{\psi_t|Y_\zeta}(\Lambda^v) = e^{-4\pi i t}
\alpha_{\phi_{g_t}}(\Lambda)$. By considering the behaviour of this
over the time $t \in [0,2\pi d]$ one finds that the canonical grading
of the monodromy satisfies
\begin{equation} \label{eq:relative-shift}
 \alphagr_{h_\zeta^d}(\Lambda^v) =
 \alphagr_{\phi_{g_{2\pi d}(y)}}(\Lambda)-2d.
\end{equation}

\begin{lemma} \label{th:local-model-1}
For every $d>0$ and $\epsilon_7>0$ there are $\delta_{13},
\delta_{14} >0$ such that the following holds. For every $y \in
Y_\zeta$ with $0 < |\zeta| < \delta_{13}$ and $\|y\| < \delta_{14}$,
and every Lagrangian subspace $\Lambda^v \subset T(Y_\zeta)_y$, we
have that $h_\zeta^d$ is well-defined near $y$, and its canonical
grading satisfies
\begin{equation}
\alphagr_{h_\zeta^d}(\Lambda^v) < n-2d + \epsilon_7.
\end{equation}
\end{lemma}

\proof Lemma \ref{th:phi-grading-estimate} says that one can choose
$\delta_{13},\delta_{14}$ so small that there is a $\delta_{15}>0$
with
\begin{equation} \label{eq:epsilon7}
|\alphagr_{\phi_t}(\Lambda)| < n + \epsilon_7
\end{equation}
for all $|t| < \delta_{15}/|p(y)|^2$. On the other hand, after
possibly making $\delta_{13}, \delta_{14}$ smaller, one knows from
Lemma \ref{th:minimum-speed} that $|g_{2\pi d}| <
\delta_{15}/|p(y)|^2$, so \eqref{eq:epsilon7} applies to the values
of $t$ that are relevant for our monodromy map.
\eqref{eq:relative-shift} completes the proof. \qed

This shows that the basic feature of our first local model, the
uniform negativity of the monodromy, survives even if one generalizes
from the idealized situation which we considered first (Assumptions
\ref{as:1}) to a more realistic one (Assumptions \ref{as:1prime}). In
principle, one should carry out a parallel discussion for the second
local model. This would be rather more delicate, since the grading
becomes zero along a submanifold; an arbitrary small perturbation can
potentially make it positive, thereby spoiling the entire argument.
Fortunately, in the $K3$ case the second model occurs only at
isolated points, and one can achieve that the local picture near
those points is standard (satisfies Assumptions \ref{as:2}), so that
the previous discussion suffices after all.

\subsection{}
We will now explain how to patch together the local models to
understand the ``large complex structure limit'' monodromy. This
is somewhat similar to the strategy used by Ruan to construct
(singular) Lagrangian torus fibrations \cite{ruan02b}. We consider a
Fano threefold $X$ with $o_X = \K_X^{-1}$, and a pair of sections
$\sigma_{X,0},\sigma_{X,\infty}$ which generate a quasi-Lefschetz
pencil of hypersurfaces $\{X_z\}$. From $\sigma_{X,\infty}$ we get a
rational complex volume form $\eta_X$ with a pole along $X_\infty$.
There are canonical induced complex volume forms $\eta_{X_z}$ on the
fibres $X_z$ (if that is smooth and $z \neq \infty$), which are
therefore Calabi-Yau surfaces. In fact, they must be $K3$s, since
$H^1(X) = 0$ for any Fano variety, hence $H^1(X_z) = 0$ by the
Lefschetz hyperplane theorem.

\begin{lemma} \label{th:kaehler-1}
One can find a K{\"a}hler form $\o_X'$ in the same cohomology class
as $\o_X$, with the following properties. (i) Near every point where
three components of $X_\infty$ meet, $(X,\o_X',1/\pi)$ is modelled on
$(Y,\o_Y,p)$ from Assumptions \ref{as:1}, with $n = k = 3$. (ii) Near
every point where $X_0$ meets two components of $X_\infty$,
$(X,\o_X',1/\pi)$ is modelled on $(Y,\o_Y,p)$ from Assumptions
\ref{as:2}, with $n = k = 3$.
\end{lemma}

This is easy because the set of such points (which we will call
``lowest stratum points'' in future) is finite. Near each of them,
there are local holomorphic coordinates in which $(X,1/\pi)$ becomes
$(Y,p)$. One can then locally modify the K{\"a}hler form to make it
standard, see e.g.\ \cite[Lemma 7.3]{ruan02b} or \cite[Lemma
1.7]{seidel01}.

\begin{lemma} \label{th:kaehler-2}
One can find a K{\"a}hler form $\o_X''$ in the same cohomology class
as $\o_X$, which agrees with $\o_X'$ in a neighbourhood of the lowest
stratum points, and such that in addition, any two components of
$X_\infty$ meet orthogonally. \qed
\end{lemma}

Ruan proves a much more general result in \cite[Theorem
7.1]{ruan02b}. The main idea is to apply the techniques from the
previous Lemma to each fibre of the normal bundle $\nu$ to the
intersection of two components of $X_\infty$.
%

Before proceeding, we need some more notation. Denote by $C$ the
irreducible components of $X_\infty$, and by $L_C = \O(-C)$ the
corresponding line bundles. Let $\gamma_C$ be their canonical
sections, with a zero along $C$. Take two components $C \neq C'$, and
choose a small neighbourhood $U_{C,C'}$ of $C \cap C'$. On this
consider the vector bundle
\begin{equation}
N_{C,C'} = L_C \oplus L_C^{-1},
\end{equation}
We want to equip $L_C$ and hence $N_{C,C'}$ with hermitian metrics
and the corresponding connections, with the following prescribed
behaviour near the lowest stratum points in $U_{C,C'}$:
\begin{itemize} \itemsep1em
\item
If $x$ is a point where $C \cap C'$ meets another component $C''$, we
choose local coordinates for $X$ as given by Lemma
\ref{th:kaehler-1}, so that Assumptions \ref{as:1} hold, and a local
trivialization of $L_C$ such that $\gamma_C(y) = y_1$. Take the
trivial metric on $L_C$ with respect to this trivialization.

\item
If $x$ is a point where $C \cap C'$ meets $X_0$, we choose local
coordinates for $X$ as given by Lemma \ref{th:kaehler-1}, so that
Assumptions \ref{as:2} hold, and a local trivialization of $L_C$ such
that $\gamma_C(y) = y_2$. Again, we take the trivial metric.
\end{itemize}

Take the K{\"a}hler metric associated to $\o_X''$, and use its
exponential map to define a retraction $r_{C,C'}: U_{C,C'}
\rightarrow C \cap C'$ of an open neighbourhood onto $C \cap C'$. A
corresponding use of parallel transport on $N_{C,C'}$, with respect
to our chosen connection, yields an isomorphism $R_{C,C'}:
N_{C,C'}|U_{C,C'} \rightarrow r_{C,C'}^*(N_{C,C'}|C \cap C')$. Taking
the two objects together, we have a map
\begin{equation}
 \nu_{C,C'}: U_{C,C'} \setminus X_0
  \xrightarrow{\xi_{C,C'} = (\gamma_C,\pi_M^{-1}\gamma_C^{-1})}
 N_{C,C'}|U_{C,C'} \xrightarrow{R_{C,C'}} N_{C,C'}|C \cap C'.
\end{equation}

\begin{lemma} \label{th:kaehler-3}
After making $U_{C,C'}$ smaller if necessary, it carries a unique
$T^2$-action $\sigma_{C,C'}$ which, under $\nu_{C,C'}$, corresponds
to the obvious fibrewise $T^2$-action on $N_{C,C'}|C \cap C'$. In
particular,
\begin{equation} \label{eq:rotate-base}
\pi_M(\sigma_{C,C',s}(x)) = e^{-i(s_1+s_2)}\pi_M(x).
\end{equation}
Moreover, at each point of $C \cap C'$, the action preserves $\o_X''$
and the complex structure; the same holds in a neighbourhood of the
lowest stratum points.
\end{lemma}

\proof It is sufficient to verify the existence and uniqueness
statement locally. First, take a point $x \in C \cap C'$ which is not
a lowest stratum point. $\gamma_C^{-1}$ has a simple pole along $C$,
and $\pi_M^{-1}$ has a simple zero along $C \cup C'$, hence
$\pi_M^{-1}\gamma_C^{-1}$ has a simple zero along $C'$. It follows
that $\xi_{C,C'}$ is transverse to the zero-section at $x$, so
$\nu_{C,C'}$ is a local diffeomorphism; we simply pull back the
$T^2$-action by it. By definition of the connection on $N_{C,C'}$,
there is a commutative diagram
\begin{equation}
\xymatrix{
 {U_{C,C'} \setminus X_0} \ar[r]^-{\xi_{C,C'}} \ar[dr]_-{\pi_M^{-1}} &
 {N_{C,C'}|U_{C,C'}} \ar[r]^-{R_{C,C'}} \ar[d] &
 {N_{C,C'}|C \cap C'} \ar[dl] \\
 & {\C} &
}
\end{equation}
where the unlabeled arrows are the canonical pairings $L_C \oplus
L_C^{-1} \rightarrow \C$. This shows that our pullback $T^2$-action
satisfies \eqref{eq:rotate-base}. At the point $x$ itself,
\begin{equation} \label{eq:dnu}
D\nu_{C,C'}: TX_x \longrightarrow T(C \cap C')_x \oplus N_x
\end{equation}
is given by orthogonal projection to $C \cap C'$ with respect to the
$o_{X,2}$-metric, together with $D\xi_x: TX_x \rightarrow N_x$. This
is $\C$-linear, which implies that $\sigma_{C,C',x}$ preserves the
complex structure. As for the symplectic structure, observe that by
Lemma \ref{th:kaehler-2} $TX_x$ splits into three complex
one-dimensional pieces, namely $T(C \cap C')_x$ and the normal
directions to that inside $C,C'$. On the right hand side of
\eqref{eq:dnu}, this corresponds to the splitting of $N_x$ into line
bundles. The symplectic structures on the each summand may not be
strictly compatible with \eqref{eq:dnu}, but they differ only by a
multiplicative scalar, which is enough to conclude that $\sigma_{C,C'}$ is symplectic at $x$.

It remains to deal with the situation near the lowest stratum points.
For this, we choose local coordinates and the trivialization of
$N_{C,C'}$ as arranged before. Near a point of $C \cap C' \cap C''$,
one then finds that $\nu_{C,C'}(y) = (y_3,y_1,y_2y_3)$, where the
first component is the coordinate on $C \cap C'$. Hence, our
$\sigma_{C,C'}$ near that point is the standard $T^2$-action in the
variables $y_1,y_2$. Similarly, near a point of $C \cap C' \cap X_0$,
$\nu_{C,C'}(y) = (y_1,y_2,y_3/y_1)$, so that one ends up with the
$T^2$-action in the variables $y_2,y_3$. These obviously have the
desired properties. \qed

We now pick an ordering of every pair $\{C,C'\}$ of components, and
will use only the $T^2$-actions $\sigma_{C,C'}$ in that order.

\begin{lemma} \label{th:kaehler-4}
There is a symplectic form $\o_X'''$ which tames the complex
structure, such that for each intersection $C \cap C'$ of two
components of $X_\infty$, there is a neighbourhood $U_{C,C'}$ such
that the local $T^2$-action $\sigma_{C,C'}$ leaves $\o_X'''$ invariant.
Moreover, $\o_X'''$ lies in the same cohomology class as $\o_X''$,
and agrees with it at each point $x \in C \cap C'$, as well as in a
neighbourhood of the lowest stratum points.
\end{lemma}

\proof Let $\bar{\o}_{C,C'} \in \Omega^2(U_{C,C'})$ be the result of
averaging $\o_X''$ with respect to $\sigma_{C,C'}$. This is a closed
two-form; it agrees with $\o_X''$ along $C \cap C'$, and also in a
neighbourhood of the lowest stratum points. Using the Poincar{\'e}
lemma one can write $\bar{\o}_{C,C'} - \o_X'' = d\beta_{C,C'}$, where
$\beta_{C,C'}$ is a one-form that vanishes to second order along $C
\cap C'$, and which is again zero in a neighbourhood of the lowest
stratum point. Take a compactly supported cutoff function
$\psi_{C,C'}$ on $U_{C,C'}$, which is equal to one near $C \cap C'$.
One sees easily that as one rescales the normal directions to make
the support of $\psi_{C,C'}$ smaller, $d(\psi_{C,C'} \beta_{C,C'})$
becomes arbitrarily small in $C^0$-sense. It follows that for
sufficiently small support, $\o_X'' + d(\psi_{C,C'} \beta_{C,C'})$ is
a symplectic form which tames the complex structure. This agrees with
$\o_X''$ near the lowest stratum points, and for that reason, one can
carry out the construction simultaneously for all intersections $C
\cap C'$. \qed

\begin{lemma}
Let $x \in C \cap C'$ be a point where two components of $X_\infty$
meet, which is not a lowest stratum point. Then there are local
$\sigma_{C,C'}$-equivariant Darboux coordinates near $x$, in which $X$
with the map $1/\pi$, the symplectic structure $\o_X'''$, the given
complex structure $J_X$, and the meromorphic complex volume form
$\eta_X$, is isomorphic to $(Y,p,\o_Y,J_Y,\eta_Y)$ as in Assumptions
\ref{as:1prime}. \qed
\end{lemma}

This is straightforward to check, and we leave it to the reader. The
outcome of the whole construction is this:

\begin{prop} \label{th:negativity}
Let $\gamma_\infty: [0,2\pi] \rightarrow \C$ be a circle of large
radius $R \gg 0$. Then the graded monodromy
\begin{equation}
\hgr_{\gamma_\infty} \in \Autgr(M_R)
\end{equation}
is isotopic (within that group) to a graded symplectic automorphism
$\phigr$ with the following property. There is a finite set $\Sigma
\subset X_z \cap X_\infty$ which admits arbitrarily small open
neighbourhoods $W \subset X_z$, such that $\phi(W) = W$, and
$\phigr|(X_z \setminus W)$ is negative.
\end{prop}

\proof The main point is to replace the given K{\"a}hler form $\o_X$
by $\o_X'''$, and we need to convince ourselves that this is
permitted. For each regular $X_z$ there is a symplectic isomorphism
\begin{equation}
\phi_z: (X_z,\o_{X_z}) \longrightarrow (X_z,\o_X'''|X_z)
\end{equation}
which carries $X_{0,\infty}$ to itself. This follows from Moser's
Lemma together with the isotopy theorem for symplectic surfaces with
orthogonal normal crossings, see \cite[Theorem 6.5]{ruan02b} for a
more general statement valid in all dimensions. Next, an inspection
of the construction of monodromy maps for $\pi_M$ shows that they can
be defined using any symplectic form which tames the complex
structure, and in particular $\o_X'''$. The resulting diagram of
symplectic maps
\begin{equation}
\xymatrix{
 (X_R,\o_{X_R}) \ar[d]^-{h_{\gamma_\infty}}
 \ar[r]^-{\phi_R} &
 (X_R,\o_X'''|X_R) \ar[d]^-{h_{\gamma_\infty}'''} \\
 (X_R,\o_{X_R}) \ar[r]^-{\phi_R} & (X_R,\o_X'''|X_R)
}
\end{equation}
commutes up to symplectic isotopy rel $X_{R,\infty}$. Since $X_R$ is
a $K3$ surface, Remark \ref{th:h1} shows that $h_{\gamma_\infty}$ is
isotopic to $\phi_R^{-1} \circ h_{\gamma_\infty}''' \circ \phi_R$
within $Aut(M_R)$. Moreover, this is compatible with gradings.
%

Let $\Sigma$ be the set of those points where two components of
$X_\infty$ intersect each other and $X_0$. This lies in the base
locus of our pencil, hence in every fibre $X_z$. Around each point $x
\in \Sigma$, choose a neighbourhood $U_x \subset X$ which is an open
ball around the origin in local coordinates as in Assumptions
\ref{as:2}, and set $U = \bigcup_x U_x$. We know from our discussion
of the local model that $U$ is invariant under parallel transport
along large circles. The claim is that for sufficiently large $R$,
the grading of
\begin{equation}
\psigr_R = (\hgr_{\gamma_\infty}''')^2
\end{equation}
(with respect to $\o_X'''$, the given complex structure and complex
volume form $\eta_{X_R}$) is negative on $X_R \setminus U$. This is
proved by contradiction: suppose that we have a sequence $R_k
\rightarrow \infty$, and points $x_k \in X_{R_k} \setminus U$ such
that the grading of $\psigr_R$ is not negative on some Lagrangian
subspace in the tangent space at $x_k$. Consider the limit point of a
subsequence, $x_\infty \in X_\infty \setminus U$. If $x_\infty$ is a
smooth point of $X_\infty$, the grading of $\psigr_R$ at $x_{z_k}$
converges to $-2$ for the same reason as in the case of Lefschetz
pencils, see Lemma \ref{th:shift-factor}. If $x_\infty$ is a point
where several components of $X_\infty$ meet, Lemma
\ref{th:local-model-1} says that the grading of $\psigr_R$ at
$x_{z_k}$ will be bounded above by $-1/2$ for large $k$ (when
applying this local result, the reader should keep in mind the change
of coordinates on the base, $\zeta = 1/z$ and $dz = d\zeta/\zeta^2$).

Fix some $R$ such that the claim made above applies. Take a
neighbourhood $W \subset U \cap X_R$ of $\Sigma$. After making $W$
smaller, we may again assume that in the local coordinates of
Assumption \ref{as:2}, this is the intersection of the fibre
$Y_{1/R}$ with an open ball around the origin. Lemma
\ref{th:bad-points} tells us that $h_{\gamma_\infty}'''(W) = W$, and
that $h_{\gamma_\infty}'''$ is negative on $U \setminus W$. \qed

\section{Fukaya categories\label{sec:fukaya}}

There are Fukaya categories for affine and projective Calabi-Yaus, as
well as a relative version which interpolates between the two. All of
them are c-unital $A_\infty$-categories, but with different
coefficient rings:
\begin{center}
\begin{tabular}{l| l | l |l}
 & notation for & & \\
 & category & notation for objects & coefficient ring \\
 \hline
 affine & $\Fuk(M)$ &
 $L^\br = (L,\alphagr_L,\Spin_L) $ & $\C$ \\
 relative & $\Fuk(X,X_\infty)$ &
 $L^\mbr = (L,\alphagr_L,\Spin_L,J_L) $ & $\Lambda_\N$ \\
 projective & $\Fuk(X)$ &
 $L^\rbr= (L,\alphagr_L,\Spin_L,\lambda_L,J_L) $ & $\Lambda_\Q$ \\
\end{tabular}
\end{center}
While our definition of affine Fukaya categories is fairly general,
in the projective and relative setups only Calabi-Yau surfaces ($K3$
and tori) will be considered. This restriction allows us to avoid
using virtual fundamental chains, and it also simplifies the relation
between the different categories: up to quasi-equivalence,
$\Fuk(X,X_\infty)$ reduces to $\Fuk(M)$ if the deformation parameter
$q$ is set to zero, while tensoring with $\Lambda_\Q$ gives a full
subcategory of $\Fuk(X)$ (for the definition of $\Fuk(X)$ in higher
dimensions, see \cite{fooo}; and for its presumed relation with
$\Fuk(M)$, \cite{seidel02}). We emphasize that all properties of
Fukaya categories which we use are quite simple ones, and largely
independent of the finer details of the definition.

\subsection{}
The following is a summary of \cite[Chapter 2]{seidel04}, with minor
changes in notation and degree of generality. Let $M = X \setminus
X_\infty$ be an affine Calabi-Yau manifold. Objects of $\Fuk(M)$ are {\em
exact Lagrangian branes}
\begin{equation} \label{eq:exact-brane}
L^\br = (L,\alphagr_L,\Spin_L),
\end{equation}
where $L \subset M$ is an oriented closed exact Lagrangian submanifold,
$\alphagr_L$ a grading which is compatible with the orientation (see
Remark \ref{th:or}), and $\Spin_L$ a $Spin$ structure. To define the
morphism spaces, choose for each pair of objects a {\em Floer datum}
\begin{equation}
(H,J) = (H_{L_0^\br,L_1^\br},J_{L_0^\br,L_1^\br}).
\end{equation}
This consists of a function $H \in \smooth_c([0,1] \times M,\R)$, as
well as a family $J = J_t$ of $\o_X$-compatible almost complex
structures parametrized by $t \in [0,1]$, which for all $t$ agree
with the given complex structure in a neighbourhood of $X_\infty$.
The condition on $H$ is that the associated Hamiltonian isotopy
$(\phi^t_H)$ should satisfy $\phi^1_H(L_0) \trans L_1$. We then
define $\CC(L_0^\br,L_1^\br)$ to be the set of trajectories $x: [0,1]
\rightarrow M$, $x(t) = \phi^t_H(x(0))$, with boundary conditions
$x(k) \in L_k$. This is obviously bijective to $\phi^1_H(L_0) \cap
L_1$, hence finite. To each such $x$ one can associate a Maslov index
$I(x)$ and orientation set $or(x)$. For $I(x)$, one observes that the
grading of $L_0$ induces one of $\phi^1_H(L_0)$, and then uses the
index of an intersection point as defined in
\eqref{eq:absolute-index}. $or(x)$ is a set with two elements, the
two possible ``coherent orientations'' of $x$; we refer to
\cite{seidel04,fooo} for its definition, which uses the $Spin$
structures on $L_0,L_1$. The Floer cochain space is the graded
$\C$-vector space
\begin{equation} \label{eq:basic-floer-complex}
\mathit{CF}^j_M(L_0^\br,L_1^\br) = \bigoplus_{\substack{x \in
\CC(L_0^\br,L_1^\br)
\\ I(x) = j}} \C_x,
\end{equation}
where $\C_x$ is defined as the quotient of the two-dimensional space
$\C[or(x)]$ by the relation that the two elements of $or(x)$ must sum
to zero (it is of course isomorphic to $\C$, but not canonically). In
the Fukaya category, $\mathit{CF}^*_M$ is the space of $\mathit{hom}$s from $L_0^\br$
to $L_1^\br$.

Given $x_0,x_1 \in \CC(L_0^\br,L_1^\br)$, we now consider the moduli
spaces $\MM^1_M(x_0,x_1)$ of solutions $u: \R \times [0,1]
\rightarrow M$ to Floer's ``gradient flow'' equation with limits
$x_0,x_1$, divided by the $\R$-action by translation. As part of the
definition of Floer datum, we require that these spaces should all be
regular. Assume from now on that $I(x_1) = I(x_0)-1$. By an index
formula $\MM^1_M(x_0,x_1)$ is zero-dimensional, and a compactness
theorem says that it is a finite set. A few words about this: since
the action functional is exact \eqref{eq:energy} and there is no
bubbling, we only have to worry about sequences of solutions which
become increasingly close to $X_\infty$. To see that this does not
happen, one uses a maximum principle argument based on the function
$-\log\,\|\sigma_{X,\infty}\|^2$, which goes to $+\infty$ at
$X_\infty$. It is plurisubharmonic since
\begin{equation} \label{eq:log-sigma}
\partial \bar\partial \log\,\|\sigma_{X,\infty}\|^2 = -F_{\nabla_X} =
2\pi i \o_M.
\end{equation}
Each $u \in \MM^1_M(x_0,x_1)$ comes with a preferred identification
$or(x_1) \iso or(x_0)$, which in turn defines a map $\C_{x_1}
\rightarrow \C_{x_0}$. The sum of these over all $u$, multiplied by $(-1)^j$, yields the Floer
differential $\mu^1_M: \mathit{CF}^j_M(L_0^\br,L_1^\br) \rightarrow
\mathit{CF}^{j+1}_M(L_0^\br,L_1^\br)$, which is the first composition map in
the Fukaya category. Its cohomology is the Floer cohomology
$\mathit{HF}^*_M(L_0^\br,L_1^\br)$, which is therefore the space of morphisms
in the cohomological category $H\Fuk(M)$.

To define the rest of the $A_\infty$-structure, one proceeds as
follows. Take some $d \geq 2$. A {\em $(d+1)$-pointed disc with
Lagrangian labels}
\begin{equation} \label{eq:s}
 (S,\zeta_0,\dots,\zeta_d,L_0^\br,\dots,L_d^\br)
\end{equation}
is a Riemann surface $S$ isomorphic to the closed disc with $d+1$
boundary points removed. Moreover, these points at infinity should
have a preferred numbering $\zeta_0,\dots,\zeta_d$, compatible with
their natural cyclic order. Denote the connected components of
$\partial S$ by $I_0,\dots,I_d$, again in cyclic order and starting
with the component $I_0$ which lies between $\zeta_0$ and $\zeta_1$.
Finally, we want to have exact Lagrangian branes
$L_0^\br,\dots,L_d^\br$ attached to the boundary components. A set of
{\em strip-like ends} for $S$ consists of proper holomorphic
embeddings
\begin{align*}
 \epsilon_{S,0} & : \Rleq \times [0,1] \longrightarrow S, &&
 \textstyle\lim_{s \rightarrow -\infty} \epsilon_{S,0}(s,\cdot) = \zeta_0, \\
 \epsilon_{S,1},\dots,\epsilon_{S,d} & : \Rgeq \times [0,1] \longrightarrow S, &&
 \textstyle\lim_{s \rightarrow +\infty} \epsilon_{S,k}(s,\cdot) = \zeta_k
\end{align*}
taking $\Rleq \times \{0;1\}$, $\Rgeq \times \{0;1\}$ to $\partial S$
and which have disjoint images. A {\em perturbation datum} on $S$ is
a pair $(K_S,J_S)$ consisting of a family $\{J_{S,z}\}_{z \in S}$ of
compatible almost complex structures on $X$, which agree with the
given complex structure in a neighbourhood of $X_\infty$, and a
one-form $K_S \in \Omega^1(S,\smooth_c(M,\R))$ with values in smooth
functions, subject to the condition that $K_S(\xi)|L_k = 0$ for $\xi
\in T(I_k) \subset T(\partial S)$. The behaviour of $(K_S,J_S)$ over
the strip-like ends is fixed by the previously chosen Floer data:
\begin{equation}
 J_{S,\epsilon_{S,k}(s,t)} = \begin{cases}
 J_{L_0^\br,L_d^\br,t} & k = 0, \\
 J_{L_{k-1}^\br,L_k^\br,t} & k > 0
 \end{cases}
 \quad \text{and} \quad
 \epsilon_{S,k}^*K_S = \begin{cases}
 H_{L_0^\br,L_d^\br,t}\, dt & k = 0, \\
 H_{L_{k-1}^\br,L_k^\br,t} \, dt & k > 0.
 \end{cases}
\end{equation}
Suppose that we have made a global choice of strip-like ends and
perturbation data. This means that for every $(d+1)$-pointed disc
with Lagrangian labels we have chosen a set of strip-like ends
$(K_S,J_S)$ (the dependence on the other data in \eqref{eq:s} is
suppressed from the notation for the sake of brevity). This needs to
satisfy two additional kinds of conditions: first of all, it needs to
vary smoothly if we change the complex structure on $S$ in a smooth
way, which means that it is defined on the universal family of such
discs. Secondly, it needs to be well-behaved with respect to the
compactification of that family, which includes degenerations to
discs with nodes. This in fact establishes recursive relations
between the choices of perturbation data for different $d$, see
\cite[Section 9]{seidel04} for details.

Supposing that all these conditions have been met, take $x_0 \in
\CC(L_0^\br,L_d^\br)$ as well as $x_k \in \CC(L_{k-1}^\br,L_k^\br)$
for $k = 1,\dots,d$. Define $\MM^d_M(x_0,\dots,x_d)$ to be the space
of equivalence classes of pairs consisting of a $(d+1)$-marked disc
with labels $L_0^\br,\dots,L_d^\br$, and a solution of the following
generalization of Floer's equation for a map $u: S \rightarrow M$:
\begin{equation} \label{eq:generalized-floer}
\begin{cases}
 \!\!\! & u(I_k) \subset L_k, \\
 \!\!\! & (du(z)-Y_{S,z,u(z)})
 + J_{S,z,u(z)} \circ (du(z)-Y_{S,z,u(z)}) \circ i = 0, \\
 \!\!\! & \lim_{s \rightarrow \pm \infty} u(\epsilon_{S,k}(s,\cdot)) = x_k.
\end{cases}
\end{equation}
Here $Y_S \in \Omega^1(S,\smooth_c(TM))$ is the
Hamiltonian-vector-field valued one-form determined by $K_S$, and $i$ is the complex structure on $S$. The
equivalence relation on pairs is isomorphism of the underlying
Riemann surfaces, compatible with all the additional structure and commuting with the maps $u$. The resulting quotient
space is (locally) the zero-set of a Fredholm section of a suitably
defined Banach vector bundle. Generically, this will be transverse,
and then $\MM^d_M(x_0,\dots,x_d)$ is smooth of dimension $I(x_0) -
I(x_1) - \dots - I(x_d) + d - 2$. Moreover, the zero-dimensional
spaces are finite. It is maybe worthwhile stating the basic energy
equality which underlies the compactness argument. Choose primitives
$K_{L_k}$ for $\theta_M|L_k$. Define the action of the limits $x_k$ by
a generalization of \eqref{eq:action},
\begin{align*}
 A(x_k) & = - \textstyle \int x_k^*\theta_M + \\
 + & \begin{cases}
 \int H_{L_0^\br,L_d^\br}(x_0(t))\, dt +
 K_{L_d}(x_0(1)) - K_{L_0}(x_0(0)) & k = 0, \\
 \int H_{L_{k-1}^\br,L_k^\br}(x_k(t))\,dt +
 K_{L_k}(x_k(1)) - K_{L_{k-1}}(x_k(0)) & k = 1,\dots,d.
 \end{cases}
\end{align*}
Then for $u \in \MM^d_M(x_0,\dots,x_d)$,
\begin{equation} \label{eq:polygon-energy}
\begin{aligned}
 E(u) & \stackrel{\mathrm{def}}{=} \half \textstyle \int_S | du - Y |^2 \\ & =
 A(x_0) - \sum_{k=1}^d A(x_d) + \textstyle \int_S u^*(dK_S + \half\{K_S,K_S\}).
\end{aligned}
\end{equation}
The integral on the right hand side can be bounded by a constant
which is independent of $u$, because the curvature term $dK_S +
\half\{K_S,K_S\}$ vanishes over the strip-like ends of $S$ (and also outside a compact subset of $M$). Of course \eqref{eq:log-sigma} also plays a role in the proof of
compactness, as was the case for the Floer differential. To continue
with the main story: any element of a zero-dimensional space
$\MM^d_M(x_0,\dots,x_d)$ defines a canonical bijection
\begin{equation} \label{eq:orientation-mechanism}
\mathit{or}(x_0) \iso \mathit{or}(x_1) \times_{\Z/2} \dots \times_{\Z/2} \mathit{or}(x_d).
\end{equation}
The order $d$ composition map $\mu^d_M$ in $\Fuk(M)$ is the sum of
the corresponding maps $\C_{x_d} \otimes \dots \otimes \C_{x_1}
\rightarrow \C_{x_0}$, multiplied by $(-1)^{I(x_1) + 2I(x_2) + \cdots + dI(x_d)}$. The
additional signs involving the Maslov indices $I(x_k)$ are arranged in such a way that (if suitable conventions are used in setting up \eqref{eq:orientation-mechanism}, see \cite[Section 12]{seidel04}) the $A_\infty$-structure equations will hold with the signs used here. $\Fuk(M)$ is independent of all the choices (Floer data,
strip-like ends, perturbation data) made in its construction up to
quasi-equivalence, and indeed quasi-isomorphism acting trivially on
the objects.

\begin{remark} \label{th:class-of-j}
Suppose that $X_\infty$ does not contain rational curves (if $X$ is an algebraic surface, this means that no irreducible component of $X_\infty$ is rational). Then
one can use a wider class of almost complex structures, namely all
those $\o_X$-compatible ones such that each irreducible component $C \subset X_\infty$ is an almost complex submanifold. Similarly, in the inhomogeneous terms one can use Hamiltonian
functions on $X$ whose values and first derivatives vanish on
$X_\infty$. Instead of using \eqref{eq:log-sigma} to prevent solutions of Floer's equation and its
generalization from escaping into $X_\infty$, one now argues as
follows: if the solution stays inside $M$, its intersection number
with $C$ is zero. The limit of a sequence of solutions has the same
property, and nonnegativity of intersection multiplicities shows that
the limit is disjoint from $X_\infty$ (more precisely, one thinks of
a solution of \eqref{eq:generalized-floer} as a map $u: S \rightarrow
S \times X$ which is pseudo-holomorphic for a suitable almost complex
structure, and considers its intersection with $S \times X_\infty$).
\end{remark}

\subsection{\label{subsec:coverings}}
Let $\bar{M}$ be an affine Calabi-Yau manifold. Suppose that we have a homomorphism $\rho: \pi_1(\bar{M}) \rightarrow \Gamma$, where $\Gamma$ is a finite abelian group, and its associated finite covering $M \rightarrow \bar{M}$. For any object $\bar{L}^\br$ of $\Fuk(\bar{M})$ such that the composition
\begin{equation} \label{eq:pi-vanishes}
\pi_1(\bar L) \longrightarrow \pi_1(\bar M) \longrightarrow \Gamma
\end{equation}
is trivial, choose a lift $L \subset M$. Once one has done that, the morphism spaces split as
\begin{equation}
 \mathit{CF}^*_{\bar M}(\bar L_0^\br,\bar L_1^\br) =
 \bigoplus_{\gamma \in \Gamma} \mathit{CF}^*_{\bar M}(\bar L_0^\br,\bar L_1^\br)_\gamma,
\end{equation}
where the $\gamma$-summand consists of those $\bar x \in \CC(\bar L_0^\br,\bar L_1^\br)$ whose unique lift $x: [0,1] \rightarrow M$ with $x(0) \in L_0$ satisfies $\gamma^{-1}(x(1)) \in L_1$. Because each pseudo-holomorphic polygon can be lifted to $M$, these
splittings are compatible with the compositions. This means that if we have objects $\bar L_0^\br,\dots,\bar L_d^\br$ as before and choices of lifts, then
\begin{equation}
 \mathit{CF}^*_{\bar M}(\bar L_{d-1}^\br,\bar L_d^\br)_{\gamma_d} \otimes \dots \otimes
 \mathit{CF}^*_{\bar M}(\bar L_0^\br,\bar L_1^\br)_{\gamma_1} \xrightarrow{\mu^d_{\Fuk(\bar M)}}
 \mathit{CF}^*_{\bar M}(\bar L_0^\br,\bar L_d^\br)_{\gamma_1\dots \gamma_d}.
\end{equation}
For general nonsense reasons, this constitutes an action of the dual group $\Gamma^* = \mathit{Hom}(\Gamma,\C^*)$ on the full subcategory $\bar\Cat \subset \Fuk(\bar M)$ of objects such that \eqref{eq:pi-vanishes} is trivial. Note that even though we have not represented the cover $M$ as an affine Calabi-Yau manifold, one can define $\Fuk(M)$ as before, using a maximum principle argument (or projection to $\bar{M}$) to show that pseudo-holomorphic maps remain within a bounded subset. The objects $\gamma(L^\br)$, for all $\gamma \in \Gamma$, then form a full subcategory $\Cat \subset \Fuk(M)$, and one has
\begin{equation} \label{eq:semidirect-fukaya}
\Cat \iso \bar\Cat \semidirect \Gamma^*.
\end{equation}
We know that $\bar\Cat \semidirect \Gamma^*$ comes with an algebraically defined action of $\Gamma$, and this can be identified with the natural action of $\Gamma$ by covering transformations on the left hand side of \eqref{eq:semidirect-fukaya}.

\begin{remark} \label{th:fg}
One can also allow a finitely generated, but not necessarily finite, abelian group $\Gamma$. For instance, if $\Gamma = \Z$ then $\Gamma^* = \C^*$. In this case, it may actually be more intuitive to think of the action of $\Gamma^*$ on $\bar\Cat$ as a coaction of $\Gamma$, which means as an additional grading on morphism spaces, which is homogeneous for the $A_\infty$-compositions.
\end{remark}

One can also approach the construction from the reverse point of view. Namely, start with an affine Calabi-Yau manifold $M$ together with a free action of $\Gamma$ by K{\"a}hler isometries. These should preserve the one-form $\theta_M$ as well as the holomorphic volume form, which ensures that these structures descend to $\bar{M} = M/\Gamma$. One starts by defining $\Cat \subset \Fuk(M)$ to be the full subcategory consisting of those Lagrangian submanifolds $L$ such that
\begin{equation} \label{eq:free-l}
\gamma(L) \cap L = \emptyset \text{ for all $\gamma \neq e$,}
\end{equation}
and then takes their images to be objects of $\bar\Cat \subset \Fuk(\bar M)$, which again leads to the relationship \eqref{eq:semidirect-fukaya}.

We will need also need another slightly more tricky construction. Take $M$ and an action of $\Gamma$ as before, but where we drop the assumption that the action is free. Suppose that we have a collection of objects $(L_i^\br)$ of $\Fuk(M)$ indexed by some set $I$, which comes equipped with an action of $\Gamma$ compatible with that on our objects. Importantly, the Lagrangian submanifolds do not have to satisfy \eqref{eq:free-l}, and in fact the action on the indexing set $I$ does not have to be free. If $C$ is the cohomology level Fukaya category with objects $(L_i^\br)$, then $\Gamma$ acts on it in a canonical way. What we want to have is an underlying equivariant cochain level category $\Cat$.

\begin{remark} \label{th:spin-group-action}
It may be worthwhile to spell out part of the assumptions in more detail. What we are given is an action of $\Gamma$ on the indexing set $I$ for our collection of objects, so that $\gamma(L_i) = L_{\gamma(i)}$ and the gradings are preserved. Moreover, the diffeomorphisms $\gamma|L_i$ must come with preferred isomorphisms of {\em Spin} structures, compatible with the group law. The last-mentioned condition is definitely an additional restriction (even for a single manifold which has a group action and is {\em Spin}, there are obstructions to making it {\em equivariantly Spin}). However, one case which is unproblematic is the following one:
\begin{itemize} \itemsep1em
\item $\Gamma$ is cyclic.
\item if $\gamma(i) = i$ for some $\gamma \in \Gamma$ and $i \in I$, then $\gamma|L_i$ is the identity.
\end{itemize}
Then the indexing set can be decomposed into orbits, each of which is a free orbit for some (still cyclic) quotient of $\Gamma$. One can adjust the isomorphisms of {\em Spin} structures at one step of each orbit so as to be compatible with the group action.
\end{remark}

In principle, one can view the issue of constructing $\Cat$ as an equivariant transversality problem, to which multivalued perturbation theory can be applied. However, we can handle the necessary ``averaging'' purely algebraically, as follows. First, note that if $\Gamma$ acted freely on $I$ there would be no problem: in that case, the action of a nontrivial $\gamma$ relates one moduli space (of holomorphic maps with boundary on $L_{i_0},\dots,L_{i_d}$) to a different one (of maps with boundary on $L_{\gamma(i_0)},\dots,L_{\gamma(i_d)}$). The perturbation data for those two moduli spaces can be chosen independently, so there are no issues of equivariant transversality. Note that this holds even if the branes associated to two elements of $I$ have the same underlying Lagrangian submanifolds, as long as we distinguish between them as objects, hence allow different choices of perturbation data.

This leads to the following formal trick for the general case. Given a collection $\{L_i\}$, let's introduce the bigger indexing set $\tilde{I} = I \times \Gamma$, and the collection $L_{\tilde{i}}^\br = L_i^\br$ for $\tilde{i} = (i,\gamma)$. The action of $\Gamma$ on $\tilde{I}$ (by left multiplication on the second factor) is always free, so we can define an equivariant version of the Fukaya category $\Dat$ with objects $L_{\tilde{i}}^\br$, and carrying an action of $\Gamma$. On the cohomology level, $D = H(\Dat)$ comes with a functor $F: D \rightarrow C$ which forgets $\gamma$, and which is an equivalence compatible with the group actions. Hence, one can apply Lemma \ref{th:transfer-group-action} to produce the desired $\Cat$.

\subsection{}
Let $X$ be a projective Calabi-Yau surface (real dimension 4), and $L \subset X$ a
Lagrangian submanifold which admits a grading. An $\o_X$-compatible
almost complex structure $J$ is called {\em regular with respect to
$L$} if there are no non-constant $J$-holomorphic spheres, and no
non-constant $J$-holomorphic discs with boundary on $L$.

\begin{lemma}
Almost complex structures which are regular with respect to a given
$L$ are dense within the set of all $\o_X$-compatible almost complex
structures.
\end{lemma}

\proof Since $c_1(X) = 0$, the virtual dimension of the space of
unparametrized $J$-holomorphic spheres at a point $v$ is $4 + 2
\leftsc c_1(X), [v] \rightsc - 6 = -2$. Standard transversality
results tell us that for generic $J$, there are no $J$-holomorphic
spheres that are somewhere injective; and since any non-constant
sphere is a multiple cover of a somewhere injective one, there are no
non-constant spheres at all. The argument for discs goes along the
same lines: because of the grading, the Maslov number $\mu(w)$ of any
disc in $(X,L)$ vanishes, and the virtual dimension is
\begin{equation} \label{eq:minus-one}
2 + \mu(w) - 3 = -1.
\end{equation}
Again, this shows that for generic $J$, there are no somewhere
injective discs. The decomposition theorem of Kwon-Oh
\cite{kwon-oh96} and Lazzarini \cite{lazzarini98} says that if there
is a non-constant $J$-holomorphic disc, there is also a somewhere
injective one, which implies the desired result. \qed

Objects of $\Fuk(X)$ are {\em rational Lagrangian branes}
\begin{equation} \label{eq:rational-brane}
L^\rbr = (L,\alphagr_L,\Spin_L,\lambda_L,J_L)
\end{equation}
where $L$ is an oriented Lagrangian submanifold of $X$ with a grading
$\alphagr_L$ and $Spin$ structure $\Spin_L$; $\lambda_L$ is a
covariantly constant multisection of the circle bundle of $o_X|L$ of
some degree $d_L \geq 1$; and $J_L$ is an $\o_X$-compatible almost
complex structure which is regular with respect to $L$. For any pair
of objects, we take a Floer datum
\begin{equation}
(H,J) = (H_{L_0^\rbr,L_1^\rbr},J_{L_0^\rbr,L_1^\rbr})
\end{equation}
as before, except that the conditions on the behaviour of $H$ and $J$
near $X_\infty$ are omitted, and the following new requirements
imposed instead: $J_t$ must be equal to $J_{L_0}$, $J_{L_1}$ for $t =
0,1$; and there should be no non-constant $J_t$-holomorphic spheres
for any $t$ (this is still true generically, because the virtual
dimension of the parametrized moduli space is $-1$). Define the set
$\CC(L_0^\rbr,L_1^\rbr)$ as before; and for any $x: [0,1] \rightarrow
X$ in that set, write $\bar{A}(x) = \bar{A}(x(1))$ for the mod $\Q$
action \eqref{eq:modq-action} of $x(1)$ as an intersection point of
$\phi^1_H(L_0) \cap L_1$ (where we have used Remark \ref{th:deform-rational}
to equip $\phi^1_H(L_0)$ with a flat multisection of $o_X$).
For any number $r \in \R/\Q$, we define a one-dimensional $\Lambda_\Q$-vector
space $q^r\Lambda_{\Q}$ whose elements are formal series $f(q) =
\sum_m a_m q^m$, with $m$ running over all real numbers in $r +
\frac{1}{d}\Z$ for some natural number $d$ which depends on $f$, and
$a_m \in \C$ vanishing for sufficiently negative $m$. Write
$\Lambda_{\Q,x}$ for the quotient of
$q^{\bar{A}(x)}\Lambda_{\Q}[\mathit{or}(x)]$ by the relation that the two
elements of $\mathit{or}(x)$ sum to zero. Then the Floer cochain space is
\begin{equation}
\mathit{CF}^j_X(L_0^\rbr,L_1^\rbr) = \bigoplus_{\substack{x \in
\CC(L_0^\rbr,L_1^\rbr)
\\ I(x) = j}} \Lambda_{\Q,x}.
\end{equation}
Denote by $\MM^1_X(x_0,x_1)$ the space of solutions of Floer's
equation in $X$ with limits $x_0,x_1$ (this means we drop the
previous assumption that the solutions stay in $M$). As before, we
require that this is regular, and in particular zero-dimensional in
the case $I(x_1) = I(x_0)-1$. Take a solution $u$ in such a
zero-dimensional moduli space, and let $E(u) = \int \|\partial_s u\|^2 \, ds\, dt > 0$ be its energy. A version of \eqref{eq:modq-energy} shows that
$E(u) \in \bar{A}(x_0) - \bar{A}(x_1) + \lcm(d_{L_0},d_{L_1})^{-1}
\Z$, so that we get a well-defined $\Lambda_\Q$-module map
\begin{equation}
\pm q^{E(u)} : \Lambda_{\Q,x_1} \longrightarrow \Lambda_{\Q,x_0}.
\end{equation}
The sign of this is determined by the induced bijection of
orientation sets, as before. Summing over all $u$ and $x_0$ yields $(-1)^{I(x_1)}$ times the
Floer differential $\mu^1_X(x_1)$ in $\Fuk(X)$. The sum is no longer
finite, but it gives a well-defined $\Lambda_\Q$-module map by Gromov
compactness; the basic point being that there is no bubbling off of
holomorphic discs or spheres, by assumption on the Floer datum.
Denote the resulting Floer cohomology by $\mathit{HF}^*_X(L_0^\rbr,L_1^\rbr)$.

The definition of the higher order compositions is largely similar.
One may use only perturbation data $(K_S,J_S)$ where $J_S = J_{L_k}$
over the boundary component $I_k \subset \partial S$. It is no longer
possible to avoid $J_{S,z}$-holomorphic spheres for some values of $z
\in \mathit{int}(S)$, but the argument from \cite{hofer-salamon95} shows that
for a generic choice of $J_S$, these spheres will never occur as
bubbles in the Gromov compactification of the spaces
$\MM^d_X(x_0,\dots,x_d)$ as long as these are of dimension $\leq 1$,
which is all that we will ever use. Equation
\eqref{eq:polygon-energy} still holds if we replace all the actions
by their mod $\Q$ counterparts. We define the contribution of a point
$(S,u)$ in a zero-dimensional moduli space $\MM^d_X(x_0,\dots,x_d)$
to $\mu^d_X$ to be
\begin{equation} \label{eq:qe}
 \pm q^{E(u) - \int_S u^*(dK_S + \half\{K_S,K_S\})} :
 \Lambda_{\Q,x_d} \otimes \dots \otimes \Lambda_{\Q,x_1}
 \longrightarrow \Lambda_{\Q,x_0},
\end{equation}
with the sign again depending on \eqref{eq:orientation-mechanism}.
There is a uniform bound for the $\int_S$ term. This, together
with Gromov compactness, ensures that the sum is a well-defined
$\Lambda_\Q$-multilinear map.

\begin{remark}
The choice of $\lambda_L$ serves only to fix the action: the Floer
cohomology groups for different choices are isomorphic in an obvious
way, by multiplication with $q^m$ for some $m \in \Q$. In particular,
changing $\lambda_L$ does not affect the quasi-isomorphism class of
the resulting object in $\Fuk(X)$.
\end{remark}

\subsection{\label{subsec:relative-fukaya}}
Take $X$ as before, but now equipped with a section
$\sigma_{X,\infty}$ of $o_X$ whose zero set $X_\infty$ is a reduced
divisor with normal crossings, so that $M = X \setminus X_\infty$ is
an affine Calabi-Yau surface. Following Remark \ref{th:class-of-j},
we also assume that $X_\infty$ has no rational components. The
objects of the relative Fukaya category $\Fuk(X,X_\infty)$ are a
mixture of the two previously introduced classes,
\eqref{eq:exact-brane} and \eqref{eq:rational-brane}. One considers
{\em relative Lagrangian branes}
\begin{equation} \label{eq:relative-brane}
L^\mbr = (L,\alphagr_L,\Spin_L,J_L)
\end{equation}
where $(L,\alphagr_L,\Spin_L)$ is an exact Lagrangian brane in $M$;
and $J_L$ is an $\o_X$-compa\-tible almost complex structure on $X$
such that each component of $X_\infty$ is an almost complex
submanifold, and which is regular with respect to $L$. The restriction
on $J_L|X_\infty$ does not interfere with the regularity theory for
somewhere injective pseudo-holo\-morphic spheres $v: S^2 \rightarrow
X$, since $v^{-1}(X_\infty)$ is necessarily a finite set, while the
set of somewhere injective points is open; and the same applies to
pseudo-holomorphic discs. For any two objects, choose a Floer datum
\begin{equation} \label{eq:relative-floer-datum}
(H,J) = (H_{L_0^\mbr,L_1^\mbr},J_{L_0^\mbr,L_1^\mbr})
\end{equation}
such that $H \in \smooth([0,1] \times X,\R)$ vanishes along $[0,1]
\times X_\infty$, along with its first derivative; and $J_t|X_\infty$
makes the components of $X_\infty$ almost complex. In addition, we
assume that $J_k = J_{L_k}$ for $k = 0,1$, and that there are no
non-constant $J_t$-holomorphic spheres for any $t$.

We have the usual set $\CC(L_0^\mbr,L_1^\mbr)$ of $H$-trajectories
joining $L_0$ to $L_1$ (these are all contained in $M$, since the Hamiltonian
vector field of $H$ vanishes on $X_\infty$), and for $x$ in that set we define
$\Lambda_{\N,x}$ in the same way as $\C_x$, just using $\Lambda_\N$
instead of $\C$ as the ground ring. The Floer cochain space in
$\Fuk(X,X_\infty)$ is
\begin{equation}
 \mathit{CF}_{X,X_\infty}^j(L_0^\mbr,L_1^\mbr) =
 \bigoplus_{\substack{x \in \CC(L_0^\mbr,L_1^\mbr)
 \\ I(x) = j}} \Lambda_{\N,x}.
\end{equation}
For the Floer differential, each $u \in \MM^1_X(x_0,x_1)$ contributes
with $\pm q^{u \cdot X_\infty}: \Lambda_{\N,x_1} \rightarrow
\Lambda_{\N,x_0}$. This makes sense for a single $u$ since the number
$u \cdot X_\infty$ is nonnegative; a version of \eqref{eq:energy} and Gromov compactness tell us that the sum over all $u$ gives rise to a well-defined map
$\mu^1_{X,X_\infty}$. Denote the resulting Floer cohomology by
$\mathit{HF}^*_{X,X_\infty}(L_0^\mbr,L_1^\mbr)$. The definition of the higher
order compositions works in the same way, and we will therefore omit
the details.

\begin{remark} \label{th:class-of-j-2}
Let $X$, $X'$ be two Calabi-Yau surfaces, with sections
$\sigma_{X,\infty}$, $\sigma_{X',\infty}$ of their ample line bundles
$o_X$, $o_{X'}$ which define reduced divisors with normal crossings,
having no rational components. Suppose that there is a symplectic
isomorphism
\begin{equation}
\phi: (X,X_\infty) \longrightarrow (X',X'_\infty)
\end{equation}
with the following additional properties: $\mathit{Hom}(o_X,\phi^*o_{X'})$ is
the trivial flat line bundle, and over $M = X \setminus X_\infty$,
its covariantly constant sections lie in the same homotopy class of
nowhere zero sections as
$\phi^*\sigma_{X',\infty}/\sigma_{X,\infty}$. Moreover, under the
isomorphism (unique up to homotopy) $\K_X \iso \phi^*\K_{X'}$ induced
by deforming the almost complex structure on $X$ to the pullback of
the one from $X'$, the nowhere zero sections $\eta_X$ and
$\phi^*\eta_{X'}$ are homotopic. Then, by suitably choosing the Floer
and perturbation data, one can define an induced quasi-isomorphism of
$A_\infty$-categories over $\Lambda_\N$,
\begin{equation}
\phi_*: \Fuk(X,X_\infty) \rightarrow \Fuk(X',X'_\infty).
\end{equation}
In the case where $H_1(X) = 0$, all the additional conditions imposed
on $\phi$ are automatically satisfied, for the same reasons as in
Remark \ref{th:h1}.
\end{remark}

\subsection{}
We need to discuss the formal relations between the different
kinds of Fukaya categories introduced above. Let $X$ be a Calabi-Yau surface with a
suitable section $\sigma_{X,\infty}$, so that $\Fuk(M)$, $\Fuk(X)$
and $\Fuk(X,X_\infty)$ are all defined.

\begin{prop} \label{th:q-zero}
Consider $\C$ as a $\Lambda_\N$-module in the obvious way, with $q$
acting trivially. Then the $\C$-linear $A_\infty$-category obtained
from $\Fuk(X,X_\infty)$ by reduction of constants, $\Fuk(X,X_\infty)
\otimes_{\Lambda_\N} \C$, is quasi-equivalent to $\Fuk(M)$.
\end{prop}

\begin{prop} \label{th:invert-q}
Consider $\Lambda_\Q \supset \Lambda_\N$ as a $\Lambda_\N$-module in
the obvious way. Then the $\Lambda_\Q$-linear $A_\infty$-category obtained from
$\Fuk(X,X_\infty)$ by extension of constants, $\Fuk(X,X_\infty)
\otimes_{\Lambda_\N} \Lambda_\Q$, is quasi-equivalent to a full
$A_\infty$-subcategory of $\Fuk(X)$. That subcategory consists of all
those $L^\rbr$ such that $L \cap X_\infty = \emptyset$, $o_X|L$ is
the trivial flat bundle, and its covariantly constant sections lie in
the same homotopy class as $\sigma_{X,\infty}|L$.
\end{prop}

Both statements are essentially obvious. For Proposition
\ref{th:q-zero}, one notices that setting $q = 0$ means that we count
only those solutions of Floer's equation or its generalization which
have zero intersection number with, hence are disjoint from,
$X_\infty$. Then, the only difference between $\Fuk(M)$ and
$\Fuk(X,X_\infty) \otimes_{\Lambda_\N} \C$ is that the latter
category has more objects, consisting of the same exact Lagrangian
brane with different choices of $J_L$. However, all these objects are
quasi-isomorphic, whence the quasi-equivalence. 
For Proposition \ref{th:invert-q} one has to choose a function $K_L$, $dK_L =
\theta_M|L$, for each exact Lagrangian brane in $M$. The associated
object in $\Fuk(X)$ then carries the multisection $\lambda_L$ from
\eqref{eq:k-lambda}, and one identifies
\begin{equation} \label{eq:mq}
 \mathit{CF}^*_{X,X_\infty}(L_0^\mbr,L_1^\mbr) \otimes_{\Lambda_\N}
 \Lambda_\Q \stackrel{\iso}{\longrightarrow}
 \mathit{CF}^*_X(L_0^\rbr,L_1^\rbr)
\end{equation}
by using the map $q^{A(x)}: \Lambda_{\N,x} \otimes_{\Lambda_\N}
\Lambda_{\Q} \rightarrow \Lambda_{\Q,x}$ on each one-dimensional
subspace. The compositions in the two categories are based on the
same moduli spaces, and the fact that the weights used to count
solutions correspond under \eqref{eq:mq} is a consequence of
\eqref{eq:energy} and a suitable generalization. Clearly,
$\Fuk(X,X_\infty) \otimes_{\Lambda_\N} \Lambda_\Q$ is not all of
$\Fuk(X)$, since only the objects which are exact Lagrangian
submanifolds in $M$ occur.

\subsection{\label{subsec:pss}}
While the basic properties of affine Fukaya categories are quite simple, the relative and projective versions reserve some small surprises. These are all covered by the powerful general theory from \cite{fooo}, but we prefer a more elementary approach. The relative situation is the more important one for the arguments later on, and we will consider it primarily, and then add some discussion of the projective case. We work on a fixed affine Calabi-Yau surface $M = X \setminus X_\infty$, which is such that $X_\infty$ has no rational components.

As a concrete consequence of Proposition \ref{th:q-zero}, we have the following. Let $L_0^\mbr,L_1^\mbr$ be two objects of $\Fuk(X,X_\infty)$, and $L_0^\br,L_1^\br$ the associated objects of $\Fuk(M)$. The $q$-adic filtration gives rise to a spectral sequence converging to $\mathit{HF}^*_{X,X_\infty}(L_0^\mbr,L_1^\mbr)$, and whose starting term is
\begin{equation} \label{eq:q-spectral-e1}
E_1^{rs} = \begin{cases} \mathit{HF}^{r+s}_M(L_0^\br,L_1^\br) & r \geq 0, \\
0 & \text{otherwise.} 
\end{cases}
\end{equation}
The differentials are compatible with the action of $q \in \Lambda_\N$, which is the obvious map $E_1^{rs} \rightarrow E_1^{r+1,s}$ (zero if $r<0$, and the identity otherwise). We will be particularly interested in the case where $L_0^\br$ and $L_1^\br$ are isotopic in $M$ (as exact Lagrangian submanifolds, compatibly with the grading and {\em Spin} structure), hence quasi-isomorphic in $\Fuk(M)$. Then \eqref{eq:q-spectral-e1} specializes to
\begin{equation} \label{eq:q-spectral-e1-b}
E_1^{rs} = \begin{cases} H^{r+s}(L_0;\C) & r \geq 0, \\ 0 & \text{otherwise.}
\end{cases}
\end{equation}
(Of course, one could substitute $H^{r+s}(L_1;\C)$ instead, since the two cohomologies are identified through the isotopy).

\begin{lemma} \label{th:endomorphisms}
For any object of the relative Fukaya category, there is a canonical isomorphism of graded rings, $\mathit{HF}^*_{X,X_\infty}(L^\mbr,L^\mbr) \iso H^*(L;\Lambda_\N)$.
\end{lemma}

This can be derived from our requirement that there should be no pseudo-holomorphic discs for the almost complex structure $J_L$, by following the argument pioneered (in the case of Hamiltonian Floer cohomology) by Piunikhin-Salamon-Schwarz \cite{piunikhin-salamon-schwarz94}, or alternatively by Morse-Bott type techniques going back to \cite{pozniak}. 

\begin{lemma} \label{th:unique-sphere}
A graded Lagrangian sphere in $M$ gives rise to an object of $\Fuk(X,X_\infty)$ which is unique up to quasi-isomorphism. Moreover, spheres that are isotopic (compatibly with gradings) in $M$ yield quasi-isomorphic objects.
\end{lemma}

\proof
Let $L_0^\mbr, L_1^\mbr$ be two objects of $\Fuk(X,X_\infty)$ obtained from the same graded Lagrangian sphere $L$. The differentials in the spectral sequence starting with \eqref{eq:q-spectral-e1-b} are necessarily zero, which implies that $1 \in H^0(L;\C)$ survives to yield a (non-unique) element of $\mathit{HF}^0_{X,X_\infty}(L_0^\mbr,L_1^\mbr)$. The multiplicative structure of the spectral sequence \eqref{eq:q-spectral-e1} shows that for any other $L_2^\mbr$, the map
\begin{equation} \label{eq:pseudo-continuation}
\mathit{HF}^*_{X,X_\infty}(L_1^\mbr,L_2^\mbr) \longrightarrow \mathit{HF}^*_{X,X_\infty}(L_0^\mbr,L_2^\mbr)
\end{equation}
given by right composition with our element is an isomorphism. The same holds for composition on the other side, which implies that our element is a quasi-isomorphism. The same reasoning applies if the graded Lagrangian spheres underlying $L_0^\mbr$ and $L_1^\mbr$ are merely isotopic.
\qed

\begin{lemma} \label{th:regular-isotopy}
Consider a Lagrangian isotopy $(L_r)_{0 \leq r \leq 1}$ in $M$, such that each $L_r$ is exact. On the total space $\Lambda = \bigcup_r \{r\} \times L_r$ of the isotopy, we want to have a {\em Spin} structure and a family of gradings.
Suppose that there is a family $(J_r)_{0 \leq r \leq 1}$ of almost complex structures, within the general class from Remark \ref{th:class-of-j}, such that $J_r$ is regular for $L_r$. Together with the previous data, this turns our Lagrangian submanifolds into objects $(L_r^\mbr)$ of $\Fuk(X,X_\infty)$. Then, $L_0^\mbr$ is quasi-isomorphic to $L_1^\mbr$.
\end{lemma}

{\sc Sketch of proof.} 
The difference between this and the previous argument is that we will define the isomorphism
\begin{equation} \label{eq:continuation-element}
[e_{(L_r^\mbr)}] \in \mathit{HF}^0_{X,X_\infty}(L_0^\mbr,L_1^\mbr)
\end{equation}
geometrically, instead of deriving its existence from algebraic arguments (this has the added advantage that the outcome is actually unique).

Suppose for simplicity that $L_0,L_1$ intersect transversally. Take the closed upper half plane $S$ (with coordinate $z$). We consider it as a once-punctured disc, with strip-like end 
\begin{equation}
\begin{aligned}
& \epsilon_S: \Rleq \times [0,1] \longrightarrow S, \\ 
& \epsilon_S(s,t) = -\exp(\pi(-s-it)). 
\end{aligned}
\end{equation}
Fix a function $\psi: \R \rightarrow [0,1]$ with $\psi(z) = 0$ for $z \leq-1$, $\psi(z) = 1$ for $z \geq 1$. We consider an inhomogeneous $\bar\partial$-equation with moving boundary conditions,
\begin{equation} \label{eq:pss-equation}
\left\{
\begin{aligned}
& u: S \longrightarrow X, \\
& u(z) \in L_{\psi(z)} \text{ for $z \in \R = \partial S$}, \\
& du + J_{S,u(z),z} \circ du \circ i = Y_{S,u(z),z} + J_{S,u(z),z} \circ Y_{S,u(z),z} \circ i, \\
& \textstyle\lim_{s \rightarrow -\infty} u(\epsilon_S(s,\cdot)) = x \in L_0 \cap L_1.
\end{aligned}
\right.
\end{equation}
Here, the family $J_S$ of almost complex structures satisfies $J_{S,z} = J_{\psi(z)}$ along the boundary, and as usual $J_{S,\epsilon_S(s,t)}$ depends only on $t$. The inhomogeneous term $Y$ is obtained from a $K_S$ as in \eqref{eq:generalized-floer}, where the precise requirements are as follows: $K_S$ and its first derivatives should vanish at all points of $X_\infty$; next, $K_S$ vanishes altogether on the strip-like ends; finally, if $\partial_z \in T_z(\partial S)$ stands for the standard vector tangent to the boundary at $z \in \R = \partial S$, then $K_{S,z}(\partial_z)|L_{\psi(z)}$ is the function (unique up to a constant) that describes the infinitesimal Lagrangian deformation $\partial_z L_{\psi(z)}$; see \cite[Section 8k]{seidel04} for more details. As in \eqref{eq:polygon-energy}, the energy of a solution $u$ equals the intersection number $u \cdot X_\infty$ up a uniformly bounded error, more precisely we have
\begin{equation} \label{eq:cont-energy}
\begin{aligned}
& E(u) - \textstyle\int_S u^*(dK_S + \half\{K_S,K_S\}) \\ & \qquad = \textstyle\int_S u^*\omega_X - d(u^*K_S) \\ & \qquad = A(x) + u \cdot X_\infty.
\end{aligned}
\end{equation}
Each isolated regular solution yields a preferred element of $\mathit{or}(x)$. One adds up these elements, with weights $q^{u \cdot X_\infty}$, to define the cocycle $e_{(L_r^\mbr)}$ underlying \eqref{eq:continuation-element}.

For the rest of the proof, there are several (related but slightly different) strategies. One can consider the reverse isotopy, and show (by a gluing and deformation argument) that the resulting class in $\mathit{HF}^0_{X,X_\infty}(L_1^\mbr,L_0^\mbr)$ is the inverse of \eqref{eq:continuation-element}. Or one can consider the maps \eqref{eq:pseudo-continuation} given by multiplication with \eqref{eq:continuation-element}, and show that these are continuation maps, hence isomorphisms. Finally, it is also an option to consider the truncation of $e_{(L_r^\mbr)}$ to $q = 0$, show that that is a quasi-isomorphism in $\Fuk(M)$, and then argue as in Lemma \ref{th:unique-sphere}. \qed

The assumption made in Lemma \ref{th:regular-isotopy}, that all the $J_r$ are regular, is not generic: while a single almost complex structure will generically admit no nonconstant pseudo-holomorphic discs of Maslov index zero, these can no longer be avoided in one-parameter families, as one can see from \eqref{eq:minus-one}. For a fixed (graded) Lagrangian submanifold $L$, the space of all almost complex structures is divided into chambers by codimension one walls consisting of non-regular ones (because there is an infinite number of homotopy classes of discs, there is potentially an infinite number of walls, which may be everywhere dense). Crossing a wall generally changes the isomorphism type of the object associated to $L$, in a way which is equivalent to turning on a $\Lambda_\N$-local system. A general discussion of this phenomenon involves dealing with multiply covered discs, hence requires the methods of \cite{fooo}. However, there is a weaker result that is more easily accessible:

\begin{lemma} \label{th:quasi-regular-isotopy}
In the situation of Lemma \ref{th:regular-isotopy}, suppose that $J_0$ is regular for $L_0$, $J_1$ is regular for $L_1$, but that the other $J_r$ only have the following partial regularity property, for some positive integer $N$: any $J_r$-holomorphic sphere, or $J_r$-holomorphic disc with boundary on $L_r$, satisfies $u \cdot X_\infty \geq N$. Then $L_0^\mbr$, $L_1^\mbr$ are quasi-isomorphic in $\Fuk(X,X_\infty) \otimes_{\Lambda_\N} \C[q]/q^N$.
\end{lemma}

To prove that, one goes through the proof of Lemma \ref{th:regular-isotopy}, noting that wherever there is a breakdown of compactness due to bubbling, the error term in the associated algebraic equation (for instance, the equation which says that counting solutions of \eqref{eq:pss-equation} yields a cocycle in the Floer complex) comes with a power of at least $q^N$.

Let's turn to the Fukaya category of $X$ itself. Some of the arguments above (notably Lemmas \ref{th:endomorphisms} and \ref{th:regular-isotopy}) carry over to that context without any issues. There is no exact analogue of the filtration by intersection number used in steps like \eqref{eq:q-spectral-e1}, but filtrations by energy can often serve as replacements (requiring a little more delicacy in the details). For later use, we record one application of that technique, which involves ``local Floer cohomology'' (there are many previous uses of this in the literature, see e.g. \cite{pozniak}).

Take two objects $L_0^\rbr,L_1^\rbr$ of $\Fuk(X)$ which intersect transversally. Fix $d$ such that $\lambda_{L_0}^{\otimes d}$ is single-valued. We will assume that our two Lagrangian submanifolds are {\em sufficiently close} in the following sense. There is an $H \in \smooth(X,\R)$ whose Hamiltonian flow $(\phi^r_H)$ satisfies $L_1 = \phi^1_H(L_0)$, and such that:
\begin{itemize} \itemsep1em
\item $|H(x)| < d^{-1}/4$ for all $x$;
\item The differential $dH$ vanishes at every point of $L_0 \cap L_1$.
\end{itemize}
(For a fixed $L_0$, one can fulfil these conditions as long as $L_1$ is $C^1$-close to $L_0$, by writing $L_1$ as a graph in a symplectic tubular neighbourhood of $L_0$). Moreover, the isotopy $L_r = \phi^r_H(L_0)$ should be compatible with the choices of {\em Spin} structure, grading, and covariantly constant multisections $\lambda_{L_0}$, $\lambda_{L_1}$ (see Remark \ref{th:deform-rational}). In particular, $\lambda_{L_1}^{\otimes d}$ is again single-valued.

Each intersection point $x \in L_0 \cap L_1$ is a stationary point of $\phi_H^r$. Moreover, the action $\bar{A}(x) \in \R/d^{-1}\Z$ has a preferred real-valued lift $A(x) = -H(x)$. This means that one can write
\begin{equation} \label{eq:crunch}
\mathit{CF}^*_X(L_0^\rbr,L_1^\rbr) \iso \mathit{CF}_{X,+}^*(L_0^\rbr,L_1^\rbr) \otimes_{\Lambda_{d^{-1}\N}} \Lambda_\Q,
\end{equation}
where $\Lambda_{d^{-1}\N}$ is the ring of formal power series in $q^{1/d}$, and $\mathit{CF}^*_{X,+}(L_0^\rbr,L_1^\rbr)$ is a free graded module over that ring with one generator (up to the usual sign issues) for each $x$. The map from right to left in \eqref{eq:crunch} multiplies the generator corresponding to an intersection point $x$ with $q^{A(x)}$. 
There is a differential on $\mathit{CF}^*_{X,+}(L_0^\rbr,L_1^\rbr)$ which corresponds to the standard Floer differential under \eqref{eq:crunch}: it counts Floer trajectories $u$ asymptotic to $x_0,x_1 \in L_0 \cap L_1$ with powers
\begin{equation} \label{eq:q-crunch}
q^{E(u)+A(x_1)-A(x_0)}. 
\end{equation}
One knows that $E(u)+A(x_1)-A(x_0) \in d^{-1}\Z$, that $E(u)>0$, and that $A(x_1)-A(x_0) > -d^{-1}/2$. This is why all of the powers \eqref{eq:q-crunch} are nonnegative. We denote the cohomology of $\mathit{CF}_{X,+}^*(L_0^\rbr,L_1^\rbr)$ by $\mathit{HF}^*_{X,+}(L_0^\rbr,L_1^\rbr)$. Local Floer cohomology is obtained from this by considering only the $q^0$ term, which counts trajectories with ``small energy'' $E(u) < d^{-1}/2$: formally,
\begin{equation}
\mathit{HF}^*_{X,\mathit{local}}(L_0^\rbr,L_1^\rbr) \stackrel{\text{def}}{=} H^*(\mathit{CF}_{X,+}^*(L_0^\rbr,L_1^\rbr) \otimes_{\Lambda_{d^{-1}\N}} \C).
\end{equation}

There are also groups $\mathit{HF}^*_{X,+}(L_k^\rbr,L_k^\rbr)$ and $\mathit{HF}^*_{X,\mathit{local}}(L_k^\rbr,L_k^\rbr)$ for $k = 0,1$, which are best defined using Morse-Bott methods (it is possible to use small Hamiltonian perturbations instead, but that has to be done carefully in order to preserve the necessary energy bounds). For those, it is easy to show as in Lemma \ref{th:endomorphisms} that they are canonically isomorphic to $H^*(L_k;\Lambda_{d^{-1}\N})$ and $H^*(L_k;\C)$, respectively. Finally, we have products
\begin{equation} \label{eq:plus-products}
\begin{aligned}
& \mathit{HF}^*_{X,+}(L_1^\rbr,L_1^\rbr) \otimes \mathit{HF}^*_{X,+}(L_0^\rbr,L_1^\rbr) \longrightarrow \mathit{HF}^*_{X,+}(L_0^\rbr,L_1^\rbr), \\
& \mathit{HF}^*_{X,+}(L_0^\rbr,L_1^\rbr) \otimes \mathit{HF}^*_{X,+}(L_0^\rbr,L_0^\rbr) \longrightarrow
\mathit{HF}^*_{X,+}(L_0^\rbr,L_1^\rbr), \\
& \mathit{HF}^*_{X,+}(L_1^\rbr,L_0^\rbr) \otimes \mathit{HF}^*_{X,+}(L_0^\rbr,L_1^\rbr) \longrightarrow
\mathit{HF}^*_{X,+}(L_0^\rbr,L_0^\rbr),
\end{aligned}
\end{equation}
plus their analogues for $\mathit{HF}^*_{X,\mathit{local}}$, and the same with the roles of $L_0^\rbr$ and $L_1^\rbr$ reversed. These products are associative (note that, while the definition of \eqref{eq:crunch} only required $|H(x)| < d^{-1}/2$, associativity uses the full strength of the bound we have imposed).

\begin{lemma} \label{th:local-floer-homology}
$\mathit{HF}^*_{X,\mathit{local}}(L_0^\rbr,L_1^\rbr) \iso H^*(L_0;\C)$.
\end{lemma}

{\sc Sketch of proof.}
Consider solutions of an appropriate equation \eqref{eq:pss-equation} for the isotopy $(L_r)$, and count them with powers
\begin{equation} \label{eq:q-corrected}
q^{E(u) - \textstyle\int_S u^*(dK_S + \half\{K_S,K_S\}) - A(x)}.
\end{equation}
By a careful choice of $K_S$ based on our given $H$, one can ensure that the curvature term satisfies
\begin{equation} \label{eq:control-curvature}
\Big| {\textstyle\int_S u^*(dK_S + \half\{K_S,K_S\})} \Big| < d^{-1}/4. 
\end{equation}
In parallel with \eqref{eq:q-crunch}, this shows that the exponents \eqref{eq:q-corrected} are nonnegative. The outcome is a cochain in $\mathit{CF}^*_{X,+}(L_0^\rbr,L_1^\rbr)$. Because of bubbling, this is not a cocycle. However, its lowest energy ($q^0$) term is a cocycle mod higher powers of $q$, hence gives rise to a class
\begin{equation} \label{eq:local-e}
[e_{(L_r^\rbr),\mathit{local}}] \in \mathit{HF}^0_{X,\mathit{local}}(L_0^\rbr,L_1^\rbr).
\end{equation}
The rest follows the same strategy as in Lemma \ref{th:regular-isotopy}: one introduces the corresponding class for the reverse isotopy, and shows that their product in either order is the identity map in each group $\mathit{HF}^*_{X,\mathit{local}}(L_k^\rbr,L_k^\rbr)$. The rest then follows from associativity of the product. \qed

\begin{lemma} \label{th:unique-sphere-2}
A graded Lagrangian sphere in $X$ gives rise to an object of $\Fuk(X)$ which is unique up to quasi-isomorphism. More generally, spheres that are isotopic (compatibly with gradings) yield quasi-isomorphic objects.
\end{lemma}

\proof It is sufficient to show that, given some sphere $L_0^\rbr$, any sphere $L_1^\rbr$ which has transverse intersection with $L_0^\rbr$ and is sufficiently $C^1$-close (including grading) yields a quasi-isomorphic object. Lemma \ref{th:local-floer-homology} applies, and the spectral sequence associated to the $q^{1/d}$-adic filtration of $\mathit{CF}^*_{X,+}(L_0^\rbr,L_1^\rbr)$ then necessarily degenerates, leading to $\mathit{HF}^*_{X,+}(L_0^\rbr,L_1^\rbr) \iso H^*(L_0;\Lambda_{d^{-1}\N})$ and hence $\mathit{HF}^*_X(L_0^\rbr,L_1^\rbr) \iso H^*(L_0;\Lambda_\Q)$. An argument along the same lines shows that if we take any class in $\mathit{HF}^0_{X,+}(L_0^\rbr,L_1^\rbr)$ whose reduction to $q = 0$ is \eqref{eq:local-e}, then its image in $\mathit{HF}^0_X(L_0^\rbr,L_1^\rbr)$ is a quasi-isomorphism in $\Fuk(X)$.
\qed

This was the analogue of Lemma \ref{th:unique-sphere}, based (as promised) on filtrations by energy rather than intersection number. There is also an analogue of Lemma \ref{th:quasi-regular-isotopy}, which goes as follows:

\begin{lemma} \label{th:extended-local-floer-homology}
In the situation of Lemma \ref{th:local-floer-homology}, suppose that there is some positive $N \in d^{-1}\N$ and a family $(J_r)_{0 \leq r \leq 1}$ of almost complex structures interpolating between $J_{L_0}$ and $J_{L_1}$, such that: every $J_r$-holomorphic sphere, or $J_r$-holomorphic disc with boundary on $L_r$, has energy $\geq N$. Then
\begin{equation}
H(CF^*_{X,+}(L_0^\rbr,L_1^\rbr) \otimes_{\Lambda_{d^{-1}\N}}
\C[q^{1/d}]/q^N) \iso H^*(L_0;\C[q^{1/d}]/q^N).
\end{equation}
\qed
\end{lemma}

\subsection{\label{subsec:skip-divisor}}
Our next task is to derive a geometric criterion for the nontriviality of the deformation given by the relative Fukaya category. We will explain how to do this, at least in principle, by starting with a Lagrangian submanifold which does intersect the divisor at infinity, and then pushing it off that divisor in different ways.

Suppose that we have a Lagrangian surface $L_{1/2} \subset X$ which is disjoint from the singular locus of $X_\infty$, and such that the intersection $Z = L_{1/2} \cap X_\infty$ is a smooth one-dimensional manifold. Fix a function $H \in \smooth(X,\R)$ such that the restriction $H|L_{1/2}$ is Morse, and its gradient (formed with respect to the restriction of the K{\"a}hler metric to $L_{1/2}$) is transverse to $Z$. Note that $\nabla (H|L_{1/2})$ equips $Z$ with a co-orientation, hence associates to it a Poincar\'e dual class $z \in H^1(L_{1/2};\R)$. We will assume that $H$ is small (by multiplying it with a small nonzero constant if necessary). Write $L_0$ and $L_1$ for the surfaces obtained from $L_{1/2}$ by flowing along the flow of $H$ with time $-1/2$ (for $L_0$) and $+1/2$ (for $L_1$). The assumptions on $\nabla(H|L_{1/2})$ ensure that $L_0$ and $L_1$ are contained in $M$. Moreover,
\begin{equation} \label{eq:l0l1}
L_0 \cap L_1 = \mathit{Crit}(H|L_{1/2}).
\end{equation}

\begin{lemma} \label{th:one-pushoff}
Under the identifications $H^1(L_0;\R) \iso H^1(L_{1/2};\R) \iso H^1(L_1;\R)$, we have
\begin{equation} \label{eq:l01}
[\theta_M|L_1] - [\theta_M|L_0] = z.
\end{equation}
\end{lemma}

\proof
Take a one-cycle $c$ on $L_{1/2}$, and let $c_0,c_1$ be the corresponding cycles in the pushoff surfaces $L_0,L_1$. There is a two-cycle $C$ in $\Lambda$ (here, $\Lambda$ is the domain of the obvious isotopy from $L_0$ to $L_1$) with boundary $c_1 - c_0$, and whose intersection number with $X_{\infty}$ is $-z(c)$. On the other hand, because our isotopies are Hamiltonian, we know that
\begin{equation}
0 = \textstyle \int_C \omega_X = \int_{c_1} \theta_M - \int_{c_0} \theta_M + (C \cdot X_{\infty}). \qed
\end{equation}

\begin{remark} \label{th:two-pushoffs}
There is a slight variation which will be useful at one point later on. Namely, suppose that we have another Hamiltonian $\tilde{H}$, with its associated cohomology class $\tilde{z}$ and Lagrangian submanifolds $\tilde{L}_k$, obtained from the same $L_{1/2}$. Then
\begin{equation} \label{eq:two-pushoffs}
\begin{aligned}
& [\theta_M|\tilde{L}_1] - [\theta_M|L_1] = \half(\tilde{z} - z), \\
& [\theta_M|\tilde{L}_0] - [\theta_M|L_0] = \half(z - \tilde{z}).
\end{aligned}
\end{equation}
The formula \eqref{eq:l01} can be recovered from this by setting $\tilde{H} = -H$. Note that $\tilde{z}$ and $z$ have the same mod $2$ reduction, hence the right hand side of \eqref{eq:two-pushoffs} is integral.
\end{remark}

From now on, we impose:

\begin{assumption}
$L_0$ and $L_1$ are exact in $M$.
\end{assumption}

Since $L_0$ and $L_1$ are Hamiltonian perturbations of $L_{1/2}$, the exactness of either one requires that the flat vector bundle $o_X|L_{1/2}$ is trivial, hence in particular that $L_{1/2}$ is rational. The exactness of both $L_0$ and $L_1$ implies that $z = 0$, by Lemma \ref{th:one-pushoff}. Hence, we can choose a bounding cochain $B$ for the underlying cocycle $Z$. Concretely, $B$ is a locally constant function on $L_{1/2} \setminus Z$ which jumps by $\pm 1$ when crossing any component of $Z$ (with the sign determined by co-orientation). 

\begin{lemma} \label{th:modified-action}
There is a choice of primitives $K_{L_k}$ for $\theta_M|L_k$ such that the action of any point $x \in L_0 \cap L_1$ is
\begin{equation} \label{eq:abh}
A(x) = B(x) - H(x).
\end{equation}
\end{lemma}

\proof
Take a path $c:\R \rightarrow L_{1/2}$ whose limits $\lim_{s \rightarrow -\infty} c(s) = x_0$ and $\lim_{s \rightarrow \infty} c(s) = x_1$ are critical points of $H|L_{1/2}$. By moving that path around using the flow of $H$, we get a map $u: \R \times [0,1] \longrightarrow M$ such that $u(s,0) \in L_0$, $u(s,1) \in L_1$, and with the same asymptotic behaviour. As in the proof of Lemma \ref{th:one-pushoff}, we have $u \cdot X_\infty = -Z \cdot c$. Because of the way in which $u$ is constructed, we have
\begin{equation}
\int_{\R \times [0,1]} u^*\omega_X = \int_{\R} dH(c'(s)) \, ds = H(x_1) - H(x_0).
\end{equation}
On the other hand, by \eqref{eq:energy} the same integral can be written as
\begin{equation}
\begin{aligned}
& \int_{\R \times [0,1]} u^*\omega_X = A(x_0) - A(x_1) + (u \cdot X_\infty) \\ & = A(x_0) - A(x_1) - (Z \cdot c) =
A(x_0) - B(x_0) - A(x_1) + B(x_1).
\end{aligned}
\end{equation}
Comparing the two expressions shows that \eqref{eq:abh} holds up to a constant independent of $x$. One can modify $K_{L_0}$ or $K_{L_1}$ to make that constant equal to zero.
\qed

The specific situation we are interested in is when
\begin{equation}
L_{1/2} \iso S^1 \times S^1, \quad  Z = \{1/4,3/4\} \times S^1.
\end{equation}
Choose the function $H$ so that $\partial_p H|\{1/4\} \times S^1 < 0$ and $\partial_p H|\{3/4\} \times S^1 > 0$, where $(p,q)$ are the coordinates on $L_{1/2}$. In particular, $z = 0$, and the function $B$ can be taken to be $B(p,q) = 0$ for $p \in (1/4,3/4)$, $B(p,q) = 1$ for $p \in (-1/4,1/4)$. We also need to choose auxiliary data as follows. $L_{1/2}$ should come with a grading and {\em Spin} structure, which are then inherited by $L_0$ and $L_1$. Moreover, we choose a family $(J_r)_{0 \leq r \leq 1}$ of almost complex structures (as in Remark \ref{th:class-of-j}) such that:
\begin{itemize} \itemsep1em
\item $J_0$ is regular for $L_0$, and $J_1$ is regular for $L_1$;
\item Every $J_r$-holomorphic sphere, or $J_r$-holomorphic disc with boundary on $L_r$, has energy $\geq 2$.
\end{itemize}
These conditions can be easily satisfied by first choosing $J_{1/2}$ generically, and then taking a small perturbation of the constant family (assuming that $H$ is chosen sufficiently small as well). 

Let's first consider our two pushoffs as objects $L_0^\br, L_1^\br$ of $\Fuk(M)$. Each point of \eqref{eq:l0l1} lies either in $(1/4,3/4) \times S^1$ or in $(-1/4,1/4) \times S^1$, and we correspondingly write
\begin{equation} \label{eq:fake-splitting}
\mathit{CF}^*_{M}(L_0^\br,L_1^\br) = 
\mathit{CF}^*_{M}(L_0^\br,L_1^\br)_{(0)}
\oplus
\mathit{CF}^*_{M}(L_0^\br,L_1^\br)_{(1)},
\end{equation}
A priori, this is a splitting of graded vector spaces, not necessarily compatible with the differential. Lemma \ref{th:modified-action} shows that the critical points contributing to the first summand have action $A(x) \approx 0$, and those for the second summand have action $A(x) \approx 1$. Hence, $\mathit{CF}^*_{M}(L_0^\br,L_1^\br)_{(1)}$ is in fact a subcomplex, and $\mathit{CF}^*_{M}(L_0^\br,L_1^\br)_{(0)}$ can be equipped with the induced quotient differential.

\begin{lemma}
As graded vector spaces,
\begin{equation} \label{eq:3-h}
\begin{aligned}
& H^*(\mathit{CF}^*_M(L_0^\br,L_1^\br)_{(0)}) \iso H^*(S^1;\C)[-1], \\
& H^*(\mathit{CF}^*_M(L_0^\br,L_1^\br)_{(1)}) \iso H^*(S^1;\C), \\
& \mathit{HF}^*_M(L_0^\br,L_1^\br) \iso H^*(S^1;\C)[-1] \oplus H^*(S^1;\C).
\end{aligned}
\end{equation}
\end{lemma}

\proof
By a continuation map argument, one can show that all the cohomology groups under consideration are independent of the particular choice of $H$, within the overall class introduced above. In fact, we may slightly extend that class by allowing the function to be Morse-Bott, and then use such a function with two critical circles $C_{(0)} \subset (1/4,3/4) \times S^1$ and $C_{(1)} \subset (-1/4,1/4) \times S^1$, whose Morse index is $1$ and $0$, respectively (this choice of Morse indices is made necessary because of the previously imposed requirement on $\partial_p H|Z$). The first two parts of \eqref{eq:3-h} then reproduce the cohomology of those circles, with the degree shift coming from the Morse index. For the last part, note that by construction, we have a long exact sequence
\begin{multline} \label{eq:01-boundary}
\cdots \rightarrow H^*(\mathit{CF}^*_M(L_0^\br,L_1^\br)_{(1)}) \iso H^*(S^1;\C)[-1] \longrightarrow
\mathit{HF}^*_M(L_0^\br,L_1^\br) \\ \longrightarrow
H^*(\mathit{CF}^*_M(L_0^\br,L_1^\br)_{(0)}) \rightarrow \cdots
\end{multline}
The previous computations show that the connecting map for this sequence is necessarily zero (for degree reasons).
\qed

\begin{lemma}
The product
\begin{equation}
\mathit{HF}^0_M(L_0^\br,L_1^\br) \otimes \mathit{HF}^2_M(L_0^\br,L_0^\br) \longrightarrow \mathit{HF}^2_M(L_0^\br,L_1^\br).
\end{equation}
vanishes.
\end{lemma}

\proof
This is easiest if one thinks of this as being given by the {\em quantum cap} action of $\mathit{HF}^2_M(L_0^\br,L_0^\br) \iso H^2(L_0;\C)$ on $\mathit{HF}^*_M(L_0^\br,L_1^\br)$. This quantum cap action is defined (roughly speaking) by fixing a generic point in $L_0$, and counting Floer trajectories whose boundary goes through that point. For action reasons, any such Floer trajectory $u$ has limits $\lim_{s \rightarrow -\infty} u(s,\cdot) \in C_{(1)}$, $\lim_{s \rightarrow +\infty} u(s,\cdot) \in C_{(0)}$. However, the dimension of the space of such trajectories up to translation is $-1$, hence it is empty (the same argument has appeared before when showing that \eqref{eq:01-boundary} has zero connecting map, but here we have spelled it out in more concrete terms).
\qed

Let's look at the same Lagrangian submanifolds as objects $L_0^\mbr, L_1^\mbr$ of $\Fuk(X,X_\infty)$. In parallel with \eqref{eq:fake-splitting}, write
\begin{equation}
\mathit{CF}^*_{X,X_\infty}(L_0^\mbr,L_1^\mbr) = 
\mathit{CF}^*_{X,X_\infty}(L_0^\mbr,L_1^\mbr)_{(0)}
\oplus
\mathit{CF}^*_{X,X_\infty}(L_0^\mbr,L_1^\mbr)_{(1)},
\end{equation}
which is a splitting of graded $\Lambda_\N$-modules but not of chain complexes. A third possible viewpoint is to look at the associated objects $L_0^\rbr, L_1^\rbr$ of $\Fuk(X)$. This is important because in that sense, they are small perturbations of each other, allowing one to use the formalism from \eqref{eq:crunch} (where $d = 1$), with its chain complex of $\Lambda_\N$-modules $\mathit{CF}^*_{X,+}(L_0^\rbr,L_1^\rbr)$. In fact, that complex is almost the same as $\mathit{CF}^*_{X,X_\infty}(L_0^\mbr,L_1^\mbr)$, with the difference lying in the bookkeeping of powers of $q$, as expressed in Lemma \ref{th:modified-action}. The precise relationship is that
\begin{equation} \label{eq:q-shift}
\begin{aligned}
\mathit{CF}^*_{X,+}(L_0^\rbr,L_1^\rbr) & \iso q\mathit{CF}^*_{X,X_\infty}(L_0^\mbr,L_1^\mbr)_{(0)} \oplus
\mathit{CF}^*_{X,X_\infty}(L_0^\mbr,L_1^\mbr)_{(1)} \\ &
\subset \mathit{CF}^*_{X,X_\infty}(L_0^\mbr,L_1^\mbr),
\end{aligned}
\end{equation}
as chain complexes of $\Lambda_\N$-modules. 
As a consequence, for each $N$ we have a commutative diagram
\begin{equation} \label{eq:q-diagram}
\xymatrix{
H^*(\mathit{CF}^*_{X,+}(L_0^\rbr,L_1^\rbr) \otimes_{\Lambda_\N} \C[q]/q^{N-1})
\ar[d] \ar@/_10pc/[dd]_-{q}
\\
H^*(\mathit{CF}^*_{X,X_\infty}(L_0^\mbr,L_1^\mbr)
\otimes_{\Lambda_\N} \C[q]/q^{N-1}) 
\ar[d]
\\
H^*(\mathit{CF}^*_{X,+}(L_0^\rbr,L_1^\rbr) \otimes_{\Lambda_\N} \C[q]/q^N).
}
\end{equation}
The top $\downarrow$ is induced by the inclusion of the subcomplex \eqref{eq:q-shift}, and the second one by multiplication with $q$, which lands in the same subcomplex.

\begin{lemma} \label{th:nonzero-product-12}
The product
\begin{equation}
\xymatrix{
\parbox{19em}{
$H^0(\mathit{CF}_{X,X_\infty}(L_0^\mbr,L_1^\mbr) \otimes_{\Lambda_\N} \C[q]/q^2)$ \newline
$\otimes_{\C[q]/q^2} \, H^2(\mathit{CF}_{X,X_\infty}(L_0^\mbr,L_0^\mbr) \otimes_{\Lambda_\N} \C[q]/q^2)$
}
\ar[d]
\\
H^2(\mathit{CF}_{X,X_\infty}(L_0^\mbr,L_1^\mbr) \otimes_{\Lambda_\N} \C[q]/q^2)
}
\end{equation}
is nonzero.
\end{lemma}

\proof
As before, it is more convenient to think of this product as being given by the quantum cap action of $H^2(L_0;\C)$ on the Floer cohomology of the pair $(L_0,L_1)$. The isomorphism from Lemma \ref{th:local-floer-homology} carries the quantum cap action into the ordinary ring structure on the cohomology of $L_0$. In particular, it follows that $H^2(L_0;\C)$ acts nontrivially on $\mathit{HF}^*_{X,\mathit{local}}(L_0^\rbr,L_1^\rbr)$. Next, as a consequence of Lemma \ref{th:extended-local-floer-homology}, one sees that the map
\begin{equation}
\xymatrix{
\mathit{HF}^*_{X,\mathit{local}}(L_0^\rbr,L_1^\rbr) \ar[r]^-{=} & 
H^*(\mathit{CF}^*_{X,+}(L_0^\rbr,L_1^\rbr) \otimes_{\Lambda_\N} \C)
\ar[d]^-{q} \\ &
H^*(\mathit{CF}^*_{X,+}(L_0^\rbr,L_1^\rbr) \otimes_{\Lambda_\N} \C[q]/q^2)
}
\end{equation}
is injective. By going through the $N = 2$ case of \eqref{eq:q-diagram}, it follows that $H^2(L_0;\C)$ also acts nontrivially on $H^*(\mathit{CF}^*_{X,X_\infty}(L_0^\mbr,L_1^\mbr) \otimes_{\Lambda_\N} \C[q]/q^2)$.
\qed

We can now apply Lemma \ref{th:product-turns-on} (with $Y_0 = L_0^\br[-2]$, $Y_1 = L_0^\br$ and $Y_2 = L_1^\br$) to conclude:

\begin{cor} \label{th:wall-crossing-2}
Suppose that there is a graded Lagrangian torus $S^1 \times S^1 \iso L_{1/2} \subset X$ which is disjoint from the singularities of $X_\infty$, and which satisfies $L_{1/2} \cap X_{\infty} = \{1/4,3/4\} \times S^1$. Suppose in addition that the tori $L_0,L_1 \subset M$ obtained by pushing $L_{1/2}$ off $X_\infty$ in the way explained above are exact. Then the relative Fukaya category $\Fuk(X,X_\infty)$, seen as a deformation of $A_\infty$-categories, gives rise to a nonzero deformation class in $\mathit{HH}^2(\Fuk(M),\Fuk(M))$.
\qed
\end{cor}

\begin{remark}
To summarize the argument above, it is instructive to briefly digress to the toy model situation one dimension lower. Take $X$ to be an elliptic curve, and $X_\infty$ an ample divisor consisting of two points. Then $L_{1/2} \iso S^1$ would be a non-contractible simple closed curve going through both points of $X_\infty$, or more concretely with $L_{1/2} \cap X_\infty = \{1/4,3/4\}$. The pushoffs $L_0,L_1$ are drawn in Figure \ref{fig:toy-pushoffs}. One can arrange that $L_0 \cap L_1 = \{x_{(0)}, x_{(1)}\}$ consists of two points with Maslov indices $1$ and $0$, respectively.
\begin{figure}
\begin{centering}
\begin{picture}(0,0)%
\includegraphics{missing.pstex}%
\end{picture}%
\setlength{\unitlength}{3355sp}%
\begingroup\makeatletter\ifx\SetFigFont\undefined%
\gdef\SetFigFont#1#2#3#4#5{%
  \reset@font\fontsize{#1}{#2pt}%
  \fontfamily{#3}\fontseries{#4}\fontshape{#5}%
  \selectfont}%
\fi\endgroup%
\begin{picture}(2812,2424)(1179,-6373)
\put(3976,-4636){\makebox(0,0)[lb]{\smash{{\SetFigFont{10}{12.0}{\rmdefault}{\mddefault}{\updefault}{\color[rgb]{0,0,0}$L_{1/2}$}%
}}}}
\put(2851,-4711){\makebox(0,0)[lb]{\smash{{\SetFigFont{10}{12.0}{\rmdefault}{\mddefault}{\updefault}{\color[rgb]{0,0,0}$L_1$}%
}}}}
\put(2851,-5686){\makebox(0,0)[lb]{\smash{{\SetFigFont{10}{12.0}{\rmdefault}{\mddefault}{\updefault}{\color[rgb]{0,0,0}$L_0$}%
}}}}
\put(3676,-5236){\makebox(0,0)[lb]{\smash{{\SetFigFont{10}{12.0}{\rmdefault}{\mddefault}{\updefault}{\color[rgb]{0,0,0}$x_{(1)}$}%
}}}}
\put(1951,-6061){\makebox(0,0)[lb]{\smash{{\SetFigFont{10}{12.0}{\rmdefault}{\mddefault}{\updefault}{\color[rgb]{0,0,0}$x_{(0)}$}%
}}}}
\put(1501,-4411){\makebox(0,0)[lb]{\smash{{\SetFigFont{10}{12.0}{\rmdefault}{\mddefault}{\updefault}{\color[rgb]{0,0,0}$X_\infty$}%
}}}}
\end{picture}%
\caption{\label{fig:toy-pushoffs}}
\end{centering}
\end{figure}

There are two obvious Floer trajectories (holomorphic disc with two corners) with small energy $E(u) \approx 0$, connecting $x_{(1)}$ to $x_{(0)}$, and each of them intersects $X_\infty$ exactly once. In view of \eqref{eq:energy}, this shows that if we consider $L_0$ and $L_1$ as exact Lagrangian submanifolds in $M = X \setminus X_\infty$, then
\begin{equation}
A(x_{(0)}) - A(x_{(1)}) = E(u) - u \cdot X_\infty \approx -1.
\end{equation}
After shifting the action by a constant, we may assume that $A(x_{(0)}) \approx 0$, $A(x_{(1)}) \approx 1$, exactly as in our previous discussion. If we adopt the formalism from \eqref{eq:crunch}, our two Floer trajectories contribute $\pm q^0$ to the differential on $\mathit{CF}^*_{X,+}(L_0^\rbr,L_1^\rbr)$ (the signs are opposite, so the contributions actually cancel out). On the other hand, the contributions of the same trajectories to the differential in $\Fuk(X,X_\infty)$ are $\pm q^1$, because of the intersection number. The difference between these two weights explains \eqref{eq:01-boundary}. For the same reason, the two trajectories do not contribute to the quantum cap action of $H^1(L_0;\C)$ on $\mathit{HF}^*_M(L_0^\br,L_1^\br)$, which is in fact zero; but they do contribute to the corresponding action on $\mathit{HF}^*_{X,\mathit{local}}(L_0^\rbr,L_1^\rbr)$, which is nonzero, leading to the analogue of Lemma \ref{th:nonzero-product-12}.
\end{remark}

\section{Computations in Fukaya categories\label{sec:induction}}

To effect the link between symplectic geometry and homological
algebra, we basically rely on two tools. The first is the
correspondence between Dehn twists and algebraic twists, which is a
reformulation of the long exact sequence in Floer cohomology. The
version in the literature \cite{seidel01,seidel04} is for the exact
or affine situation, but as will be outlined below, the proof can be
easily adapted to the projective and relative cases. Our second tool
is a dimensional induction machine based on the notion of matching
cycle, developed in \cite[Chapter 3]{seidel04}; this is strictly limited to Lefschetz fibrations with base $\C$, which is the reason why affine varieties appear in this paper at all. We will treat dimensional induction as a black box, stating its properties in an axiomatic way only; ultimately it
is again an application of the long exact sequence, but the way in
which this works is too complicated to be summarized briefly.

\subsection{}
Let $M$ be an affine Calabi-Yau of dimension $n>1$, and $L \subset M$
a Lagrangian sphere. Equip it with a grading and the unique $Spin$
structure, making it into an object $L^\br$ of $\Fuk(M)$. Let
$\tau_L$ be the Dehn twist along $L$, and $\taugr_L$ its canonical
grading (again, there is the small matter of framing, which we
ignore). Consider also the algebraic twist $T_{L^\br}$ along $L^\br$, acting on objects of $D^b\Fuk(M)$. Let $L_1^\br$ be any exact
Lagrangian brane, with underlying Lagrangian submanifold $L_1$.
Clearly $\tau_L^{-1}(L_1)$ carries an induced brane structure, making
it into an object $\taugr_L^{-1}(L_1^\br)$ of the Fukaya category.
Then

\begin{prop} \label{th:exact-sequence-1}
$T_{L^\br}(\taugr_L^{-1}(L_1^\br))$ is isomorphic to $L_1^\br$ in
$D^b\Fuk(M)$.
\end{prop}

On the algebraic side, this is an application of Lemma
\ref{th:recognize-twist}, where $W$ is an (arbitrary) other exact Lagrangian
brane $L_0^\br$. We need to explain how a suitable moduli space setup
provides the necessary pair $(k,h)$. This is a variation of the
construction in \cite{seidel01}, based on pseudo-holomorphic sections
of exact symplectic Lefschetz fibrations with Lagrangian boundary
conditions. Specifically, taking $S$ to be the closed upper half
plane, there is such a fibration $E_L \rightarrow S$, which has a
single critical point somewhere over $\mathit{int}(S)$. The regular fibres of
$E_L$ are isomorphic to $M$; the monodromy around the singular fibre is
$\tau_L$; and the fibration is trivial outside a compact subset. Given any exact Lagrangian brane in $M$ such
as $L_1^\br$, one can define a subbundle $F_{L_1} \subset
E_L|\partial S$ which, in the canonical trivialization, satisfies
\begin{equation}
 F_{L_1} \cap E_z = \begin{cases}
 \tau_L^{-1}(L_1) & z \ll 0, \\
 L_1 & z \gg 0.
 \end{cases}
\end{equation}
By considering an equation similar to \eqref{eq:generalized-floer}
for sections of $E_L$ with boundary values on $F_{L_1}$, one then
obtains a ``relative invariant'' in the style of
\cite{piunikhin-salamon-schwarz94} taking values in
$\mathit{HF}^0(\taugr_L^{-1}(L_1^\br),L_1^\br)$. Take $k$ to be the cocycle
representing it, which of course depends on the specific choices of
almost complex structures etc.\ used in the definition.
\FIGURE{fig:homotopy}{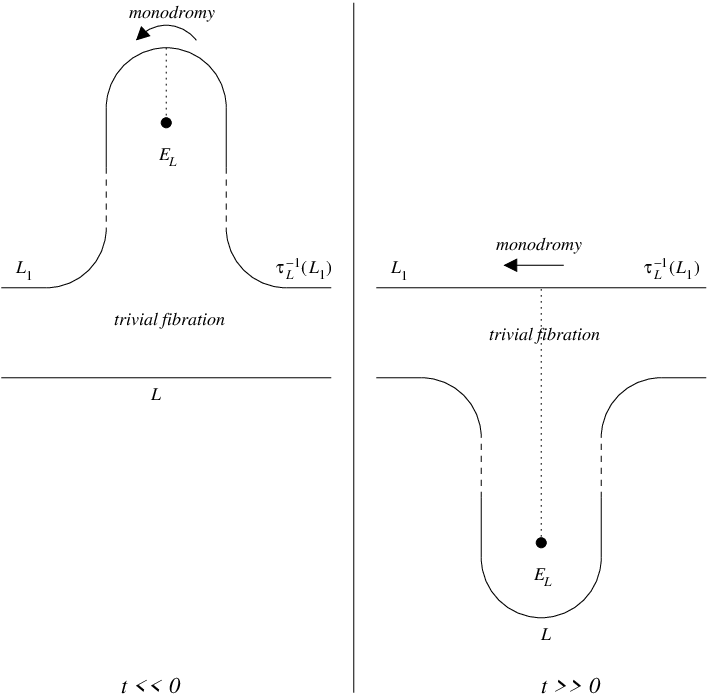}{hb}%

For $h$ one uses a parametrized version of the same construction,
applied to a family of Lefschetz fibrations $E^t$ with boundary
conditions $F^t$, $t \in \R$. Figure \ref{fig:homotopy} shows the
asymptotic behaviour of this family: for $t \ll 0$, $(E^t,F^t)$ is
the fibre connected sum of $(E_L,F_{L_1})$ and a trivial fibration
over the three-punctured disc, with boundary conditions over the
different boundary components given by the Lagrangian submanifolds
$L$, $\tau_L^{-1}(L_1)$ and $L_1$. The invariant associated to the
degenerate limit $t = -\infty$ is therefore the composition
\begin{equation} \label{eq:t-plus}
\begin{aligned}
 \mathit{CF}_M^*(L^\br,\taugr_L^{-1}(L_1^\br)) & \longrightarrow \mathit{CF}_M^*(L^\br,L_1^\br), \\
 a & \longmapsto \mu^2_M(a,k).
\end{aligned}
\end{equation}
For $t \gg 0$, $(E^t,F^t)$ again decomposes as a connected sum. Its
first component is $E_L$ as before, but now with boundary condition
$F_L$; and the second component is trivial, with boundary conditions
given by $L$, $\tau_L^{-1}(L_1)$ and again $L$. In principle, the
associated invariant is
\begin{equation} \label{eq:t-minus}
\begin{aligned}
 \mathit{CF}_M^*(L^\br,\taugr_L^{-1}(L_1^\br)) \iso
 \mathit{CF}_M^*(\taugr_L(L^\br),L_1^\br) &
 \longrightarrow \mathit{CF}_M^*(L^\br,L_1^\br), \\
 a & \longmapsto \mu^2_M(l,a),
\end{aligned}
\end{equation}
where $l \in \mathit{CF}^0_M(L^\br,\taugr_L(L^\br)) \iso
\mathit{CF}^0_M(\taugr^{-1}(L^\br),L^\br)$ is obtained from $(E_L,F_L)$.
However, since $\taugr_L(\Lgr) = \Lgr[1-n]$ with $n>1$, 
$\mathit{HF}^0_M(L^\br,\taugr_L(L^\br)) = \mathit{HF}^{1-n}_M(L^\br,L^\br) \iso
H^{1-n}(S^n;\C) = 0$, and by suitably choosing the Floer data, one
can achieve that $CF^0(L^\br,\taugr_L(L^\br))$ vanishes too. Hence
$l$ and the expression \eqref{eq:t-minus} are both zero, and from
this one deduces that the parametrized moduli space of sections
associated to the deformation from $t = -\infty$ to $t = +\infty$
yields a chain homotopy $h$ between \eqref{eq:t-plus} and zero, as
desired (this vanishing argument, which relies on the grading, is
different from the one in \cite[Proposition 2.13]{seidel01}; we use
it because it extends easily to the projective and relative cases).

Since all Lagrangian submanifolds involved are exact, there are
real-valued action functionals for all the Floer cochain groups
involved, so that in the relevant complex \eqref{eq:total-cone} one
can introduce an $\R$-filtration by the values of the action. Its
acyclicity is then proved by a spectral sequence argument. Briefly,
one replaces the $\R$-filtration by a $\Z$-subfiltration, which
divides the action values into intervals of some suitably chosen
small size $\epsilon>0$. After passing to the associated graded space
of the filtration, one finds that the Floer differentials $\mu^1$ and
the terms $\mu^2(h(a_2),a_1)$, $\mu^3(k,a_2,a_1)$ disappear, while
the remaining terms $\mu^2(a_2,a_1)$ and $\mu^2(k,b)$ make the
complex exact, hence acyclic. For this to work, some care needs to be
observed when choosing the map $\tau_L$ and the cochain level
representatives for $k$ and $h$. A more conceptual but slightly more complicated version of this argument is to first introduce another chain complex quasi-isomorphic to \eqref{eq:total-cone}, and then prove that that one is acyclic; this is carried out in detail in \cite[Section 17]{seidel04}, leading to \cite[Corollary 17.17]{seidel04} which is equivalent to our Proposition \ref{th:exact-sequence-1}.

\subsection{}
We now turn to the consequences of Proposition
\ref{th:exact-sequence-1}, combining arguments from \cite{seidel04}
with the notion of negative graded symplectic automorphisms
introduced in Section \ref{sec:negativity}.

\begin{lemma} \label{th:negative-exact}
Suppose that $m_X = 0$, so $M$ is the affine part of a projective
Calabi-Yau $X$. Let $L_1,\dots,L_m \subset M$ be Lagrangian spheres,
and suppose that
\begin{equation}
\taugr_{L_1}\dots \taugr_{L_m} \in \Autgr(M)
\end{equation}
is isotopic, within that group, to a graded symplectic automorphism
$\phigr$ whose extension to $X$ is negative (see Remark
\ref{th:extend-grading}). Then $L_1^\br,\dots,L_m^\br$ are
split-generators for $D^\pi\Fuk(M)$.
\end{lemma}

\proof Let $L^\br$ be any exact Lagrangian brane. From Lemma
\ref{th:maslov-growth} we know that if $d$ is sufficiently large, a
small exact perturbation $L'$ of $L$ will cause the Maslov index of
any intersection point $x \in L' \cap \phi^d(L)$ to be strictly
negative. By definition of the Floer cochain space
\eqref{eq:basic-floer-complex}, it follows that
\begin{equation}
\mathit{Hom}_{H^0\Fuk(M)}(L^\br,\phigr^d(L^\br)) =
\mathit{HF}^0_M(L^\br,\phigr^d(L^\br)) = 0.
\end{equation}
On the other hand, $\phigr^d(L^\br)$ is isomorphic in $D^b\Fuk(M)$ to
$(T_{L_1^\br} \dots T_{L_m^\br})^d(L^\br)$. Lemma
\ref{th:twist-generators} says that $L^\br$ lies in the subcategory
of $D^\pi\Fuk(M)$ split-generated by the $L_k^\br$; by construction
of that category, this implies that these objects are
split-generators. \qed

\begin{remark} \label{th:codimension-n}
In Lemma \ref{th:negative-exact}, one can replace the negativity of
$\phigr$ by the following slightly weaker assumption. There is a
closed subset $\Sigma \subset X$ which is a finite union of real
submanifolds of dimension $<n$, which admits arbitrarily small open
neighbourhoods $W$ satisfying $\phi(W) = W$, and such that $\phigr$
is negative on $X \setminus W$. In that situation, a generic perturbation $L'$ will
be disjoint from $\Sigma$. It is then also disjoint from a sufficiently small $W$, 
and the same will hold for all $\phi^d(L')$, so that the breakdown of negativity
on $\Sigma$ can be disregarded.
\end{remark}

\subsection{}
To adapt Proposition \ref{th:exact-sequence-1} to the case of a
projective Calabi-Yau surface $X$, with $L$ a Lagrangian sphere and
$L_1^\rbr$ a rational brane, two points have to be taken into
account. Firstly, more care needs to be exercised when defining the
invariants derived from pseudo-holomorphic sections, to avoid
bubbling off of holomorphic discs. This is not really new, since we
have dealt with the same problem in the construction of $\Fuk(X)$.
Note that up to shifts, $L$ gives rise to a unique object $L^\rbr$ of
the Fukaya category, by Lemma \ref{th:unique-sphere-2}. For
$\taugr_L^{-1}(L_1^\rbr)$, one needs to choose the regular almost
complex structure which is the $\tau_L$-pullback of that associated
to $L_1^\rbr$. The second modification concerns the filtration
argument. Our ground ring is now $\Lambda_\Q$, and each
pseudo-holomorphic section contributes $\pm q^e$ as in \eqref{eq:qe}.
In the $\R$-filtration of the Floer cochain groups, one now uses
powers of $q$ as a replacement for values of the action functional.
Since only a finite number of rational branes are involved ($L^\rbr$,
$L_0^\rbr$, $L_1^\rbr$), one can actually work over a subfield
$\Lambda_{(1/d)\Z} \subset \Lambda_\Q$ containing only roots of $q$
of some fixed order $d$. An essentially equivalent and more familiar
formulation is to say that all action functionals together have only
a discrete set of periods. This allows one to carry over the previous
argument, with its use of a $\Z$-subfiltration: it suffices to make
sure that $\epsilon \ll 1/d$ (actually, this is the only point in this
entire paper that would have to be modified substantially if one
wished to include non-rational branes in the Fukaya category of $X$). The outcome is

\begin{prop} \label{th:exact-sequence-2}
$T_{L^\br}(\taugr_L^{-1}(L_1^\br))$ is isomorphic to $L_1^\br$ in
$D^b\Fuk(X)$. \qed
\end{prop}

There is also a version for the relative Fukaya category, where $L$
is again a Lagrangian sphere in $M$ and $L_1^\mbr$ an object of that
category. The definition of the relevant version of $(k,h)$ is as in
the projective case, except sections $u$ now contribute with integer
powers of $q$ according to their intersection with $X_\infty$. As for
the acyclicity of \eqref{eq:total-cone}, one filters it by powers of
$q$, and then taking the associated graded space brings one back to
the affine situation, where we already know that acyclicity holds.

\begin{prop} \label{th:exact-equence-3}
$T_{L^\mbr}(\taugr_L^{-1}(L_1^\mbr))$ is isomorphic to $L_1^\mbr$ in
$D^b\Fuk(X,X_\infty)$. \qed
\end{prop}

The argument from Lemma \ref{th:negative-exact} carries over to the
projective and relative cases without any modifications, and so does
Remark \ref{th:codimension-n}. By combining this with Proposition
\ref{th:negativity}, one arrives at the following conclusion:

\begin{cor} \label{th:generates-everything}
Let $X$ be a projective K{\"a}hler threefold which is Fano, with $o_X
= \K_X^{-1}$, equipped with sections $\sigma_{X,0}$,
$\sigma_{X,\infty}$ which generate an almost Lefschetz pencil of
Calabi-Yau surfaces $\{X_z\}$. Assume that $X_0 \cap X_\infty$ has no
rational components. Pick some base point $z_* \in \C \setminus
Critv(\pi_M)$, base path $c_*$, and a distinguished basis of
vanishing cycles $V_{c_1},\dots,V_{c_r} \subset M_{z_*}$. Then the
$V_{c_i}$, equipped with some gradings, are split-generators for all
versions of the split-closed derived Fukaya category:
$D^\pi\Fuk(M_{z_*})$, $D^\pi\Fuk(X_{z_*})$ and $D^\pi\Fuk(X_{z_*},
X_{z_*,\infty})$. \qed
\end{cor}

\subsection{\label{subsec:induction}}
Let $M = X \setminus X_\infty$ be an affine Calabi-Yau of dimension
$n > 1$, equipped with a quasi-Lefschetz pencil $\{X_z\}$ and the
associated holomorphic function $\pi_M: M \rightarrow \C$. Choose a
base point $z_*$ and base path $c_*$. Let $c_1,\dots,c_r$ be a basis
of vanishing paths, and $V_{c_1},\dots,V_{c_r} \subset M_{z_*}$ the
associated vanishing cycles. Make them into objects $V_{c_i}^\br$ of
$\Fuk(M_{z_*})$; any choice of grading is allowed, and the $Spin$
structure is unique except for $n = 2$, where one needs to take the
{\em nontrivial} one (the one that bounds a $Spin$ structure on a disc). 

\begin{remark}
It may be worthwhile to take a closer look at the $n = 2$ situation. The
choice of the nontrivial $Spin$ structure is important in two
ways. First, note that the vanishing cycles do naturally bound discs (Lefschetz 
thimbles) in $M$. Hence, if one wants to define the Fukaya category of
$\pi_M$ and relate it to the Fukaya category of $M_{z_*}$, as in \cite[Section 18]{seidel04},
the nontrivial $Spin$ structure appears automatically. Second, the long exact 
sequence for Dehn twists \cite{seidel01} works
(over coefficient fields of characteristic $\neq 2$) only for this particular choice
of $Spin$ structure on the curve defining the Dehn twist, because of a cancellation
argument explained in \cite[Example 17.15]{seidel04}.
\end{remark}

Let
\begin{equation} \label{eq:directed}
\F^\rightarrow = \Fuk^\rightarrow(V_{c_1}^\br,\dots,V_{c_r}^\br)
\end{equation}
be the directed $A_\infty$-subcategory of $\Fuk(M_{z_*})$ consisting
of these objects. We associate to each embedded vanishing path $c$,
subject to the conditions $im(c) \cap im(c_*) = \{z_*\}$ and $\Rgeq
c'(0) \neq \Rgeq c_*'(0)$, an object $D_c$ of $D^b(\F^\rightarrow)$,
unique up to isomorphism and shifts, according to the following
rules:
\begin{itemize}
\item
to $c_i$ belongs the object $D_{c_i} = V_{c_i}^\br$ of
$\F^\rightarrow$ itself;
\item
an isotopy of $c$ within the class of such paths does not affect
$D_c$;
\item if we have vanishing paths $c,c',c''$ as in Figure
\ref{fig:hurwitz}, then $D_{c'} \iso T_{D_c}(D_{c''})$ and $D_{c''}
\iso T^\vee_{D_c}(D_{c'})$.
\end{itemize}
From this and the elements of the theory of mutations
\cite{rudakov90}, one sees that all the $D_c$ are exceptional
objects, which means that $\mathit{Hom}^*_{D^b(\Fuk^\rightarrow)}(D_c,D_c)$ is one-dimensional in
degree zero, and zero in all other degrees.

Next, we associate objects $S_d$ of $D^b(\Fuk^\rightarrow)$ to
matching paths which are disjoint from $c_*$. For this, take paths
$c',c'',d$ as in Figure \ref{fig:triangle}. Suppose that the
vanishing cycles associated to $c',c''$ are isotopic, so that $d$ is
a matching path. The objects associated to $c',c''$ will satisfy
\begin{equation} \label{eq:matching-objects}
 \mathit{Hom}^k_{D^b(\Fuk^\rightarrow)}(D_{c'},D_{c''}) =
 \begin{cases} \C & k = l,l+n-1, \\
 0 & \text{otherwise} \end{cases}
\end{equation}
for some $l \in \Z$. After a shift, we may assume that $l = 0$, and
then the rule is:
\begin{itemize}
\item
$S_d$ is the mapping cone $\{D_{c'} \rightarrow D_{c''}\}$ over some
nontrivial morphism.
\end{itemize}
An easy computation shows that $\mathit{Hom}^*_{D^b(\Fuk^\rightarrow)}(S_d,S_d)$ is one-dimensional in
degrees $0$ and $n$, and zero otherwise. This resemblance to the
cohomology of an $n$-sphere is not accidental:
\FIGURE{fig:hurwitz}{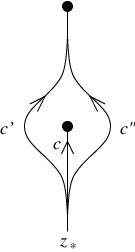}{ht}
\FIGURE{fig:triangle}{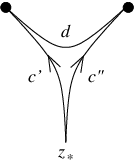}{ht} %

\begin{theorem} \label{th:induction}
Let $d_1,\dots,d_m$ be matching paths which are disjoint from our
base path $c_*$. Assign to them objects $S_{d_1},\dots,S_{d_m}$ in
$D^b(\Fuk^\rightarrow)$ as explained above. Then the full
$A_\infty$-subcategory of $\mathit{Tw}\Fuk^\rightarrow$ consisting of these
objects is quasi-isomorphic to the full subcategory of the Fukaya
category of the total space, $\Fuk(M)$, whose objects are the
matching cycles $\Sigma_{d_1},\dots,\Sigma_{d_m}$, with some choice
of gradings. \qed
\end{theorem}

This way of stating the result leaves several things implicit. Since
one can get any vanishing path by starting with $c_1,\dots,c_r$ and
applying the moves shown in Figure \ref{fig:triangle}, our rules
determine the objects $D_c$. However, it is not a priori clear that
the assignment $c \mapsto D_c$ is unambiguously defined, since one
can get the same path by different sequences of moves. Similarly,
\eqref{eq:matching-objects} and the well-definedness of $S_d$ are not
obvious. The fact that all of this works is a nontrivial part of the
statement. For proofs see \cite[Section 18]{seidel04}.

\subsection{\label{subsec:combinatorial}}
The usefulness of Theorem \ref{th:induction} is particularly obvious
for $n = 2$, when the vanishing cycles are curves on a surface, since
then $\Fuk^\rightarrow$ can be computed in a purely combinatorial
way, up to quasi-isomorphism (versions of this fact have been
discovered independently in various contexts, see e.g.\
\cite{desilva98,polishchuk-zaslow98,chekanov99}). We will now outline
how to do this computation, referring to the papers just quoted and
\cite[Section 13]{seidel04} for further details and justification.

{\em Exactness:} For simplicity, we use a compact piece $N \subset
M_{z_*}$ which is a surface with boundary, capturing all of the
topology of $M_{z_*}$ (the inclusion is a homotopy equivalence). Draw
simple closed curves $\nu_k$ on $N$ which are isotopic to $V_{c_k}$
and in general position (any two intersect transversally, and there
are no triple intersections). One should then verify that these
curves are exact Lagrangian submanifolds for some choice of one-form
$\theta_N$ with $d\theta_N > 0$. For this, consider the decomposition
of $N$ given by the graph $\Gamma = \bigcup_i \nu_i$. Suppose that
one can associate to each component $Z \subset N \setminus \Gamma$
which is disjoint from $\partial N$ a positive weight $w(Z) > 0$,
with the following property: if $\sum_k m_k \bar{Z}_k$, $m_k \in \Z$,
is a two-chain whose boundary is a linear combination of the $\nu_k$,
then
\begin{equation} \label{eq:zero-energy}
\sum_k m_k w(Z_k) = 0.
\end{equation}
One needs to check this only for a finite number of chains, which
form a basis for the homology relations between the $\nu_k$. If this
works out, one can find a symplectic form $\o_N$ with $\int_Z \o_N =
w(Z)$, and that can be written as $\o_N = d\theta_N$ with
$\int_{\nu_k} \theta_N = 0$ for all $k$.

\begin{remark} \label{th:painting-by-numbers}
This ``painting by numbers'' step is rather tedious, so we will
mention a way of reducing the amount of work. Relations in homology
between the vanishing cycles correspond to elements of $H_2(M)$
(because $H_2(M_{z_*}) = 0$, and $M$ is obtained topologically by
attaching two-cells to the vanishing cycles). It is sufficient to
check \eqref{eq:zero-energy} for relations corresponding to elements
in $im(\pi_2(M) \rightarrow H_2(M))$. Equivalently, these are
relations obtained from maps $\Sigma \rightarrow M_{z_*}$ where
$\Sigma$ is a genus zero surface with boundary, each boundary
component getting mapped to some $\nu_i$. This weaker condition does
not guarantee that the vanishing cycles can be made exact, but it
allows one to define the directed Fukaya category, and moving the
$\nu_i$ within that class does not affect the quasi-isomorphism type
of the category. Hence, after a (possibly non-Hamiltonian) isotopy,
one sees that the directed Fukaya category is the same as for the
exact choice of vanishing cycles. The appearance of $\pi_2(M)$ should
remind the reader of the assumption $\o|\pi_2(M) = 0$ used frequently
in symplectic topology, and indeed it addresses the same issue
(getting a bound on the energy of pseudo-holomorphic polygons).
\end{remark}

{\em Gradings:} Recall that the holomorphic one-form $\eta_{M_{z_*}}$
defines a phase function $\alpha_{M_{z_*}}$. One can associate to
this an unoriented foliation of $M_{z_*}$, consisting of those
tangent lines $\Lambda$ for which $\alpha_{M_{z_*}}(\Lambda) = 1$. We
draw some unoriented foliation $\Fol_N$ on $N$ which is homotopic to
the one determined by $\alpha_{M_{z_*}}$. To equip the $\nu_k$ with
gradings, one has to choose a point $b_k \in \nu_k$ which should not
be an intersection point with any other $\nu_l$, and a homotopy
between $\Fol_N$ and $T\nu_k$ at that point (within the $\RP{1}$ of
lines in $TN_{b_k}$). This determines the Maslov indices $I(x)$ of
all intersection points $x \in \nu_i \cap \nu_j$, $i < j$, and also
orientations of the $\nu_i$. Note that $(-1)^{I(x)+1}$ is the local
intersection number $(\nu_i \cdot \nu_j)_x$.

{\em $Spin$ structures:} For every $x$ as before with $(\nu_i \cdot
\nu_j)_x = +1$, we choose a nonzero tangent vector $\xi_x$ to $\nu_j$
at $x$. Of course, one can take these to be the vectors pointing
along the orientation, but that is not strictly speaking necessary.
Next, recall that we have to equip the $\nu_i$ with nontrivial $Spin$
structures, or what is the same, with nontrivial double covers.
Having already chosen points $b_k$, the simplest way to do this is to
take the double cover which is trivial on $\nu_k \setminus b_k$, and
whose two sheets get exchanged when passing over that point.

{\em Counting polygons:} We now introduce a directed
$A_\infty$-category $\Cat^\rightarrow$ with objects
$\nu_1,\dots,\nu_r$. The nontrivial morphism spaces are defined as usual
\begin{equation}
 \mathit{hom}_{\Cat^\rightarrow}^k(\nu_{i_0},\nu_{i_1}) = \bigoplus_{I(x) = k} \C_x,
 \quad i_0<i_1,
\end{equation}
but where now $\C_x$ is identified with $\C$ in a fixed way. Given
$i_0 < \dots < i_d$ in $\{1,\dots,r\}$ and intersection points $x_0
\in \nu_{i_0} \cap \nu_{i_d}$, $x_1 \in \nu_{i_0} \cap \nu_{i_1}$,
\dots, $x_d \in \nu_{i_{d-1}} \cap \nu_{i_d}$, one considers oriented
immersions $u: \Delta \rightarrow N$, where $\Delta$ is some model convex
$(d+1)$-gon in the plane with numbered sides and corners. Under $u$,
the sides should go to the $\nu_{i_k}$, and the corners to the $x_k$;
moreover, the immersion should be a local diffeomorphism up to the
corners, so that the image has locally an angle $<\pi$. One considers
two such $u$ to be equivalent if they differ by a diffeomorphism of
$\Delta$. Each equivalence class contributes to
$\mu^k_{\Cat^\rightarrow}$ by
\begin{equation}
 \pm 1: \C_{x_d} \otimes \dots \otimes \C_{x_1} \rightarrow \C_{x_0},
\end{equation}
and the sign is the product of the following factors: (i) for each $k
\geq 1$ such that $(\nu_{i_{k-1}} \cdot \nu_{i_k})_{x_k} = 1$, we get
a $-1$ iff the tangent vector $\xi_{x_k}$ points away from
$u(\Delta)$. (ii) We get another $-1$ if $(\nu_{i_0} \cap
\nu_{i_d})_{x_0} = 1$ and the tangent vector $\xi_{x_0}$ points
towards $u(\Delta)$. (iii) Every time that $u|\partial\Delta$ passes
over a point $b_k$, this adds a further $-1$.

\begin{lemma}
$\Cat^{\rightarrow}$ is quasi-isomorphic to the directed subcategory
$\Fuk^{\rightarrow}$ from \eqref{eq:directed}, at least up to a shift
in the gradings of the $V_{c_i}^\br$. \qed
\end{lemma}

Starting from this, Theorem \ref{th:induction} allows one to
determine the full $A_\infty$-sub\-ca\-te\-gory of $\Fuk(M)$ consisting of
any finite number of matching cycles, and in particular the Floer
cohomology of any two such cycles. Note that matching cycles are
Lagrangian two-spheres in a four-manifold, so their Floer cohomology
is not amenable to computation in any direct way.

\subsection{\label{subsec:addition}}
We will need a minor addition to the previously explained
computational method, which takes into account the presence of a
homomorphism $\rho: \pi_1(M) \rightarrow \Gamma$ into a finitely
generated abelian group $\Gamma$. Restriction yields a homomorphism
from the fundamental group of the fibre $M_{z_*}$, and all vanishing
cycles $V_{c_i}$ bound balls in $M$, hence the composition
$\pi_1(V_{c_i}) \rightarrow \pi_1(M_{z_*}) \rightarrow \Gamma$ is
zero. As discussed in Section \ref{subsec:coverings}, after choosing
lifts of the vanishing cycles to the associated $\Gamma$-cover, one has a $\Gamma^*$-action on \eqref{eq:directed}. We can use the same rules as before to associate to each vanishing
path $c$ an object $D_c$ of the equivariant derived category
$D^b_{\Gamma^*}(\Fuk^\rightarrow)$, see Remark \ref{th:twist-is-equivariant}. For
$c',c'',d$ are as in Figure \ref{fig:triangle}, the objects $D_{c'}$ and $D_{c''}$ will now satisfy
\begin{equation}
 \mathit{Hom}^k_{D^b_{\Gamma^*}(\Fuk^\rightarrow)}(D_{c'},D_{c''}) \iso \begin{cases} W & k = l,l+n-1, \\
 0 & \text{otherwise,} \end{cases}
\end{equation}
where $W$ is some one-dimensional $\Gamma^*$-representation. After
shifting $D_{c'}$ and tensoring it with $W^*$, one may assume that $l
= 0$ and $W$ is trivial. Then the cone $S_d = \{D_{c'} \rightarrow
D_{c''}\}$ over a nontrivial degree zero morphism is again an object
of the equivariant derived category. The equivariant version of
Theorem \ref{th:induction} says that the full $A_\infty$-subcategory
of $\mathit{Tw}_{\Gamma^*}\Fuk^\rightarrow$ with objects
$S_{d_1},\dots,S_{d_m}$ is $\Gamma^*$-equivariantly quasi-isomorphic
to the full subcategory of $\Fuk(M)$ with objects
$\Sigma_{d_1},\dots,\Sigma_{d_m}$, where we assume that gradings and
lifts to $M_\Gamma$ have been chosen in some appropriate way for
these matching cycles. For the proof, which is not particularly
difficult, there are two essentially equivalent viewpoints: one
either goes through the relevant material in \cite{seidel04}, making
sure that everything can be made compatible with the
$\Gamma^*$-actions; or one works directly with the associated $\Gamma$-covering of $M$ and pulls back the holomorphic function to it.

\section{The algebras $Q_4$ and $Q_{64}$\label{sec:4-64}}

Using split-generators, the categories which occur in many instances of homological
mirror symmetry can be described in terms of $A_\infty$-structures on
finite-dimensional algebras. This section introduces the algebra
relevant for our computation, as well as a simpler version which
serves as a stepping-stone. The main point is to use the general
theory from Section \ref{sec:deformations} to partially classify
$A_\infty$-structures on those algebras. We end by applying the
results to coherent sheaves on the mirror of the quartic surface:
this finishes work on the algebro-geometric side of mirror symmetry,
and it also helps to clarify the amount of things still to be proved
on the symplectic side.

\subsection{\label{subsec:q-preliminary}}
We begin by recalling a general ``doubling'' construction. Let
$C^\rightarrow$ be a directed ($\C$-linear graded) category, with
objects $\{X_1,\dots,X_m\}$. The {\em trivial extension category} of
degree $d > 0$ is a category $C$ with the same objects, and it
satisfies:
\begin{itemize} \itemsep1em
\item
The directed subcategory of $C$ is the given $C^\rightarrow$.
\item
Each object of $C$ is spherical, meaning that $\mathit{Hom}_C^*(X_j,X_j)$ has total dimension $2$.
\item
$C$ is a Frobenius category of degree $d$, in the sense that there
are linear maps $\int_{X_j}: \mathit{Hom}_C^d(X_j,X_j) \rightarrow \C$, such
that
\begin{equation} \label{eq:int-pairing}
\leftsc a,b \rightsc = \textstyle \int_{X_j} ab
\end{equation}
for $a \in \mathit{Hom}_C^*(X_k,X_j)$, $b \in \mathit{Hom}_C^*(X_j,X_k)$ is a
nondegenerate pairing. These pairings are graded symmetric, $\leftsc
a,b \rightsc = (-1)^{\mathrm{deg}(a)\mathrm{deg}(b)} \leftsc b,a \rightsc$.
\item
One can split $\mathit{Hom}_{C}^*(X_j,X_k) = \mathit{Hom}_{C^\rightarrow}^*(X_j,X_k)
\oplus J_{kj}$ in such a way that the product of two morphisms
belonging to the $J$ spaces is zero (the last-mentioned property is automatic if the morphisms in $C^\rightarrow$ are all concentrated in degree $0$).
\end{itemize}
These properties characterize $C$ up to isomorphism. Alternatively,
one can define it explicitly by
\begin{equation} \label{eq:double}
 \mathit{Hom}_C(X_j,X_k) = \mathit{Hom}_{C^\rightarrow}(X_j,X_k)
 \oplus \mathit{Hom}_{C^\rightarrow}(X_k,X_j)^\vee[-d]
\end{equation}
and $(a,a^\vee)(b,b^\vee) = (ab,a^\vee(b\,\cdot) +
(-1)^{\mathrm{deg}(a)(\mathrm{deg}(b)+\mathrm{deg}(b^\vee))} b^\vee(\cdot\, a))$.

The example relevant to us is defined as follows: let $V$ be a four-dimensional vector space. Consider the (linear graded) directed category $C^{\rightarrow}_4$ having four objects $X_1,\dots,X_4$ and morphisms
\begin{equation}
 \mathit{Hom}_{C^\rightarrow_4}(X_j,X_k) = \begin{cases}
 \Lambda^{k-j}V & j \leq k, \\
 0 & \text{otherwise,} \end{cases}
\end{equation}
with the obvious grading and (wedge product) composition. We define
$C_4$ to be the trivial extension category of $C_4^\rightarrow$ in
degree $d = 2$. Its total morphism algebra, linear over $R_4 \iso
\C^4$, will be denoted by $Q_4$. Once one has fixed a choice of
volume element in $\Lambda^4 V^\vee$, the duals in the general
formula \eqref{eq:double} can be removed, so that then
\begin{equation} \label{eq:c4-morphism-revised}
 \mathit{Hom}_{C_4}(X_j,X_k) = \begin{cases}
 \Lambda^{k-j}V & j < k, \\
 \Lambda^0V \oplus (\Lambda^4V)[2], & j = k, \\
 (\Lambda^{k-j+4}V)[2] & j > k.
 \end{cases}
\end{equation}

From a different viewpoint, take the exterior algebra $\Lambda V$,
and consider the action of the center $\Gamma_4 \subset \mathit{SL}(V)$,
which is generated by $\gamma = i \cdot \mathit{id}_V$. From our general discussion of semidirect products, we know that $\Lambda V \semidirect \Gamma_4$ is an algebra over $\C
\Gamma_4 \iso R_4$, with a basis of idempotents given by $e_j =
\quarter\chi_j = \quarter(e + i^{-j} \gamma + i^{-2j} \gamma^2 +
i^{-3j} \gamma^3)$ for $j = 1,\dots,4$. Concretely,
\begin{equation} \label{eq:c4-chi}
 e_k (\Lambda V \semidirect \Gamma_4) e_j =
 \begin{cases}
 \Lambda^{k-j} V \otimes \C e_j & j<k, \\
 (\Lambda^0 V \oplus \Lambda^4 V) \otimes \C e_j & j = k, \\
 \Lambda^{k-j+4} V \otimes \C e_j & j>k. \\
 \end{cases}
\end{equation}
So $\Lambda V \semidirect \Gamma_4$ is the total morphism algebra of
a category with four objects. If we view it in this way, the directed
subcategory is isomorphic to $C_4^\rightarrow$. Moreover, after
choosing a volume element on $V$, the resulting integration maps make
$\Lambda V \semidirect \Gamma_4$ into a Frobenius category of degree
4. One can split the morphism spaces into those of the directed
subcategory and a complementary summand $J$, consisting of those
$\Lambda^jV \otimes \C e_k$ with $j+k > 4$. This satisfies the
properties stated above, so we conclude that $\Lambda V \semidirect
\Gamma_4$ is the trivial extension category of $C_4^\rightarrow$ in
degree 4. Hence there is an isomorphism $\Lambda V \semidirect
\Gamma_4 \iso Q_4$ of algebras, which is $R_4$-linear but only
$\Z/2$-graded. To make the matter of gradings slightly more precise,
introduce another version $\tilde{Q}_4 = \Lambda V \semidirect \Gamma_4$
which is bigraded by $\half \Z \times \Z/2$, so that the
$r$-th exterior power has grading $(r/2,r)$. The idea is that we use the $\Z/2$ part as parity, determining the Koszul signs, and the $\half\Z$ part for the actual gradings.
Then:

\begin{lemma}
$Q_4$ is obtained from $\tilde{Q}_4$ by changing the grading of $e_k \tilde{Q}_4 e_j$ by $((k-j)/2,0)$. \qed
\end{lemma}

More precisely, after this change in the bigrading of $\tilde{Q}_4$, the $\half\Z$ part becomes integral, and the $\Z/2$ part is its mod $2$ reduction, so we are indeed left with an ordinary graded algebra, which as such is isomorphic to $Q_4$.

In particular, the Hochschild cohomology of $Q_4$ and $\tilde{Q}_4$ are isomorphic. It is straightforward to adapt the proof of Proposition \ref{th:twisted-hkr} to this situation, which yields
\begin{equation} \label{eq:fractional-hkr}
\begin{aligned}
 & \mathit{HH}^{s+t}(Q_4,Q_4)^t \iso
 \bigoplus_{\gamma \in \Gamma_4} \Big( \mathit{Sym}^s(V^{\gamma,\vee}) \otimes
 \Lambda^{s+2t-\mathrm{codim}(V^\gamma)}(V^\gamma) \\[-1.5em] & \hspace{20em} \otimes \, \Lambda^{\mathrm{codim}(V^\gamma)}(V/V^\gamma) \Big)^{\Gamma_4}.
\end{aligned}
\end{equation}
Note that the summands with nontrivial $\gamma$ only contribute to $(s,t) = (0,2)$. In particular,
\begin{equation} \label{eq:written}
 \mathit{HH}^2(Q_4,Q_4)^t = \begin{cases} \mathit{Sym}^4(V^\vee) & t = -2, \\
 0 & \text{all other $t<0$.}
 \end{cases}
\end{equation}

The group of (graded $R_4$-algebra) automorphisms of $Q_4$ is
\begin{equation} \label{eq:aut-q4}
\mathit{Aut}(Q_4) \iso I_4 \times_{\Gamma_4} \mathit{GL}(V).
\end{equation}
Here, $I_4 = (\C^*)^4/\C^*$, where the quotient is by the diagonal $\C^*$. We embed $\Gamma_4$ into $I_4$ by taking $i \cdot \mathit{id}_V$ to $[1,i,i^2,i^3]$. The notation in \eqref{eq:aut-q4} means that we divide the product $I_4 \times \mathit{GL}(V)$ by the normal subgroup generated by $([1,i^{-1},i^{-2},i^{-3}], i \cdot \mathit{id}_V)$. This describes it as an abstract group. As for its action on $Q_4$, $I_4$ is a group of inner automorphisms, with $\alpha = [\alpha_1,\dots,\alpha_4] \in I_4$ acting by multiplying $e_k Q_4 e_j$ by $\alpha_k/\alpha_j$. As explained before, it is a general feature of semidirect products that the action of $\Gamma_4 \subset \mathit{GL}(V)$ becomes inner on $Q_4$. Indeed, the diagonal matrix $i \cdot \mathit{id}_V$ acts in the same way as $\alpha = [1,i,i^2,i^3]$, which is why we get an overall action of \eqref{eq:aut-q4}.
From now on we fix an isomorphism $V \iso \C^4$, which at the same time singles out the maximal torus of diagonal matrices $H \subset \mathit{SL}(V)$. Let $T = H/\Gamma_4$. The inclusion $T \rightarrow \mathit{GL}(V)/\Gamma_4$ can be lifted to a homomorphism $\rho_T: T \rightarrow \mathit{Aut}(Q_4)$, by letting $[t_0,t_1,t_2,t_3]$ act by the combination of $\alpha = [1,t_0^{-1},t_0^{-2},t_0^{-3}] \in I_4$ and $\mathit{diag}(t_0,t_1,t_2,t_3) \in \mathit{GL}(V)$. In view of \eqref{eq:written}, the induced action on \eqref{eq:fractional-hkr} satisfies
\begin{equation} \label{eq:written-2}
 \big( \mathit{HH}^2(Q_4,Q_4)^{-2} \big)^T = \C \cdot y_0y_1y_2y_3.
\end{equation}
Now suppose that $\QQ_4$ is an $A_\infty$-algebra structure on $Q_4$, satisfying \eqref{eq:mu1zero} and which is compatible with the structure of an $R_4$-bimodule as well as with the action of $T$. Finally, assume that $\QQ_4$ is not formal (not isomorphic to the trivial $A_\infty$-structure). We have no need of proving the existence of $\QQ_4$ algebraically, since it will be constructed geometrically later on from coherent sheaves, and also again as a Fukaya category. The uniqueness statement is:

\begin{lemma} \label{th:recognize-q4}
$\QQ_4$ is unique up to ($R_4$-linear, $T$-equivariant) $A_\infty$-isomorphism.
\end{lemma}

\proof We will use the deformation theory from Section \ref{sec:deformations} in a version which is equivariant with respect to $T$. This is unproblematic, since we can consider the action as an additional grading by $T^* \iso \Z^3$, and equivariance amounts to considering $A_\infty$-structures which are homogeneous for that grading (as already mentioned in Remark \ref{th:fg}). The equivariant version of Lemma \ref{th:versal-1} and \eqref{eq:written-2} imply that any other $A_\infty$-structure with the same properties as $\QQ_4$ is gauge equivalent to $\epsilon^*\QQ_4$ for some $\epsilon \in \C^*$, and in particular isomorphic to $\QQ_4$. \qed

\begin{remark} \label{th:more-recognize-q4}
We will need a concrete criterion for recognizing non-formality in this context. Suppose that we have an $A_\infty$-algebra $\A$ linear over $R_4$ and carrying a $T$-action, whose cohomology algebra is
equivariantly isomorphic to $Q_4$. Choose $\lambda$, $\pi$ and $h$ as in Remark \ref{th:explicit-hpt} (and $R_4$-linearly), so that we can run the explicit form of the Homological Perturbation Lemma to
produce an $A_\infty$-structure $\tilde{\A}$ on $H(\A)$. Since $e_k Q_4 e_j$ is nonzero only in degrees $\equiv k-j$ mod 2, any odd degree $R_4$-bimodule map
\begin{equation}
Q_4 \otimes_{R_4} \cdots \otimes_{R_4} Q_4 \longrightarrow Q_4
\end{equation}
is necessarily zero. That applies to all the odd order $\mu^d_{\tilde{\A}}$ and in particular to $d = 3$, so that we get a well-defined class
\begin{equation} \label{eq:mu4-class}
[\mu^4_{\tilde{\A}}] \in \mathit{HH}^2(Q_4,Q_4)^{-2} \iso \mathit{Sym}^4V^\vee.
\end{equation}
Concretely, take nonzero elements $\bar{a}_4 \in e_1 \tilde{\A} e_4$, $\bar{a}_3 \in e_4 \tilde{\A} e_3$, $\bar{a}_2 \in e_3 \tilde{\A} e_2$, $\bar{a}_1 \in e_2 \tilde{\A} e_1$, such that the product of any two successive ones is zero. Plug the $\bar{a}_k$ into the expression \eqref{eq:mu4-hpt} for $\mu^4_{\tilde{\A}}$; if the outcome is nonzero, then $\A$ is ($T$-equivariantly) quasi-isomorphic to $\QQ_4$. To see this, note that by the explicit formula \eqref{eq:reverse-hkr} for the HKR isomorphism, or rather a slight variation of it for semidirect products, $\mu^4_{\tilde{\A}}(\bar{a}_4,\bar{a}_3,\bar{a}_2,\bar{a}_1) = p(v,v,v,v) \otimes e_1$, where $p \in \mathit{Sym}^4V^\vee$ is the polynomial corresponding to the obstruction class \eqref{eq:mu4-class}, and $v \in V$ is a nonzero vector reflecting the choice of the $\bar{a}_k$.
\end{remark}

\subsection{}
Let $\Gamma_{64} \subset H \subset \mathit{SL}(V)$ be the subgroup of diagonal matrices which have order $4$, and $\Gamma_{16} = \Gamma_{64}/\Gamma_4$ its image in $T$. We use the restriction of $\rho_T$ to that finite group to form the semidirect product
\begin{equation}
 Q_{64} \stackrel{\text{def}}{=} Q_4 \semidirect \Gamma_{16}.
\end{equation}
This is an algebra linear over $\C\Gamma_4 \otimes \C\Gamma_{16} = \C\Gamma_{64}$, so one can also view it as a category with 64 objects.

\begin{remark}
The definition of $\rho_T$ involved a choice of lift $T \rightarrow \mathit{Aut}(Q_4)$. If we take another choice differing from the given one by a homomorphism $T \rightarrow Q_4^\times \iso (\C^*)^4$, the result of restricting that to $\Gamma_{16}$ is an action differing from the given one as in \eqref{eq:conjugate-action}. Hence, $\Gamma_{64}$ would remain the same.
\end{remark}

Explicitly, $\rho_T|\Gamma_{16}$ lets $[t_0,t_1,t_2,t_3]$ act by $\mathit{diag}(1,t_0^{-1}t_1,t_0^{-1}t_2,t_0^{-1}t_3) \in \mathit{SL}(V) \subset \mathit{Aut}(Q_4)$, which is one of the possible splittings of the projection $\Gamma_{64} \rightarrow \Gamma_{16}$. By applying \eqref{eq:in-stages} one finds that there is a $\Z/2$-graded isomorphism
\begin{equation} \label{eq:z2-graded-64}
 Q_{64} \iso (\Lambda V \semidirect \Gamma_4) \semidirect \Gamma_{16}
 \iso \Lambda V \semidirect \Gamma_{64}.
\end{equation}
As before, we can refine this to take gradings better into account, and then apply a version of the Hochschild-Kostant-Rosenberg isomorphism to conclude that
\begin{equation} \label{eq:fractional-hkr-2}
\begin{aligned}
 & \mathit{HH}^{s+t}(Q_{64},Q_{64})^t \iso
 \bigoplus_{\gamma \in \Gamma_{64}} \Big( \mathit{Sym}^s(V^{\gamma,\vee}) \otimes
 \Lambda^{s+2t-\mathrm{codim}(V^\gamma)}(V^\gamma) \\[-1.5em] & \hspace{20em} \otimes \, \Lambda^{\mathrm{codim}(V^\gamma)}(V/V^\gamma) \Big)^{\Gamma_{64}}.
\end{aligned}
\end{equation}
Case-by-case inspection shows that the $\gamma \neq e$ summands do not contribute to the terms in \eqref{eq:fractional-hkr-2} which have $s+t = 2$, with the exception of $(s = 0, t = 2)$.

Let $\QQ_4$ be as in Lemma \ref{th:recognize-q4}, which in particular means that it is equivariant for the action of $T$. By restricting that action to $\Gamma_{16}$, one defines an $A_\infty$-structure
\begin{equation}
\QQ_{64} \stackrel{\text{def}}{=} \QQ_4 \semidirect \Gamma_{16},
\end{equation}
whose underlying cohomology level algebra is $Q_{64}$. This carries an action of the dual group $\Gamma_{16}^*$, which can be extended to an action of $\Gamma_{64}^*$ if one is willing to reduce the grading to $\Z/2$.

\begin{lemma} \label{th:qq64-trunc}
The truncated Hochschild cohomology of $\QQ_{64}$ satisfies
\begin{equation} \label{eq:hh12}
 \mathit{HH}^2(\QQ_{64},\QQ_{64})^{\scriptscriptstyle\leq 0} \iso \C^7.
\end{equation}
\end{lemma}

\proof We use the spectral sequence \eqref{eq:truncated-ss}. Since the differentials in the spectral sequence preserve the action of $\Gamma_{64}^*$, we can concentrate on the trivial $\gamma$ summand in \eqref{eq:fractional-hkr-2}. Here is the relevant part of the $E_2$ term ($S^j\Lambda^k$ is shorthand for $(\mathit{Sym}^j V^\vee \otimes \Lambda^k V )^{\Gamma_{64}}$, and the crossed out summands vanish because the $\Gamma_4$-fixed point sets are already trivial):
\begin{equation}
\begin{array}{r|ccccccc}
 & s = 0 & s = 1 & s = 2 & s = 3 & s = 4 & s = 5 & s = 6
 \\[0.25em]
 \hline
 t = 1 &
 0
 & \xout{S^1\Lambda^3}
 & \xout{S^2\Lambda^4}
 & 0
 & 0
 & 0
 & 0
 \\[0.25em]
 t = 0 &
 0
 & S^1\Lambda^1
 & S^2\Lambda^2
 & S^3\Lambda^3
 & S^4\Lambda^4
 & 0
 & 0
 \\[0.25em]
 t = -1 &
 0
 & 0
 & \xout{S^2\Lambda^0}
 & \xout{S^3\Lambda^1}
 & \xout{S^4\Lambda^2}
 & \xout{S^5\Lambda^3}
 & \xout{S^6\Lambda^4}
 \\[0.25em]
 t = -2 &
 0
 & 0
 & 0
 & 0
 & S^4\Lambda^0
 & S^5\Lambda^1
 & S^6\Lambda^2
 \\[0.25em]
 t = -3 &
 0
 & 0
 & 0
 & 0
 & 0
 & 0
 & \xout{S^6\Lambda^0}
\end{array}
\end{equation}
The first nonzero differential is $\delta_3$, which is the Schouten
bracket with the obstruction class \eqref{eq:md-def-class}. We identify $V = \C^4$
as before, so that the obstruction class is $y_0 y_1 y_2 y_3$ (up to
a multiple, which we may safely ignore since it can be absorbed into
a rescaling of $V$). In total degree $s+t = 1$, we have the
$\Gamma_{64}$-invariant part of $V^\vee \otimes V$ which is spanned
by elements $y_k \otimes v_k$, and these satisfy
\begin{equation} \label{eq:d311}
 \delta_3^{1,0}(y_k \otimes v_k) = y_0 y_1 y_2 y_3
\end{equation}
for all $k$. Note that $\half \sum_k y_k \otimes v_k$ is what we
called the Euler element in Section \ref{sec:deformations}, and so
the fact that its image under $\delta_3$ is twice the obstruction
class is an instance of the general property
\eqref{eq:euler-element}. For $s+t = 2$ we have firstly the
$\Gamma_{64}$-invariant part of $S^2V^\vee \otimes \Lambda^2V$,
generated by the six elements $y_j y_k \otimes v_j \wedge v_k$ which
satisfy
\begin{equation}
 \delta_3^{2,0}(y_j y_k \otimes v_j \wedge v_k) =
 (y_0 y_1 y_2 y_3) y_k \otimes v_k -
 (y_0 y_1 y_2 y_3) y_j \otimes v_j;
\end{equation}
and secondly the invariant part of $S^4V^*$, whose basis consists of
the $y_k^4$ together with $y_0y_1y_2y_3$, the latter being
annihilated by being in the image of \eqref{eq:d311}. The structure
of the $E_2$ term shows that the differentials of order $d>3$ all
vanish, so our computation is complete. To summarize, explicit
generators on the $E_2$ page which represent \eqref{eq:hh12} are: $y_0y_1 \otimes v_0
\wedge v_1 + y_1y_2 \otimes v_1 \wedge v_2 + y_2y_0 \otimes v_2
\wedge v_0$ and its cyclic permutations (the sum of all being zero),
together with the monomials $y_k^4$, $k = 0,\dots,3$. \qed

\subsection{\label{subsec:u-symmetry}}
Let $U_4 \in \mathit{GL}(V)$ be the matrix which cyclically permutes the coordinates to the right. This generates an additional action of $\Z/4$ on $Q_4$. Moreover, it can be lifted to a $\Z/4$-action on $Q_{64}$, which is generated by
\begin{equation} \label{eq:u64}
U_{64}( a \otimes [t_0,t_1,t_2,t_3] ) = U_4(a) \otimes [t_3,t_0,t_1,t_2]
\end{equation}
for $a \in Q_4$, $[t_0,t_1,t_2,t_3] \in \Gamma_{16}$.

\begin{remark} \label{th:fake-u4}
Occasionally, we will face the situation where we have an additional automorphism $\tilde{U}_4$ of $Q_4$ which is not explicitly given, but is known to have the following properties:
\begin{itemize} \itemsep1em
\item $\rho_T([t_0,t_1,t_2,t_3]) \circ \tilde{U}_4 = \tilde{U}_4 \circ \rho_T([t_1,t_2,t_3,t_0])$.
\item for any $k$, $\tilde{U}_4$ acts as $-1$ on the degree $2$ part $(e_k Q_4 e_k)^2 \iso \Lambda^4 V$;
\item $\tilde{U}_4^4$ is the identity.
\end{itemize}
By \eqref{eq:aut-q4}, $\tilde{U}_4$ could be the combination of some $\alpha = [\alpha_0,\dots,\alpha_3] \in I_4$ with a matrix
\begin{equation} \label{eq:delta-matrix}
\textstyle \left(\begin{smallmatrix} &&& \delta_0 \\ \delta_1 &&& \\ & \delta_2 && \\ && \delta_3 & \end{smallmatrix}\right) \in \mathit{GL}(V),
\end{equation}
whose determinant (because of the action on $\Lambda^4V$) is $-\delta_0 \delta_1 \delta_2 \delta_3 = -1$.
This implies that the $\alpha_k$ are fourth roots of unity. We can conjugate $\tilde{U}_4$ with
\begin{equation}
[\delta_0^{-1/4}\delta_1^{1/2}\delta_2^{1/4},
\delta_0^{-1/4}\delta_1^{-1/2}\delta_2^{1/4},
\delta_0^{-1/4}\delta_1^{-1/2}\delta_2^{-3/4},
\delta_0^{3/4}\delta_1^{1/2}\delta_2^{1/4}] \in T
\end{equation}
to reduce the matrix \eqref{eq:delta-matrix} to the case where all $u_i$ are equal to $1$. After this conjugation, composition with $[\alpha_0^{-1},\dots,\alpha_3^{-1}] \in I_4$ yields $U_4$.
\end{remark}

Lemma \ref{th:recognize-q4} holds even if we add the requirement of symmetry with respect to $U_4$. This means that if $\QQ_4$ is invariant under $T \semidirect \Z/4$, then it is unique in a way which preserves that additional symmetry. We assume from now on that this is the case. As a consequence, $\QQ_{64}$ is invariant under $U_{64}$.

\begin{lemma} \label{th:recognize-q64-2}
Suppose that $\QQ_{64,q}$ is a one-parameter deformation of $\QQ_{64}$ which is invariant under $U_{64}$, and whose class in $\mathit{HH}^2(\QQ_{64},\QQ_{64})^{\scriptscriptstyle\leq 0}$ is nontrivial. Such a deformation is unique up to equivalence of $A_\infty$-deformations combined with the action of $\mathit{End}(\Lambda_\N)^\times$.
\end{lemma}

\proof We proceed as in Lemma \ref{th:qq64-trunc} but using only the part of Hochschild cohomology which is fixed by the additional $\Z/4$-action. The action of $U_{64}$ on $\mathit{HH}^*(Q_{64},Q_{64})$ cyclically permutes $y_0y_1 \otimes v_0 \wedge v_1 + y_1y_2 \otimes v_1 \wedge v_2 + y_2y_0 \otimes v_2$ and its three analogues. However, the sum of these generators is zero in $\mathit{HH}^*(\QQ_{64},\QQ_{64})^{\scriptscriptstyle\leq 0}$, hence they do not contribute to the $\Z/4$-invariant part. The remaining four generators are $y_0^4$ and its analogues, which are again cyclically permuted, leaving a one-dimensional fixed part
\begin{equation}
\mathit{HH}^2(\QQ_{64},\QQ_{64})^{\leq 0,\Z/4} \iso \C \cdot (y_0^4 + y_1^4 + y_2^4 + y_3^4).
\end{equation}
A version of Lemma \ref{th:versal-3} which takes finite group actions into account completes the argument. \qed

Finally, here is the way in which we will apply this classification result:

\begin{proposition} \label{th:recognize-q64-3}
Let $\SS_4$ be an $R_4$-linear $A_\infty$-algebra, such that $H(\SS_4) \iso Q_4$. This should carry an action of $T \semidirect \Z/4$, which on the cohomology level agrees with the one on $Q_4$ described above. Moreover, $\SS_4$ should not be formal. Then $\SS_4$ is quasi-isomorphic to $\QQ_4$ in a way which is compatible with all these symmetries.

Now use the restriction of the given $T$-action to form $\SS_{64} = \SS_4 \semidirect \Gamma_{16}$, and lift the $\Z/4$-action to $\SS_{64}$ as in \eqref{eq:u64}. Suppose that $\SS_{64,q}$ is a one-parameter deformation of $\SS_{64}$, which is $\Z/4$-equivariant, and whose deformation class in $\mathit{HH}^2(\SS_{64},\SS_{64})$ is nontrivial. Then $\SS_{64,q}$ is quasi-isomorphic to $\psi^*\QQ_{64,q}$, for some $\psi \in \mathit{End}(\Lambda_\N)^\times$.
\end{proposition}

\proof
Applying the Homological Perturbation Lemma (equivariantly), we can assume without loss of generality that $\SS_4$ is minimal. Moreover, since each of the graded pieces $e_k Q_{64} e_j$ is concentrated in either odd or even degree, $H(\SS_{64,q})$ is necessarily torsion-free, hence we can carry over the deformation to the minimal model.

After these preparations, Lemma \ref{th:recognize-q4} yields the desired $A_\infty$-isomorphism $\SS_4 \iso \QQ_4$, and also $\SS_{64} \iso \QQ_{64}$. Next, Lemma \ref{th:product-turns-on-2} shows that the deformation class is also nontrivial in $\mathit{HH}^2(\SS_{64},\SS_{64})^{\scriptscriptstyle\leq 0}$. Finally, Lemma \ref{th:recognize-q64-2} applies and completes the argument.
\qed

\subsection{}
Take the Beilinson basis of the derived category of coherent sheaves on $P = {\mathbb P}(V)$, which consists of $F_k = \Omega^{4-k}(4-k)[4-k]$ for $k = 1,\dots,4$. The $\Omega^k$ are naturally $\mathit{PGL}(V)$-equivariant sheaves, and if we give $\O(-1)$ the $\mathit{GL}(V)$-action obtained by embedding it into $\O \otimes V$, the $\Omega^{4-k}(4-k) = \Omega^{4-k} \otimes \O(-1)^{\otimes k-4}$
become $\mathit{GL}(V)$-equivariant. A straightforward computation, which is the core of the argument in \cite{beilinson78}, shows that

\begin{lemma} \label{th:beilinson}
The linear graded category with objects $F_1,\dots,F_4$, and where
the morphisms are homomorphisms of all degrees in the derived
category, is $\mathit{GL}(V)$-equiva\-ri\-antly isomorphic to
$C_4^\rightarrow$. \qed
\end{lemma}

Let $\iota_0: Y_0 \hookrightarrow P$ be the
reduced quartic hypersurface given by $Y_0 = \{y_0y_1y_2y_3 = 0\}$. Write $E_{0,k} = \iota_0^*F_k$. Repeating an
argument from \cite{seidel-thomas99}, we observe that
\begin{align*}
 \mathit{Hom}^*_{Y_0}(E_{0,j},E_{0,k})
 & = \mathit{Hom}^*_P(F_j,(\iota_0)_*\iota_0^*F_k)
 = \mathit{Hom}^*_P(F_j,\O_{Y_0} \stackrel{{\bf L}}{\otimes} F_k) \\
 & = \mathit{Hom}^*_P(F_j,\{\K_P \otimes F_k \rightarrow \O_P \otimes F_k\})
\end{align*}
sits in a long exact sequence which one can write, using Serre
duality on $P$, as
\begin{equation}
 \cdots \rightarrow \mathit{Hom}^d_P(F_j,F_k) \stackrel{\iota_0^*}{\longrightarrow}
 \mathit{Hom}^d_{Y_0}(E_{0,j},E_{0,k}) \stackrel{(\iota_0^*)^\vee}\longrightarrow
 \mathit{Hom}^{2-d}_P(F_k,F_j)^\vee \rightarrow \cdots
\end{equation}
Since the $F_k$ form an exceptional collection, we have unique
splittings
\begin{equation} \label{eq:doubling}
 \mathit{Hom}^*_{Y_0}(E_{0,j},E_{0,k}) = \mathit{Hom}^*_P(F_j,F_k) \oplus \mathit{Hom}^*_P(F_k,F_j)^\vee[-2].
\end{equation}
$Y_0$ has trivial dualizing sheaf, so Serre duality makes
$D^b\mathit{Coh}(Y_0)$ into a Frobenius category of degree $d = 2$. If we
denote by $J_{jk}$ the second summand in \eqref{eq:doubling}, then
the composition of two morphisms in $J$ spaces is zero for degree
reasons (these morphisms are $\mathit{Ext}^2$'s of locally free sheaves on
$Y_0$). We have now checked off the properties of a trivial extension
category. Moreover, the induced $\mathit{GL}(V)$-action on
$\mathit{Ext}^2_{Y_0}(E_{0,k},E_{0,k}) \iso \C$ equals the determinant, and we therefore conclude that:

\begin{lemma}
The linear graded category with objects $E_{0,k}$, and where the
morphisms are maps of arbitrary degree in $D^b\mathit{Coh}(Y_0)$, is
$\mathit{GL}(V)$-equiva\-riantly isomorphic to $C_4$. \qed
\end{lemma}

Denote by $S_4$ the total morphism algebra of that subcategory, so
$S_4 \iso Q_4$. This is the cohomology algebra of a dg algebra
$\SS_4$, canonical up to quasi-isomorphism, which can be defined as
total dg algebra of the same objects in the {\v C}ech dg category
$\cech{\SS}(Y_0)$. The next step comes from
\cite{douglas-govindarajan-jayaraman-tomasiello01}:

\begin{lemma} \label{th:massey}
$\SS_4$ is not formal.
\end{lemma}

\proof The skyscraper sheaf at any point $x = \C v \in P$ admits a
locally free Koszul resolution
\begin{equation}
 \big\{\Omega^3(3) \longrightarrow \Omega^2(2)
 \longrightarrow \Omega^1(1) \longrightarrow \O\big\} \iso \O_x
\end{equation}
where each differential corresponds to $v \in V$ under the
isomorphism from Lemma \ref{th:beilinson}. After restricting, it
follows that there is a spectral sequence converging to
$\mathit{Hom}^*_{Y_0}(\Omega^3(3)[3], L\iota_0^*\O_x)$ whose $E_1$ term is
\begin{equation}
\begin{array}{r|ccccccc}
 & s = 0 & s = 1 & s = 2 & s = 3 & s = 4
 \\[0.25em]
 \hline
 t = 2 & \C & 0 & 0 & 0 & 0 \\
 t = 1 & 0  & 0 & 0 & 0 & 0 \\
 t = 0 & \C & V & \Lambda^2V & \Lambda^3V & 0
\end{array}
\end{equation}
The differential $\delta_1: E_1^{s,t} \rightarrow E_1^{s+1,t}$ is
wedge product with $v$ in the bottom row, which leaves just two
nonzero terms $E_2^{0,2} \iso \C$, $E_2^{3,0} \iso \C$. If $x \notin
Y_0$ then $L\iota_0^*\O_x = 0$, so these two terms must kill each
other through the only remaining differential $\delta_3^{0,2}$. On
the other hand, if $x \in Y_0$ then $L\iota_0^*\O_x = \O_x \oplus
\O_x[1]$ has nonzero sections, so the spectral sequence degenerates
at $E_2$.

The connection with the question we are interested in comes about as
follows. If $\SS_4$ was formal, the triangulated subcategory of
$D^b\mathit{Coh}(Y_0)$ generated by the $E_{0,k}$ could be modelled by twisted
complexes defined on the category $C_4$ with its trivial
$A_\infty$-structure. We know that $L\iota_0^*\O_x$ is a twisted
complex of the form $(E_{0,1} \oplus E_{0,2} \oplus E_{0,3} \oplus
E_{0,4},\delta)$, and the behaviour of first differential in the
spectral sequence means that $\delta$ must be (up to irrelevant
scalar multiples)
\begin{equation}
\delta = \begin{pmatrix}
 0 & & & \\
 v & 0 & & \\
 \ast & v & 0 & \\
 \ast & \ast & v & 0
\end{pmatrix},
\end{equation}
But for degree reasons the terms $\ast$ are necessarily equal to
zero. Having identified the object, we can check that in $D^b(C_4)$
the analogue of the spectral sequence above degenerates at $E_2$ for
any nonzero $v$, which is a contradiction. \qed

As before, let $H \subset \mathit{SL}(V)$ be the maximal torus, and $T = H/\Gamma_4$.
By choosing the {\v C}ech covering to be $H$-invariant, one can arrange that $\SS_4$ comes with a natural $H$-action. Following the same strategy as in Section \ref{subsec:q-preliminary}, one can use inner automorphisms to define an induced action of $T$ on $\SS_4$.

On the geometric side, we take a splitting $\Gamma_{16} \rightarrow \Gamma_{64} \subset H$ as before, and consider each $E_{0,k}$ as a $\Gamma_{16}$-equivariant sheaf in all 16 possible ways (differing from each other by twisting with a character of $\Gamma_{16}$). This gives a total of 64 objects in the derived category $D^b\mathit{Coh}_{\Gamma_{16}}(Y_0)$. Remark \ref{th:equivariant-sheaves} tells us that the dg algebra underlying their total morphism algebra is
\begin{equation}
\SS_{64} = \SS_4 \semidirect \Gamma_{16}.
\end{equation}
At this point, we can apply Lemma \ref{th:recognize-q4} to a minimal model of $\SS_4$, deducing that $\SS_4$ is $T$-equivariantly quasi-isomorphic to $\QQ_4$, and hence $\SS_{64}$ is quasi-isomorphic to $\QQ_{64}$. In particular, we know that $\mathit{HH}^*(\SS_{64},\SS_{64}) \iso \mathit{HH}^*(\QQ_{64},\QQ_{64})$ (which alternatively could also be determined using more geometric means).

The formal deformation $Y_q = \{p_q(v) = 0\}$ of $Y_0$ given by \eqref{eq:fermat} still admits a $\Lambda_\N$-linear $\Gamma_{16}$-action. Since the $E_{0,k}$ are restricted from projective space, they have obvious equivariant extensions $E_{q,k}$ to locally free coherent sheaves over $Y_q$. By the discussion in Section \ref{subsec:formal-schemes}, or rather its equivariant analogue, this gives rise to a $\Gamma_{16}$-equivariant one-parameter deformation $\SS_{64,q}$ of $\SS_{64}$.

\begin{lemma} \label{th:64-nontrivial}
${\SS}_{64,q}$ has nontrivial deformation class in $\mathit{HH}^2(\SS_{64},\SS_{64})$.
\end{lemma}

\proof
The argument is similar to the previous one, and will be only sketched. Let $x_q$ be a
$\Lambda_\N$-point of $P \times_\C \Lambda_{\N}$, which is now given by a family $v_q$ of nonzero vectors in $V$ varying with the formal parameter $q$. We have the same spectral sequence as before for the
derived restriction of its structure sheaf $\O_{x_q}$ to $\iota_q: Y_q \hookrightarrow P \times_\C \Lambda_\N$, and therefore $\mathit{Hom}^*_{Y_q}(\Omega^3(3)[3],L\iota_q^*\O_{x_q})$ is the cohomology of
\begin{equation}
 \Lambda_\N \xrightarrow{q \mapsto p_q(v_q)} \Lambda_\N.
\end{equation}
$L\iota_q^*\O_{x_q}$ is not an equivariant sheaf, but one can remedy this by summing over the $\Gamma_{16}$-orbit of $x_q$. This shows that the deformation class corresponds to $\partial_q p_q|_{q = 0} \in \mathit{Sym}^4(V)$. By a computation similar to that in Lemma \ref{th:qq64-trunc}, this is nonzero in $\mathit{HH}^2(\SS_{64},\SS_{64})$. \qed

Take the linear automorphism $U_4 \in \mathit{GL}_4(V)$, let it act on $P$ as well as our collection of objects. The induced action on $S_4$ matches that on $Q_4$ described previously, and there is an underlying action on $\SS_4$ if one chooses the Cech covers appropriately. Then, we also get an induced action of $\SS_{64}$. More precisely, when defining that induced action, there is a freedom of choice of an element of $\Gamma_{16}^*$, and one can use that freedom to ensure that the result reproduces $U_{64}$ on the cohomology level. Finally, this $\Z/4$-action persists when we deform to $\SS_{64,q}$.

Together with the previous Lemma, this allows us to apply Proposition \ref{th:recognize-q64-3} which shows that $\SS_{64,q}$ is quasi-isomorphic to $\QQ_{64,q}$, at least in the sense in which the
latter object is unique, which means up to reparametrization by $\psi \in \mathit{End}(\Lambda_\N)^\times$.

Let $Y_q^* = Y_q \times_{\Lambda_\N} \Lambda_\Q$ be the general fibre of our deformation, and $E_{q,k}^*$ the restrictions of the $E_{q,k}$ to it. Let $Z_q^*$ be the minimal resolution of the quotient $Y_q^*/\Gamma_{16}$. We need two more basic facts:

\begin{lemma}
The 64 equivariant versions of the $E_{q,k}^*$ are split-generators for $D^b\mathit{Coh}_{\Gamma_{16}}(Y_q^*)$.
\end{lemma}

\proof Beilinson's resolution of the diagonal \cite{beilinson78} shows that the $F_{q,k}^*$ are generators for the derived category of $P \times_\C \Lambda_\Q$. His argument carries over to the equivariant case, which means that the 64 objects obtained by making the $F_{q,k}^*$ equivariant in all possible ways generate the derived category of $\Gamma_{16}$-equivariant sheaves. The equivariant analogue of Lemma \ref{th:kontsevich} shows that the restrictions of these objects to $Y_q^*$ split-generate its derived category of equivariant sheaves. \qed

\begin{lemma} \label{th:kapranov-vasserot00}
$D^b\mathit{Coh}(Z_q^*) \iso D^b\mathit{Coh}_{\Gamma_{16}}(Y_q^*)$, as triangulated categories linear over $\Lambda_\Q$.
\end{lemma}

This is a theorem of Kapranov and Vasserot \cite{kapranov-vasserot00}, valid for any finite group action on a smooth $K3$ surface which preserves the holomorphic two-form. We now have the situation where we know that
\begin{equation}
D^b\mathit{Coh}_{\Gamma_{16}}(Y_q^*) \iso D^\pi(\SS(Y_q^*))
\end{equation}
is split-generated by certain objects, and where we know that the dg subcategory of $\SS_{\Gamma_{16}}(Y_q^*)$ consisting of those objects, denoted by $\SS_{64,q}$ above, is quasi-isomorphic to $\QQ_{64,q}$. Lemma \ref{th:derived-equivalence} then says that

\begin{cor} \label{th:mirror-1}
We have
$D^b\mathit{Coh}(Z_q^*) \iso D^\pi(\psi^*\QQ_{64,q} \otimes_{\Lambda_\N}
\Lambda_\Q)$ for some $\psi \in \mathit{End}(\Lambda_\N)^\times$. \qed
\end{cor}

\section{Counting polygons\label{sec:computation}}

We begin by showing that, for a quartic surface in $\CP{3}$ whose part at infinity is the intersection with the ``simplex'' $x_0x_1x_2x_3 = 0$, the deformation from affine to relative Fukaya category is
nontrivial at first order. This is based on Corollary \ref{th:wall-crossing-2}, with the relevant Lagrangian torus constructed by a degeneration argument.
%
Having done that, the next step is to apply Corollary
\ref{th:generates-everything} to prove that the 64 vanishing cycles
of the Fermat pencil of quartics are split-generators for the derived
Fukaya category. This reduces us to studying the full
$A_\infty$-subcategory formed by these cycles, and finally, the
dimensional induction machinery from Section
\ref{subsec:combinatorial} allows us to compute that subcategory
combinatorially.

\subsection{\label{subsec:nontrivial-geometric}}
Let $X = \CP{3}$ with $o_X = \K_X^{-1} = \O(4)$, equipped with a rescaled version of the Fubini-Study metric (so a line has area $4$ for the associated symplectic form $\omega_X$). Choose some section $\sigma_{X,0}$ which with $\sigma_{X,\infty} = x_0x_1x_2x_3$ generates a quasi-Lefschetz pencil. Fix some $z \in \C$ such that $X_z$ is smooth, and define $M_z = X_z \setminus X_{0,\infty}$. We will prove:

\begin{proposition} \label{th:quartic-deformation-2}
The deformation of $A_\infty$-categories $\Fuk(X_z,X_{0,\infty})$ yields, at first order, a nontrivial class in $\mathit{HH}^2(\Fuk(M_z),\Fuk(M_z))$.
\end{proposition}

A first remark is that the validity of the statement is independent of the choice of $z$ and $\sigma_{X,0}$. This is because the space of choices is path-connected. Note that the components of $X_{0,\infty}$ intersect orthogonally with respect to the restriction of our metric to $X_z$. Hence, one can apply a relative Moser argument, to show that for any two choices there is a symplectic isomorphism between the resulting $X_z$, which preserves the open subsets $M_z$. Since quartic surfaces are simply-connected, this is compatible with gradings, hence (see Remark \ref{th:class-of-j-2}) it induces an equivalence of relative Fukaya categories.

Take the irreducible component $C = \{x_3 = 0\} \iso \CP{2}$ of $X_\infty$, with its line bundle $o_C = o_X|C$ and corresponding symplectic form $\omega_C$. Let $C_\infty = \{x_0x_1x_2 = 0\} \subset C$ be its intersection with the other components, and $C_0 = X_0 \cap C = X_{0,\infty} \cap C$ the part of the base locus lying on $C$. On $C \setminus C_0$ we have a one-form $\theta_{C \setminus C_0}$, obtained by pulling back the connection via $\sigma_{X,0}|C$, which satisfies $d\theta_{C \setminus C_0} = \omega_C$.

Consider a moment map fibre
\begin{equation}
T_\lambda = \{ |x_k|^2 = \lambda_k \} \subset C \setminus C_\infty
\end{equation}
for $\lambda = (\lambda_0,\lambda_1,\lambda_2)$ with $\lambda_k >0$ and $\lambda_0 + \lambda_1 + \lambda_2 = 1$. Since we have written $T_\lambda = U(1)^3/U(1)$, it is natural to identify $H_1(T_\lambda) = \Z^3/(1,1,1)$.

\begin{lemma} \label{th:compute-theta}
Suppose that $T_\lambda$ is disjoint from $C_0$. Let $w: (D,\partial D) \rightarrow (C,T_{\lambda})$
be a disc whose boundary circle is in class $(b_0,b_1,b_2)$. Then
\begin{equation} \label{eq:c0-infty}
\textstyle \int_{\partial D} w^*\theta_{C \setminus C_0} = -\textstyle w \cdot C_0 + \frac{4}{3} w \cdot C_\infty + b_0(4\lambda_0 - \frac{4}{3}) + b_1(4\lambda_1 - \frac{4}{3}) + b_2(4\lambda_2 - \frac{4}{3}).
\end{equation}
\end{lemma}

To make that formula a little more plausible, note that $w \cdot C_\infty \equiv b_0 + b_1 + b_2$ mod $3$. Hence, the right hand side of \eqref{eq:c0-infty} reduces to $4\lambda_0 b_0 + 4\lambda_1 b_1 + 4 \lambda_2 b_2$ mod $\Z$. This of course reflects the fact that the flat bundle $o_C|T_{\lambda}$ gives rise to a class in $H^1(T_{\lambda};\R/\Z)$, which is precisely $(4\lambda_0,4\lambda_1,4\lambda_2)$.

\proof As usual we have
\begin{equation}
\textstyle \int_D w^*\omega_C = \int_{\partial D} w^*\theta_{C \setminus C_0} + w \cdot C_0.
\end{equation}
On the other hand, the area of discs with boundary on $T_{\lambda}$ is well-known (there is a disc whose boundary class is the $k$-th unit vector ($k = 0,1,2$), which intersects the divisor $C_\infty$ once, and has area $4\lambda_k$).
\qed

Let's temporarily consider the degenerate choice
\begin{equation} \label{eq:conic-line}
\sigma_{X,0}|C = x_0^2(x_1^2 + x_2^2).
\end{equation}
This of course does not satisfy the conditions for a quasi-Lefschetz pencil, but we will remedy that shortly. With that choice, $C_0 \setminus C_\infty \subset C \setminus C_\infty$ is the union of two lines $\{x_1 = \pm i x_2\}$. Using Lemma \ref{th:compute-theta} one sees that in $H^1(T_{\lambda};\R) \iso H_1(T_{\lambda};\R)^\vee$ we have
\begin{equation} \label{eq:theta-fibre}
[\theta_{C \setminus C_0}|T_{\lambda}] = \begin{cases} (4\lambda_0-2,4\lambda_1,4\lambda_2-2) & \lambda_1 < \lambda_2, \\
(4\lambda_0-2,4\lambda_1-2,4\lambda_2) & \lambda_1 > \lambda_2.
\end{cases}
\end{equation}
$K_{1/2} = T_{1/2,1/4,1/4}$ is a Lagrangian torus intersecting $C_0$ in two parallel circles, whose homology class is $(0,1,1)$, hence in our notation Poincar{\'e} dual to $(0,-1,1) \in H^1(K_{1/2})$. Suppose first that we choose a small Hamiltonian function $\tilde{H}$ supported in a neighbourhood of $K_{1/2}$, and such that the gradient vector field of $\tilde{H}|K_{1/2}$ is transverse to both circles, and points in the same direction along each. One can use the backwards and forwards flow of that function to produce two perturbed Lagrangian tori $\tilde{K}_0, \tilde{K}_1 \subset C$, which are disjoint from $C_0 \cup C_\infty$. These are Hamiltonian isotopic to $K_{1/2}$, but at the same time they are Lagrangian isotopic to $T_{1/2,1/4-\epsilon,1/4+\epsilon}$ and $T_{1/2,1/4+\epsilon,1/4-\epsilon}$ (for small nonzero $\epsilon$, let's fix the signs by saying $\epsilon>0$) in the complement of $C_0 \cup C_\infty$ (the last-mentioned fact follows by interpolating between $d\tilde{H}$ and a torus-invariant closed one-form). In particular, they behave like the limiting case of \eqref{eq:theta-fibre} when approaching the point $(1/2,1/4,1/4)$ from both sides:
\begin{equation} \label{eq:limit-theta}
\begin{aligned}
& [\theta_{C \setminus C_0}|\tilde{K}_1] = (0,-1,1), \\
& [\theta_{C \setminus C_0}|\tilde{K}_0] = (0,1,-1).
\end{aligned}
\end{equation}
Now consider a different Hamiltonian function $H$, such that the gradient vector field of $H|K_{1/2}$ points in opposite directions transversally to our two circles. This again gives rise to Lagrangian tori $K_0$ and $K_1$ in the complement of $C_0 \cup C_\infty$, but now \eqref{eq:limit-theta} and \eqref{eq:two-pushoffs} tell us that
\begin{equation} \label{eq:pm-push}
\begin{aligned}
& [\theta_{C \setminus C_0}|K_1] = [\theta_{C \setminus C_0}|\tilde{K}_1] - (0,-1,1) = (0,0,0), \\
& [\theta_{C \setminus C_0}|K_0] = [\theta_{C \setminus C_0}|\tilde{K}_0] + (0,-1,1) = (0,0,0).
\end{aligned}
\end{equation}
(Note that both classes have to be equal, by Lemma \ref{th:one-pushoff}, which is a useful check on the signs in this computation). Suppose now that one slightly perturbs \eqref{eq:conic-line} so as to satisfy the requirements for a quasi-pencil. By using the isotopy theorem for symplectic submanifolds locally, one can move $K_{1/2}$ by a Hamiltonian isotopy, so that its intersection with the perturbed version of $C_0$ again consists of two circles, and the entire argument leading to \eqref{eq:pm-push} carries over. We will assume from now on (without changing notation) that we are in this more generic situation.

The final step is to move away from the fibre at infinity $X_\infty$ in our pencil. Consider part of the graph of the pencil of quartic surfaces, and the base locus in it:
\begin{equation}
\begin{aligned}
& \hat{X} = \{(y,x) \in \C \times X \,:\, \sigma_{X,\infty}(x) = y\sigma_{X,0}(x)\} \xrightarrow{\text{projection}} \C, \\
& \hat{X}_\infty = \C \times X_{0,\infty} \subset \hat{X}.
\end{aligned}
\end{equation}
The fibre of $(\hat{X},\hat{X}_\infty)$ over a point $y$ is $(X_{1/y},X_{0,\infty})$. We equip $\hat{X}$ with the product symplectic form. Using parallel transport, one can find a locally closed Lagrangian submanifold $\hat{L}_{1/2} \subset \hat{X}$, which is fibered over a small interval $(-\delta,\delta) \subset \C$, and whose fibre over $0$ is $K_{1/2}$. Moreover, the parallel transport vector field can be made compatible with $\hat{X}_{\infty}$ near $K_{1/2}$, which means that in each fibre $\hat{L}_{1/2}$ intersects $\hat{X}_{\infty}$ in two circles. Now specialize to some small nonzero $y = 1/z \in (-\delta,\delta)$, and let $L_{1/2}$ be the corresponding fibre of $\hat{L}_{1/2}$, which is a Lagrangian submanifold in $X_z$.

\begin{lemma}
$L_{1/2} \subset X_z$ admits a grading.
\end{lemma}

\proof We have a short exact sequence
\begin{equation} \label{eq:tangent-sequence}
 0 \rightarrow TC \longrightarrow
 T\hat{X}|C \longrightarrow
 \O_C(-C_\infty) \rightarrow 0,
\end{equation}
and if we restrict that to $K_{1/2}$, a subsequence of real vector bundles
\begin{equation}
 0 \longrightarrow TK_{1/2} \longrightarrow T\hat{L}_{1/2}|K_{1/2} \longrightarrow \R \rightarrow 0.
\end{equation}
Take a disc $w: (D,\partial D) \rightarrow (\hat{X},\hat{L}_{1/2})$ whose image is contained in $C$. The two sequences above show that its Maslov index in $\hat{X}$ differs from its Maslov index in $C$ by the intersection number (taken in $C$) $-2(w \cdot C_\infty)$. This precisely cancels out the familiar formula for the Maslov index in $C$. We have now shown that $\hat{L}_{1/2} \subset \hat{X}$ has zero Maslov index. Specializing to the fibre at $z$ yields the desired conclusion.
\qed

A similar (in fact even simpler) argument shows that the flat line bundle $o_{X_z}|L_{1/2}$ is trivial, and that the Lagrangian submanifolds $L_0,L_1 \subset M_z = X_z \setminus X_{0,\infty}$ obtained by a pushoff construction applied to $L_{1/2} \cap X_{0,\infty}$ are exact. Corollary \ref{th:wall-crossing-2} then applies and concludes the proof of Proposition \ref{th:quartic-deformation-2}.

\subsection{}
As before take $X = \CP{3}$ with $\sigma_{X,\infty}(x) = x_0x_1x_2x_3$, so that the complement of $X_\infty$ can be written as a quotient by the action of $\Gamma_4^* = \{\pm \One, \pm i \One\}
\subset SL_4(\C)$,
\begin{equation} \label{eq:m-coordinates}
 M = \{x \in \C^4 \suchthat x_0 x_1 x_2 x_3 = 1\}/\Gamma_4^*.
\end{equation}
We now specify $X_0$ to be the Fermat quartic, $\sigma_{X,0} =
\quarter( x_0^4+x_1^4+x_2^4+x_3^4)$. $X_0$ intersects the strata of
$X_\infty$ transversally, and $\pi_M =
\sigma_{X,0}/\sigma_{X,\infty}$ has 64 nondegenerate critical points,
coming in groups of 16 which lie in the same fibre (this is not a
problem, see Remark \ref{th:multiple-fibres}). We take $z_* = 0$ as
our base point, and use the same embedded vanishing path for all
critical points which lie in a given fibre. This yields a collection
of 64 vanishing cycles in $M_0 = X_0 \setminus X_{0,\infty}$, divided
into four groups of 16 mutually disjoint ones. We make them into
branes by choosing gradings and, for the projective and relative
versions of that notion, the other additional data. Notation:
$\Fuk_{64}$ will be the full $A_\infty$-subcategory of $\Fuk(M_0)$
whose objects are these vanishing cycles; and $\Fuk_{64,q}$,
$\Fuk_{64,q}^*$ the corresponding $A_\infty$-subcategories of
$\Fuk(X_0,X_{0,\infty})$ and $\Fuk(X_0)$, respectively. From Propositions
\ref{th:q-zero} and \ref{th:invert-q} we know that $\Fuk_{64,q}$ is a
one-parameter deformation of $\Fuk_{64}$, and $\Fuk_{64,q}^* \iso
\Fuk_{64,q} \otimes_{\Lambda_\N} \Lambda_\Q$.

\begin{lemma} \label{th:f64-generates}
$D^\pi\Fuk(X_0) \iso D^\pi(\Fuk_{64,q}^*)$.
\end{lemma}

This follows from Corollary \ref{th:generates-everything}, which shows that the
vanishing cycles are split-generators for $D^\pi\Fuk(X_0)$, and the general nonsense Lemma
\ref{th:generating-subcategory}.

\begin{lemma} \label{th:f64-nontrivial}
$\Fuk_{64,q}$ has nontrivial deformation class in $\mathit{HH}^2(\Fuk_{64},\Fuk_{64})$.
\end{lemma}

This is a consequence of Proposition \ref{th:quartic-deformation-2} and the general derived invariance of Hochschild cohomology, see the proof of Lemma \ref{th:product-turns-on-2}.

To simplify the geometry of our pencil, consider the action of the
group $\Gamma_{16}^* \subset \mathit{PSL}_4(\C)$ of diagonal matrices $A$
satisfying $A^4 = \mathit{id}$ on $X$. This action is free on $M$; it can be
lifted to $o_X$ and preserves the Fubini-Study connection, as well as
$\sigma_{X,0},\sigma_{X,\infty}$; and it is compatible with the
isomorphism $o_X \iso \K_X^{-1}$. Hence, on the quotient $\bar{M} =
M/\Gamma_{16}^*$ we have an induced symplectic form, a one-form
primitive for it, a holomorphic volume form, and a function
$\pi_{\bar{M}}: \bar{M} \rightarrow \C$. Explicitly,
\begin{equation}
 \bar{M} = \{ u_0 u_1 u_2 u_3 = 1\}, \quad
 \pi_{\bar{M}}(u) = \quarter(u_0 + u_1 + u_2 + u_3)
\end{equation}
where the coordinates are related to the previous ones by $u_k =
x_k^4$. Since $\Gamma_{16}^*$ permutes the critical points in any
given fibre transitively, $\pi_{\bar{M}}$ has one nondegenerate
critical point in each each of the four singular fibres. In other
words, the 64 vanishing cycles in $M_0$ are all possible lifts of the
corresponding 4 vanishing cycles in the quotient $\bar{M}_0$. We need
some basic information about the topology of this space:
\begin{equation} \label{eq:pi12}
\begin{split}
 & \pi_1(\bar{M}_0) = \{ h \in \Z^4 \suchthat h_0 + \dots + h_3 = 0\}, \\
 & im(\pi_2(\bar{M}_0) \otimes \Q \rightarrow H_2(\bar{M}_0;\Q)) \iso \Q^3
\end{split}
\end{equation}
where the first group maps isomorphically to $\pi_1(\bar{M})$, with
the generators $h_k$ becoming loops around the coordinates axes $u_k
= 0$ at infinity; and the second is generated by the 4 vanishing
cycles. Both things follow from elementary Morse and Lefschetz
theory, since $\bar{M} \htp T^3$ is obtained by attaching four
3-cells to $\bar{M}_0$, which is itself homotopy equivalent to a
2-complex. As explained in Section \ref{subsec:coverings}, this means
that we have an action of $\pi_1(\bar{M}_0)^*$ on a subcategory of
$\Fuk(\bar{M}_0)$ which consists of branes with a certain condition
on their fundamental groups, and in particular contains our vanishing
cycles since they are spheres. Recall that in Section
\ref{subsec:q-preliminary} we considered the group $T = H/\Gamma_4
\subset PSL(V)$. This can also be written as $T = (\C^*)^4/\C^*$
(the quotient of the group of all diagonal matrices by the multiples
of the identity) and hence can be identified with $\pi_1(\bar{M}_0)^*$.
The following is our main computational result, and will be proved later:

\begin{prop} \label{th:computation}
For a suitable choice of the lifts to $M_0$ which define the
$T$-action, the full $A_\infty$-subcategory $\Fuk_4 \subset
\Fuk(\bar{M}_0)$ whose objects are the 4 vanishing cycles of
$\pi_{\bar{M}}$ is $T$-equivariantly quasi-isomorphic to $\QQ_4$.
\end{prop}

The covering $M_0 \rightarrow \bar{M}_0$ provides an epimorphism
$\pi_1(\bar{M}_0) \rightarrow \Gamma_{16}^*$, and the dual of that is
the inclusion of the subgroup $\Gamma_{16} \hookrightarrow T$
considered in Section \ref{subsec:q-preliminary}. Combining the
statement above with \eqref{eq:semidirect-fukaya}, one finds that
\begin{equation}
 \Fuk_{64} \iso \Fuk_4 \semidirect \Gamma_{16} \iso
 \QQ_4 \semidirect \Gamma_{16} \iso \QQ_{64}.
\end{equation}

Now consider the symplectic automorphism $\bar\phi$ of $\bar{M}$
given by $\bar\phi(u_0,u_1,u_2,u_3) = (u_1,u_2,u_3,u_0)$.
This acts fibrewise with respect to $\pi_{\bar{M}}$, hence takes each
vanishing cycle to itself, but it reverses their orientation since it
acts with determinant $-1$ on the tangent space to each critical
point. If we equip $\bar\phi_0 = \bar\phi|\bar{M}_0$ with a suitable
``odd'' grading $\alphagr_{\bar\phi_0}$, it will map each vanishing
cycle to itself as a graded Lagrangian submanifold. The resulting $\Z/4$-action can also be made compatible with the {\em Spin} structures, as explained in Remark \ref{th:spin-group-action}. Therefore, one can define $\Fuk_4$ in such a way that it carries a $\Z/4$-action, which moreover is also compatible with $T$-action (meaning that the two combine to yield an action of $T \semidirect \Z/4$). We know that the action of the generator of this $\Z/4$-action on $H(\Fuk_4) \iso Q_4$, temporarily denoted by $\tilde{U}_4$, satisfies the conditions from Remark \ref{th:fake-u4}. We are certainly free to change the quasi-isomorphism $\Fuk_4 \htp \QQ_4$ by some element of $T$, and by so doing we can achieve that the matrix in $\mathit{GL}(V)$ associated to $\tilde{U}_4$ is of the form \eqref{eq:delta-matrix} with all nonzero coefficients equal to $1$. Moreover, the action of any element $[\alpha_1,\dots,\alpha_4] \in I_4$ on $Q_4$ lifts to an action on $\Fuk_4$, which simply multiplies each space $\mathit{hom}_{\Fuk_4}(X_j,X_k)$ with $\alpha_k\alpha_j^{-1}$. In particular, by using $\alpha_i$ which are fourth roots of unity, we can get a modified $\Z/4$-action whose associated generator is the automorphism $U_4$ (all of this follows Remark \ref{th:fake-u4}).

We can lift $\bar\phi_0$ to an automorphism $\phi_0$ of automorphism $\phi_0$ of $M_0$, in a way which is unique up to the action of the covering group $\Gamma_{16}^*$. The interplay of this with \eqref{eq:semidirect-fukaya} is not hard to understand, and a suitable choice of lift yields a $\Z/4$-action on $\Fuk_{64}$ whose generator on the cohomology level is the automorphism $U_{64}$. Moreover, the deformation $\Fuk_{64,q}$ can be made equivariant with respect to this (this requires an analogue of the discussion at the end of Section \ref{subsec:coverings} for relative Fukaya categories, which is entirely parallel to the affine case). This means that we have verified all the conditions of Proposition \ref{th:recognize-q64-2}, and therefore $\Fuk_{64,q}$ is quasi-isomorphic to $\QQ_{64,q}$ up to some change of parameter in $\mathit{End}(\Lambda_\N)^\times$. Using Lemma
\ref{th:f64-generates}, we obtain:

\begin{cor} \label{th:mirror-2}
There is a $\psi \in \mathit{End}(\Lambda_\N)^\times$ such that
\begin{equation} D^\pi\Fuk(X_0) \iso D^\pi(\Fuk_{64,q}^*) \iso
\psi^*D^\pi(\QQ_{64,q} \otimes_{\Lambda_\N} \Lambda_\Q).\qed \end{equation}
\end{cor}

This is the mirror statement to Corollary \ref{th:mirror-1}, and by
comparing the two, Theorem \ref{th:main} follows.

\subsection{\label{subsec:computation}}
It remains to explain Proposition \ref{th:computation}. Following the
strategy from Section \ref{subsec:braid-monodromy}, the first step in
understanding the geometry of the vanishing cycles is to represent
them as matching cycles, which involves the choice of a generic
auxiliary section. Let us momentarily work ``upstairs'' in $M$. An
ansatz for the section is
\begin{equation} \label{eq:ansatz}
 \sigma_{X,0}' = \quarter(\xi_0 x_0^4 + \xi_1 x_1^4 + \xi_2 x_2^4 + \xi_3 x_3^4)
\end{equation}
with pairwise distinct constants $\xi_k \in \C^*$. This always
satisfies the first of the genericity assumptions, namely
$(\sigma_{X,0}')^{-1}(0)$ intersects each stratum of $X_{0,\infty}$
transversally. The critical point set of the associated map
\eqref{eq:bm-map} is always smooth,
\begin{equation}
 \mathit{Crit}(b_M) = \Big\{ x_k^4 = (s+\xi_k t)^{-1},
 \text{ where } \prod_k (s+t \xi_k) = 1 \Big\}.
\end{equation}
The parameter value $t = 0$ corresponds precisely to the 64 critical
points of $\pi_M$. At those points, $b_M: \mathit{Crit}(b_M) \rightarrow \C^2$
is an immersion, and its second component
$\sigma'_{X,0}/\sigma_{X,\infty} : \mathit{Crit}(b_M) \rightarrow \C$ a simple
branched cover, iff the coefficients $\xi_k$ satisfy
\begin{equation} \label{eq:nondegeneracy}
 \sum_k \xi_k^2 \neq \quarter \Big(\sum_k \xi_k\Big)^2.
\end{equation}
Next, $b_M|\mathit{Crit}(b_M)$ cannot be a generic embedding because it is
$\Gamma_{16}^*$-invariant, but after passing to the quotient we do
get a generic embedding $b_{\bar{M}}: \mathit{Crit}(b_{\bar{M}}) \rightarrow
\C^2$, for elementary reasons: the image $\bar{C}$ of this map is an
irreducible algebraic curve in $\C^2$ having four distinct points at
infinity, and since $\mathit{Crit}(b_{\bar{M}})$ has Euler characteristic -2,
the map cannot be a multiple branch cover. The upshot is that the
dimensional induction technique from Section \ref{subsec:induction}
can be applied ``downstairs'', meaning to the vanishing cycles of
$\pi_{\bar{M}}$.

\begin{figure}[H]
\begin{centering}
\begin{picture}(0,0)%
\includegraphics{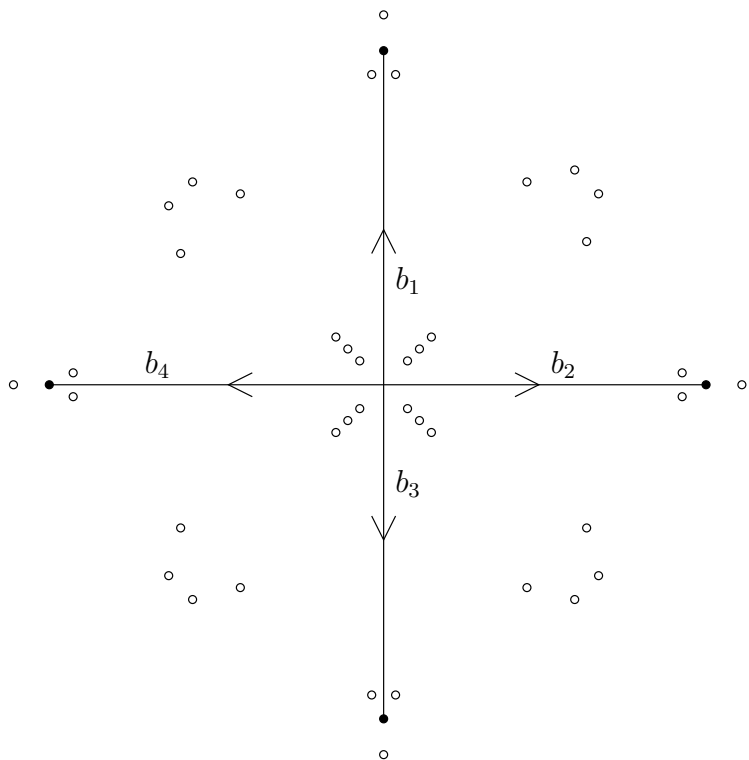}%
\end{picture}%
\setlength{\unitlength}{3947sp}%
\begingroup\makeatletter\ifx\SetFigFont\undefined%
\gdef\SetFigFont#1#2#3#4#5{%
  \reset@font\fontsize{#1}{#2pt}%
  \fontfamily{#3}\fontseries{#4}\fontshape{#5}%
  \selectfont}%
\fi\endgroup%
\begin{picture}(4641,4716)(268,-3694)
\put(2701,-736){\makebox(0,0)[lb]{\smash{\SetFigFont{12}{14.4}{\familydefault}{\mddefault}{\updefault}{\color[rgb]{0,0,0}$b_1$}%
}}}
\put(3676,-1261){\makebox(0,0)[lb]{\smash{\SetFigFont{12}{14.4}{\familydefault}{\mddefault}{\updefault}{\color[rgb]{0,0,0}$b_2$}%
}}}
\put(2701,-2011){\makebox(0,0)[lb]{\smash{\SetFigFont{12}{14.4}{\familydefault}{\mddefault}{\updefault}{\color[rgb]{0,0,0}$b_3$}%
}}}
\put(1126,-1261){\makebox(0,0)[lb]{\smash{\SetFigFont{12}{14.4}{\familydefault}{\mddefault}{\updefault}{\color[rgb]{0,0,0}$b_4$}%
}}}
\end{picture}
\end{centering}
\caption{\label{fig:paths1} The actual ($\bullet$) and fake ($\circ$)
critical points of $\pi_{\bar{M}}$, together with our choice of
vanishing paths $b_1,\dots,b_4$. The image has been distorted for
better legibility.}
\end{figure}

\begin{figure}[H]
\begin{center}
\begin{picture}(0,0)%
\includegraphics{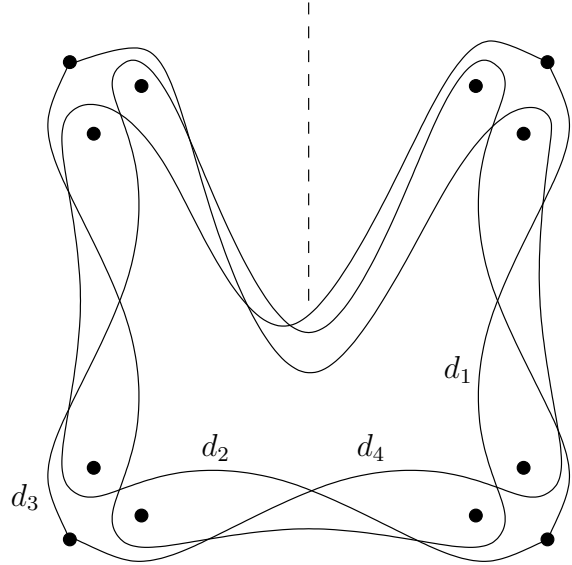}%
\end{picture}%
\setlength{\unitlength}{3947sp}%
\begingroup\makeatletter\ifx\SetFigFont\undefined%
\gdef\SetFigFont#1#2#3#4#5{%
  \reset@font\fontsize{#1}{#2pt}%
  \fontfamily{#3}\fontseries{#4}\fontshape{#5}%
  \selectfont}%
\fi\endgroup%
\begin{picture}(3750,3562)(301,-3211)
\put(1726,-2536){\makebox(0,0)[lb]{\smash{\SetFigFont{12}{14.4}{\familydefault}{\mddefault}{\updefault}{\color[rgb]{0,0,0}$d_2$}%
}}}
\put(530,-2850){\makebox(0,0)[lb]{\smash{\SetFigFont{12}{14.4}{\familydefault}{\mddefault}{\updefault}{\color[rgb]{0,0,0}$d_3$}%
}}}
\put(3250,-2050){\makebox(0,0)[lb]{\smash{\SetFigFont{12}{14.4}{\familydefault}{\mddefault}{\updefault}{\color[rgb]{0,0,0}$d_1$}%
}}}
\put(2701,-2536){\makebox(0,0)[lb]{\smash{\SetFigFont{12}{14.4}{\familydefault}{\mddefault}{\updefault}{\color[rgb]{0,0,0}$d_4$}%
}}}
\end{picture}
\end{center}

\caption{\label{fig:paths2}The four matching paths obtained from the
vanishing paths in Figure \ref{fig:paths1} by braid monodromy. These
were drawn using a computer to watch how the points of $\bar{C} \cap
\{z = b_k(t)\}$ move in $\C$ as one goes along the vanishing paths.
The picture would be $\Z/4$-symmetric if we hadn't moved the paths by
an isotopy, so as to avoid the future base path (the dashed line).}
\end{figure}

Let's get down to concrete numbers. Our choice of coefficients is $\xi_0 =
5/4$, $\xi_1 = i$, $\xi_2 = -1$, $\xi_3 = -i$, which satisfies
\eqref{eq:nondegeneracy}. The branch curve (determined by computer
using elimination theory) is
\begin{equation} \label{eq:elimination-theory}
\begin{split}
 \bar{C} = \{ & 128000 z^{12} - 76800 z^{11}w
 + 92160 z^{10}w^2 - 108544 z^9w^3 - 258048 z^8w^4 \\ & + 86016z^7w^5
 - 112640 z^6w^6 + 141312 z^5w^7 + 211968 z^4w^8 - 25600z^3w^9 \\
 & + 36864 z^2 w^{10} - 49152 zw^{11} - 65536 w^{12}
 - 391307 z^8 + 180632 z^7w \\ & + 1940 z^6w^2 - 503112 z^5w^3
 - 2102866 z^4w^4 + 457512 z^3w^5 - 122764 z^2w^6 \\ & - 102904 zw^7
 - 251483 w^8 + 399430 z^4 - 106312 z^3w - 155028z^2w^2 \\ & - 53752 zw^3
 - 320858 w^4 - 136161 = 0 \}.
\end{split}
\end{equation}
Figure \ref{fig:paths1} shows the $z$-plane with the 4 points of
$\mathit{Critv}(\pi_{\bar{M}}) = \{\pm 1, \pm i\}$, as well the 40 points of
$\mathit{Fakev}(\pi_{\bar{M}})$. Luckily, none of the latter lie on $[-1,1]$
or $i[-1,1]$, so we can take the basis of embedded critical paths
which are straight lines $b_k(t) = i^{2-k}t$, $k = 1,\dots,4$. Denote
by $V_{b_k} \subset \bar{M}_0$ the resulting vanishing cycles. We now
turn to the function
\begin{equation} \label{eq:q0}
q_{\bar{M}_0} = \quarter \sum_k \xi_k u_k: \bar{M}_0 =
\pi_{\bar{M}}^{-1}(0) \longrightarrow \C.
\end{equation}
This is not quite a Lefschetz fibration, but it is the free
$\Gamma_{16}^*$-quotient of such a fibration on $M_0$, so the same techniques can be applied. Figure
\ref{fig:paths2} shows its base $\C$ with the 12 critical values,
which are also the points of $\bar{C} \cap \{z = 0\}$, and the
matching paths $d_k$ obtained from the braid monodromy of $\bar{C}$
over $b_k$. The associated matching cycles are our vanishing cycles:
$V_{b_k} \htp \Sigma_{d_k}$.

\begin{figure}[H]
\begin{center}
 \epsfig{file=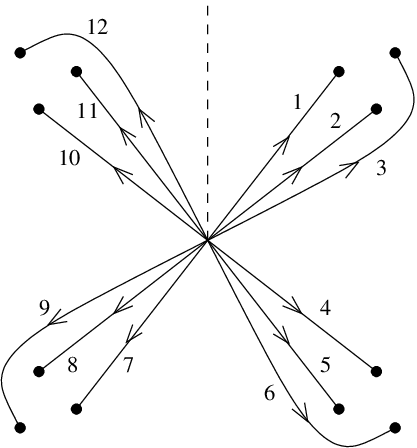}

 \caption{\label{fig:paths3}
 The base path and the 12 vanishing paths $c_1,\dots,c_{12}$
 (ordered as they must be, clockwise with respect to
 the derivatives at the base point).
 }
\end{center}
\end{figure}

Taking $0$ as a base point, choose a base path $c_*$ and a basis of
vanishing paths $c_1,\dots,c_{12}$ for \eqref{eq:q0} as indicated in
Figure \ref{fig:paths3}. This gives rise to vanishing cycles $V_{c_k}
\subset \bar{M}_{0,0} = q_{\bar{M}_0}^{-1}(0)$, which are simple
closed loops on a four-pointed genus three Riemann surface.

\begin{remark}
If one took $\xi_k = i^k$ instead, condition \eqref{eq:nondegeneracy}
would be violated but the structure of $q_{\bar{M}_0}$ would be much
simpler: it would have 4 degenerate critical points of type $A_3$.
The actually chosen $\xi_k$ are a small perturbation of these values,
in the sense that one can still recognize the groups of 3 critical
points obtained by ``morsifying'' the degenerate singularities. This
allows one to choose the vanishing paths $c_k$ in such a way as to
minimize the intersections between the vanishing cycles in each
group.
\end{remark}

Let $\Fuk^\rightarrow_{12}$ be the directed $A_\infty$-subcategory of
the Fukaya category of $\bar{M}_{0,0}$ whose objects are the
$V_{c_k}$. From the restriction of the universal cover of
$\bar{M}_0$, we get an action of $\pi_1(\bar{M}_0)^* \iso T$ on this
$A_\infty$-category. Let $\mathit{Tw}_T(\Fuk_{12}^\rightarrow)$ be the equivariant version of the category of twisted complexes (this is straightforward to define; the forgetful map $\mathit{Tw}_T(\Fuk_{12}^\rightarrow) \rightarrow \mathit{Tw}(\Fuk_{12}^\rightarrow)$ is injective, but not surjective, on morphism spaces). The equivariant version of Theorem
\ref{th:induction} shows that $\Fuk_4$ is equivariantly
quasi-isomorphic to the full $A_\infty$-subcategory of
$\mathit{Tw}_T(\Fuk_{12}^\rightarrow)$ whose objects are the following four
equivariant twisted complexes, obtained by applying the rules from
Section \ref{subsec:induction} to the paths from Figures
\ref{fig:paths2}, \ref{fig:paths3}:
\begin{equation} \label{eq:s-objects}
\begin{split}
 &
 S_1 = \big\{ T'_{V_{c_7}}T'_{V_{c_5}}(V_{c_3})
 \longrightarrow T_{V_{c_{10}}}T_{V_{c_{11}}}(V_{c_{12}}) \big\}, \\
 &
 S_2 = \big\{ T_{V_{c_1}}T_{V_{c_2}}(V_{c_3})
 \longrightarrow T'_{V_{c_{10}}}T'_{V_{c_8}}(V_{c_6}) \big\}, \\
 &
 S_3 = \big\{ T_{V_{c_1}}T_{V_{c_4}}T_{V_{c_5}}(V_{c_6})
 \longrightarrow T_{V_{c_{11}}}'(V_{c_9}) \big\}, \\
 &
 S_4 = \big\{ T_{V_{c_2}}T_{V_{c_4}}T_{V_{c_7}}T_{V_{c_8}}(V_{c_9})
 \longrightarrow V_{c_{12}} \big\}
\end{split}
\end{equation}
where $\{\cdots \rightarrow \cdots\}$ is the cone over the lowest
degree nonzero morphism, which is unique up to a $\C^*$ multiple.
Actually, that morphism may not have degree zero, or it may not be
$T$-invariant, so one needs to change one of the two objects forming
the cone by a translation and tensoring with a suitable
$T$-character. In that sense, the formulae given \eqref{eq:s-objects}
are not entirely precise, but the appropriate corrections can not yet
be determined at this stage in the computation.

\begin{figure}[H]
\begin{center}
 \epsfig{file=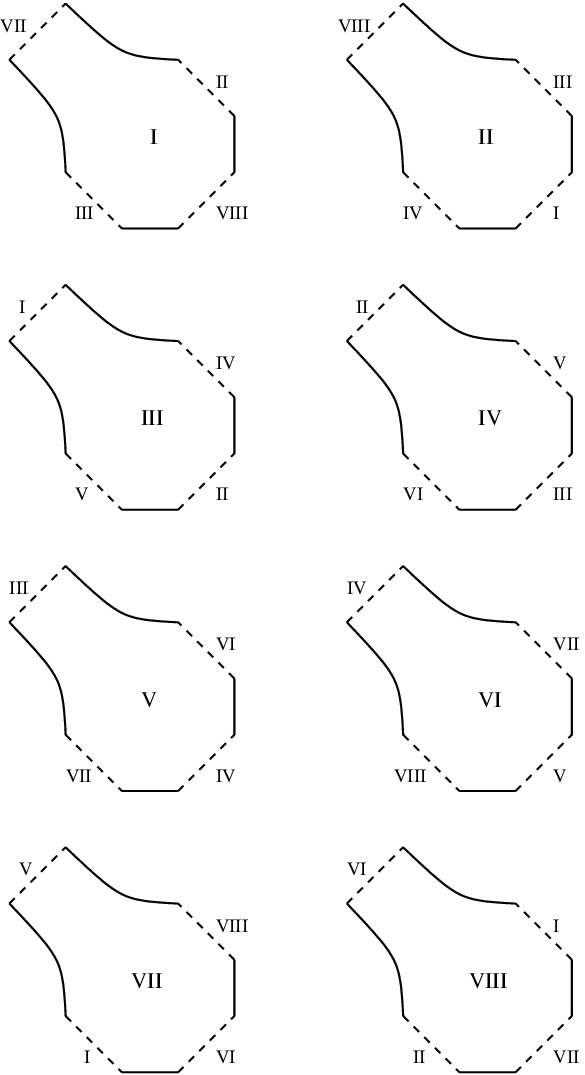}

 \caption{\label{fig:curve1}}
\end{center}
\end{figure}

\begin{figure}[H]
\begin{center}
 \epsfig{file=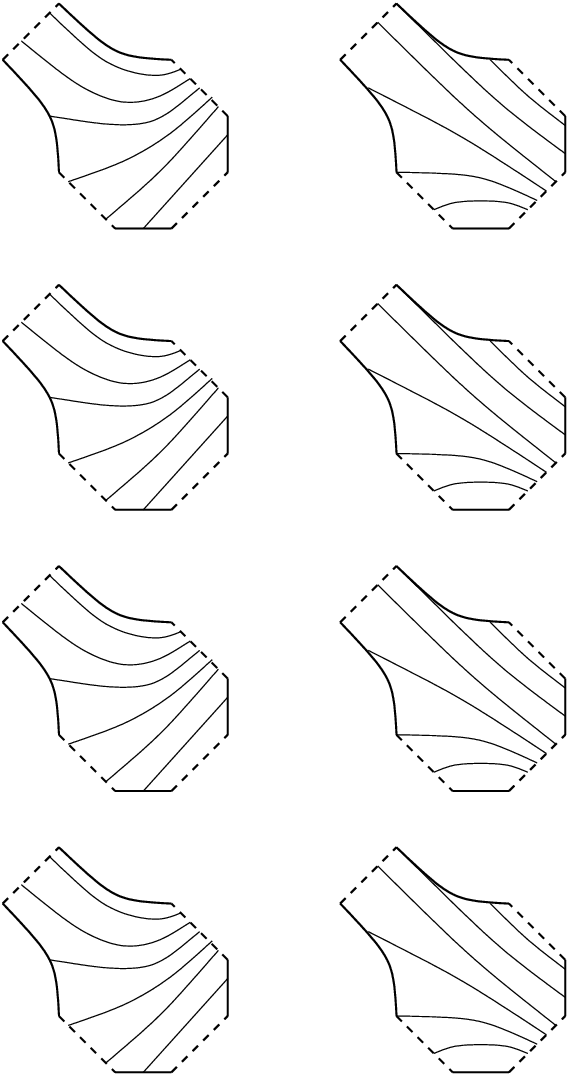}

 \caption{\label{fig:curve2}}
\end{center}
\end{figure}

\begin{figure}[H]
\begin{center}
 \epsfig{file=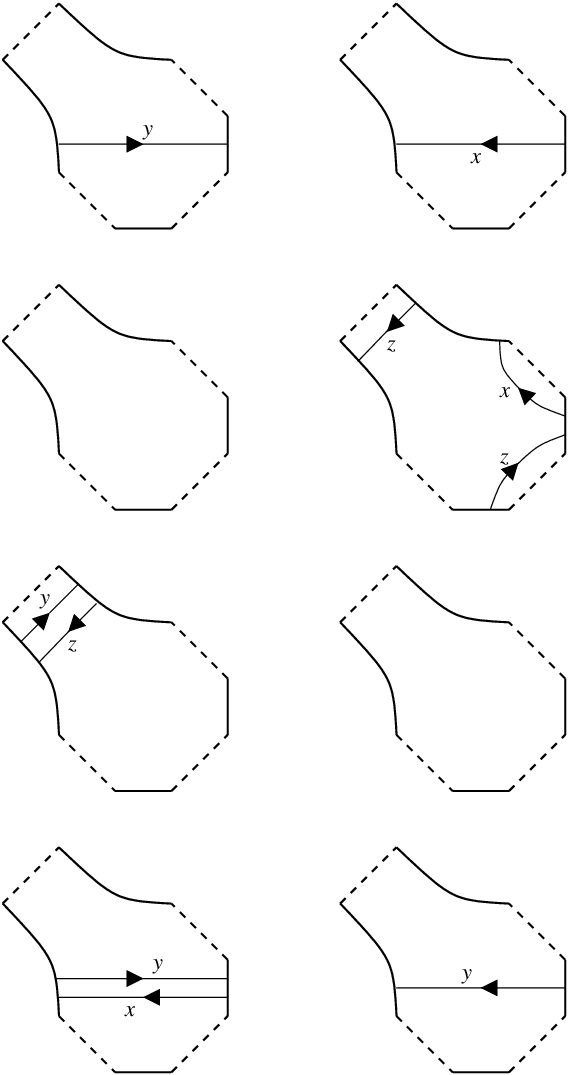}
 \caption{\label{fig:curve3}}
\end{center}
\end{figure}

The next step is to compute $\Fuk^\rightarrow_{12}$, or rather a
directed $A_\infty$-category $\Cat^\rightarrow_{12}$ quasi-isomorphic
to it, using the combinatorial description of directed Fukaya
categories on Riemann surfaces. Following the procedure described in
Section \ref{subsec:combinatorial}, we take a compact piece $N
\subset \bar{M}_{0,0}$ which is a surface with boundary and a
deformation retract. Figure \ref{fig:curve1} shows a decomposition of
$N$ into eight polygons (I)--(VIII), with the solid lines the actual
$\partial N$, and the dashed lines to be glued to the boundaries of
other polygons as indicated. $\bar{M}_{0,0}$ inherits the structure
of an affine Calabi-Yau from $\bar{M}$, and Figure \ref{fig:curve2}
shows an unoriented foliation on $N$ lying in the same homotopy class
as the distinguished trivialization of the squared canonical bundle.
Moreover, restriction
\begin{equation}
\Z^4/\Z(1,1,1,1) \iso H^1(\bar{M}) \longrightarrow H^1(\bar{M}_{0,0})
\end{equation}
gives four canonical elements of $H^1(N)$ whose sum is zero; the
Poincar{\'e} duals of the first three elements are represented by the
one-cycles $x,y,z$ in $(N,\partial N)$ shown in Figure
\ref{fig:curve3}. Figures \ref{fig:vanishing1}--\ref{fig:vanishing8}
show simple closed curves $\nu_k$ isotopic to the vanishing cycles
$V_{c_k}$, drawn on each of the eight polygons which make up $N$.

\begin{remark}
These pictures were arrived at by an application of braid monodromy
(in one dimension lower than before), which means that we introduce
another auxiliary holomorphic function on $\bar{M}_0$, whose
restriction to any generic fibre of $q_{\bar{M}_0}$ represents that
surface as a branched cover with only simple branch points. The
function we use is simply a projection
\begin{equation}
 r_{\bar{M}_0}(u) = u_2 \in \C^*,
\end{equation}
making $\bar{M}_{0,0}$ into a fourfold branched cover of $\C^*$ with
eight simple branch points. Cutting the base of that into suitable
pieces gives rise to our decomposition of $N$ into polygons.
Elimination theory is used to compute the branch curve of
\begin{equation}
(q_{\bar{M}_0},r_{\bar{M}_0}): \bar{M}_{0,0} \longrightarrow \C
\times \C^*,
\end{equation}
and from that one obtains the vanishing cycles
$V_{c_1},\dots,V_{c_{12}}$. We have omitted the details of this step,
since it is essentially topological, and completely parallel to the
procedure which led to the paths $b_1,\dots,b_4$ in Figure
\ref{fig:paths2}.
\end{remark}

When drawing the $\nu_k$, there is some freedom in their relative
position with respect to each other. The only constraint, which comes
from exactness, says that we should be able to assign positive areas
to the connected components of $N \setminus \bigcup_k \nu_k$ such
that the two-chains which represent relations between the $\nu_k$
have signed area zero. In fact, as explained in Remark
\ref{th:painting-by-numbers}, not all relations need to be
considered; in our case \eqref{eq:pi12} shows that there are three
essential ones. These can be read off from \eqref{eq:s-objects} and
the intersection numbers of the $\nu_k$, and are
\begin{equation}
\begin{split}
 & \nu_5 + \nu_8 + \nu_{10} - \nu_1 - \nu_2 - \nu_3 \sim 0, \\
 & \nu_3 + \nu_5 + \nu_7 - \nu_{10} - \nu_{11} - \nu_{12} \sim 0, \\
 & \nu_8 + \nu_2 + \nu_4 - \nu_7 - \nu_8 - \nu_0 \sim 0.
\end{split}
\end{equation}
After having drawn the relevant two-chains explicitly in Figures
\ref{fig:vanishing1}--\ref{fig:vanishing8}, one realizes that for
each of them, there are pieces of $N \setminus \bigcup_k \nu_k$ which
appear in it with either sign, and which do not appear in the other
two-chains. This means that the configuration of curves $\nu_k$ can
indeed be used as basis for the computation.

Choose gradings and $Spin$ structures for the $\nu_k$, lifts to the
abelian covering of $N$ induced from $\pi_1(N) \iso
\pi_1(\bar{M}_{0,0}) \rightarrow \pi_1(\bar{M}_0)$, and
identifications $\C_x \iso \C$ for each intersection point of $\nu_k
\cap \nu_l$ with $k<l$. One can then list all intersection points
with their (Maslov index) degrees as well as their ($T$-action)
weights. This is done in Tables \ref{table:intersections1} and
\ref{table:intersections2}, where the notation is that
\begin{equation} \label{eq:abcd}
 m^d x^a y^b z^c
\end{equation}
is an intersection point with Maslov index $d$, and on which
$[x,y,z,1] \in T$ acts as $x^ay^bz^c$, respectively. By inspecting
the table, one sees that different intersection points of the same
$\nu_k \cap \nu_l$ are always distinguished by $(a,b,c,d)$. We may
therefore denote the generator of
$\mathit{hom}_{\Cat^\rightarrow_{12}}(\nu_k,\nu_l)$ given by an intersection
point \eqref{eq:abcd} by
\begin{equation}
 \nu_k \cap \nu_l: m^d x^a y^b z^c.
\end{equation}
The next step is to count all immersed holomorphic polygons which
have sides on the $\nu_k$ (in the correct order), with the
appropriate signs. The outcome of that are the composition maps
$\mu^d$ in $\Cat^{\rightarrow}_{12}$. It turns out that only
$\mu^2,\mu^3$ are nonzero; their coefficients are listed in Tables
\ref{table:mu2} and \ref{table:mu3}.

\begin{remark}
The $T$-action is actually a big bonus in this computation, since the
fact that the $\mu^d$ are homogeneous means that a lot of
their coefficients have to vanish. For instance, the fact that $\mu^d
= 0$ for $d > 3$ follows by inspection of the list of $T$-weights of
the intersection points, so one does not need to verify geometrically
that there are no immersed pentagons, hexagons, and so forth. Another important self-check on the computation is provided by the $A_\infty$-relations themselves.
\end{remark}

At this point, we can write down the objects $S_k$ from
\eqref{eq:s-objects}, or rather the corresponding objects $C_k$ in
$\mathit{Tw}_T(\Cat^{\rightarrow}_{12})$, more explicitly. They are
\begin{align*}
 &
 C_1 = \nu_3[1]\twist{x} \oplus \nu_5\twist{y} \oplus \nu_7\twist{yz}
 \oplus \nu_{10}\twist{x^{-1}y} \oplus \nu_{11}[-1]\twist{y} \oplus
 \nu_{12}[-1], \\
 & \delta_1 \!=\!
 {\arraycolsep2.5pt \begin{pmatrix}
 0 & 0 & 0 & 0 & 0 & 0 \\
 \nu_3 \cap \nu_5\!:\!x^{-1}y \!\! & 0 & 0 & 0 & 0 & 0 \\
 0 & \!\! \nu_5 \cap \nu_7\!:\!mz & 0 & 0 & 0 & 0 \\
 0 & \!\! \nu_5 \cap \nu_{10}\!:\!mx^{-1} & 0 & 0 & 0 & 0 \\
 0 & 0 & \nu_7 \cap \nu_{11}\!:\!z^{-1} & \nu_{10} \cap \nu_{11}\!:\!x & 0 & 0 \\
 0 & 0 & 0 & 0 & \!\!\!\! \nu_{11} \cap \nu_{12}\!:\!my^{-1} & 0
 \end{pmatrix}},
 \\
 &
 C_2 = \nu_1[1]\twist{x^{-1}y} \oplus \nu_2\twist{y} \oplus \nu_3
 \oplus \nu_6\twist{z^{-1}} \oplus \nu_8[-1]\twist{z^{-1}}
 \oplus \nu_{10}[-1]\twist{x^{-1}}, \\
 & \delta_2 =
 {\arraycolsep2.5pt \begin{pmatrix}
 0 & 0 & 0 & 0 & 0 & 0 \\
 \nu_1 \cap \nu_2\!:\!x\!\! & 0 & 0 & 0 & 0 & 0 \\
 0 & \nu_2 \cap \nu_3\!:\!my^{-1} & 0 & 0 & 0 & 0 \\
 0 & \nu_2 \cap \nu_6\!:\!my^{-1}z^{-1} & 0 & 0 & 0 & 0 \\
 0 & 0 & \!\!\!\!\!\! -(\nu_3 \cap \nu_8\!:\!z^{-1}) & \nu_6 \cap \nu_8\!:\!1 & 0 & 0 \\
 0 & 0 & 0 & 0 & \!\! \nu_{8} \cap \nu_{10}\!:\!mx^{-1}z & 0
 \end{pmatrix}},
 \\
 &
 C_3 = \nu_1\twist{x^{-1}} \oplus \nu_4\twist{x^{-1}y} \oplus
 \nu_5[-1]\twist{x^{-1}y} \oplus \nu_6[-1]\twist{x^{-1}z^{-1}} \\
 & \qquad \qquad
 \oplus \nu_9[-1]\twist{y^{-1}z^{-1}} \oplus \nu_{11}[-2]\twist{z^{-1}}, \\
 & \delta_3 =
 {\arraycolsep2.5pt \begin{pmatrix}
 0 & 0 & 0 & 0 & 0 & 0 \\
 0 & 0 & 0 & 0 & 0 & 0 \\
 \nu_1 \cap \nu_5\!:\!y & \nu_4 \cap \nu_5\!:\!1 & 0 & 0 & 0 & 0 \\
 0 & 0 & \nu_5 \cap \nu_6\!:\!my^{-1}z^{-1} & 0 & 0 & 0 \\
 0 & 0 & \!\!\nu_5 \cap \nu_9\!:\!mxy^{-2}z^{-1}\!\! & 0 & 0 & 0 \\
 0 & 0 & 0 & \nu_6 \cap \nu_{11}\!:\!x & \nu_9 \cap \nu_{11}\!:\!y & 0
 \end{pmatrix}}, \displaybreak[0]
 \\
 &
 C_4 = \nu_2[1]\twist{z^{-1}} \oplus \nu_4[1]\twist{x^{-1}} \oplus \nu_7[1]\twist{x^{-1}y}
 \oplus \nu_8\twist{x^{-1}z^{-1}} \\ & \qquad \qquad \oplus \nu_9\twist{x^{-1}y^{-1}z^{-1}}
 \oplus \nu_{12}\twist{y^{-1}z^{-1}}, \displaybreak[0] \\
 & \delta_4 =
 {\arraycolsep2.5pt \begin{pmatrix}
 0 & 0 & 0 & 0 & 0 & 0 \\
 \nu_2 \cap \nu_4 \!:\! mx^{-1}z & 0 & 0 & 0 & 0 & 0 \\
 \nu_2 \cap \nu_7 \!:\! mx^{-1}yz & 0 & 0 & 0 & 0 & 0 \\
 0 & \nu_4 \cap \nu_8\!:\!z^{-1} & \nu_7 \cap \nu_8\!:\!y^{-1}z^{-1} & 0 & 0 & 0 \\
 0 & 0 & 0 & \nu_8 \cap \nu_9\!:\!my^{-1} & 0 & 0 \\
 0 & 0 & 0 & \nu_8 \cap \nu_{12}\!:\!mxy^{-1} & 0 & 0
 \end{pmatrix}}.
\end{align*}
where $\nu_k[d]\twist{x^ay^bz^c}$ means the object $\nu_k$ shifted
down by $d$, and tensored with the one-dimensional representation of
$T$ on which $[x,y,z,1]$ acts as $x^{-a}y^{-b}z^{-c}$.

It is now easy to compute explicitly the full $A_\infty$-subcategory
of $\mathit{Tw}_T\Cat^{\rightarrow}_{12}$ consisting of these four equivariant
twisted complexes $C_k$. More precisely, one verifies (by computer)
first that the cohomology of this subcategory is $T$-equivariantly
isomorphic to $Q_4$ and then, following the indications in Remark
\ref{th:more-recognize-q4}, that a suitable fourth order Massey
product is nonzero, which implies that our full subcategory is
$T$-equivariantly quasi-isomorphic to $\QQ_4$. By Theorem
\ref{th:induction}, or rather the equivariant version of it explained
in Section \ref{subsec:addition}, we have a $T$-equivariant
quasi-isomorphism of our subcategory with $\Fuk_4$, and this
completes the proof of Proposition \ref{th:computation}.

\begin{figure}[H]
\begin{picture}(0,0)%
\includegraphics{vanishing1.pstex}%
\end{picture}%
\setlength{\unitlength}{3947sp}%
\begingroup\makeatletter\ifx\SetFigFont\undefined%
\gdef\SetFigFont#1#2#3#4#5{%
  \reset@font\fontsize{#1}{#2pt}%
  \fontfamily{#3}\fontseries{#4}\fontshape{#5}%
  \selectfont}%
\fi\endgroup%
\begin{picture}(4049,3690)(151,-3206)
\put(475,-330){\makebox(0,0)[lb]{\smash{\SetFigFont{10}{12.0}{\familydefault}{\mddefault}{\updefault}{\color[rgb]{0,0,0}$\nu_9$}%
}}}
\put(850, 14){\makebox(0,0)[lb]{\smash{\SetFigFont{10}{12.0}{\familydefault}{\mddefault}{\updefault}{\color[rgb]{0,0,0}$\nu_5$}%
}}}
\put(3526,-586){\makebox(0,0)[lb]{\smash{\SetFigFont{10}{12.0}{\familydefault}{\mddefault}{\updefault}{\color[rgb]{0,0,0}$\nu_3$}%
}}}
\put(4100,-2500){\makebox(0,0)[lb]{\smash{\SetFigFont{10}{12.0}{\familydefault}{\mddefault}{\updefault}{\color[rgb]{0,0,0}$\nu_8$}%
}}}
\put(1375,-2500){\makebox(0,0)[lb]{\smash{\SetFigFont{10}{12.0}{\familydefault}{\mddefault}{\updefault}{\color[rgb]{0,0,0}$\nu_{12}$}%
}}}
\put(2050,-3100){\makebox(0,0)[lb]{\smash{\SetFigFont{10}{12.0}{\familydefault}{\mddefault}{\updefault}{\color[rgb]{0,0,0}$\nu_5$}%
}}}
\put(3475,-3100){\makebox(0,0)[lb]{\smash{\SetFigFont{10}{12.0}{\familydefault}{\mddefault}{\updefault}{\color[rgb]{0,0,0}$\nu_6$}%
}}}
\put(3901,-2700){\makebox(0,0)[lb]{\smash{\SetFigFont{10}{12.0}{\familydefault}{\mddefault}{\updefault}{\color[rgb]{0,0,0}$\nu_1$}%
}}}
\put(3650,-750){\makebox(0,0)[lb]{\smash{\SetFigFont{10}{12.0}{\familydefault}{\mddefault}{\updefault}{\color[rgb]{0,0,0}$\nu_2$}%
}}}
\put(3800,-886){\makebox(0,0)[lb]{\smash{\SetFigFont{10}{12.0}{\familydefault}{\mddefault}{\updefault}{\color[rgb]{0,0,0}$\nu_{10}$}%
}}}
\end{picture}

\caption{\label{fig:vanishing1}}
\end{figure}

\begin{figure}[H]
\begin{picture}(0,0)%
\includegraphics{vanishing2.pstex}%
\end{picture}%
\setlength{\unitlength}{3947sp}%
\begingroup\makeatletter\ifx\SetFigFont\undefined%
\gdef\SetFigFont#1#2#3#4#5{%
  \reset@font\fontsize{#1}{#2pt}%
  \fontfamily{#3}\fontseries{#4}\fontshape{#5}%
  \selectfont}%
\fi\endgroup%
\begin{picture}(3974,3753)(226,-3269)
\put(4100,-2490){\makebox(0,0)[lb]{\smash{\SetFigFont{12}{14.4}{\familydefault}{\mddefault}{\updefault}{\color[rgb]{0,0,0}$\nu_5$}%
}}}
\put(3876,-2711){\makebox(0,0)[lb]{\smash{\SetFigFont{12}{14.4}{\familydefault}{\mddefault}{\updefault}{\color[rgb]{0,0,0}$\nu_1$}%
}}}
\put(1800,-2921){\makebox(0,0)[lb]{\smash{\SetFigFont{12}{14.4}{\familydefault}{\mddefault}{\updefault}{\color[rgb]{0,0,0}$\nu_{10}$}%
}}}
\put(2070,-3120){\makebox(0,0)[lb]{\smash{\SetFigFont{12}{14.4}{\familydefault}{\mddefault}{\updefault}{\color[rgb]{0,0,0}$\nu_3$}%
}}}
\put(3440,-3150){\makebox(0,0)[lb]{\smash{\SetFigFont{12}{14.4}{\familydefault}{\mddefault}{\updefault}{\color[rgb]{0,0,0}$\nu_{12}$}%
}}}
\put(550,-180){\makebox(0,0)[lb]{\smash{\SetFigFont{12}{14.4}{\familydefault}{\mddefault}{\updefault}{\color[rgb]{0,0,0}$\nu_{10}$}%
}}}
\put(1050,280){\makebox(0,0)[lb]{\smash{\SetFigFont{12}{14.4}{\familydefault}{\mddefault}{\updefault}{\color[rgb]{0,0,0}$\nu_3$}%
}}}
\put(3350,-466){\makebox(0,0)[lb]{\smash{\SetFigFont{12}{14.4}{\familydefault}{\mddefault}{\updefault}{\color[rgb]{0,0,0}$\nu_6$}%
}}}
\put(3676,-800){\makebox(0,0)[lb]{\smash{\SetFigFont{12}{14.4}{\familydefault}{\mddefault}{\updefault}{\color[rgb]{0,0,0}$\nu_7$}%
}}}
\put(3526,-630){\makebox(0,0)[lb]{\smash{\SetFigFont{12}{14.4}{\familydefault}{\mddefault}{\updefault}{\color[rgb]{0,0,0}$\nu_2$}%
}}}
\end{picture}

\caption{\label{fig:vanishing2}}
\end{figure}

\begin{figure}[H]
\begin{picture}(0,0)%
\includegraphics{vanishing3.pstex}%
\end{picture}%
\setlength{\unitlength}{3947sp}%
\begingroup\makeatletter\ifx\SetFigFont\undefined%
\gdef\SetFigFont#1#2#3#4#5{%
  \reset@font\fontsize{#1}{#2pt}%
  \fontfamily{#3}\fontseries{#4}\fontshape{#5}%
  \selectfont}%
\fi\endgroup%
\begin{picture}(3674,3690)(526,-3206)
\put(3601,-700){\makebox(0,0)[lb]{\smash{\SetFigFont{12}{14.4}{\familydefault}{\mddefault}{\updefault}{\color[rgb]{0,0,0}$\nu_{11}$}%
}}}
\put(3376,-480){\makebox(0,0)[lb]{\smash{\SetFigFont{12}{14.4}{\familydefault}{\mddefault}{\updefault}{\color[rgb]{0,0,0}$\nu_{12}$}%
}}}
\put(3770,-890){\makebox(0,0)[lb]{\smash{\SetFigFont{12}{14.4}{\familydefault}{\mddefault}{\updefault}{\color[rgb]{0,0,0}$\nu_7$}%
}}}
\put(3901,-2686){\makebox(0,0)[lb]{\smash{\SetFigFont{12}{14.4}{\familydefault}{\mddefault}{\updefault}{\color[rgb]{0,0,0}$\nu_{10}$}%
}}}
\put(4100,-2461){\makebox(0,0)[lb]{\smash{\SetFigFont{12}{14.4}{\familydefault}{\mddefault}{\updefault}{\color[rgb]{0,0,0}$\nu_5$}%
}}}
\put(3526,-3061){\makebox(0,0)[lb]{\smash{\SetFigFont{12}{14.4}{\familydefault}{\mddefault}{\updefault}{\color[rgb]{0,0,0}$\nu_3$}%
}}}
\put(2000,-3100){\makebox(0,0)[lb]{\smash{\SetFigFont{12}{14.4}{\familydefault}{\mddefault}{\updefault}{\color[rgb]{0,0,0}$\nu_2$}%
}}}
\put(1410,-2500){\makebox(0,0)[lb]{\smash{\SetFigFont{12}{14.4}{\familydefault}{\mddefault}{\updefault}{\color[rgb]{0,0,0}$\nu_9$}%
}}}
\put(576,-246){\makebox(0,0)[lb]{\smash{\SetFigFont{12}{14.4}{\familydefault}{\mddefault}{\updefault}{\color[rgb]{0,0,0}$\nu_6$}%
}}}
\put(826,0){\makebox(0,0)[lb]{\smash{\SetFigFont{12}{14.4}{\familydefault}{\mddefault}{\updefault}{\color[rgb]{0,0,0}$\nu_2$}%
}}}
\end{picture}

\caption{\label{fig:vanishing3}}
\end{figure}

\begin{figure}[H]
\begin{picture}(0,0)%
\includegraphics{vanishing4.pstex}%
\end{picture}%
\setlength{\unitlength}{3947sp}%
\begingroup\makeatletter\ifx\SetFigFont\undefined%
\gdef\SetFigFont#1#2#3#4#5{%
  \reset@font\fontsize{#1}{#2pt}%
  \fontfamily{#3}\fontseries{#4}\fontshape{#5}%
  \selectfont}%
\fi\endgroup%
\begin{picture}(3674,3753)(526,-3269)
\put(1100,375){\makebox(0,0)[lb]{\smash{\SetFigFont{12}{14.4}{\familydefault}{\mddefault}{\updefault}{\color[rgb]{0,0,0}$\nu_{12}$}%
}}}
\put(3440,-520){\makebox(0,0)[lb]{\smash{\SetFigFont{12}{14.4}{\familydefault}{\mddefault}{\updefault}{\color[rgb]{0,0,0}$\nu_3$}%
}}}
\put(3725,-820){\makebox(0,0)[lb]{\smash{\SetFigFont{12}{14.4}{\familydefault}{\mddefault}{\updefault}{\color[rgb]{0,0,0}$\nu_4$}%
}}}
\put(3575,-670){\makebox(0,0)[lb]{\smash{\SetFigFont{12}{14.4}{\familydefault}{\mddefault}{\updefault}{\color[rgb]{0,0,0}$\nu_{11}$}%
}}}
\put(4051,-2536){\makebox(0,0)[lb]{\smash{\SetFigFont{12}{14.4}{\familydefault}{\mddefault}{\updefault}{\color[rgb]{0,0,0}$\nu_2$}%
}}}
\put(3901,-2686){\makebox(0,0)[lb]{\smash{\SetFigFont{12}{14.4}{\familydefault}{\mddefault}{\updefault}{\color[rgb]{0,0,0}$\nu_{10}$}%
}}}
\put(3500,-3100){\makebox(0,0)[lb]{\smash{\SetFigFont{12}{14.4}{\familydefault}{\mddefault}{\updefault}{\color[rgb]{0,0,0}$\nu_9$}%
}}}
\put(1821,-2936){\makebox(0,0)[lb]{\smash{\SetFigFont{12}{14.4}{\familydefault}{\mddefault}{\updefault}{\color[rgb]{0,0,0}$\nu_7$}%
}}}
\put(570,-200){\makebox(0,0)[lb]{\smash{\SetFigFont{12}{14.4}{\familydefault}{\mddefault}{\updefault}{\color[rgb]{0,0,0}$\nu_7$}%
}}}
\put(2101,-3251){\makebox(0,0)[lb]{\smash{\SetFigFont{12}{14.4}{\familydefault}{\mddefault}{\updefault}{\color[rgb]{0,0,0}$\nu_{12}$}%
}}}
\end{picture}

\caption{\label{fig:vanishing4}}
\end{figure}

\begin{figure}[H]
\begin{picture}(0,0)%
\includegraphics{vanishing5.pstex}%
\end{picture}%
\setlength{\unitlength}{3947sp}%
\begingroup\makeatletter\ifx\SetFigFont\undefined%
\gdef\SetFigFont#1#2#3#4#5{%
  \reset@font\fontsize{#1}{#2pt}%
  \fontfamily{#3}\fontseries{#4}\fontshape{#5}%
  \selectfont}%
\fi\endgroup%
\begin{picture}(3622,3753)(578,-3269)
\put(3451,-536){\makebox(0,0)[lb]{\smash{\SetFigFont{12}{14.4}{\familydefault}{\mddefault}{\updefault}{\color[rgb]{0,0,0}$\nu_9$}%
}}}
\put(3601,-686){\makebox(0,0)[lb]{\smash{\SetFigFont{12}{14.4}{\familydefault}{\mddefault}{\updefault}{\color[rgb]{0,0,0}$\nu_8$}%
}}}
\put(3701,-836){\makebox(0,0)[lb]{\smash{\SetFigFont{12}{14.4}{\familydefault}{\mddefault}{\updefault}{\color[rgb]{0,0,0}$\nu_4$}%
}}}
\put(606,-136){\makebox(0,0)[lb]{\smash{\SetFigFont{12}{14.4}{\familydefault}{\mddefault}{\updefault}{\color[rgb]{0,0,0}$\nu_{11}$}%
}}}
\put(531,-286){\makebox(0,0)[lb]{\smash{\SetFigFont{12}{14.4}{\familydefault}{\mddefault}{\updefault}{\color[rgb]{0,0,0}$\nu_3$}%
}}}
\put(1600,-2636){\makebox(0,0)[lb]{\smash{\SetFigFont{12}{14.4}{\familydefault}{\mddefault}{\updefault}{\color[rgb]{0,0,0}$\nu_6$}%
}}}
\put(2000,-3136){\makebox(0,0)[lb]{\smash{\SetFigFont{12}{14.4}{\familydefault}{\mddefault}{\updefault}{\color[rgb]{0,0,0}$\nu_{11}$}%
}}}
\put(3376,-3211){\makebox(0,0)[lb]{\smash{\SetFigFont{12}{14.4}{\familydefault}{\mddefault}{\updefault}{\color[rgb]{0,0,0}$\nu_{12}$}%
}}}
\put(3826,-2761){\makebox(0,0)[lb]{\smash{\SetFigFont{12}{14.4}{\familydefault}{\mddefault}{\updefault}{\color[rgb]{0,0,0}$\nu_7$}%
}}}
\put(4050,-2536){\makebox(0,0)[lb]{\smash{\SetFigFont{12}{14.4}{\familydefault}{\mddefault}{\updefault}{\color[rgb]{0,0,0}$\nu_2$}%
}}}
\end{picture}

\caption{\label{fig:vanishing5}}
\end{figure}

\begin{figure}[H]
\begin{picture}(0,0)%
\includegraphics{vanishing6.pstex}%
\end{picture}%
\setlength{\unitlength}{3947sp}%
\begingroup\makeatletter\ifx\SetFigFont\undefined%
\gdef\SetFigFont#1#2#3#4#5{%
  \reset@font\fontsize{#1}{#2pt}%
  \fontfamily{#3}\fontseries{#4}\fontshape{#5}%
  \selectfont}%
\fi\endgroup%
\begin{picture}(3622,3753)(578,-3269)
\put(1626,-2711){\makebox(0,0)[lb]{\smash{\SetFigFont{12}{14.4}{\familydefault}{\mddefault}{\updefault}{\color[rgb]{0,0,0}$\nu_4$}%
}}}
\put(3485,-3080){\makebox(0,0)[lb]{\smash{\SetFigFont{12}{14.4}{\familydefault}{\mddefault}{\updefault}{\color[rgb]{0,0,0}$\nu_6$}%
}}}
\put(2121,-3211){\makebox(0,0)[lb]{\smash{\SetFigFont{12}{14.4}{\familydefault}{\mddefault}{\updefault}{\color[rgb]{0,0,0}$\nu_9$}%
}}}
\put(3901,-2686){\makebox(0,0)[lb]{\smash{\SetFigFont{12}{14.4}{\familydefault}{\mddefault}{\updefault}{\color[rgb]{0,0,0}$\nu_7$}%
}}}
\put(4100,-2461){\makebox(0,0)[lb]{\smash{\SetFigFont{12}{14.4}{\familydefault}{\mddefault}{\updefault}{\color[rgb]{0,0,0}$\nu_{11}$}%
}}}
\put(3700,-811){\makebox(0,0)[lb]{\smash{\SetFigFont{12}{14.4}{\familydefault}{\mddefault}{\updefault}{\color[rgb]{0,0,0}$\nu_1$}%
}}}
\put(3566,-661){\makebox(0,0)[lb]{\smash{\SetFigFont{12}{14.4}{\familydefault}{\mddefault}{\updefault}{\color[rgb]{0,0,0}$\nu_5$}%
}}}
\put(3361,-466){\makebox(0,0)[lb]{\smash{\SetFigFont{12}{14.4}{\familydefault}{\mddefault}{\updefault}{\color[rgb]{0,0,0}$\nu_{12}$}%
}}}
\put(1100,300){\makebox(0,0)[lb]{\smash{\SetFigFont{12}{14.4}{\familydefault}{\mddefault}{\updefault}{\color[rgb]{0,0,0}$\nu_9$}%
}}}
\put(701,-106){\makebox(0,0)[lb]{\smash{\SetFigFont{12}{14.4}{\familydefault}{\mddefault}{\updefault}{\color[rgb]{0,0,0}$\nu_4$}%
}}}
\end{picture}

\caption{\label{fig:vanishing6}}
\end{figure}

\begin{figure}[H]
\begin{picture}(0,0)%
\includegraphics{vanishing7.pstex}%
\end{picture}%
\setlength{\unitlength}{3947sp}%
\begingroup\makeatletter\ifx\SetFigFont\undefined%
\gdef\SetFigFont#1#2#3#4#5{%
  \reset@font\fontsize{#1}{#2pt}%
  \fontfamily{#3}\fontseries{#4}\fontshape{#5}%
  \selectfont}%
\fi\endgroup%
\begin{picture}(3974,3690)(226,-3206)
\put(700,-75){\makebox(0,0)[lb]{\smash{\SetFigFont{12}{14.4}{\familydefault}{\mddefault}{\updefault}{\color[rgb]{0,0,0}$\nu_8$}%
}}}
\put(460,-250){\makebox(0,0)[lb]{\smash{\SetFigFont{12}{14.4}{\familydefault}{\mddefault}{\updefault}{\color[rgb]{0,0,0}$\nu_{12}$}%
}}}
\put(1576,-2636){\makebox(0,0)[lb]{\smash{\SetFigFont{12}{14.4}{\familydefault}{\mddefault}{\updefault}{\color[rgb]{0,0,0}$\nu_3$}%
}}}
\put(2056,-3136){\makebox(0,0)[lb]{\smash{\SetFigFont{12}{14.4}{\familydefault}{\mddefault}{\updefault}{\color[rgb]{0,0,0}$\nu_8$}%
}}}
\put(3450,-536){\makebox(0,0)[lb]{\smash{\SetFigFont{12}{14.4}{\familydefault}{\mddefault}{\updefault}{\color[rgb]{0,0,0}$\nu_6$}%
}}}
\put(3600,-686){\makebox(0,0)[lb]{\smash{\SetFigFont{12}{14.4}{\familydefault}{\mddefault}{\updefault}{\color[rgb]{0,0,0}$\nu_5$}%
}}}
\put(3740,-836){\makebox(0,0)[lb]{\smash{\SetFigFont{12}{14.4}{\familydefault}{\mddefault}{\updefault}{\color[rgb]{0,0,0}$\nu_1$}%
}}}
\put(4050,-2550){\makebox(0,0)[lb]{\smash{\SetFigFont{12}{14.4}{\familydefault}{\mddefault}{\updefault}{\color[rgb]{0,0,0}$\nu_{11}$}%
}}}
\put(3866,-2761){\makebox(0,0)[lb]{\smash{\SetFigFont{12}{14.4}{\familydefault}{\mddefault}{\updefault}{\color[rgb]{0,0,0}$\nu_4$}%
}}}
\put(3500,-3100){\makebox(0,0)[lb]{\smash{\SetFigFont{12}{14.4}{\familydefault}{\mddefault}{\updefault}{\color[rgb]{0,0,0}$\nu_9$}%
}}}
\end{picture}

\caption{\label{fig:vanishing7}}
\end{figure}

\begin{figure}[H]
\begin{picture}(0,0)%
\includegraphics{vanishing8.pstex}%
\end{picture}%
\setlength{\unitlength}{3947sp}%
\begingroup\makeatletter\ifx\SetFigFont\undefined%
\gdef\SetFigFont#1#2#3#4#5{%
  \reset@font\fontsize{#1}{#2pt}%
  \fontfamily{#3}\fontseries{#4}\fontshape{#5}%
  \selectfont}%
\fi\endgroup%
\begin{picture}(3622,3753)(578,-3269)
\put(1640,-2721){\makebox(0,0)[lb]{\smash{\SetFigFont{12}{14.4}{\familydefault}{\mddefault}{\updefault}{\color[rgb]{0,0,0}$\nu_1$}%
}}}
\put(2101,-3181){\makebox(0,0)[lb]{\smash{\SetFigFont{12}{14.4}{\familydefault}{\mddefault}{\updefault}{\color[rgb]{0,0,0}$\nu_6$}%
}}}
\put(4101,-2461){\makebox(0,0)[lb]{\smash{\SetFigFont{12}{14.4}{\familydefault}{\mddefault}{\updefault}{\color[rgb]{0,0,0}$\nu_8$}%
}}}
\put(3556,-661){\makebox(0,0)[lb]{\smash{\SetFigFont{12}{14.4}{\familydefault}{\mddefault}{\updefault}{\color[rgb]{0,0,0}$\nu_2$}%
}}}
\put(3361,-466){\makebox(0,0)[lb]{\smash{\SetFigFont{12}{14.4}{\familydefault}{\mddefault}{\updefault}{\color[rgb]{0,0,0}$\nu_9$}%
}}}
\put(3720,-811){\makebox(0,0)[lb]{\smash{\SetFigFont{12}{14.4}{\familydefault}{\mddefault}{\updefault}{\color[rgb]{0,0,0}$\nu_{10}$}%
}}}
\put(1080,280){\makebox(0,0)[lb]{\smash{\SetFigFont{12}{14.4}{\familydefault}{\mddefault}{\updefault}{\color[rgb]{0,0,0}$\nu_6$}%
}}}
\put(721,-80){\makebox(0,0)[lb]{\smash{\SetFigFont{12}{14.4}{\familydefault}{\mddefault}{\updefault}{\color[rgb]{0,0,0}$\nu_1$}%
}}}
\put(3500,-3066){\makebox(0,0)[lb]{\smash{\SetFigFont{12}{14.4}{\familydefault}{\mddefault}{\updefault}{\color[rgb]{0,0,0}$\nu_3$}%
}}}
\put(3921,-2656){\makebox(0,0)[lb]{\smash{\SetFigFont{12}{14.4}{\familydefault}{\mddefault}{\updefault}{\color[rgb]{0,0,0}$\nu_4$}%
}}}
\end{picture}

\caption{\label{fig:vanishing8}}
\end{figure}

\begin{table}[H] \extrarowheight0.2em
\[
\begin{array}{l||l|l|l|l|l|l|}
 & 1 & 2 & 3 & 4 & 5 & 6 \\
 \hhline{=#======} 1  \\
 \hhline{--} 2 & x  \\
 \hhline{---} 3 & & my^{-1}  \\
 \hhline{----} 4 & 1 & mx^{-1}z & \\
 & x && \\
 \hhline{-----} 5 & x & x^{-1}y & x^{-1}y & 1  \\
 & 1      & 1 && \\
 & y &&& \\
 \hhline{------} 6 & & my^{-1}z^{-1} & x^{-1}yz^{-1} && my^{-1}z^{-1}  \\
 & & mx^{-1}z^{-1} & x^{-1} && \\
 & & mx^{-1}y^{-1} &&& \\
 \hline 7 & myz & mx^{-1}z && 1 & mz & \\
 & x & mx^{-1}yz && x && \\
 && mz &&&& \\
 \hline 8 & xz^{-1} & z^{-1} & x^{-1} & xy^{-1}z^{-1} & z^{-1} & 1 \\
 & xy^{-1}z^{-1} & mx^{-1}y^{-1} & z^{-1} & y^{-1}z^{-1} & xy^{-1}z^{-1} & \\
 & z^{-1} && x^{-1}yz^{-1} & z^{-1} && \\
 \hline 9 && my^{-1}z^{-1} & x^{-1} && mxy^{-2}z^{-1} & 1 \\
 && mx^{-1}y^{-1} & my^{-1}z^{-1} && my^{-1}z^{-1} & y^{-1}z \\
 && my^{-2} &&& my^{-2} & \\
 \hline 10 & m & mx^{-2} && mx^{-1} & mx^{-1}y^{-1} & \\
 & mx^{-1} & mx^{-1} && y^{-1}z^{-1} & mx^{-1} & \\
 & & mx^{-1}y^{-1} & & & my^{-1} & \\
 \hline 11 & xz^{-1} & my^{-1} & x^{-1}y & xz^{-1} & xz^{-1} & z \\
 && mx^{-1} & 1 & xy^{-1}z^{-1} & my^{-1} & x \\
 && & yz^{-1} & z^{-1} && y \\
 \hline 12 && my^{-1} & my^{-1} && mxy^{-1}z^{-1} & z \\
 &&& mz^{-1} && my^{-1} & mxy^{-1} \\
 &&&&& mxy^{-2} & \\ \hline
\end{array}
\]

\caption{\label{table:intersections1} The points of $\nu_k \cap
\nu_l$ for $1 \leq k \leq 6$ and $l > k$, with their Maslov indices
and $T$-action weights.}
\end{table}

\begin{table}[H] \extrarowheight0.2em
\[
\begin{array}{l||l|l|l|l|l|l|}
 & 7 & 8 & 9 & 10 & 11 & 12 \\
 \hhline{=#======} 7  \\
 \hhline{--} 8 & y^{-1}z^{-1}  \\
 \hhline{---} 9 & & my^{-1}  \\
 \hhline{----} 10 & x^{-1}y^{-1}z^{-1} & mx^{-1}z & \\
 & y^{-1}z^{-1} && \\
 \hhline{-----} 11 & xy^{-1}z^{-1} & y & y & x \\
 & y^{-1}z^{-1}      & x && \\
 & z^{-1} &&& \\
 \hhline{------} 12 & & mxy^{-1} & y && my^{-1}  \\
 & & m & z && \\
 & & my^{-1}z &&& \\
 \hline
\end{array}
\]

\caption{\label{table:intersections2} The points of $\nu_k \cap
\nu_l$ for $7 \leq k \leq 12$ and $l > k$, with their Maslov indices
and $T$-action weights.}
\end{table}

{ \smaller\smaller
\begin{longtable}[c]{ll|l|l}
 $a_1$ & $a_2$ & $b$ & $\pm$ \\
 \hline
 $\nu_1 \cap \nu_2:x$
 & $\nu_2 \cap \nu_5:1$
 & $\nu_1 \cap \nu_5:x$
 & $+$ \\
$\nu_1 \cap \nu_2:x$
 & $\nu_2 \cap \nu_5:x^{-1}y$
 & $\nu_1 \cap \nu_5:y$
 & $-$ \\
$\nu_1 \cap \nu_4:x$
 & $\nu_4 \cap \nu_5:1$
 & $\nu_1 \cap \nu_5:x$
 & $+$ \\
$\nu_1 \cap \nu_4:1$
 & $\nu_4 \cap \nu_5:1$
 & $\nu_1 \cap \nu_5:1$
 & $-$ \\
$\nu_2 \cap \nu_3:my^{-1}$
 & $\nu_3 \cap \nu_6:x^{-1}yz^{-1}$
 & $\nu_2 \cap \nu_6:mx^{-1}z^{-1}$
 & $-$ \\
$\nu_2 \cap \nu_3:my^{-1}$
 & $\nu_3 \cap \nu_6:x^{-1}$
 & $\nu_2 \cap \nu_6:mx^{-1}y^{-1}$
 & $-$ \\
$\nu_2 \cap \nu_5:x^{-1}y$
 & $\nu_5 \cap \nu_6:my^{-1}z^{-1}$
 & $\nu_2 \cap \nu_6:mx^{-1}z^{-1}$
 & $-$ \\
$\nu_2 \cap \nu_5:1$
 & $\nu_5 \cap \nu_6:my^{-1}z^{-1}$
 & $\nu_2 \cap \nu_6:my^{-1}z^{-1}$
 & $+$ \\
$\nu_4 \cap \nu_5:1$
 & $\nu_5 \cap \nu_8:xy^{-1}z^{-1}$
 & $\nu_4 \cap \nu_8:xy^{-1}z^{-1}$
 & $+$ \\
$\nu_4 \cap \nu_5:1$
 & $\nu_5 \cap \nu_8:z^{-1}$
 & $\nu_4 \cap \nu_8:z^{-1}$
 & $-$ \\
$\nu_4 \cap \nu_7:x$
 & $\nu_7 \cap \nu_8:y^{-1}z^{-1}$
 & $\nu_4 \cap \nu_8:xy^{-1}z^{-1}$
 & $+$ \\
$\nu_4 \cap \nu_7:1$
 & $\nu_7 \cap \nu_8:y^{-1}z^{-1}$
 & $\nu_4 \cap \nu_8:y^{-1}z^{-1}$
 & $-$ \\
$\nu_5 \cap \nu_6:my^{-1}z^{-1}$
 & $\nu_6 \cap \nu_9:1$
 & $\nu_5 \cap \nu_9:my^{-1}z^{-1}$
 & $-$ \\
$\nu_5 \cap \nu_6:my^{-1}z^{-1}$
 & $\nu_6 \cap \nu_9:y^{-1}z$
 & $\nu_5 \cap \nu_9:my^{-2}$
 & $-$ \\
$\nu_5 \cap \nu_8:z^{-1}$
 & $\nu_8 \cap \nu_9:my^{-1}$
 & $\nu_5 \cap \nu_9:my^{-1}z^{-1}$
 & $-$ \\
$\nu_5 \cap \nu_8:xy^{-1}z^{-1}$
 & $\nu_8 \cap \nu_9:my^{-1}$
 & $\nu_5 \cap \nu_9:mxy^{-2}z^{-1}$
 & $+$ \\
$\nu_7 \cap \nu_8:y^{-1}z^{-1}$
 & $\nu_8 \cap \nu_{11}:x$
 & $\nu_7 \cap \nu_{11}:xy^{-1}z^{-1}$
 & $+$ \\
$\nu_7 \cap \nu_8:y^{-1}z^{-1}$
 & $\nu_8 \cap \nu_{11}:y$
 & $\nu_7 \cap \nu_{11}:z^{-1}$
 & $-$ \\
$\nu_7 \cap \nu_{10}:y^{-1}z^{-1}$
 & $\nu_{10} \cap \nu_{11}:x$
 & $\nu_7 \cap \nu_{11}:xy^{-1}z^{-1}$
 & $+$ \\
$\nu_7 \cap \nu_{10}:x^{-1}y^{-1}z^{-1}$
 & $\nu_{10} \cap \nu_{11}:x$
 & $\nu_7 \cap \nu_{11}:y^{-1}z^{-1}$
 & $-$ \\
$\nu_8 \cap \nu_9:my^{-1}$
 & $\nu_9 \cap \nu_{12}:y$
 & $\nu_8 \cap \nu_{12}:m$
 & $-$ \\
$\nu_8 \cap \nu_9:my^{-1}$
 & $\nu_9 \cap \nu_{12}:z$
 & $\nu_8 \cap \nu_{12}:my^{-1}z$
 & $-$ \\
$\nu_8 \cap \nu_{11}:y$
 & $\nu_{11} \cap \nu_{12}:my^{-1}$
 & $\nu_8 \cap \nu_{12}:m$
 & $-$ \\
$\nu_8 \cap \nu_{11}:x$
 & $\nu_{11} \cap \nu_{12}:my^{-1}$
 & $\nu_8 \cap \nu_{12}:mxy^{-1}$
 & $+$ \\
$\nu_2 \cap \nu_4:mx^{-1}z$
 & $\nu_4 \cap \nu_7:1$
 & $\nu_2 \cap \nu_7:mx^{-1}z$
 & $-$ \\
$\nu_2 \cap \nu_4:mx^{-1}z$
 & $\nu_4 \cap \nu_7:x$
 & $\nu_2 \cap \nu_7:mz$
 & $+$ \\
$\nu_2 \cap \nu_5:1$
 & $\nu_5 \cap \nu_7:mz$
 & $\nu_2 \cap \nu_7:mz$
 & $+$ \\
$\nu_2 \cap \nu_5:x^{-1}y$
 & $\nu_5 \cap \nu_7:mz$
 & $\nu_2 \cap \nu_7:mx^{-1}yz$
 & $-$ \\
$\nu_3 \cap \nu_5:x^{-1}y$
 & $\nu_5 \cap \nu_8:xy^{-1}z^{-1}$
 & $\nu_3 \cap \nu_8:z^{-1}$
 & $-$ \\
$\nu_3 \cap \nu_5:x^{-1}y$
 & $\nu_5 \cap \nu_8:z^{-1}$
 & $\nu_3 \cap \nu_8:x^{-1}yz^{-1}$
 & $-$ \\
$\nu_3 \cap \nu_6:x^{-1}yz^{-1}$
 & $\nu_6 \cap \nu_8:1$
 & $\nu_3 \cap \nu_8:x^{-1}yz^{-1}$
 & $-$ \\
$\nu_3 \cap \nu_6:x^{-1}$
 & $\nu_6 \cap \nu_8:1$
 & $\nu_3 \cap \nu_8:x^{-1}$
 & $+$ \\
$\nu_1 \cap \nu_2:x$
 & $\nu_2 \cap \nu_7:mx^{-1}yz$
 & $\nu_1 \cap \nu_7:myz$
 & $+$ \\
$\nu_1 \cap \nu_4:x$
 & $\nu_4 \cap \nu_7:1$
 & $\nu_1 \cap \nu_7:x$
 & $-$ \\
$\nu_1 \cap \nu_4:1$
 & $\nu_4 \cap \nu_7:x$
 & $\nu_1 \cap \nu_7:x$
 & $+$ \\
$\nu_1 \cap \nu_5:y$
 & $\nu_5 \cap \nu_7:mz$
 & $\nu_1 \cap \nu_7:myz$
 & $+$ \\
$\nu_2 \cap \nu_3:my^{-1}$
 & $\nu_3 \cap \nu_8:x^{-1}$
 & $\nu_2 \cap \nu_8:mx^{-1}y^{-1}$
 & $+$ \\
$\nu_2 \cap \nu_4:mx^{-1}z$
 & $\nu_4 \cap \nu_8:y^{-1}z^{-1}$
 & $\nu_2 \cap \nu_8:mx^{-1}y^{-1}$
 & $+$ \\
$\nu_2 \cap \nu_5:x^{-1}y$
 & $\nu_5 \cap \nu_8:xy^{-1}z^{-1}$
 & $\nu_2 \cap \nu_8:z^{-1}$
 & $+$ \\
$\nu_2 \cap \nu_5:1$
 & $\nu_5 \cap \nu_8:z^{-1}$
 & $\nu_2 \cap \nu_8:z^{-1}$
 & $-$ \\
$\nu_2 \cap \nu_6:mx^{-1}y^{-1}$
 & $\nu_6 \cap \nu_8:1$
 & $\nu_2 \cap \nu_8:mx^{-1}y^{-1}$
 & $+$ \\
$\nu_2 \cap \nu_7:mx^{-1}z$
 & $\nu_7 \cap \nu_8:y^{-1}z^{-1}$
 & $\nu_2 \cap \nu_8:mx^{-1}y^{-1}$
 & $-$ \\
$\nu_3 \cap \nu_5:x^{-1}y$
 & $\nu_5 \cap \nu_9:mxy^{-2}z^{-1}$
 & $\nu_3 \cap \nu_9:my^{-1}z^{-1}$
 & $+$ \\
$\nu_3 \cap \nu_6:x^{-1}$
 & $\nu_6 \cap \nu_9:1$
 & $\nu_3 \cap \nu_9:x^{-1}$
 & $-$ \\
$\nu_3 \cap \nu_6:x^{-1}yz^{-1}$
 & $\nu_6 \cap \nu_9:y^{-1}z$
 & $\nu_3 \cap \nu_9:x^{-1}$
 & $+$ \\
$\nu_3 \cap \nu_8:z^{-1}$
 & $\nu_8 \cap \nu_9:my^{-1}$
 & $\nu_3 \cap \nu_9:my^{-1}z^{-1}$
 & $-$ \\
$\nu_1 \cap \nu_2:x$
 & $\nu_2 \cap \nu_8:z^{-1}$
 & $\nu_1 \cap \nu_8:xz^{-1}$
 & $-$ \\
$\nu_1 \cap \nu_4:x$
 & $\nu_4 \cap \nu_8:z^{-1}$
 & $\nu_1 \cap \nu_8:xz^{-1}$
 & $-$ \\
$\nu_1 \cap \nu_4:x$
 & $\nu_4 \cap \nu_8:y^{-1}z^{-1}$
 & $\nu_1 \cap \nu_8:xy^{-1}z^{-1}$
 & $+$ \\
$\nu_1 \cap \nu_4:1$
 & $\nu_4 \cap \nu_8:z^{-1}$
 & $\nu_1 \cap \nu_8:z^{-1}$
 & $+$ \\
$\nu_1 \cap \nu_4:1$
 & $\nu_4 \cap \nu_8:xy^{-1}z^{-1}$
 & $\nu_1 \cap \nu_8:xy^{-1}z^{-1}$
 & $+$ \\
$\nu_1 \cap \nu_5:x$
 & $\nu_5 \cap \nu_8:z^{-1}$
 & $\nu_1 \cap \nu_8:xz^{-1}$
 & $+$ \\
$\nu_1 \cap \nu_5:y$
 & $\nu_5 \cap \nu_8:xy^{-1}z^{-1}$
 & $\nu_1 \cap \nu_8:xz^{-1}$
 & $+$ \\
$\nu_1 \cap \nu_5:1$
 & $\nu_5 \cap \nu_8:z^{-1}$
 & $\nu_1 \cap \nu_8:z^{-1}$
 & $+$ \\
$\nu_1 \cap \nu_5:1$
 & $\nu_5 \cap \nu_8:xy^{-1}z^{-1}$
 & $\nu_1 \cap \nu_8:xy^{-1}z^{-1}$
 & $-$ \\
$\nu_1 \cap \nu_7:x$
 & $\nu_7 \cap \nu_8:y^{-1}z^{-1}$
 & $\nu_1 \cap \nu_8:xy^{-1}z^{-1}$
 & $+$ \\
$\nu_2 \cap \nu_3:my^{-1}$
 & $\nu_3 \cap \nu_9:x^{-1}$
 & $\nu_2 \cap \nu_9:mx^{-1}y^{-1}$
 & $+$ \\
$\nu_2 \cap \nu_5:1$
 & $\nu_5 \cap \nu_9:my^{-2}$
 & $\nu_2 \cap \nu_9:my^{-2}$
 & $+$ \\
$\nu_2 \cap \nu_5:1$
 & $\nu_5 \cap \nu_9:my^{-1}z^{-1}$
 & $\nu_2 \cap \nu_9:my^{-1}z^{-1}$
 & $+$ \\
$\nu_2 \cap \nu_5:x^{-1}y$
 & $\nu_5 \cap \nu_9:mxy^{-2}z^{-1}$
 & $\nu_2 \cap \nu_9:my^{-1}z^{-1}$
 & $+$ \\
$\nu_2 \cap \nu_5:x^{-1}y$
 & $\nu_5 \cap \nu_9:my^{-2}$
 & $\nu_2 \cap \nu_9:mx^{-1}y^{-1}$
 & $+$ \\
$\nu_2 \cap \nu_6:my^{-1}z^{-1}$
 & $\nu_6 \cap \nu_9:1$
 & $\nu_2 \cap \nu_9:my^{-1}z^{-1}$
 & $-$ \\
$\nu_2 \cap \nu_6:mx^{-1}z^{-1}$
 & $\nu_6 \cap \nu_9:y^{-1}z$
 & $\nu_2 \cap \nu_9:mx^{-1}y^{-1}$
 & $+$ \\
$\nu_2 \cap \nu_6:mx^{-1}y^{-1}$
 & $\nu_6 \cap \nu_9:1$
 & $\nu_2 \cap \nu_9:mx^{-1}y^{-1}$
 & $-$ \\
$\nu_2 \cap \nu_6:my^{-1}z^{-1}$
 & $\nu_6 \cap \nu_9:y^{-1}z$
 & $\nu_2 \cap \nu_9:my^{-2}$
 & $-$ \\
$\nu_2 \cap \nu_8:z^{-1}$
 & $\nu_8 \cap \nu_9:my^{-1}$
 & $\nu_2 \cap \nu_9:my^{-1}z^{-1}$
 & $+$ \\
$\nu_5 \cap \nu_7:mz$
 & $\nu_7 \cap \nu_{10}:x^{-1}y^{-1}z^{-1}$
 & $\nu_5 \cap \nu_{10}:mx^{-1}y^{-1}$
 & $-$ \\
$\nu_5 \cap \nu_7:mz$
 & $\nu_7 \cap \nu_{10}:y^{-1}z^{-1}$
 & $\nu_5 \cap \nu_{10}:my^{-1}$
 & $+$ \\
$\nu_5 \cap \nu_8:xy^{-1}z^{-1}$
 & $\nu_8 \cap \nu_{10}:mx^{-1}z$
 & $\nu_5 \cap \nu_{10}:my^{-1}$
 & $+$ \\
$\nu_5 \cap \nu_8:z^{-1}$
 & $\nu_8 \cap \nu_{10}:mx^{-1}z$
 & $\nu_5 \cap \nu_{10}:mx^{-1}$
 & $-$ \\
$\nu_6 \cap \nu_8:1$
 & $\nu_8 \cap \nu_{11}:x$
 & $\nu_6 \cap \nu_{11}:x$
 & $-$ \\
$\nu_6 \cap \nu_8:1$
 & $\nu_8 \cap \nu_{11}:y$
 & $\nu_6 \cap \nu_{11}:y$
 & $-$ \\
$\nu_6 \cap \nu_9:1$
 & $\nu_9 \cap \nu_{11}:y$
 & $\nu_6 \cap \nu_{11}:y$
 & $-$ \\
$\nu_6 \cap \nu_9:y^{-1}z$
 & $\nu_9 \cap \nu_{11}:y$
 & $\nu_6 \cap \nu_{11}:z$
 & $+$ \\
$\nu_4 \cap \nu_5:1$
 & $\nu_5 \cap \nu_{10}:mx^{-1}$
 & $\nu_4 \cap \nu_{10}:mx^{-1}$
 & $+$ \\
$\nu_4 \cap \nu_7:x$
 & $\nu_7 \cap \nu_{10}:x^{-1}y^{-1}z^{-1}$
 & $\nu_4 \cap \nu_{10}:y^{-1}z^{-1}$
 & $-$ \\
$\nu_4 \cap \nu_7:1$
 & $\nu_7 \cap \nu_{10}:y^{-1}z^{-1}$
 & $\nu_4 \cap \nu_{10}:y^{-1}z^{-1}$
 & $+$ \\
$\nu_4 \cap \nu_8:z^{-1}$
 & $\nu_8 \cap \nu_{10}:mx^{-1}z$
 & $\nu_4 \cap \nu_{10}:mx^{-1}$
 & $+$ \\
$\nu_5 \cap \nu_6:my^{-1}z^{-1}$
 & $\nu_6 \cap \nu_{11}:z$
 & $\nu_5 \cap \nu_{11}:my^{-1}$
 & $+$ \\
$\nu_5 \cap \nu_7:mz$
 & $\nu_7 \cap \nu_{11}:y^{-1}z^{-1}$
 & $\nu_5 \cap \nu_{11}:my^{-1}$
 & $+$ \\
$\nu_5 \cap \nu_8:z^{-1}$
 & $\nu_8 \cap \nu_{11}:x$
 & $\nu_5 \cap \nu_{11}:xz^{-1}$
 & $+$ \\
$\nu_5 \cap \nu_8:xy^{-1}z^{-1}$
 & $\nu_8 \cap \nu_{11}:y$
 & $\nu_5 \cap \nu_{11}:xz^{-1}$
 & $-$ \\
$\nu_5 \cap \nu_9:my^{-2}$
 & $\nu_9 \cap \nu_{11}:y$
 & $\nu_5 \cap \nu_{11}:my^{-1}$
 & $+$ \\
$\nu_5 \cap \nu_{10}:mx^{-1}y^{-1}$
 & $\nu_{10} \cap \nu_{11}:x$
 & $\nu_5 \cap \nu_{11}:my^{-1}$
 & $-$ \\
$\nu_6 \cap \nu_8:1$
 & $\nu_8 \cap \nu_{12}:mxy^{-1}$
 & $\nu_6 \cap \nu_{12}:mxy^{-1}$
 & $+$ \\
$\nu_6 \cap \nu_9:y^{-1}z$
 & $\nu_9 \cap \nu_{12}:y$
 & $\nu_6 \cap \nu_{12}:z$
 & $-$ \\
$\nu_6 \cap \nu_9:1$
 & $\nu_9 \cap \nu_{12}:z$
 & $\nu_6 \cap \nu_{12}:z$
 & $+$ \\
$\nu_6 \cap \nu_{11}:x$
 & $\nu_{11} \cap \nu_{12}:my^{-1}$
 & $\nu_6 \cap \nu_{12}:mxy^{-1}$
 & $-$ \\
$\nu_4 \cap \nu_5:1$
 & $\nu_5 \cap \nu_{11}:xz^{-1}$
 & $\nu_4 \cap \nu_{11}:xz^{-1}$
 & $-$ \\
$\nu_4 \cap \nu_7:x$
 & $\nu_7 \cap \nu_{11}:z^{-1}$
 & $\nu_4 \cap \nu_{11}:xz^{-1}$
 & $-$ \\
$\nu_4 \cap \nu_7:x$
 & $\nu_7 \cap \nu_{11}:y^{-1}z^{-1}$
 & $\nu_4 \cap \nu_{11}:xy^{-1}z^{-1}$
 & $+$ \\
$\nu_4 \cap \nu_7:1$
 & $\nu_7 \cap \nu_{11}:z^{-1}$
 & $\nu_4 \cap \nu_{11}:z^{-1}$
 & $+$ \\
$\nu_4 \cap \nu_7:1$
 & $\nu_7 \cap \nu_{11}:xy^{-1}z^{-1}$
 & $\nu_4 \cap \nu_{11}:xy^{-1}z^{-1}$
 & $+$ \\
$\nu_4 \cap \nu_8:xy^{-1}z^{-1}$
 & $\nu_8 \cap \nu_{11}:y$
 & $\nu_4 \cap \nu_{11}:xz^{-1}$
 & $+$ \\
$\nu_4 \cap \nu_8:z^{-1}$
 & $\nu_8 \cap \nu_{11}:x$
 & $\nu_4 \cap \nu_{11}:xz^{-1}$
 & $+$ \\
$\nu_4 \cap \nu_8:y^{-1}z^{-1}$
 & $\nu_8 \cap \nu_{11}:y$
 & $\nu_4 \cap \nu_{11}:z^{-1}$
 & $+$ \\
$\nu_4 \cap \nu_8:y^{-1}z^{-1}$
 & $\nu_8 \cap \nu_{11}:x$
 & $\nu_4 \cap \nu_{11}:xy^{-1}z^{-1}$
 & $-$ \\
$\nu_4 \cap \nu_{10}:y^{-1}z^{-1}$
 & $\nu_{10} \cap \nu_{11}:x$
 & $\nu_4 \cap \nu_{11}:xy^{-1}z^{-1}$
 & $+$ \\
$\nu_5 \cap \nu_6:my^{-1}z^{-1}$
 & $\nu_6 \cap \nu_{12}:z$
 & $\nu_5 \cap \nu_{12}:my^{-1}$
 & $+$ \\
$\nu_5 \cap \nu_8:xy^{-1}z^{-1}$
 & $\nu_8 \cap \nu_{12}:my^{-1}z$
 & $\nu_5 \cap \nu_{12}:mxy^{-2}$
 & $+$ \\
$\nu_5 \cap \nu_8:xy^{-1}z^{-1}$
 & $\nu_8 \cap \nu_{12}:m$
 & $\nu_5 \cap \nu_{12}:mxy^{-1}z^{-1}$
 & $+$ \\
$\nu_5 \cap \nu_8:z^{-1}$
 & $\nu_8 \cap \nu_{12}:mxy^{-1}$
 & $\nu_5 \cap \nu_{12}:mxy^{-1}z^{-1}$
 & $+$ \\
$\nu_5 \cap \nu_8:z^{-1}$
 & $\nu_8 \cap \nu_{12}:my^{-1}z$
 & $\nu_5 \cap \nu_{12}:my^{-1}$
 & $+$ \\
$\nu_5 \cap \nu_9:mxy^{-2}z^{-1}$
 & $\nu_9 \cap \nu_{12}:y$
 & $\nu_5 \cap \nu_{12}:mxy^{-1}z^{-1}$
 & $-$ \\
$\nu_5 \cap \nu_9:my^{-1}z^{-1}$
 & $\nu_9 \cap \nu_{12}:z$
 & $\nu_5 \cap \nu_{12}:my^{-1}$
 & $+$ \\
$\nu_5 \cap \nu_9:my^{-2}$
 & $\nu_9 \cap \nu_{12}:y$
 & $\nu_5 \cap \nu_{12}:my^{-1}$
 & $-$ \\
$\nu_5 \cap \nu_9:mxy^{-2}z^{-1}$
 & $\nu_9 \cap \nu_{12}:z$
 & $\nu_5 \cap \nu_{12}:mxy^{-2}$
 & $-$ \\
$\nu_5 \cap \nu_{11}:xz^{-1}$
 & $\nu_{11} \cap \nu_{12}:my^{-1}$
 & $\nu_5 \cap \nu_{12}:mxy^{-1}z^{-1}$
 & $+$ \\
$\nu_2 \cap \nu_4:mx^{-1}z$
 & $\nu_4 \cap \nu_{10}:y^{-1}z^{-1}$
 & $\nu_2 \cap \nu_{10}:mx^{-1}y^{-1}$
 & $-$ \\
$\nu_2 \cap \nu_5:x^{-1}y$
 & $\nu_5 \cap \nu_{10}:my^{-1}$
 & $\nu_2 \cap \nu_{10}:mx^{-1}$
 & $+$ \\
$\nu_2 \cap \nu_5:x^{-1}y$
 & $\nu_5 \cap \nu_{10}:mx^{-1}y^{-1}$
 & $\nu_2 \cap \nu_{10}:mx^{-2}$
 & $-$ \\
$\nu_2 \cap \nu_5:1$
 & $\nu_5 \cap \nu_{10}:mx^{-1}y^{-1}$
 & $\nu_2 \cap \nu_{10}:mx^{-1}y^{-1}$
 & $+$ \\
$\nu_2 \cap \nu_5:1$
 & $\nu_5 \cap \nu_{10}:mx^{-1}$
 & $\nu_2 \cap \nu_{10}:mx^{-1}$
 & $+$ \\
$\nu_2 \cap \nu_7:mx^{-1}yz$
 & $\nu_7 \cap \nu_{10}:y^{-1}z^{-1}$
 & $\nu_2 \cap \nu_{10}:mx^{-1}$
 & $-$ \\
$\nu_2 \cap \nu_7:mx^{-1}yz$
 & $\nu_7 \cap \nu_{10}:x^{-1}y^{-1}z^{-1}$
 & $\nu_2 \cap \nu_{10}:mx^{-2}$
 & $-$ \\
$\nu_2 \cap \nu_7:mz$
 & $\nu_7 \cap \nu_{10}:x^{-1}y^{-1}z^{-1}$
 & $\nu_2 \cap \nu_{10}:mx^{-1}y^{-1}$
 & $-$ \\
$\nu_2 \cap \nu_7:mx^{-1}z$
 & $\nu_7 \cap \nu_{10}:y^{-1}z^{-1}$
 & $\nu_2 \cap \nu_{10}:mx^{-1}y^{-1}$
 & $-$ \\
$\nu_2 \cap \nu_8:z^{-1}$
 & $\nu_8 \cap \nu_{10}:mx^{-1}z$
 & $\nu_2 \cap \nu_{10}:mx^{-1}$
 & $+$ \\
$\nu_3 \cap \nu_5:x^{-1}y$
 & $\nu_5 \cap \nu_{11}:xz^{-1}$
 & $\nu_3 \cap \nu_{11}:yz^{-1}$
 & $-$ \\
$\nu_3 \cap \nu_6:x^{-1}yz^{-1}$
 & $\nu_6 \cap \nu_{11}:z$
 & $\nu_3 \cap \nu_{11}:x^{-1}y$
 & $+$ \\
$\nu_3 \cap \nu_6:x^{-1}yz^{-1}$
 & $\nu_6 \cap \nu_{11}:x$
 & $\nu_3 \cap \nu_{11}:yz^{-1}$
 & $+$ \\
$\nu_3 \cap \nu_6:x^{-1}$
 & $\nu_6 \cap \nu_{11}:y$
 & $\nu_3 \cap \nu_{11}:x^{-1}y$
 & $+$ \\
$\nu_3 \cap \nu_6:x^{-1}$
 & $\nu_6 \cap \nu_{11}:x$
 & $\nu_3 \cap \nu_{11}:1$
 & $-$ \\
$\nu_3 \cap \nu_8:x^{-1}yz^{-1}$
 & $\nu_8 \cap \nu_{11}:x$
 & $\nu_3 \cap \nu_{11}:yz^{-1}$
 & $+$ \\
$\nu_3 \cap \nu_8:x^{-1}$
 & $\nu_8 \cap \nu_{11}:x$
 & $\nu_3 \cap \nu_{11}:1$
 & $+$ \\
$\nu_3 \cap \nu_8:z^{-1}$
 & $\nu_8 \cap \nu_{11}:y$
 & $\nu_3 \cap \nu_{11}:yz^{-1}$
 & $-$ \\
$\nu_3 \cap \nu_8:x^{-1}$
 & $\nu_8 \cap \nu_{11}:y$
 & $\nu_3 \cap \nu_{11}:x^{-1}y$
 & $-$ \\
$\nu_3 \cap \nu_9:x^{-1}$
 & $\nu_9 \cap \nu_{11}:y$
 & $\nu_3 \cap \nu_{11}:x^{-1}y$
 & $+$ \\
$\nu_1 \cap \nu_2:x$
 & $\nu_2 \cap \nu_{10}:mx^{-1}$
 & $\nu_1 \cap \nu_{10}:m$
 & $+$ \\
$\nu_1 \cap \nu_2:x$
 & $\nu_2 \cap \nu_{10}:mx^{-2}$
 & $\nu_1 \cap \nu_{10}:mx^{-1}$
 & $+$ \\
$\nu_1 \cap \nu_4:x$
 & $\nu_4 \cap \nu_{10}:mx^{-1}$
 & $\nu_1 \cap \nu_{10}:m$
 & $+$ \\
$\nu_1 \cap \nu_4:1$
 & $\nu_4 \cap \nu_{10}:mx^{-1}$
 & $\nu_1 \cap \nu_{10}:mx^{-1}$
 & $+$ \\
$\nu_1 \cap \nu_5:y$
 & $\nu_5 \cap \nu_{10}:mx^{-1}y^{-1}$
 & $\nu_1 \cap \nu_{10}:mx^{-1}$
 & $+$ \\
$\nu_1 \cap \nu_5:1$
 & $\nu_5 \cap \nu_{10}:mx^{-1}$
 & $\nu_1 \cap \nu_{10}:mx^{-1}$
 & $-$ \\
$\nu_1 \cap \nu_5:y$
 & $\nu_5 \cap \nu_{10}:my^{-1}$
 & $\nu_1 \cap \nu_{10}:m$
 & $-$ \\
$\nu_1 \cap \nu_5:x$
 & $\nu_5 \cap \nu_{10}:mx^{-1}$
 & $\nu_1 \cap \nu_{10}:m$
 & $+$ \\
$\nu_1 \cap \nu_7:myz$
 & $\nu_7 \cap \nu_{10}:x^{-1}y^{-1}z^{-1}$
 & $\nu_1 \cap \nu_{10}:mx^{-1}$
 & $-$ \\
$\nu_1 \cap \nu_7:myz$
 & $\nu_7 \cap \nu_{10}:y^{-1}z^{-1}$
 & $\nu_1 \cap \nu_{10}:m$
 & $-$ \\
$\nu_1 \cap \nu_8:xz^{-1}$
 & $\nu_8 \cap \nu_{10}:mx^{-1}z$
 & $\nu_1 \cap \nu_{10}:m$
 & $-$ \\
$\nu_1 \cap \nu_8:z^{-1}$
 & $\nu_8 \cap \nu_{10}:mx^{-1}z$
 & $\nu_1 \cap \nu_{10}:mx^{-1}$
 & $+$ \\
$\nu_2 \cap \nu_3:my^{-1}$
 & $\nu_3 \cap \nu_{11}:x^{-1}y$
 & $\nu_2 \cap \nu_{11}:mx^{-1}$
 & $-$ \\
$\nu_2 \cap \nu_3:my^{-1}$
 & $\nu_3 \cap \nu_{11}:1$
 & $\nu_2 \cap \nu_{11}:my^{-1}$
 & $+$ \\
$\nu_2 \cap \nu_4:mx^{-1}z$
 & $\nu_4 \cap \nu_{11}:z^{-1}$
 & $\nu_2 \cap \nu_{11}:mx^{-1}$
 & $+$ \\
$\nu_2 \cap \nu_4:mx^{-1}z$
 & $\nu_4 \cap \nu_{11}:xy^{-1}z^{-1}$
 & $\nu_2 \cap \nu_{11}:my^{-1}$
 & $-$ \\
$\nu_2 \cap \nu_5:1$
 & $\nu_5 \cap \nu_{11}:my^{-1}$
 & $\nu_2 \cap \nu_{11}:my^{-1}$
 & $+$ \\
$\nu_2 \cap \nu_5:x^{-1}y$
 & $\nu_5 \cap \nu_{11}:my^{-1}$
 & $\nu_2 \cap \nu_{11}:mx^{-1}$
 & $+$ \\
$\nu_2 \cap \nu_6:my^{-1}z^{-1}$
 & $\nu_6 \cap \nu_{11}:z$
 & $\nu_2 \cap \nu_{11}:my^{-1}$
 & $+$ \\
$\nu_2 \cap \nu_6:mx^{-1}y^{-1}$
 & $\nu_6 \cap \nu_{11}:x$
 & $\nu_2 \cap \nu_{11}:my^{-1}$
 & $-$ \\
$\nu_2 \cap \nu_6:mx^{-1}z^{-1}$
 & $\nu_6 \cap \nu_{11}:z$
 & $\nu_2 \cap \nu_{11}:mx^{-1}$
 & $-$ \\
$\nu_2 \cap \nu_6:mx^{-1}y^{-1}$
 & $\nu_6 \cap \nu_{11}:y$
 & $\nu_2 \cap \nu_{11}:mx^{-1}$
 & $-$ \\
$\nu_2 \cap \nu_7:mx^{-1}z$
 & $\nu_7 \cap \nu_{11}:xy^{-1}z^{-1}$
 & $\nu_2 \cap \nu_{11}:my^{-1}$
 & $-$ \\
$\nu_2 \cap \nu_7:mz$
 & $\nu_7 \cap \nu_{11}:y^{-1}z^{-1}$
 & $\nu_2 \cap \nu_{11}:my^{-1}$
 & $+$ \\
$\nu_2 \cap \nu_7:mx^{-1}yz$
 & $\nu_7 \cap \nu_{11}:y^{-1}z^{-1}$
 & $\nu_2 \cap \nu_{11}:mx^{-1}$
 & $-$ \\
$\nu_2 \cap \nu_7:mx^{-1}z$
 & $\nu_7 \cap \nu_{11}:z^{-1}$
 & $\nu_2 \cap \nu_{11}:mx^{-1}$
 & $+$ \\
$\nu_2 \cap \nu_8:mx^{-1}y^{-1}$
 & $\nu_8 \cap \nu_{11}:y$
 & $\nu_2 \cap \nu_{11}:mx^{-1}$
 & $-$ \\
$\nu_2 \cap \nu_8:mx^{-1}y^{-1}$
 & $\nu_8 \cap \nu_{11}:x$
 & $\nu_2 \cap \nu_{11}:my^{-1}$
 & $-$ \\
$\nu_2 \cap \nu_9:my^{-2}$
 & $\nu_9 \cap \nu_{11}:y$
 & $\nu_2 \cap \nu_{11}:my^{-1}$
 & $+$ \\
$\nu_2 \cap \nu_9:mx^{-1}y^{-1}$
 & $\nu_9 \cap \nu_{11}:y$
 & $\nu_2 \cap \nu_{11}:mx^{-1}$
 & $+$ \\
$\nu_2 \cap \nu_{10}:mx^{-1}y^{-1}$
 & $\nu_{10} \cap \nu_{11}:x$
 & $\nu_2 \cap \nu_{11}:my^{-1}$
 & $-$ \\
$\nu_2 \cap \nu_{10}:mx^{-2}$
 & $\nu_{10} \cap \nu_{11}:x$
 & $\nu_2 \cap \nu_{11}:mx^{-1}$
 & $+$ \\
$\nu_3 \cap \nu_5:x^{-1}y$
 & $\nu_5 \cap \nu_{12}:mxy^{-1}z^{-1}$
 & $\nu_3 \cap \nu_{12}:mz^{-1}$
 & $+$ \\
$\nu_3 \cap \nu_5:x^{-1}y$
 & $\nu_5 \cap \nu_{12}:mxy^{-2}$
 & $\nu_3 \cap \nu_{12}:my^{-1}$
 & $-$ \\
$\nu_3 \cap \nu_6:x^{-1}$
 & $\nu_6 \cap \nu_{12}:mxy^{-1}$
 & $\nu_3 \cap \nu_{12}:my^{-1}$
 & $-$ \\
$\nu_3 \cap \nu_6:x^{-1}yz^{-1}$
 & $\nu_6 \cap \nu_{12}:mxy^{-1}$
 & $\nu_3 \cap \nu_{12}:mz^{-1}$
 & $+$ \\
$\nu_3 \cap \nu_8:z^{-1}$
 & $\nu_8 \cap \nu_{12}:my^{-1}z$
 & $\nu_3 \cap \nu_{12}:my^{-1}$
 & $+$ \\
$\nu_3 \cap \nu_8:x^{-1}$
 & $\nu_8 \cap \nu_{12}:mxy^{-1}$
 & $\nu_3 \cap \nu_{12}:my^{-1}$
 & $-$ \\
$\nu_3 \cap \nu_8:z^{-1}$
 & $\nu_8 \cap \nu_{12}:m$
 & $\nu_3 \cap \nu_{12}:mz^{-1}$
 & $-$ \\
$\nu_3 \cap \nu_8:x^{-1}yz^{-1}$
 & $\nu_8 \cap \nu_{12}:mxy^{-1}$
 & $\nu_3 \cap \nu_{12}:mz^{-1}$
 & $-$ \\
$\nu_3 \cap \nu_9:my^{-1}z^{-1}$
 & $\nu_9 \cap \nu_{12}:y$
 & $\nu_3 \cap \nu_{12}:mz^{-1}$
 & $-$ \\
$\nu_3 \cap \nu_9:my^{-1}z^{-1}$
 & $\nu_9 \cap \nu_{12}:z$
 & $\nu_3 \cap \nu_{12}:my^{-1}$
 & $+$ \\
$\nu_3 \cap \nu_{11}:yz^{-1}$
 & $\nu_{11} \cap \nu_{12}:my^{-1}$
 & $\nu_3 \cap \nu_{12}:mz^{-1}$
 & $-$ \\
$\nu_3 \cap \nu_{11}:1$
 & $\nu_{11} \cap \nu_{12}:my^{-1}$
 & $\nu_3 \cap \nu_{12}:my^{-1}$
 & $-$ \\
$\nu_1 \cap \nu_4:1$
 & $\nu_4 \cap \nu_{11}:xz^{-1}$
 & $\nu_1 \cap \nu_{11}:xz^{-1}$
 & $+$ \\
$\nu_1 \cap \nu_4:x$
 & $\nu_4 \cap \nu_{11}:z^{-1}$
 & $\nu_1 \cap \nu_{11}:xz^{-1}$
 & $+$ \\
$\nu_1 \cap \nu_5:1$
 & $\nu_5 \cap \nu_{11}:xz^{-1}$
 & $\nu_1 \cap \nu_{11}:xz^{-1}$
 & $+$ \\
$\nu_1 \cap \nu_7:x$
 & $\nu_7 \cap \nu_{11}:z^{-1}$
 & $\nu_1 \cap \nu_{11}:xz^{-1}$
 & $-$ \\
$\nu_1 \cap \nu_8:z^{-1}$
 & $\nu_8 \cap \nu_{11}:x$
 & $\nu_1 \cap \nu_{11}:xz^{-1}$
 & $+$ \\
$\nu_1 \cap \nu_8:xy^{-1}z^{-1}$
 & $\nu_8 \cap \nu_{11}:y$
 & $\nu_1 \cap \nu_{11}:xz^{-1}$
 & $+$ \\
$\nu_2 \cap \nu_5:x^{-1}y$
 & $\nu_5 \cap \nu_{12}:mxy^{-2}$
 & $\nu_2 \cap \nu_{12}:my^{-1}$
 & $+$ \\
$\nu_2 \cap \nu_5:1$
 & $\nu_5 \cap \nu_{12}:my^{-1}$
 & $\nu_2 \cap \nu_{12}:my^{-1}$
 & $-$ \\
$\nu_2 \cap \nu_6:my^{-1}z^{-1}$
 & $\nu_6 \cap \nu_{12}:z$
 & $\nu_2 \cap \nu_{12}:my^{-1}$
 & $-$ \\
$\nu_2 \cap \nu_8:z^{-1}$
 & $\nu_8 \cap \nu_{12}:my^{-1}z$
 & $\nu_2 \cap \nu_{12}:my^{-1}$
 & $+$ \\
$\nu_2 \cap \nu_9:my^{-2}$
 & $\nu_9 \cap \nu_{12}:y$
 & $\nu_2 \cap \nu_{12}:my^{-1}$
 & $+$ \\
$\nu_2 \cap \nu_9:my^{-1}z^{-1}$
 & $\nu_9 \cap \nu_{12}:z$
 & $\nu_2 \cap \nu_{12}:my^{-1}$
 & $-$ \\[1em]

\caption*{{\normalsize {\sc Table} \thetable. Products $\mu^2$ in the
directed $A_\infty$-category $\Cat_{12}^\rightarrow$. The notation
means that $\mu^2(a_2,a_1) = \pm b$. \label{table:mu2}}}
\\
\end{longtable}
}

{ \smaller\smaller
\begin{longtable}[l]{lll|l|l}
 $a_1$ & $a_2$ & $a_3$ & $b$ & $\pm$ \\
 \hline
 $\nu_1 \cap \nu_2:x$
 & $\nu_2 \cap \nu_3:my^{-1}$
 & $\nu_3 \cap \nu_5:x^{-1}y$
 & $\nu_1 \cap \nu_5:1$
 & $+$ \\
$\nu_2 \cap \nu_4:mx^{-1}z$
 & $\nu_4 \cap \nu_5:1$
 & $\nu_5 \cap \nu_6:my^{-1}z^{-1}$
 & $\nu_2 \cap \nu_6:mx^{-1}y^{-1}$
 & $-$ \\
$\nu_4 \cap \nu_5:1$
 & $\nu_5 \cap \nu_6:my^{-1}z^{-1}$
 & $\nu_6 \cap \nu_8:1$
 & $\nu_4 \cap \nu_8:y^{-1}z^{-1}$
 & $+$ \\
$\nu_5 \cap \nu_7:mz$
 & $\nu_7 \cap \nu_8:y^{-1}z^{-1}$
 & $\nu_8 \cap \nu_9:my^{-1}$
 & $\nu_5 \cap \nu_9:my^{-2}$
 & $-$ \\
$\nu_7 \cap \nu_8:y^{-1}z^{-1}$
 & $\nu_8 \cap \nu_9:my^{-1}$
 & $\nu_9 \cap \nu_{11}:y$
 & $\nu_7 \cap \nu_{11}:y^{-1}z^{-1}$
 & $+$ \\
$\nu_8 \cap \nu_{10}:mx^{-1}z$
 & $\nu_{10} \cap \nu_{11}:x$
 & $\nu_{11} \cap \nu_{12}:my^{-1}$
 & $\nu_8 \cap \nu_{12}:my^{-1}z$
 & $-$ \\
$\nu_1 \cap \nu_2:x$
 & $\nu_2 \cap \nu_3:my^{-1}$
 & $\nu_3 \cap \nu_8:x^{-1}yz^{-1}$
 & $\nu_1 \cap \nu_8:z^{-1}$
 & $-$ \\
$\nu_1 \cap \nu_2:x$
 & $\nu_2 \cap \nu_3:my^{-1}$
 & $\nu_3 \cap \nu_8:z^{-1}$
 & $\nu_1 \cap \nu_8:xy^{-1}z^{-1}$
 & $+$ \\
$\nu_2 \cap \nu_3:my^{-1}$
 & $\nu_3 \cap \nu_5:x^{-1}y$
 & $\nu_5 \cap \nu_7:mz$
 & $\nu_2 \cap \nu_7:mx^{-1}z$
 & $+$ \\
$\nu_2 \cap \nu_7:mz$
 & $\nu_7 \cap \nu_8:y^{-1}z^{-1}$
 & $\nu_8 \cap \nu_9:my^{-1}$
 & $\nu_2 \cap \nu_9:my^{-2}$
 & $-$ \\
$\nu_2 \cap \nu_7:mx^{-1}yz$
 & $\nu_7 \cap \nu_8:y^{-1}z^{-1}$
 & $\nu_8 \cap \nu_9:my^{-1}$
 & $\nu_2 \cap \nu_9:mx^{-1}y^{-1}$
 & $+$ \\
$\nu_1 \cap \nu_2:x$
 & $\nu_2 \cap \nu_6:my^{-1}z^{-1}$
 & $\nu_6 \cap \nu_8:1$
 & $\nu_1 \cap \nu_8:xy^{-1}z^{-1}$
 & $+$ \\
$\nu_1 \cap \nu_2:x$
 & $\nu_2 \cap \nu_6:mx^{-1}z^{-1}$
 & $\nu_6 \cap \nu_8:1$
 & $\nu_1 \cap \nu_8:z^{-1}$
 & $+$ \\
$\nu_1 \cap \nu_5:x$
 & $\nu_5 \cap \nu_6:my^{-1}z^{-1}$
 & $\nu_6 \cap \nu_8:1$
 & $\nu_1 \cap \nu_8:xy^{-1}z^{-1}$
 & $+$ \\
$\nu_1 \cap \nu_5:y$
 & $\nu_5 \cap \nu_6:my^{-1}z^{-1}$
 & $\nu_6 \cap \nu_8:1$
 & $\nu_1 \cap \nu_8:z^{-1}$
 & $+$ \\
$\nu_2 \cap \nu_4:mx^{-1}z$
 & $\nu_4 \cap \nu_5:1$
 & $\nu_5 \cap \nu_9:my^{-1}z^{-1}$
 & $\nu_2 \cap \nu_9:mx^{-1}y^{-1}$
 & $-$ \\
$\nu_2 \cap \nu_4:mx^{-1}z$
 & $\nu_4 \cap \nu_5:1$
 & $\nu_5 \cap \nu_9:mxy^{-2}z^{-1}$
 & $\nu_2 \cap \nu_9:my^{-2}$
 & $+$ \\
$\nu_2 \cap \nu_4:mx^{-1}z$
 & $\nu_4 \cap \nu_8:xy^{-1}z^{-1}$
 & $\nu_8 \cap \nu_9:my^{-1}$
 & $\nu_2 \cap \nu_9:my^{-2}$
 & $+$ \\
$\nu_2 \cap \nu_4:mx^{-1}z$
 & $\nu_4 \cap \nu_8:z^{-1}$
 & $\nu_8 \cap \nu_9:my^{-1}$
 & $\nu_2 \cap \nu_9:mx^{-1}y^{-1}$
 & $-$ \\
$\nu_3 \cap \nu_5:x^{-1}y$
 & $\nu_5 \cap \nu_7:mz$
 & $\nu_7 \cap \nu_8:y^{-1}z^{-1}$
 & $\nu_3 \cap \nu_8:x^{-1}$
 & $+$ \\
$\nu_4 \cap \nu_5:1$
 & $\nu_5 \cap \nu_6:my^{-1}z^{-1}$
 & $\nu_6 \cap \nu_{11}:y$
 & $\nu_4 \cap \nu_{11}:z^{-1}$
 & $-$ \\
$\nu_4 \cap \nu_5:1$
 & $\nu_5 \cap \nu_6:my^{-1}z^{-1}$
 & $\nu_6 \cap \nu_{11}:x$
 & $\nu_4 \cap \nu_{11}:xy^{-1}z^{-1}$
 & $+$ \\
$\nu_5 \cap \nu_6:my^{-1}z^{-1}$
 & $\nu_6 \cap \nu_8:1$
 & $\nu_8 \cap \nu_{10}:mx^{-1}z$
 & $\nu_5 \cap \nu_{10}:mx^{-1}y^{-1}$
 & $+$ \\
$\nu_5 \cap \nu_{10}:my^{-1}$
 & $\nu_{10} \cap \nu_{11}:x$
 & $\nu_{11} \cap \nu_{12}:my^{-1}$
 & $\nu_5 \cap \nu_{12}:mxy^{-2}$
 & $-$ \\
$\nu_5 \cap \nu_{10}:mx^{-1}$
 & $\nu_{10} \cap \nu_{11}:x$
 & $\nu_{11} \cap \nu_{12}:my^{-1}$
 & $\nu_5 \cap \nu_{12}:my^{-1}$
 & $+$ \\
$\nu_4 \cap \nu_5:1$
 & $\nu_5 \cap \nu_9:mxy^{-2}z^{-1}$
 & $\nu_9 \cap \nu_{11}:y$
 & $\nu_4 \cap \nu_{11}:xy^{-1}z^{-1}$
 & $+$ \\
$\nu_4 \cap \nu_5:1$
 & $\nu_5 \cap \nu_9:my^{-1}z^{-1}$
 & $\nu_9 \cap \nu_{11}:y$
 & $\nu_4 \cap \nu_{11}:z^{-1}$
 & $+$ \\
$\nu_4 \cap \nu_8:xy^{-1}z^{-1}$
 & $\nu_8 \cap \nu_9:my^{-1}$
 & $\nu_9 \cap \nu_{11}:y$
 & $\nu_4 \cap \nu_{11}:xy^{-1}z^{-1}$
 & $+$ \\
$\nu_4 \cap \nu_8:z^{-1}$
 & $\nu_8 \cap \nu_9:my^{-1}$
 & $\nu_9 \cap \nu_{11}:y$
 & $\nu_4 \cap \nu_{11}:z^{-1}$
 & $+$ \\
$\nu_5 \cap \nu_7:mz$
 & $\nu_7 \cap \nu_8:y^{-1}z^{-1}$
 & $\nu_8 \cap \nu_{12}:m$
 & $\nu_5 \cap \nu_{12}:my^{-1}$
 & $-$ \\
$\nu_5 \cap \nu_7:mz$
 & $\nu_7 \cap \nu_8:y^{-1}z^{-1}$
 & $\nu_8 \cap \nu_{12}:mxy^{-1}$
 & $\nu_5 \cap \nu_{12}:mxy^{-2}$
 & $+$ \\
$\nu_5 \cap \nu_7:mz$
 & $\nu_7 \cap \nu_{11}:xy^{-1}z^{-1}$
 & $\nu_{11} \cap \nu_{12}:my^{-1}$
 & $\nu_5 \cap \nu_{12}:mxy^{-2}$
 & $+$ \\
$\nu_5 \cap \nu_7:mz$
 & $\nu_7 \cap \nu_{11}:z^{-1}$
 & $\nu_{11} \cap \nu_{12}:my^{-1}$
 & $\nu_5 \cap \nu_{12}:my^{-1}$
 & $-$ \\
$\nu_6 \cap \nu_8:1$
 & $\nu_8 \cap \nu_{10}:mx^{-1}z$
 & $\nu_{10} \cap \nu_{11}:x$
 & $\nu_6 \cap \nu_{11}:z$
 & $+$ \\
$\nu_1 \cap \nu_2:x$
 & $\nu_2 \cap \nu_3:my^{-1}$
 & $\nu_3 \cap \nu_{11}:yz^{-1}$
 & $\nu_1 \cap \nu_{11}:xz^{-1}$
 & $-$ \\
$\nu_2 \cap \nu_{10}:mx^{-1}$
 & $\nu_{10} \cap \nu_{11}:x$
 & $\nu_{11} \cap \nu_{12}:my^{-1}$
 & $\nu_2 \cap \nu_{12}:my^{-1}$
 & $-$ \\
$\nu_1 \cap \nu_2:x$
 & $\nu_2 \cap \nu_6:mx^{-1}z^{-1}$
 & $\nu_6 \cap \nu_{11}:x$
 & $\nu_1 \cap \nu_{11}:xz^{-1}$
 & $-$ \\
$\nu_1 \cap \nu_2:x$
 & $\nu_2 \cap \nu_6:my^{-1}z^{-1}$
 & $\nu_6 \cap \nu_{11}:y$
 & $\nu_1 \cap \nu_{11}:xz^{-1}$
 & $-$ \\
$\nu_1 \cap \nu_2:x$
 & $\nu_2 \cap \nu_9:my^{-1}z^{-1}$
 & $\nu_9 \cap \nu_{11}:y$
 & $\nu_1 \cap \nu_{11}:xz^{-1}$
 & $+$ \\
$\nu_1 \cap \nu_5:y$
 & $\nu_5 \cap \nu_6:my^{-1}z^{-1}$
 & $\nu_6 \cap \nu_{11}:x$
 & $\nu_1 \cap \nu_{11}:xz^{-1}$
 & $-$ \\
$\nu_1 \cap \nu_5:x$
 & $\nu_5 \cap \nu_6:my^{-1}z^{-1}$
 & $\nu_6 \cap \nu_{11}:y$
 & $\nu_1 \cap \nu_{11}:xz^{-1}$
 & $-$ \\
$\nu_2 \cap \nu_3:my^{-1}$
 & $\nu_3 \cap \nu_5:x^{-1}y$
 & $\nu_5 \cap \nu_{10}:mx^{-1}$
 & $\nu_2 \cap \nu_{10}:mx^{-2}$
 & $-$ \\
$\nu_2 \cap \nu_3:my^{-1}$
 & $\nu_3 \cap \nu_5:x^{-1}y$
 & $\nu_5 \cap \nu_{10}:my^{-1}$
 & $\nu_2 \cap \nu_{10}:mx^{-1}y^{-1}$
 & $-$ \\
$\nu_2 \cap \nu_3:my^{-1}$
 & $\nu_3 \cap \nu_8:z^{-1}$
 & $\nu_8 \cap \nu_{10}:mx^{-1}z$
 & $\nu_2 \cap \nu_{10}:mx^{-1}y^{-1}$
 & $+$ \\
$\nu_2 \cap \nu_3:my^{-1}$
 & $\nu_3 \cap \nu_8:x^{-1}yz^{-1}$
 & $\nu_8 \cap \nu_{10}:mx^{-1}z$
 & $\nu_2 \cap \nu_{10}:mx^{-2}$
 & $-$ \\
$\nu_2 \cap \nu_6:my^{-1}z^{-1}$
 & $\nu_6 \cap \nu_8:1$
 & $\nu_8 \cap \nu_{10}:mx^{-1}z$
 & $\nu_2 \cap \nu_{10}:mx^{-1}y^{-1}$
 & $+$ \\
$\nu_2 \cap \nu_6:mx^{-1}z^{-1}$
 & $\nu_6 \cap \nu_8:1$
 & $\nu_8 \cap \nu_{10}:mx^{-1}z$
 & $\nu_2 \cap \nu_{10}:mx^{-2}$
 & $+$ \\
$\nu_2 \cap \nu_4:mx^{-1}z$
 & $\nu_4 \cap \nu_8:xy^{-1}z^{-1}$
 & $\nu_8 \cap \nu_{12}:m$
 & $\nu_2 \cap \nu_{12}:my^{-1}$
 & $-$ \\
$\nu_2 \cap \nu_4:mx^{-1}z$
 & $\nu_4 \cap \nu_8:z^{-1}$
 & $\nu_8 \cap \nu_{12}:mxy^{-1}$
 & $\nu_2 \cap \nu_{12}:my^{-1}$
 & $+$ \\
$\nu_3 \cap \nu_5:x^{-1}y$
 & $\nu_5 \cap \nu_7:mz$
 & $\nu_7 \cap \nu_{11}:z^{-1}$
 & $\nu_3 \cap \nu_{11}:x^{-1}y$
 & $+$ \\
$\nu_3 \cap \nu_5:x^{-1}y$
 & $\nu_5 \cap \nu_7:mz$
 & $\nu_7 \cap \nu_{11}:xy^{-1}z^{-1}$
 & $\nu_3 \cap \nu_{11}:1$
 & $+$ \\
$\nu_1 \cap \nu_8:xz^{-1}$
 & $\nu_8 \cap \nu_9:my^{-1}$
 & $\nu_9 \cap \nu_{11}:y$
 & $\nu_1 \cap \nu_{11}:xz^{-1}$
 & $-$ \\
$\nu_3 \cap \nu_5:x^{-1}y$
 & $\nu_5 \cap \nu_{10}:mx^{-1}$
 & $\nu_{10} \cap \nu_{11}:x$
 & $\nu_3 \cap \nu_{11}:x^{-1}y$
 & $-$ \\
$\nu_3 \cap \nu_5:x^{-1}y$
 & $\nu_5 \cap \nu_{10}:my^{-1}$
 & $\nu_{10} \cap \nu_{11}:x$
 & $\nu_3 \cap \nu_{11}:1$
 & $-$ \\
$\nu_3 \cap \nu_8:x^{-1}yz^{-1}$
 & $\nu_8 \cap \nu_{10}:mx^{-1}z$
 & $\nu_{10} \cap \nu_{11}:x$
 & $\nu_3 \cap \nu_{11}:x^{-1}y$
 & $-$ \\
$\nu_3 \cap \nu_8:z^{-1}$
 & $\nu_8 \cap \nu_{10}:mx^{-1}z$
 & $\nu_{10} \cap \nu_{11}:x$
 & $\nu_3 \cap \nu_{11}:1$
 & $+$ \\
$\nu_1 \cap \nu_5:x$
 & $\nu_5 \cap \nu_9:my^{-1}z^{-1}$
 & $\nu_9 \cap \nu_{11}:y$
 & $\nu_1 \cap \nu_{11}:xz^{-1}$
 & $+$ \\
$\nu_1 \cap \nu_5:y$
 & $\nu_5 \cap \nu_9:mxy^{-2}z^{-1}$
 & $\nu_9 \cap \nu_{11}:y$
 & $\nu_1 \cap \nu_{11}:xz^{-1}$
 & $-$ \\
$\nu_2 \cap \nu_4:mx^{-1}z$
 & $\nu_4 \cap \nu_5:1$
 & $\nu_5 \cap \nu_{12}:mxy^{-1}z^{-1}$
 & $\nu_2 \cap \nu_{12}:my^{-1}$
 & $-$ \\
$\nu_2 \cap \nu_4:mx^{-1}z$
 & $\nu_4 \cap \nu_{11}:xz^{-1}$
 & $\nu_{11} \cap \nu_{12}:my^{-1}$
 & $\nu_2 \cap \nu_{12}:my^{-1}$
 & $+$ \\
$\nu_2 \cap \nu_7:mx^{-1}yz$
 & $\nu_7 \cap \nu_8:y^{-1}z^{-1}$
 & $\nu_8 \cap \nu_{12}:mxy^{-1}$
 & $\nu_2 \cap \nu_{12}:my^{-1}$
 & $-$ \\
$\nu_2 \cap \nu_7:mz$
 & $\nu_7 \cap \nu_8:y^{-1}z^{-1}$
 & $\nu_8 \cap \nu_{12}:m$
 & $\nu_2 \cap \nu_{12}:my^{-1}$
 & $+$ \\
$\nu_2 \cap \nu_7:mz$
 & $\nu_7 \cap \nu_{11}:z^{-1}$
 & $\nu_{11} \cap \nu_{12}:my^{-1}$
 & $\nu_2 \cap \nu_{12}:my^{-1}$
 & $+$ \\
$\nu_2 \cap \nu_7:mx^{-1}yz$
 & $\nu_7 \cap \nu_{11}:xy^{-1}z^{-1}$
 & $\nu_{11} \cap \nu_{12}:my^{-1}$
 & $\nu_2 \cap \nu_{12}:my^{-1}$
 & $-$ \\[1em]

\caption*{{\normalsize {\sc Table} \thetable. Products $\mu^3$ in
$\Cat_{12}^\rightarrow$, written in the same way as in the previous
table. \label{table:mu3}}}
\\
\end{longtable}
}

\backmatter

\end{document}